\documentclass[11pt,letterpaper]{article}
\usepackage[margin=1in]{geometry}
\usepackage{parskip}
\usepackage[utf8]{inputenc}
\usepackage[T1]{fontenc} %
\usepackage{mathpazo}
\usepackage{microtype}
\usepackage{upgreek}
\renewcommand{\gamma}{\upgamma}
\renewcommand{\pi}{\uppi}

\usepackage{derivative}
\usepackage{amsmath,amssymb,amsthm,mathtools}
\usepackage{empheq}
\usepackage{graphicx}
\usepackage{subcaption}
\usepackage{bm}
\usepackage{alignnotes}
\usepackage{nicefrac}
\usepackage{aliascnt}
\usepackage[ruled,vlined,linesnumbered]{algorithm2e}

\allowdisplaybreaks[1]

\theoremstyle{plain}
\newtheorem{theorem}{Theorem}[section]

\newaliascnt{proposition}{theorem}
\newtheorem{proposition}[proposition]{Proposition}
\aliascntresetthe{proposition}

\newaliascnt{lemma}{theorem}
\newtheorem{lemma}[lemma]{Lemma}
\aliascntresetthe{lemma}

\newaliascnt{corollary}{theorem}
\newtheorem{corollary}[corollary]{Corollary}
\aliascntresetthe{corollary}

\theoremstyle{definition}
\newaliascnt{definition}{theorem}
\newtheorem{definition}[definition]{Definition}
\aliascntresetthe{definition}

\newaliascnt{assumption}{theorem}
\newtheorem{assumption}[assumption]{Assumption}
\aliascntresetthe{assumption}

\newaliascnt{hyperparameters}{theorem}
\newtheorem{hyperparameters}[hyperparameters]
  {Hyperparameters}
\aliascntresetthe{hyperparameters}

\theoremstyle{remark}
\newaliascnt{remark}{theorem}

\aliascntresetthe{remark}

\usepackage{thm-restate}

\SetKwInput{KwIn}{Input}
\SetKwInput{KwOut}{Output}
\SetKw{Return}{return}
\SetKwFunction{TrimmedMean}{TrimmedMean}
\SetKwFunction{ResidualMean}{ResidualMean}
\SetKwFunction{ProjectedMomentFamily}{EstimateProjectedMoments}
\SetKwFunction{SelectAndNormalize}{ConstructObservableOperators}
\SetKwFunction{DiagonalizeAndRead}{RecoverDirections}
\SetKwFunction{ImplicitTensorDecomposition}{ImplicitTensorDecomposition}
\SetKwFunction{Refine}{RestrictedLeastSquares}
\SetKwFunction{RecoverAll}{RecoverAllRegressors}

\makeatletter
\@addtoreset{algocf}{section}

\makeatother

\usepackage{booktabs}

\usepackage{float}

\usepackage[shortlabels]{enumitem}

\usepackage{cancel}

\usepackage{manfnt} %

\usepackage{figchild}
\usepackage{fontawesome5}

\usepackage{accents}

\usepackage[dvipsnames]{xcolor}
\usepackage{color}
\definecolor{niceRed}{RGB}{190,38,38}
\definecolor{niceYellow}{HTML}{f5b400}
\definecolor{blueGrotto}{HTML}{059DC0}
\definecolor{royalBlue}{HTML}{057DCD}
\definecolor{navyBlue}{HTML}{0B579C}
\definecolor{yaleBlue}{HTML}{00356b}
\definecolor{limeGreen}{HTML}{81B622}
\definecolor{nicePurple}{HTML}{9c27b0}
\definecolor{lightRoyalBlue}{HTML}{def2ff}
\definecolor{gold}{HTML}{ffa300}
\definecolor{claudeBlue}{HTML}{4C72B0}
\definecolor{claudeRed}{HTML}{C44E52}

\SetKwFunction{PeeledMomentMean}{EstimateMomentsByPeeling}

\newcommand*{\red}[1]{\textcolor{red}{#1}}

\usepackage[most]{tcolorbox}

\usepackage{xfp} %

\newcommand{\R}{\mathbb R}
\newcommand{\E}{\mathbb E}
\newcommand{\Prb}{\mathbb P}
\newcommand{\ip}[2]{\left\langle #1,#2\right\rangle}
\newcommand{\norm}[1]{\left\lVert #1\right\rVert}
\newcommand{\abs}[1]{\left\lvert #1\right\rvert}
\newcommand{\Mat}{\operatorname{Mat}}
\newcommand{\Span}{\operatorname{span}}
\newcommand{\He}{\operatorname{He}}
\newcommand{\PhiN}{\Phi}
\newcommand{\Ot}{\widetilde O}
\newcommand{\1}{\mathbf 1}
\newcommand{\Sym}{\operatorname{Sym}}
\newcommand{\iu}{\mathrm i}
\newcommand{\wh}{\widehat}
\newcommand{\eps}{\varepsilon}

\newcommand{\Regressor}{\bm{w}}
\newcommand{\NormalizedRegressor}{\bm{u}}
\newcommand{\ScaledRegressor}{\bm{a}}
\newcommand{\Covariate}{\bm{X}}
\newcommand{\NoiseVector}{\bm{\xi}}
\newcommand{\HermiteTensor}[2]{\mathsf H_{#1}(#2)}
\newcommand{\RestrictedMomentTensor}[2]{\mathsf T_{#1}^{#2}}
\newcommand{\MomentRemainder}[1]{\mathsf R_{#1}}
\newcommand{\TensorPower}[2]{{#1}^{\otimes #2}}
\newcommand{\MomentCoefficient}[2]{a_{#1}(#2)}
\newcommand{\ComponentWeight}[2]{\lambda_{#2}^{(#1)}}
\newcommand{\ProjectionMap}[1]{C_{#1}}
\newcommand{\EvenMomentMatrix}[2]{A_{#2}^{(#1)}}
\newcommand{\OddMomentMatrix}[3]{M_{#2}^{(#1)}(#3)}
\newcommand{\RefinementEvent}[1]{\mathcal E_{#1}}
\newcommand{\RestrictedSecondMoment}[1]{\Sigma_{#1}}
\newcommand{\RestrictedResponse}[1]{\bm c_{#1}}
\newcommand{\RefinementBias}[1]{\bm e_{#1}}

\newcommand{\TailIndicator}[1]{%
  \1\{T\leq #1\leq T+1\}}
\newcommand{\HermitePairing}[3]{%
  \ip{\HermiteTensor{#1}{#2}}{#3}}
\newcommand{\HermiteContraction}[3]{%
  \HermiteTensor{#1}{#2}\mathbin{\lrcorner}#3}

\newcommand{\match}{\mathrm{match}}
\newcommand{\op}{\mathrm{op}}

\newcommand{\refine}{\mathrm{ref}}
\newcommand{\odd}{\mathrm{odd}}
\newcommand{\even}{\mathrm{even}}

\DeclareMathOperator{\RealPart}{Re}
\DeclareMathOperator{\Range}{range}

\DeclareMathOperator{\Diag}{diag}
\DeclareMathOperator{\tr}{tr}

\DeclareMathOperator{\poly}{poly}

\DeclareMathOperator{\argmax}{argmax}

\usepackage{hyperref}
\hypersetup{
  colorlinks = true,
  urlcolor = {blueGrotto},
  linkcolor = {royalBlue},
  citecolor = {navyBlue},
  pdftitle = {Linear Regression with Self-Selection Bias in Polynomial Time},
  pdfauthor = {Felix Zhou}
}

\usepackage[capitalise,noabbrev,sort,nameinlink]{cleveref}

\crefname{proposition}{Proposition}{Propositions}
\crefname{lemma}{Lemma}{Lemmas}
\crefname{corollary}{Corollary}{Corollaries}
\crefname{definition}{Definition}{Definitions}
\crefname{assumption}{Assumption}{Assumptions}
\crefname{hyperparameters}{Hyperparameters}{Hyperparameters}
\Crefname{hyperparameters}{Hyperparameters}{Hyperparameters}
\crefname{remark}{Remark}{Remarks}
\crefname{algocf}{Algorithm}{Algorithms}
\Crefname{algocf}{Algorithm}{Algorithms}

\usepackage[
  backref=true,
  backend=biber,
  natbib=true,
  style=alphabetic,
  sorting=alphabeticlabel,
  sortcites=true,
  minbibnames=3,
  maxbibnames=999,
  mincitenames=4,
  maxcitenames=4,
  minalphanames=4,
  maxalphanames=4,
  url=false,
  doi=false,
]{biblatex}

\DeclareSortingTemplate{alphabeticlabel}{
  \sort[final]{%
    \field{labelalpha} %
  }
  \sort{%
    \field{year}   %
  }
  \sort{%
    \field{title}  %
  }
}

\AtBeginRefsection{\GenRefcontextData{sorting=ynt}}
\AtEveryCite{\localrefcontext[sorting=ynt]}
\title{Linear Regression with Self-Selection Bias\\in Polynomial Time}
\author{\begin{tabular}{c}
  Felix Zhou\\
  Yale University\\
  \texttt{felix.zhou@yale.edu}
\end{tabular}}
\date{}

\begin{document}
\pagenumbering{roman}
\pagestyle{empty}
\maketitle
\thispagestyle{empty}

\begin{abstract}
\noindent We revisit the problem of estimating $k$ linear regressors with self-selection bias in $n$ dimensions up to error $\eps$ with the maximum self-selection criterion,
as introduced by \citet[STOC]{CDIZ23}.
Our main result is a $\poly(n, k, \nicefrac1\eps)$ time algorithm for the unknown-index variant of this problem,
thereby improving upon the running times of the best known algorithms due to \citet{GM24} and \citet[COLT]{KMZ26}
with exponential dependence on $k$.
Our algorithm succeeds in the presence of correlated noise
between regressors, which captures the related setting of max-affine
regression. Consequently, we also give the first polynomial-time
algorithm for this problem, improving on the prior works of
\citet[IEEE Transactions on Information Theory]{GPGR22} and \citet{GM24}
with exponential dependence on $k$.

Prior works exhaustively search over a $k$-dimensional subspace spanned by the regressors,
using conditional moments of the response to test candidates.
In contrast, our algorithm is driven by an improved understanding of the moments of the covariates.
We show that, restricting to observations with large response values, 
the degree-$(\log n)$ Hermite moment of the covariates
approximately decomposes as a sum of
tensor powers of the regressors.
The first step of our algorithm can be viewed as a tensor decomposition algorithm for these implicitly represented high-order tensors,
which outputs unit vectors parallel to the regressors.
A secondary step based on a simple variant of ordinary least squares extracts estimates for the regressors given these unit vectors.
We believe our analytic insight and algorithmic components may have applications to designing efficient algorithms for other latent-variable problems.
\end{abstract}

\clearpage
\tableofcontents

\clearpage
\pagestyle{plain}
\pagenumbering{arabic}
\section{Introduction}

\begin{figure}[ht]
  \centering
  \begin{subfigure}{0.5\textwidth}
    \centering
    \includegraphics[width=\linewidth]
      {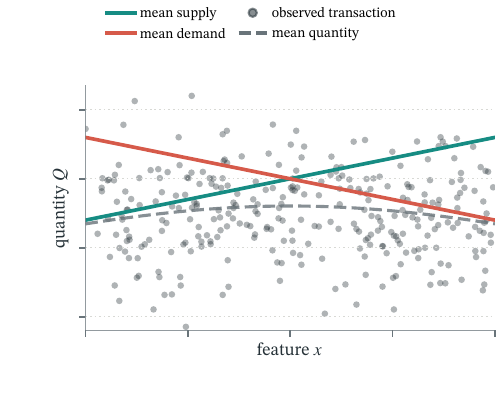}
    \caption{\small Minimizer of supply/demand is unobserved.}
    \label{subfig:intro:fair-jaffee-unknown-index}
  \end{subfigure}\hfill
  \begin{subfigure}{0.5\textwidth}
    \centering
    \includegraphics[width=\linewidth]
      {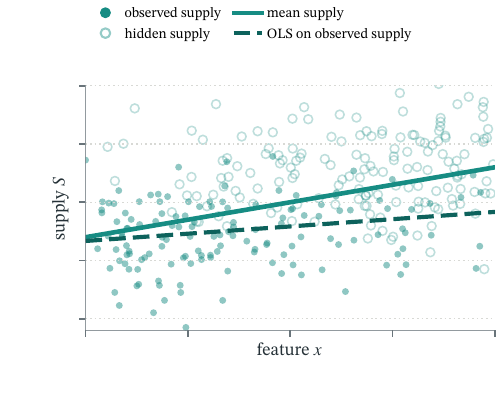}
    \caption{\small OLS on partitioned data is biased.}
    \label{subfig:intro:fair-jaffee-known-index}
  \end{subfigure}
  \caption{\small For a one-dimensional \citet{FJ72} market with independent
  Gaussian noises and no intercepts,
  the observed transactions can obscure the mean supply and demand
  curves when the minimizer is hidden (a). Selection bias occurs
  even with a known minimizer and vanilla OLS regression on the
  observed supply can lead to biased estimates (b).}
\end{figure}

In their seminal work on disequilibrium econometrics,
\citet{FJ72} studied a linear model of the housing market where there is a supply $S(\bm x, \xi_S) = \Regressor_S^\top \bm x + \xi_S$ and a demand $D(\bm x, \xi_D) = \Regressor_D^\top \bm x + \xi_D$ of houses with feature $\bm x$,
but supply does not equal demand.
The actual quantity transacted is $Q(\bm x) = \min\{S(\bm x, \xi_S), D(\bm x, \xi_D)\}$ where $(\xi_S, \xi_D)$ are random noise.
A statistician aiming to estimate the regressors $\Regressor_D, \Regressor_S$ from observations $(\bm x, Q(\bm x))$ faces difficulties with vanilla linear regression methods:
\begin{enumerate}[(1)]
  \item We do not observe whether a given transaction is constrained by supply or demand,
    so observations cannot be separated into supply and demand regression samples~(see \Cref{subfig:intro:fair-jaffee-unknown-index});
  \item Even if this is known, partitioning the data and running an ordinary least squares (OLS) estimator on all the supply-constrained data to estimate $\Regressor_S$ will likely yield biased estimates~(see \Cref{subfig:intro:fair-jaffee-known-index}).
    Similar concerns hold for the demand.
\end{enumerate}
To see the second point, note that the statistician will never observe the supply $S(\bm x, \xi_S)$ for houses with low demand $S(\bm x, \xi_S)\gg D(\bm x, \xi_D)$ and vice versa.
This leads to biased observations in the partitioned datasets which in turn induce biases in the OLS estimators on those datasets.

\paragraph{Unknown-Index Maximum Self-Selection.}
The example above is an instance of linear regression with \emph{self-selection} bias.
In the general setting, we observe features $\Covariate$ accompanied by an outcome $Y\in \{Y_1,\dots,Y_k\}$ selected from the responses of $k$ underlying affine models:
$
  Y_i = \Regressor_i^\top \Covariate + b_i + \xi_i
$.
In this work, we focus on the maximum self-selection mechanism, for which the observed response is $Y_{\max}\coloneqq\max_{i\in[k]}Y_i$.\footnote{The minimum in the model of \citet{FJ72} becomes
a maximum after negating the observed response, parameter
vectors, biases, and noises.}
Importantly, to capture the setting of \citet{FJ72},
we do \emph{not} assume the index of the selected model is known.
This is also known as the \emph{unknown-index} setting.

More broadly, self-selection also arises in occupational choice,
where workers choose among occupations with different potential
earnings and only the earnings from the chosen occupation are
observed~\citep{Roy51}.
In observational medicine, treatment choices depend on prognostic
information, which can bias comparisons of outcomes
\citep{GreenlandNeutra80,RosenbaumRubin83}.
A related inverse problem arises in inverse reinforcement learning,
where one seeks a reward function from an expert's observed
behavior~\citep{NgRussell00}.
Other examples arise in auctions and voluntary disclosure.
Auction data may reveal only the transaction price, leaving the
winning bidder and losing bids unobserved~\citep{AH02}.
Finally, models of voluntary disclosure study how informed sellers choose
what information about product quality to reveal
\citep{Grossman81,Milgrom81}.

\paragraph{Computational and Statistical Efficiency.}
Despite the long history and use of self-selection models,
algorithms with provable guarantees for finite samples have only recently been studied.
\citet{CDIZ23} initiated the study of sample- and time-efficient algorithms for linear regression under the maximum self-selection bias with $k$ regressors in the unknown-index setting.
If $n$ is the ambient dimension, they gave an algorithm with
$\poly(n)\exp(\poly(k))$ sample complexity and running time
that achieves $\poly(\nicefrac1k)$ error.

The algorithm of \citet{GM24} uses
$\poly(n,k,\nicefrac1\eps)$ samples and runs in time
$\poly(n,k)+(\nicefrac1\eps)^{\Ot(k)}$.
Using this algorithm to supply a warm start to the local-convergence algorithm
of \citet{KMZ26} further improves the running time to $\poly(n,k,\nicefrac1\eps)+2^{\Ot(k)}$.

Henceforth, \textbf{we refer to unknown-index linear regression with maximum self-selection mechanism as the \emph{self-selection problem}}.
Our main result is the first algorithm for the self-selection problem with both $\poly(n, k, \nicefrac1\eps)$ sample complexity and running time,
improving on the running time of all prior works.

\subsection{Our Contributions}
Before stating our results, let us introduce the formal model of the self-selection problem.
\begin{definition}[Maximum Self-Selection; \cite{CDIZ23}]
\label{def:gaussian-self-selection}
Fix integers $n\geq1$ and $2\leq k\leq n$.\footnote{The case $k=1$ is just ordinary linear regression. By padding the covariate if necessary, we may assume without loss of generality that $k\leq n$.} 
We wish to recover unknown regressors
$\Regressor_1,\ldots,\Regressor_k\in\R^n$ and biases $b_1,\ldots,b_k\in\R$,
given independent samples $(\Covariate,Y_{\max})$,
such that
\begin{equation}
\begin{aligned}
 \Covariate&\sim N(\bm0,I_n),\qquad
 &\xi_i&\sim N(0,1)\qquad(i\in[k]),\\
 Y_i&=\Regressor_i^\top\Covariate+b_i+\xi_i,\qquad
 &Y_{\max}&=\max_{i\in[k]}Y_i.
\end{aligned}
\end{equation}
The noises are jointly Gaussian, with \emph{known} correlation
$\operatorname{Cov}(\xi_i,\xi_j)=\rho\in[0,1]$ for $i\ne j$.
Moreover, $\NoiseVector=(\xi_i)_{i\in[k]}$ is
independent of the Gaussian covariate $\Covariate$ and the maximizing index $\argmax_{i\in [k]} Y_i$ is \emph{not} observed.
\end{definition}
When $\rho=0$ and all $b_i=0$, we recover the original
self-selection model~\citep{CDIZ23}.
When $\rho=1$, all noises coincide with a single
$\xi\sim N(0,1)$, so
\[
 Y_{\max}
 =\max_{i\in[k]}\left\{ \Regressor_i^\top\Covariate+b_i \right\}+\xi.
\]
This is known as \emph{max-affine regression} with Gaussian noise~\citep{MagnaniBoyd09}.

In order to design algorithms with finite-sample guarantees,
we impose some standard geometric assumptions introduced in \cite{CDIZ23} and shown to be necessary in \cite{GM24}.
\begin{assumption}[Boundedness and Uncoveredness]
\label{ass:boundedness-uncoveredness}
Let $\norm{\cdot}_2$ denote the Euclidean norm.
For some absolute constants $B\geq1$ and $\Delta\in(0,1]$, the regressors and biases satisfy
\begin{equation}
 \Delta\leq\norm{\Regressor_i}_2\leq B,\qquad |b_i|\leq B,
 \qquad
 \abs{\ip{\Regressor_i}{\Regressor_j}}
 \leq
 \min\{\norm{\Regressor_i}_2^2,\norm{\Regressor_j}_2^2\}-\Delta^2
 \quad(i\neq j).
\end{equation}
\end{assumption}
A simple consequence of \Cref{ass:boundedness-uncoveredness} that we use
repeatedly is the angular separation bound
\begin{equation}
 \frac{\abs{\ip{\Regressor_i}{\Regressor_j}}}
      {\norm{\Regressor_i}_2\norm{\Regressor_j}_2}
 \leq1-\frac{\Delta^2}{B^2}
 \qquad(i\ne j).
\end{equation}

We are now ready to formally state our main recovery result.
As the self-selection model is invariant under permutation,
we recover the regressors up to some permutation.
\begin{theorem}[Polynomial-Time Self-Selection]\label{thm:intro:polynomial-time-recovery}
There is a randomized
algorithm with the following guarantee.
Given $k\geq 2$, known correlation $\rho\in[0,1]$, desired error $\varepsilon\in(0,1)$, failure probability $\delta\in(0,\nicefrac13)$,
and sample access to an instance of the self-selection model (\Cref{def:gaussian-self-selection}) satisfying boundedness and uncoveredness conditions (\Cref{ass:boundedness-uncoveredness}) with constants $\Delta, B$,
the algorithm returns estimates $(\widehat{\Regressor}_i,\widehat b_i)_{i\in[k]}$ such that, with
probability at least $1-\delta$,
there is a permutation $\pi$ of the indices under which
\[
 \forall i\in [k]\,,\quad
 \norm{\widehat\Regressor_i-\Regressor_{\pi(i)}}_2
 +\left| \widehat b_i-b_{\pi(i)} \right| \leq\varepsilon\,.
\]
Its sample complexity and running time are both
\[
 \Ot\!\left(
 \left(\frac{nk}{\varepsilon}\right)^{\poly(\nicefrac{B}{\Delta})}
 \log(\nicefrac1\delta)
 \right)\,.
\]
\end{theorem}
As we mentioned, this yields the first polynomial-time algorithm for the self-selection problem.

Taking $\rho=1$ also gives the first polynomial-time algorithm
for max-affine regression, with $\poly(n,k,\nicefrac1\eps)$ sample
complexity and running time.
The max-affine recovery algorithm of \citet{GPGR22} runs in
$\poly(n,k,\nicefrac1\eps)\exp\!\bigl(\Ot(k^2)\bigr)$ time.
For max-linear regression (bias is $\bm0$), \citet{GM24} improve
the running time to $\poly(n,k,\nicefrac1\eps)+k^{O(k)}$,
with polynomial sample complexity.

\subsubsection{Technical Contributions}
Our algorithmic improvements are driven by an improved understanding of the moments of the self-selection problem.
The key insight is that, conditional on large responses,
the conditional distribution of the \emph{covariate} is essentially a $k$-mixture distribution over the $k$ regressors (see \Cref{fig:overview:self-selection}).
Mathematically, the higher-order conditional Hermite moments approximately decompose as a weighted sum of nearly orthogonal rank-one tensors.
This departs from the conditional moment-testing paradigm from prior work, which exhaustively searches over a $k$-dimensional subspace spanned by the $k$ regressors using conditional moments of the \emph{response} $Y_{\max}$.

{Unlike the exact Hermite decompositions used for spherical
Gaussian mixtures~\citep{LL24,DK25}, our decomposition contains
irreducible error even with exact moments, which can wash out the
contribution of regressors with relatively small norms.
We therefore aim to recover a subset with the largest norms
at each step.
We carefully combine implicit projections~\citep{LL24,DK25} with
and spectral methods~\citep{AHK12,HK13} to extract unit vectors
approximately parallel to these regressors.
The main difficulty is to show that these ideas are robust to the
irreducible error.
Post-processing by a simple variant of ordinary least squares then
yields accurate estimates of the corresponding regressors and biases.
This can also be viewed as a local convergence algorithm, improving
on the initialization radius of \citet{KMZ26}.
To recover the rest, we approximately simulate the restricted
Hermite moments of the self-selection problem induced by the
remaining regressors, allowing us to recurse on the smaller instance.}

\section{Technical Overview}
We now present the key technical ingredients of our algorithm.
For the sake of clarity, it is easiest to explain our algorithm in terms of true population expectations and ignore both sampling error and small bias terms. 
We will also treat $B, \Delta$ as absolute constants that disappears in the big-$O$ notation. 

\paragraph{Notation.}
In this overview, we use bold letters for vectors, such as
$\Covariate,\Regressor_i,\bm{\zeta}$;
ordinary capital letters for matrices, such as
$A,\ProjectionMap{r},M(\bm{\zeta})$; and
sans-serif capital letters for tensors, such as
$\HermiteTensor{d}{\Covariate}$ and $\MomentRemainder{d}$.
Write $\NoiseVector=(\xi_i)_{i\in[k]}$ and
$\NoiseVector_I=(\xi_i)_{i\in I}$ for $I\subseteq[k]$.
For a positive integer $m$, write $[m]=\{1,\ldots,m\}$.
We use $\ip{\cdot}{\cdot}$ for the Euclidean inner product,
and $\norm{\cdot}_2,\norm{\cdot}_F,\norm{\cdot}_{\op}$
for the Euclidean, Frobenius, and operator norms.
The symbol $\1\{\mathcal E\}$ denotes the indicator of an
event $\mathcal E$.

We use $\otimes$ for tensor products and
$\bm v^{\otimes d}$ for the $d$-fold tensor product of $\bm v$.
The space $\Sym^d(\R^n)$ consists of symmetric order-$d$
tensors on $\R^n$.
The notation $\mathsf A\mathbin{\lrcorner}\bm{\zeta}$ denotes
contraction against $\bm{\zeta}$ in the last tensor factor.
We write $\iu=\sqrt{-1}$ and use $\RealPart$ for the
entrywise real part.
Finally, $\Phi$ denotes the standard Gaussian CDF and
$\phi_{\sigma^2}$ the density of $N(0,\sigma^2)$.

\subsection{Prior Techniques}
Before we begin, it is instructive to recall prior approaches.

\paragraph{Moments of Responses.} \citet{CDIZ23} 
test candidate directions using higher-order moments of the response $Y_{\max}$.
They restrict to samples whose projected covariates lie in a thin tube around the candidate line.
This gives rise to an exhaustive search algorithm that tests every member of an $\eps$-net using conditional moments.
A straightforward implementation would require $(\nicefrac1\eps)^{\Omega(n)}$ samples and running time,
but the authors were able to reduce this to $\poly(n)\exp(\poly(k))$ samples and running time
for $\eps = \poly(\nicefrac1k)$ accuracy, 
by restricting the search space to a $k$-dimensional subspace approximately spanned by the regressors.

Building on these ideas, \citet{GM24} test a candidate $\wh{\Regressor}$ by conditioning on a large projection of $\Covariate$ along the candidate direction.
This event has inverse-polynomial probability, while the tube in the earlier test can have exponentially small probability in $k$.
They only need conditional moments up to degree two.
As a consequence, empirical estimates of these moments concentrate and they can reuse samples across tests.
This reduces the sample complexity to $\poly(n, k, \nicefrac1\eps)$,
but the running time remains $\poly(n,k,\nicefrac1\eps)+(\nicefrac1\eps)^{\Ot(k)}$.

\paragraph{Local Convergence.}
\citet{KMZ26} 
showed that stochastic gradient descent on the natural log-likelihood function of observations converges,
assuming the initialization point is within $\poly(\nicefrac1k)$-distance of the true regressors.
Combined with the previous algorithms for the warm start,
their algorithm maintains the $\poly(n,k,\nicefrac1\eps)$ sample complexity
and improves the running time to $\poly(n,k,\nicefrac1\eps)+2^{\Ot(k)}$.

In this work, we design an algorithm with $\poly(n, k, \nicefrac1\eps)$ sample complexity and running time,
improving on the running time of all prior works.
As part of our algorithm we will design a local convergence subroutine that only requires the initialization point to be within $O(1)$-distance of the true regressors,
improving on the convergence radius of \cite{KMZ26}.

\subsection{Our Approach}
We now present an overview of our techniques.
As mentioned, it is easiest to state our algorithm in terms of exact expectations and ignore small bias terms.
The actual analysis uses samples to estimate this expectation and carefully track the bias.
We defer these details to the respective technical sections.

We first show that restricting the maximum response to a tail
slice exposes the regressor directions in the Hermite moments of the
covariates (\Cref{subsubsec:overview:hermite-moments}).
We then explain how to project these moments without forming the full
tensors (\Cref{subsubsec:overview:tensor-projections}) and use
observable-operator learning to recover a subset of approximate
regressor directions (\Cref{subsubsec:overview:observable-operator-learning}).
Post-processing these unit vectors using a variant of ordinary least squares regression yields estimates of
the full regressors and biases (\Cref{subsubsec:overview:restricted-least-squares}).
Finally, we show that we can simulate the Hermite moments of the remaining
regressors so that we can repeat the procedure on the reduced instance
(\Cref{subsubsec:overview:moment-peeling}).

\subsubsection{Hermite Moments of Covariates}\label{subsubsec:overview:hermite-moments}
The focus of prior works is to study the conditional distribution of the \emph{response} $Y_{\max}$.
In particular, conditioned on the covariate $\Regressor_i^\top \Covariate = t \geq \Omega(\sqrt{\log n})$ being highly correlated with a regressor $\Regressor_i$,
we expect the maximum response to be attained by $Y_i$,
\mbox{$(Y_{\max} \mid \Regressor_i^\top \Covariate = t) = (Y_i\mid \Regressor_i^\top \Covariate = t)$},
except with $\poly(\nicefrac1n)$ probability.

\begin{figure}[ht]
  \centering
  \begin{subfigure}{0.5\textwidth}
    \centering
    \includegraphics[width=\linewidth]{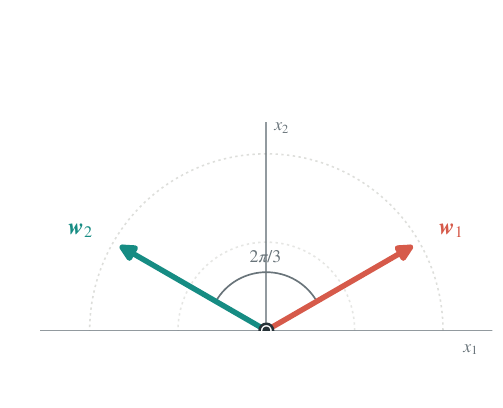}
    \caption{\small Regressors}
  \end{subfigure}\hfill
  \begin{subfigure}{0.5\textwidth}
    \centering
    \includegraphics[width=\linewidth]{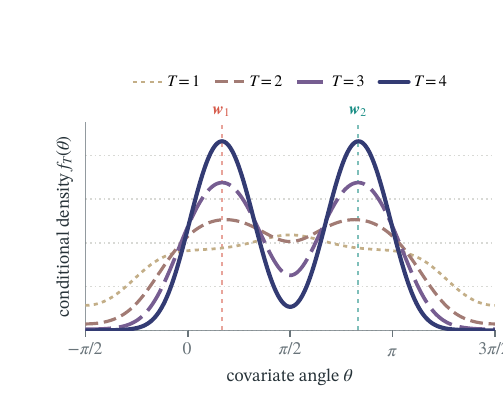}
    \caption{\small Conditional angle density}
  \end{subfigure}
  \caption{\small For two regressors $\Regressor_1,\Regressor_2\in\R^2$ separated by angle $\nicefrac{2\pi}{3}$ (a), the density of the covariate angle conditioned on $Y_{\max}\in[T,T+1]$ concentrates around the regressor directions as $T$ grows (b).
  For simplicity, the intercepts are zero and the noises
  are independent.}
  \label{fig:overview:self-selection}
\end{figure}

Our key observation is that the converse intuition also holds.
That is, conditioned on the response $Y_{\max} = T \geq \Omega(\sqrt{\log n})$ being large,
we expect the \emph{covariates} to be highly correlated with some regressor.
In particular, the conditional distribution of the covariates $\Covariate\mid Y_{\max}=T$ should be a $k$ modal distribution about the regressors $\Regressor_1, \dots, \Regressor_k$ (see \Cref{fig:overview:self-selection}).
We formalize this intuition through the restricted Hermite
moments in \Cref{sec:hermite:main} and highlight the key ideas below.

\paragraph{Ordinary Linear Regression ($k=1$).} To make this more precise,
let us study the conditional moments of $\Covariate\mid Y_{\max}=T$.
More specifically, we study the degree-$d$ conditional \emph{Hermite tensors} $\HermiteTensor{d}{\Covariate}\in \Sym^d(\R^{n})$ of the covariate $\Covariate\sim N(\bm{0}, I_n)$ (see \Cref{subsec:preliminaries:hermite-polynomials}).
Recall that the (univariate) \emph{Hermite polynomials} $\He_d(x): \R\to \R, d\geq 0$ are the unique family of monic polynomials that are orthogonal under Gaussian measure.
Their higher-dimensional analogues are the Hermite tensors $\HermiteTensor{d}{\cdot}:\R^n\to\Sym^d(\R^n)$, $d\geq0$, taking values in a space of dimension $\dim\Sym^d(\R^n)=\binom{n+d-1}{d}$.

We first consider ordinary linear regression (OLS) as a special case of $k=1$ with covariate $\Covariate\sim N(\bm{0}, I_n)$, 
response $Y = \Regressor^\top \Covariate + b + \xi$,
and independent noise $\xi\sim N(0, 1)$.
Let $\phi_{\sigma^2}(t)$ denote the density of $N(0, \sigma^2)$ and $\NormalizedRegressor\coloneqq \Regressor/\norm{\Regressor}_2$ denote the unit vector parallel to the regressor.
From a classic identity for the restricted Hermite tensors of the OLS covariate~\citep{Lancaster57,WN14,Severini26},
\begin{align*}
  \E_{\Covariate,\xi}[\HermiteTensor{d}{\Covariate}\TailIndicator{Y}]
  = \lambda^{(d)}(b,\norm{\Regressor}_2) \NormalizedRegressor^{\otimes d}\,,\qquad
  \text{where}\quad
  \lambda^{(d)}
  = \lambda^{(d)}(b,\norm{\Regressor}_2) > 0\,.
\end{align*}
Thus, up to a coefficient $\lambda^{(d)} = \lambda^{(d)}(b,\norm{\Regressor}_2)$ that depends on the bias and the norm of $\Regressor$,
the $d$-th restricted moment of the covariate is precisely the degree-$d$ tensor power of $\NormalizedRegressor$.
This identity can be proved using the orthogonal decomposition of $\Covariate$ along $\NormalizedRegressor$ and the fact that every positive-degree Hermite polynomial has mean zero under a standard Gaussian.

\paragraph{Regression with Self-Selection ($k\geq 2$).}
Let us now consider the general case of $k\geq 2$ and replace $Y$ with the maximum $Y_{\max}$.
Let $I = I(Y_1, \dots, Y_k)\in [k]$ denote the index of the maximizer response.
By decomposing the indicator, we see that
\begin{align*}
  \TailIndicator{Y_{\max}}
  &= \sum_{i\in [k]} \1\{I=i\}\TailIndicator{Y_i}
  = \sum_{i\in [k]} (1-\1\{I\neq i\})\TailIndicator{Y_i}\,.
\end{align*}
Write $\ComponentWeight{d}{i}=\lambda^{(d)}(b_i,\norm{\Regressor_i}_2)$
for the degree-$d$ weight of the $i$-th regressor.
Using this decomposition and the special case of $k=1$, we see that
\begin{align*}
  &\E_{\Covariate,\NoiseVector}[\HermiteTensor{d}{\Covariate}\TailIndicator{Y_{\max}}] \\
  &\qquad= \sum_{i\in [k]} \lambda_i^{(d)}\NormalizedRegressor_i^{\otimes d}
  - \underbrace{\sum_{i\in [k]}\E_{\Covariate,\NoiseVector}[\HermiteTensor{d}{\Covariate}\1\{I\neq i\}\TailIndicator{Y_i}]}_{\eqqcolon \MomentRemainder{d}} \,.
\end{align*}
In other words, up to a remainder term $\MomentRemainder{d}$,
the degree-$d$ restricted Hermite moment is precisely a weighted sum of degree-$d$ tensor powers of the normalized regressors!
Let us analyze some properties of this decomposition.

\paragraph{Near Orthogonality.}
A key property of this decomposition is that the rank-one tensors are nearly orthogonal.
Indeed, by the separation conditions (\Cref{ass:boundedness-uncoveredness}),
\[
    \abs{\ip{\NormalizedRegressor_i^{\otimes s}}{\NormalizedRegressor_j^{\otimes s}}}
    = \abs{\ip{\frac{\Regressor_i}{\norm{\Regressor_i}_2}}{\frac{\Regressor_j}{\norm{\Regressor_j}_2}}}^s
    \leq \left( 1 - \frac{\Delta^2}{B^2} \right)^s
    = e^{-\Omega(s)}\,.
\]
Choosing $s=\Theta(\log n)$ makes this inner product at most $\poly(\nicefrac1n)$.

\paragraph{Small Bias.}
Another important feature is that the bias term $\MomentRemainder{d}$ is small.
Fix an index $i\in [k]$.
The remainder event $\1\{I\neq i\}\TailIndicator{Y_i}$ occurs when $Y_i\geq T = \Theta(\sqrt{\log n})$ but does not attain the maximum.
Intuitively, if every other $Y_j$ was independent of $Y_i$,
then by a union bound,
the remainder event occurs with probability at most $\sum_{j\neq i}\Pr[Y_j > T \mid Y_i=T]\leq ke^{-\Omega(T^2)} = \poly(\nicefrac1n)$.
While this is too good to be true,
we can apply the separation conditions to show that $\E[Y_j\mid Y_i=T] \leq cT$ for some absolute constant $c\in (0, 1)$ and that the conditional variance remains constant.
Hence we can still prove a similar tail bound on the remainder event.
Combined with a careful bound on the variance of the (unrestricted) degree-$s$ Hermite moment,
this shows that the remainder $\norm{\MomentRemainder{d}}_F = \poly(\nicefrac1n)$ is indeed negligible.

For the rest of the overview, we drop the remainder terms
and treat the nearly orthogonal lifted vectors as orthonormal,
with the understanding that both errors are accounted for
in the formal proofs.
Next, we will see how to exploit the structural properties above to design efficient algorithms for self-selection.

\subsubsection{Projections of Hermite Tensors}
\label{subsubsec:overview:tensor-projections}
In order to exploit the structure of the high-order tensors to design a polynomial time algorithm,
we must avoid materializing the tensors explicitly.
Perhaps a natural attempt is to project the tensors into a lower-dimensional space and perform the necessary computation there.
In \Cref{subsec:tensor-decomposition:iterative-projections}, we analyze an iterative projection method that avoids explicitly forming the higher-order tensors.

\paragraph{Implicitly Accessing Higher-Order Tensors.}
Our algorithm builds on the works of \citet{LL24,DK25} that use iterative projections to compute projections of high-order tensors.
This enables them to evaluate inner products of these tensors without explicitly computing their entries.
We were unable to directly apply their framework for a few technical reasons.
One issue is that our restricted Hermite tensors do not cleanly decompose as a sum of rank-one tensors and have some irreducible error $\MomentRemainder{d}$.
Another issue is that the coefficients $\ComponentWeight{d}{i}$
depend on the degree $d$,
whereas the framework of \citet{DK25} assumes a sequence of moment
tensors with the same weights across degrees.
Most importantly, we were unable to reduce the problem of recovering the regressors to evaluating inner products.
Nonetheless, we can borrow the idea of iterative projections for dimensionality reduction in our setting,
although we need to go through a new analysis that handles the irreducible bias $\MomentRemainder{d}$.

\paragraph{Complex Hermite Outer Products.}
The iterative projection algorithm is easiest to state if we think of our Hermite tensors as flattened matrices.
For this purpose, it is convenient to introduce an alternative representation of the Hermite tensors by lifting the covariates to $\mathbb{C}^n$.
Indeed, for each covariate $\Covariate$, independently draw an auxiliary Gaussian vector $\Covariate'\sim N(\bm{0}, I_n)$ and form the complexified covariate $\widetilde{\Covariate}\coloneqq \Covariate + \iu \Covariate'$.
By the Gaussian integral representation of Hermite polynomials (see \Cref{subsec:preliminaries:hermite-polynomials}),
we have $\HermiteTensor{d}{\Covariate} = \frac1{\sqrt{d!}}\E_{\Covariate'} [\widetilde\Covariate^{\otimes d}\mid \Covariate]$.
Let $\Mat_d: (\R^n)^{\otimes 2d}\to \R^{n^d\times n^d}$ denote the flattening map.
This gives the identity
\begin{align*}
 \Mat_r\!\bigl(\HermiteTensor{2r}{\Covariate}\bigr)
 &=
 \E_{\Covariate'}\!\left[
 \frac{\RealPart[
  \widetilde{\Covariate}^{\otimes r}
  (\widetilde{\Covariate}^{\otimes r})^\top]}
 {\sqrt{(2r)!}}
 \;\middle|\;\Covariate\right]\,,
\end{align*}
where $\RealPart[\cdot]$ denotes the real part of the complex argument.
Note the identity holds even when fixing the real part $\Covariate$ so that it still holds under restrictions of $\Covariate$.

Thus for any linear projection $ C:\R^{n^r}\to\R^{p}$,
we can compute projections of the flattened restricted Hermite tensors as
\begin{align}
  C\, \E_{\Covariate,\NoiseVector} \left[ \Mat_r(\HermiteTensor{2r}{\Covariate}) \TailIndicator{Y_{\max}} \right]  C^\top
  &= \E_{\Covariate,\NoiseVector,\Covariate'}\!\left[
 \frac{\RealPart[
   C\widetilde{\Covariate}^{\otimes r}
  ( C\widetilde{\Covariate}^{\otimes r})^\top]}
 {\sqrt{(2r)!}}
 \TailIndicator{Y_{\max}}\right]
 \label{eq:overview:projected-even-moment}\\
 &\approx
 \sum_{i\in[k]}\ComponentWeight{2r}{i}
 \left[
  ( C\NormalizedRegressor_i^{\otimes r})
 \right]
 \left[
  ( C\NormalizedRegressor_i^{\otimes r})
 \right]^\top\,.
\end{align}
In other words, if we can compute $ C \widetilde{\Covariate}^{\otimes r}$ in polynomial time,
we can compute the lower-dimensional projections of the flattened matrices without explicitly forming them in high dimensions beforehand!

\paragraph{Iterative Projections.}
We now present the iterative projection method due to \citet{LL24} that implicitly computes projections of flattened Hermite tensors.

Our goal is to implicitly compute $\ProjectionMap{s}\widetilde{\Covariate}^{\otimes s}$ such that $\Span\{\NormalizedRegressor_i^{\otimes s}:i\in[k]\}\subseteq \Range(\ProjectionMap{s}^\top)$.
We begin with the trivial base case of the degree $r=1$.
We can simply take the degree-$1$ projection matrix to be $\ProjectionMap{1}=I_n$.
Now, assume that we have computed $\ProjectionMap{r}\widetilde{\Covariate}^{\otimes r}$ for some $\ProjectionMap{r}: (\R^n)^{\otimes r}\to \R^p$ satisfying $\Span\{\NormalizedRegressor_i^{\otimes r}:i\in[k]\}\subseteq \Range(\ProjectionMap{r}^\top)$,
where $p\leq n$.
Extend the projected covariate by one tensor degree $(\ProjectionMap{r}\widetilde{\Covariate}^{\otimes r})\otimes \widetilde\Covariate\in \mathbb C^{pn}$ and form the extended projected Hermite matrix
\begin{align*}
 &\E_{\Covariate,\NoiseVector,\Covariate'}\!\left[
  \frac{\RealPart\!\left[
   \bigl(
    ( \ProjectionMap{r}\widetilde{\Covariate}^{\otimes r})\otimes\widetilde{\Covariate}
   \bigr)
   \bigl(
    ( \ProjectionMap{r}\widetilde{\Covariate}^{\otimes r})\otimes\widetilde{\Covariate}
   \bigr)^\top
  \right]}
  {\sqrt{(2r+2)!}}
  \TailIndicator{Y_{\max}}
 \right] \\
 &\qquad\approx
 \sum_{i\in[k]}\ComponentWeight{2r+2   }{i}
 \left[
  ( \ProjectionMap{r}\NormalizedRegressor_i^{\otimes r})\otimes \NormalizedRegressor_i
 \right]
 \left[
  ( \ProjectionMap{r}\NormalizedRegressor_i^{\otimes r})\otimes \NormalizedRegressor_i
 \right]^\top\,.
\end{align*}
The $\approx$ here is due to the remainder term $\MomentRemainder{2r+2}$.
Let $V_{r+1}\in \R^{k\times pn}$ be the matrix whose rows are the $k$ leading eigenvectors of the extended projected Hermite matrix on the left (padding with arbitrary orthonormal vectors if necessary).
This allows us to extend the degree-$r$ projection map to degree $r+1$ as
$
  \ProjectionMap{r+1}\coloneqq V_{r+1}(\ProjectionMap{r}\otimes I_n)
$.

Using this technique, we can iteratively extend the projection up to degree $s$.
Note that it suffices to store the matrices $V_2, \dots, V_s$ with at most $knp\leq kn^2$ entries each and we can implicitly evaluate projections of the covariate via
\begin{align*}
  \ProjectionMap{1}\widetilde{\Covariate}
 =\widetilde{\Covariate}\,,\qquad
\ProjectionMap{r}\widetilde{\Covariate}^{\otimes r}
 =V_r\!\left[
  ( \ProjectionMap{r-1}\widetilde{\Covariate}^{\otimes(r-1)})\otimes\widetilde{\Covariate}
 \right]\,,
 \quad 2\leq r\leq s\,.
\end{align*}
By construction, we always have $\Span\{\NormalizedRegressor_i^{\otimes r}:i\in[k]\}\subseteq \Range(\ProjectionMap{r}^\top)$.
The first intermediate matrix has size $n^2\times n^2$,
and the later matrices have size at most $kn\times kn$.
Hence we can indeed compute \eqref{eq:overview:projected-even-moment} in polynomial time.

We have successfully reduced dimensions of the higher-order Hermite tensors to $O(k^2)$.
However, the projected vectors
$\ProjectionMap{r}\NormalizedRegressor_i^{\otimes r}$
do not directly reveal the original regressors.
Next, we will see how to circumvent this issue using odd-degree Hermite tensors in addition to the even-degree Hermite tensors we have seen so far.
This requires borrowing some tools from the observable-operator learning literature.

\subsubsection{Observable-Operator Learning}
\label{subsubsec:overview:observable-operator-learning}
Abstractly speaking, observable-operator learning is a technique for learning latent-variable models by constructing operators whose eigenvectors or eigenvalues encode the latent variables~\citep{LRA93,AHK12}.
We will design a way to encode the regressors in the spectrum of the projected Hermite tensors.
This allows us to apply existing tools from the observable-operator learning literature due to \citet{HK13} for learning spherical GMMs.
Our algorithms and exposition found in \Cref{subsec:tensor-decomposition:observable-operator-learning} follow \cite{HK13} but the analysis is our own,
as we need to carefully handle the irreducible remainder terms $\MomentRemainder{2s}$ and $\MomentRemainder{2s+1}$.

\paragraph{Reduction from Projected Hermite Tensors.}
Assume we are able to efficiently compute projections of the complexified covariates $\ProjectionMap{s}\widetilde{\Covariate}^{\otimes s}$ as above.
This enables us to compute projections of the flattened degree-$2s$ Hermite matrices
\begin{align*}
  \ProjectionMap{s} \E_{\Covariate,\NoiseVector}\left[ \Mat_s(\HermiteTensor{2s}{\Covariate}) \TailIndicator{Y_{\max}} \right] \ProjectionMap{s}^\top
  &\approx \underbrace{\sum_{i\in [k]} \ComponentWeight{2s}{i} ({\ProjectionMap{s} \NormalizedRegressor_i^{\otimes s}}) (\ProjectionMap{s} \NormalizedRegressor_i^{\otimes s})^\top}_{\eqqcolon A}\,.
\end{align*}
We can similarly consider odd-degree Hermite tensors where the last tensor component is contracted with some vector $\bm{\zeta}\in \R^n$.
As it turns out, for some $\ScaledRegressor_i \coloneqq \alpha_i \Regressor_i$ where $\alpha_i > 0$ is a positive scalar,
\begin{align*}
  &\ProjectionMap{s}\, \E_{\Covariate,\NoiseVector} \left[ \Mat_s(\HermiteContraction{2s+1}{\Covariate}{\bm{\zeta}}) \TailIndicator{Y_{\max}} \right]  \ProjectionMap{s}^\top \\
  &\quad= \E_{\Covariate,\NoiseVector,\Covariate'}\!\left[
 \frac{\RealPart[
   \ProjectionMap{s}\widetilde{\Covariate}^{\otimes s}
  ( \ProjectionMap{s}\widetilde{\Covariate}^{\otimes s})^\top \ip{\bm{\zeta}}{\widetilde{\Covariate}}]}
 {\sqrt{(2s+1)!}}
 \TailIndicator{Y_{\max}}\right] \\
&\quad\approx
 \underbrace{\sum_{i\in [k]} \ComponentWeight{2s}{i} \ip{\bm{\zeta}}{\ScaledRegressor_i}
  ( \ProjectionMap{s}\NormalizedRegressor_i^{\otimes s})
  ( \ProjectionMap{s}\NormalizedRegressor_i^{\otimes s})^\top}_{\eqqcolon M(\bm{\zeta})} \,.
\end{align*}
Here the $\approx$ is only due to the remainder term $\MomentRemainder{2s+1}$.
Note that the projected even and odd matrices $A$ and $M(\bm{\zeta})$ differ only in the spectrum and that $M(\bm{\zeta})$ encodes a scaled version of the regressors within its spectrum.
Moreover, the vector $\bm{\zeta}\in \R^n$ was arbitrary.
Hence we can hope to recover the scaled regressors $\bm a_i$ using $A$ and several instantiations of $M(\bm{\zeta})$!

\paragraph{Observable-Operator Learning.}
Compute $A^{-1/2}$ from an eigendecomposition
of $A$ and remove the component weights from $M(\bm{\zeta})$:
\[
  D(\bm{\zeta})
  \coloneqq A^{-\frac12} M(\bm{\zeta}) A^{-\frac12}
  = \sum_{i\in [k]} \ip{\bm{\zeta}}{\ScaledRegressor_i}
  ( \ProjectionMap{s}\NormalizedRegressor_i^{\otimes s})
  ( \ProjectionMap{s}\NormalizedRegressor_i^{\otimes s})^\top\,.
\]
The collection $\{D(\bm{\zeta})\}_{\bm{\zeta}\in \R^n}$ is sometimes known as an \emph{observable-operator family}.
For a random unit vector $\bm{\zeta}$, the eigenvalues
$\ip{\bm{\zeta}}{\ScaledRegressor_i}$ are distinct with probability one,
so an eigendecomposition of $D(\bm{\zeta})$ recovers the common
eigenvectors $\bm v_i$.
Then, by evaluating $D(\bm{\zeta})$ at the standard basis vectors $\bm e_1, \dots, \bm e_n$,
we can extract the scaled regressors $\ScaledRegressor_i$ using the eigenvectors.
That is, the $j$-th coordinate of $\ScaledRegressor_i$ can be computed as
$
  (\ScaledRegressor_i)_j 
  = \bm v_i^\top D(\bm e_j) \bm v_i
$.

\paragraph{Additional Complications.}
This leaves us with two additional complications.
The first is that we have only computed some positive scalings of the regressors $\bm a_i = \alpha_i \Regressor_i$ and need to somehow recover the lengths of the regressors.
In \Cref{subsubsec:overview:restricted-least-squares},
we will explain a simple refinement procedure that recovers a regressor $\Regressor_i$ and its bias $b_i$ given a unit-norm estimate $\widehat\NormalizedRegressor_i\approx \NormalizedRegressor_i = \Regressor_i/\norm{\Regressor_i}_2$.
We can then run this procedure with $\widehat\NormalizedRegressor_i = \ScaledRegressor_i/\norm{\ScaledRegressor_i}_2$.

The second more subtle issue is that the magnitude of the remainder term $\MomentRemainder{2s}$ can exceed the minimum eigenvalue $\min_i \ComponentWeight{2s}{i}$ of $A$,
so that inverting estimates of $A$ to estimate its inverse square root is an ill-conditioned operation and may exponentially amplify estimation error.
This can be handled by truncating the spectrum of $A$ at its largest spectral gap to obtain $\bar A$.
This ensures that inverting $\bar A$ will not blow up the estimation error.
However, carrying through the observable-operator learning and subsequent refinement procedures using $\bar A$ in place of $A$ only recovers a subset of the regressors.

This still leaves some unrecovered regressors.
In \Cref{subsubsec:overview:moment-peeling},
we explain how to simulate the Hermite moments of the reduced self-selection problem induced by the unrecovered regressors using only observations from the original problem as well as the recovered regressors.
Thus we can simply recurse on this smaller instance and eventually recover all regressors.

\subsubsection{Restricted Least-Squares}\label{subsubsec:overview:restricted-least-squares}
Consider the classical linear regression setting,
where $n$-dimensional covariates $\Covariate\sim\mathcal D$ are sampled from an arbitrary distribution
and noisy responses $Y = \Regressor^\top \Covariate + b + \xi$ are observed, with independent noise $\xi\sim N(0, 1)$.
As long as the augmented covariates have an invertible second moment,
$\E_{\Covariate}[\begin{bsmallmatrix}\Covariate\\1\end{bsmallmatrix}
\begin{bsmallmatrix}\Covariate^\top&1\end{bsmallmatrix}]\succ0$,
then we can recover $\Regressor$ and $b$ via ordinary least squares (OLS):
\begin{align}
 \begin{bmatrix}\Regressor\\b\end{bmatrix}
 &=\E_{\Covariate}\!\left[
 \begin{bmatrix}\Covariate\\1\end{bmatrix}
 \begin{bmatrix}\Covariate^\top&1\end{bmatrix}
 \right]^{-1}
 \E_{\Covariate,\xi}\!\left[
 \begin{bmatrix}\Covariate\\1\end{bmatrix}Y\right]\,.
 \label{eq:overview:affine-ols}
\end{align}
Below, we will explain the key steps to reduce the problem of recovering $\Regressor_i$ and $b_i$ to OLS,
given only directional estimates of the regressors,
and defer further details to \Cref{sec:refinement:main}.

\paragraph{Reduction to Restricted Least-Squares.}
Let us return to the self-selection setting and assume now that we have an estimate $\widehat\NormalizedRegressor_i\approx \NormalizedRegressor_i = \Regressor_i/\norm{\Regressor_i}_2$ for some $i\in [k]$.
Our goal is to estimate $\Regressor_i$ and $b_i$ up to some desired error $\varepsilon_{\mathrm{int}}$.
Note that we do not need to know the precise regressor index $i$,
just that we have a good estimate to some normalized regressor.
As the discussion for ordinary least-squares above holds for an arbitrary distribution $\mathcal D$,
our strategy is to restrict to an event
$\RefinementEvent{\widehat\NormalizedRegressor_i}$
on which $Y_{\max}=Y_i$ with high probability.
Then, we can simply use a plug-in estimator for \eqref{eq:overview:affine-ols} with
$\begin{bsmallmatrix}\Covariate\\1\end{bsmallmatrix}
 \begin{bsmallmatrix}\Covariate^\top&1\end{bsmallmatrix}
 \1_{\RefinementEvent{\widehat{\NormalizedRegressor}_i}}$
and
$\begin{bsmallmatrix}\Covariate\\1\end{bsmallmatrix}Y_{\max}
 \1_{\RefinementEvent{\widehat{\NormalizedRegressor}_i}}$.

\begin{figure}[ht]
  \centering
  \begin{subfigure}{0.5\textwidth}
    \centering
    \includegraphics[width=\linewidth]{figures/self_selection_regressors.pdf}
    \caption{\small Regressors}
  \end{subfigure}\hfill
  \begin{subfigure}{0.5\textwidth}
    \centering
    \includegraphics[width=\linewidth]{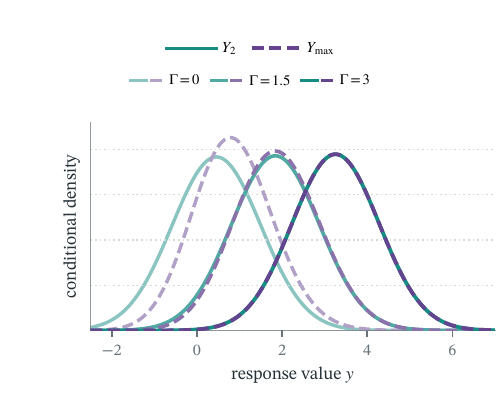}
    \caption{\small Conditional response densities}
  \end{subfigure}
  \caption{\small For the regressors in (a), conditioned on $\Gamma\leq\ip{\NormalizedRegressor_2}{\Covariate}\leq\Gamma+1$, the density of the maximum response $Y_{\max}$ becomes indistinguishable from that of $Y_2$ as $\Gamma$ grows (b).
  For simplicity, the intercepts are zero and the noises
  are independent.}
  \label{fig:overview:refinement}
\end{figure}

\citet{GM24} showed that if a candidate is sufficiently close to $\Regressor_i$, then $Y_i$ attains the maximum with high probability conditional on a large projection along the candidate direction.
We use the same idea with the event
\[
 \RefinementEvent{\widehat{\NormalizedRegressor}_i}
 \coloneqq\{\Gamma\leq\ip{\widehat{\NormalizedRegressor}_i}{\Covariate}\leq\Gamma+1\}\,.
\]
Intuitively, since the regressors have constant angular separation,
if a covariate is highly correlated with a single regressor $\Regressor_i$,
it cannot also be highly correlated with other regressors $\Regressor_j, j\neq i$,
so the maximizer response should be $Y_i$ with high probability (see \Cref{fig:overview:refinement}).
To see this, we first remark that the separation condition (\Cref{ass:boundedness-uncoveredness}) ensures that
$
  \ip{\Regressor_i-\Regressor_j}{\widehat\NormalizedRegressor_i}
  = \frac{\norm{\Regressor_i}_2^2 - \ip{\Regressor_i}{\Regressor_j}}{\norm{\Regressor_i}_2}
  \geq \Omega(1)
$.
The biases change the conditional mean by at most $2B$,
which is absorbed by a sufficiently large $\Gamma$.
Hence the conditional distribution satisfies
\[
  (Y_i - Y_j \mid \ip{\widehat\NormalizedRegressor}{\Covariate} = \Gamma) \sim N(\mu, \sigma^2)\,,\qquad
  \mu\geq \Omega(\Gamma)\,, \sigma^2=O(1)\,.
\]
We can thus choose $\Gamma=\Theta(\sqrt{\log(\nicefrac{n}{\varepsilon_{\mathrm{int}}})})$
so that by a union bound,
the probability that $Y_i$ does not attain the max is at most $ke^{-\Omega(\Gamma^2)} = \poly(\nicefrac{\varepsilon_{\mathrm{int}}}{n})$.

\paragraph{Local Convergence.}
As a corollary of our restricted least squares algorithm,
we get a simple polynomial-time local convergence algorithm for self-selection improving on prior works.
\citet{KMZ26} designed a polynomial-time local convergence algorithm for self-selection that required estimates $\norm{\widehat\Regressor_i - \Regressor_i}_2\leq \poly(\nicefrac1k)$ with inverse polynomial initialization distance.
We note that $\norm{\widehat\Regressor_i - \Regressor_i}_2 = O(1)$ implies that $\norm{\widehat\Regressor_i/\norm{\widehat\Regressor_i}_2 - \NormalizedRegressor_i}_2 = O(1)$,
possibly with bigger constants.
Hence our restricted least squares algorithm immediately implies a local convergence algorithm with $O(1)$ convergence radius,
improving on \cite{KMZ26}.

\subsubsection{Moment Peeling}\label{subsubsec:overview:moment-peeling}
As previously mentioned, our observable-operator learning subroutine extracts estimates for a block of the normalized regressors
which can be refined to actual regressor and bias estimates using restricted least squares.
The final piece of the puzzle is to show that we can recursively run our algorithm to recover the rest of the regressors.
This is not an issue for the restricted least squares algorithm as it only requires an accurate angular estimate of some regressor.
For our implicit tensor decomposition subroutine, our key observation is that we only need to access samples when computing the observable-operator family.
In particular, we only need to estimate restricted Hermite tensors of the covariate.
Below, we describe the intuition for how to simulate these tensors for the unrecovered regressors which is all we need to run our algorithm and defer further details to \Cref{sec:peeling:main}.
This idea of iteratively recovering and removing components
has been used in self-selection~\citep{CDIZ23},
Gaussian mixture clustering~\citep{LL24}, and general
tensor decomposition~\citep{AGHKT14}.

\paragraph{Moment Peeling.}
Let $R\subseteq [k], Q\coloneqq [k]\setminus R$ and assume we recovered $(\Regressor_i,b_i), i\in R$.
Our goal is to be able to simulate restricted Hermite tensors for the (unobserved) maximum $Y_{\max,Q} = \max_{i\in Q} Y_i$,
\begin{align*}
  \E_{\Covariate, \NoiseVector} \left[ \HermiteTensor{d}{\Covariate} \TailIndicator{Y_{\max,Q}} \right]\,.
\end{align*}
Using only the observations and the
recovered regressors,
we will construct a weight function $\omega_R(\Covariate,Y_{\max})$
whose conditional expectation approximates that of
$\TailIndicator{Y_{\max,Q}}$, 
\[
  \E_{\NoiseVector}[\omega_R(\Covariate, Y_{\max})\mid \Covariate]\approx \E_{\NoiseVector}[\TailIndicator{Y_{\max, Q}}\mid \Covariate]\,.
\]
We emphasize that we do not observe $Y_{\max,Q}$ directly.

Let $Y_{\max,R} \coloneqq \max_{i\in R} Y_i$ and note that $Y_{\max} = \max(Y_{\max,R}, Y_{\max,Q})$.
We begin with the pointwise identity
\begin{align}
&\TailIndicator{Y_{\max}} + \TailIndicator{\min(Y_{\max,R}, Y_{\max,Q})} \notag\\
&\qquad= \TailIndicator{Y_{\max,R}} + \TailIndicator{Y_{\max,Q}}\,, \notag \\
\notag\\
&\TailIndicator{Y_{\max}} - \TailIndicator{Y_{\max,R}} \notag\\
&\qquad= \TailIndicator{Y_{\max,Q}} - \TailIndicator{\min(Y_{\max,R}, Y_{\max,Q})}\,.
\label{eq:overview:peeling-identity}
\end{align}
Now, the RHS is nearly the indicator variable we want modulo a bias term.
Intuitively, this bias term is a much rarer event since it roughly requires at least two Gaussians to exceed $T = \Theta(\sqrt{\log n})$,
while all other events roughly only require one.

We would like to define $\omega_R(\Covariate, Y_{\max})$ to be the LHS.
Although $\TailIndicator{Y_{\max,R}}$ is unobserved, 
we only need the weight function to be a good approximation under conditional expectation.
It suffices to take
\begin{align*}
  \omega_R(\Covariate, Y_{\max})
  &\coloneqq \TailIndicator{Y_{\max}}
  - \Prb_{\NoiseVector_R}\!\left[
 T\leq Y_{\max,R}\leq T+1\mid\Covariate\right] \\
  \E_{\NoiseVector}[\omega_R(\Covariate, Y_{\max})\mid \Covariate]
  &=\E_{\NoiseVector}\left[ \TailIndicator{Y_{\max}} - \TailIndicator{Y_{\max,R}} \mid \Covariate \right]\,.
\end{align*}
Since the recovered regressors, their biases, and $\rho$ are known,
we can estimate this conditional probability by simulating noise for the recovered regressors.

Now, introducing the Hermite tensor $\HermiteTensor{d}{\Covariate}$,
taking the expectation over the randomness $(\Covariate, \NoiseVector)$,
and using the identity \eqref{eq:overview:peeling-identity},
we have that
\begin{align*}
  &\E_{\Covariate, \NoiseVector} \left[ \HermiteTensor{d}{\Covariate} \omega_R(\Covariate, Y_{\max}) \right] \\
  &\quad = \E_{\Covariate, \NoiseVector} \left[ \HermiteTensor{d}{\Covariate} \TailIndicator{Y_{\max,Q}} \right]
  - \E_{\Covariate, \NoiseVector} \left[ \HermiteTensor{d}{\Covariate} \TailIndicator{\min(Y_{\max,R}, Y_{\max,Q})} \right]\,.
\end{align*}
Thus we can indeed produce estimates of the desired restricted Hermite moment up to some bias term.

This concludes the overview of our techniques.

\subsection{Organization}
We organize the rest of this work as follows.
First, related work is discussed in \Cref{sec:related-work:main} and some preliminaries are collected in \Cref{sec:preliminaries:main}.
In \Cref{sec:hermite:main}, we analyze the structure of high-order restricted Hermite moments.
In \Cref{sec:tensor-decomposition:main}, we use this structure to design an implicit tensor decomposition procedure that recovers a subset of normalized regressors without explicitly forming the moment tensors.
\Cref{sec:refinement:main} uses restricted least squares to estimate the full regressors from these unit vector estimates.
\Cref{sec:peeling:main} implements an important estimation subroutine used in our implicit tensor decomposition procedure.
We combine these steps in \Cref{sec:recovery:main} to prove \Cref{thm:intro:polynomial-time-recovery}.

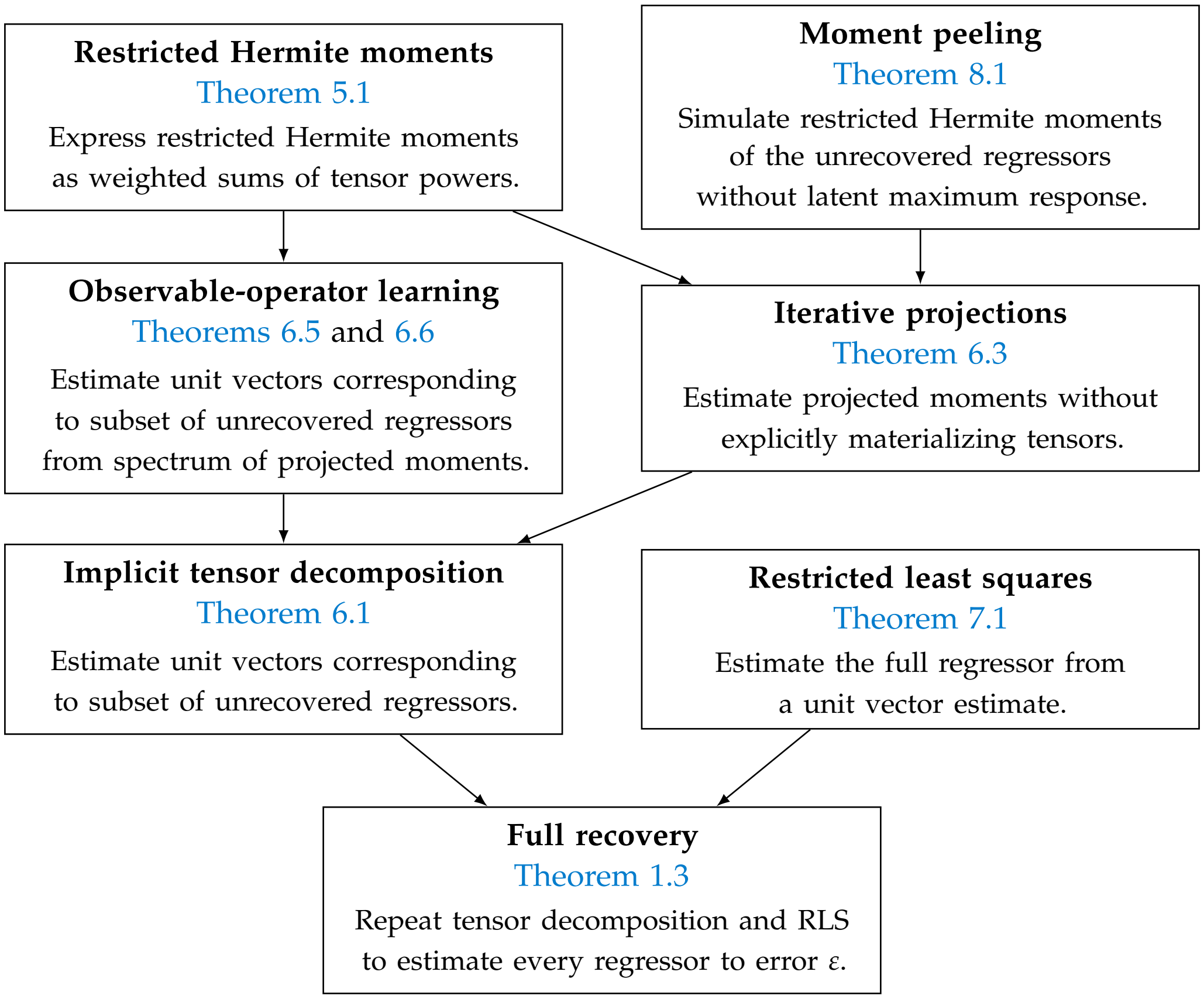
\begin{figure}[ht]
\centering
\begin{tikzpicture}[
  result/.style={draw,align=center,text width=5.95cm,
                 minimum height=1.95cm,inner sep=5pt,font=\small},
  >=latex]
\node[result] (structure) at (-3.6,0) {
  \textbf{Restricted Hermite moments}\\
  \Cref{thm:hermite:restricted-moment-structure}\\[2pt]
  {\footnotesize Express restricted Hermite moments\\as weighted sums of tensor powers.}};
\node[result] (peel) at (3.6,0) {
  \textbf{Moment peeling}\\
  \Cref{thm:peeling:hermite-moment-estimation}\\[2pt]
  {\footnotesize Simulate restricted Hermite moments\\of the unrecovered regressors\\without latent maximum response.}};

\node[result] (operators) at (-3.6,-2.95) {
  \textbf{Observable-operator learning}\\
  \Cref{thm:tensor-decomposition:observable-operator-construction,thm:tensor-decomposition:normalized-regressor-recovery}\\[2pt]
  {\footnotesize Estimate unit vectors corresponding to subset of unrecovered regressors from spectrum of projected moments.}};
\node[result] (project) at (3.6,-2.95) {
  \textbf{Iterative projections}\\
  \Cref{thm:tensor-decomposition:projected-moment-accuracy}\\[2pt]
  {\footnotesize Estimate projected moments without\\explicitly materializing tensors.}};

\node[result] (block) at (-3.6,-5.90) {
  \textbf{Implicit tensor decomposition}\\
  \Cref{thm:tensor-decomposition:implicit-tensor-decomposition}\\[2pt]
  {\footnotesize Estimate unit vectors corresponding to subset of unrecovered regressors.}};
\node[result] (refine) at (3.6,-5.90) {
  \textbf{Restricted least squares}\\
  \Cref{thm:refinement:restricted-least-squares}\\[2pt]
  {\footnotesize Estimate the full regressor from\\a unit vector estimate.}};

\node[result] (main) at (0,-8.85) {
  \textbf{Full recovery}\\
  \Cref{thm:intro:polynomial-time-recovery}\\[2pt]
  {\footnotesize Repeat tensor decomposition and RLS\\to estimate every regressor to error $\varepsilon$.}};

\draw[->] (structure) edge (operators) edge (project);
\draw[->] (peel) -- (project);
\draw[->] (operators) -- (block);
\draw[->] (project) -- (block);
\draw[->] (block) -- (main);
\draw[->] (refine) -- (main);
\end{tikzpicture}
\caption{Main dependencies in the proof of \Cref{thm:intro:polynomial-time-recovery}.
Each arrow points from a result to a downstream application.}
\label{fig:overview:proof-dependencies}
\end{figure}

For convenience, we indicate the main dependencies in the proof of \Cref{thm:intro:polynomial-time-recovery} in \Cref{fig:overview:proof-dependencies}.

\section{Related Works}\label{sec:related-work:main}

\paragraph{Regression with Self-Selection.}
As discussed, \citet{CDIZ23} initiated the study of time- and
sample-efficient algorithms for regression with self-selection.
Subsequent work improved the sample complexity and running time,
but the best known end-to-end algorithms still require time
exponential in the number of regressors $k$ \cite{GM24,KMZ26}.
Our main contribution is the first polynomial-time algorithm for regression with self-selection.

\paragraph{Max-Affine Regression.}
Max-affine regression estimates a maximum of affine functions
from noisy observations. The noise is added after taking the maximum,
so it does not affect the selection of the affine function.
This is a parametric form of convex regression, whose study goes back
to \citet{Hildreth54,Holloway79}.
\citet{MagnaniBoyd09} studied fitting a prescribed number of affine
pieces, while \citet{HannahDunson13,BalazsGS15} developed max-affine
estimators for general convex regression functions.

For recovery of the individual parameters, \citet{GPGR22} analyze
alternating minimization with spectral initialization.
Suppressing their additional spectral conditioning factors,
their algorithm runs in
$\poly(n,k,\nicefrac1\eps)\exp\!\bigl(\Ot(k^2)\bigr)$ time.
Subsequent work gives local guarantees for first-order and sparse
methods~\citep{KimLee24,KanjKimLee26}.
As a corollary of our general algorithm, we give the first
polynomial-time algorithm for max-affine regression.

\paragraph{Max-Linear Regression.}
Max-linear regression is a special case of max-affine regression with biases $\bm0$.
For this problem, \citet{KimBL24} give a convex estimator under appropriate initialization.
\citet{GM24} combine their algorithm with the local convergence
guarantee of \citet{GPGR22} to obtain polynomial sample complexity
but $\poly(n,k,\nicefrac1\eps)+k^{O(k)}$ running time.
Our algorithm improves the running time to
$\poly(n,k,\nicefrac1\eps)$.

\paragraph{Gaussian Mixture Models.}
Gaussian mixture learning seeks to estimate a mixture's density
or parameters, or to cluster its samples.
Learning and clustering guarantees have been studied under
various structural assumptions and in the presence of corrupted
samples \cite{Dasgupta99,AroraKannan05,VempalaWang04,
AchlioptasMcSherry05,KalaiMoitraValiant10,MoitraValiant10,
BelkinSinha10,RegevVijayaraghavan17,HK13,GeHuangKakade15,
HopkinsLi18,KothariSteinhardtSteurer18,BakshiEtAl22,LiuMoitra22}.
Our iterative projection construction follows work on clustering
separated spherical Gaussian mixtures using implicit high-order
Hermite moments \cite{LL24}.

\paragraph{Mixtures of Linear Regressions.}
In mixed linear regression, each response is generated by an
unobserved regressor sampled independently of the covariate.
Parameter recovery in linear and generalized-linear mixtures has
been studied under different assumptions on the covariates, noise,
and initialization
\cite{DeVeaux89,ChagantyLiang13,YiCS14,SedghiJA16,ZhongJD16,
LiLiang18,BWY17,KwonC20,CLS20,DK20}.
More directly related to our approach, implicit high-order
moments yield polynomial-time algorithms for density estimation
\cite{DK25}.

\paragraph{Multi-View Mixtures and Hidden Markov Models.}
In multi-view mixtures, conditionally independent observations
share a latent component; in hidden Markov models, the latent
states form a Markov chain.
Both admit learning algorithms under suitable rank and
identifiability assumptions \cite{Jaeger00,HsuKakadeZhang12,AGHKT14}.
Our observable-operator family uses the same common-eigenbasis
construction developed for multi-view mixtures \cite{AHK12}.

\paragraph{Gaussian Multi-Index Models.}
In Gaussian multi-index models, the response depends on the
covariate through an unknown low-dimensional projection.
Algorithms and statistical limitations have been studied for
recovering this subspace or learning the response function
\cite{Li91,Li92,Vempala10,CM20,DPLB24}.
Our recursive projection construction is related to sequential
recovery of directions from Hermite or conditional moments,
although recovering the subspace does not identify the individual
regressors in our model \cite{DLB25,DIKZ25,DIKR25}.

\paragraph{Learning Shallow Neural Networks.}
Learning a shallow network asks to learn a sum of activations
of linear forms or recover its weights from labeled examples.
Learning guarantees under Gaussian covariates have been studied
for several activations and network classes
\cite{JanzaminNNLIFT15,ZhongSJBD17,DKKZ20,ChenN24,DLS22}.
Our use of random contractions and comparisons between Hermite
degrees is related to algorithms for learning ReLU networks
\cite{ChenDGKM23,DK24}.

\paragraph{Learning from Coarse Observations.}
Coarse observations reveal only a set containing the latent sample.
Results include learning algorithms for coarse labels and Gaussian
mean estimation under identifiable convex partitions
\cite{FKKT21,KMMZZ26}.
This general framework captures the self-selection problem and through this connection,
local convergence algorithms for self-selection have been developed~\cite{KMZ26}.

\paragraph{Truncated and Censored Data.}
In truncated or censored estimation, observations are either
discarded or only partially revealed \cite{Tobin58}.
Algorithms for known truncation estimate regression and Gaussian
parameters under suitable survival and design assumptions,
with applications to differentially private estimation
\cite{DGTZ18,DGTZ19,IZD20,DRZ20,ZampetakisZhou25}.
Unknown truncation has also been studied for regression and
exponential-family estimation under structural assumptions on
the survival set and distribution
\cite{KTZ19,LMZ24,LMZ26,KMKC26,LNSW26}.
Other related work studies general selection filters, computational
lower bounds, and truncation detection \cite{KLMZ26,DKPZ24,DLNS24}.

\section{Preliminaries}
\label{sec:preliminaries:main}
In this section,
we introduce the notational conventions as well as some preliminary tools used throughout this work.

\subsection{Hyperparameters}
To simplify the presentation, we choose $K=K(B, \Delta)$ to be a sufficiently large constant depending only on the geometric quantities as below.
We choose the tensor order $s$ to make the lifted vectors nearly orthogonal, 
and the threshold $T$ so that the bias in the decomposition of restricted Hermite moments is negligible.  
The internal accuracy
$\varepsilon_{\mathrm{int}}$ is chosen to be sufficiently small to account for errors that accumulate across rounds can be absorbed.
Finally, the refinement threshold ${\Gamma}$ is chosen to control the bias of the restricted least squares estimate.

\begin{hyperparameters}\label{hyp:algorithm}
Throughout this work, we fix the following hyperparameters
\begin{align}
 K&\coloneqq\frac{16(1+B^2)^3}{\Delta^4},
 \\
 s&\coloneqq K^4\log n,
 \\
 T&\coloneqq2K^2\sqrt{s},
 \\
 \varepsilon_{\mathrm{int}}
 &\coloneqq\min\{\varepsilon,e^{-K^6s}\},
 \\
 {\Gamma}&\coloneqq
 2K\sqrt{\log\!\left(\frac{n}{\varepsilon_{\mathrm{int}}}\right)}.
\end{align}
\end{hyperparameters}

\subsection{Notation}

For a positive integer $m$, write $[m]=\{1,\ldots,m\}$.
We use bold letters for vectors, ordinary capital letters for matrices,
and sans-serif capital letters for tensors.
We use $\ip{\cdot}{\cdot}$ for the Euclidean inner product
on vectors and the Frobenius inner product on real matrices
and tensors, and $\norm{\cdot}_2,\norm{\cdot}_F,\norm{\cdot}_{\op}$
for the Euclidean, Frobenius, and operator norms.
The symbol $\1\{\mathcal E\}$ denotes the indicator of an event $\mathcal E$.

We use $\otimes$ for tensor products and
$\bm v^{\otimes d}$ for the $d$-fold tensor product of $\bm v$.
The space $\Sym^d(\R^n)$ consists of symmetric order-$d$ tensors on $\R^n$.
The notation $\mathsf A\mathbin{\lrcorner}\bm{\zeta}$ denotes
contraction against $\bm{\zeta}$ in the last tensor factor.
For an order-$2r$ tensor on $\R^n$,
$\Mat_r(\cdot)$ flattens its argument into an $n^r\times n^r$ matrix using the first $r$ multi-indices as the row index and the last $r$ as the column index.
We write $\iu=\sqrt{-1}$ and use $\RealPart$ for the entrywise real part.
Finally, $\Phi$ denotes the standard Gaussian CDF and
$\phi_{\sigma^2}$ the density of $N(0,\sigma^2)$.
We write $\phi=\phi_1$ for the standard Gaussian density.

\subsection{Properties of Gaussians}

\begin{proposition}[Gaussian Tail Bound;
 {\cite[Equation~(2.7)]{Wainwright19}}]
\label{prop:gaussian-tail}
Let $G\sim N(\mu,\sigma^2)$ with $\sigma>0$.
For every $t\geq0$,
\[
 \Prb_G(G\geq\mu+t)
 \leq\exp\!\left(-\frac{t^2}{2\sigma^2}\right).
\]
By symmetry, the same bound holds for $\Prb_G(G\leq\mu-t)$.
\end{proposition}

\begin{proposition}[Conditional Gaussian Distribution;
 {\cite[Equation~(A.6)]{RW06}}]
\label{prop:gaussian-conditioning}
Let $\bm{G}\in\R^n$ and $Y\in\R$ be centered jointly Gaussian.
Write $\Sigma=\E_{\bm{G}}[\bm{G}\bm{G}^\top]$, $\bm{v}=\E_{\bm{G},Y}[\bm{G}Y]$, and
$\sigma^2=\E_Y[Y^2]>0$.
For every $t\in\R$,
\[
 \bm{G}\mid Y=t
 \sim N\!\left(
  \frac{\bm{v}t}{\sigma^2},
  \Sigma-\frac{\bm{v}\bm{v}^\top}{\sigma^2}
 \right).
\]
\end{proposition}

\subsection{Hermite Polynomials}
\label{subsec:preliminaries:hermite-polynomials}

We define the Hermite features and fix their normalization.
Hermite expansions have been used in algorithms and statistical query
lower bounds for learning neural networks under Gaussian
inputs~\cite{DKKZ20,DK24}.

Let $\Covariate'\sim N(\bm{0},I_n)$ be an auxiliary Gaussian vector.

For every nonnegative integer $d$, the probabilists' Hermite
polynomial has the Gaussian integral representation
\cite[Equations~18.10.7 and 18.7.11]{DLMF}
\[
 \He_d(t)\coloneqq\E_g[(t+\iu g)^d],
 \qquad g\sim N(0,1),\quad t\in\R.
\]
This expectation is real, since $g$ and $-g$ have the same law.
For $\bm{x}\in\R^n$, define the normalized Hermite tensor
\begin{equation}
 \HermiteTensor{d}{\bm{x}}\coloneqq
 \frac1{\sqrt{d!}}\,
 \E_{\Covariate'}[\RealPart((\bm{x}+\iu \Covariate')^{\otimes d})].
\end{equation}
These tensors are real and symmetric.  For example,
\[
 \HermiteTensor{0}{\bm{x}}=1,\qquad
 \HermiteTensor{1}{\bm{x}}=\bm{x},\qquad
 \HermiteTensor{2}{\bm{x}}=\frac{\bm{x}\bm{x}^\top-I_n}{\sqrt2}.
\]
In one dimension, $\HermiteTensor{d}{t}=\He_d(t)/\sqrt{d!}$.

The normalized Hermite polynomials are orthonormal under a
standard Gaussian~\cite[Proposition~11.33]{ODonnell14}.
For $G\sim N(0,1)$ and nonnegative integers $d_1,d_2$,
\begin{equation}\label{eq:hermite-orthogonality}
 \E_G\!\left[
  \frac{\He_{d_1}(G)}{\sqrt{d_1!}}
  \frac{\He_{d_2}(G)}{\sqrt{d_2!}}
 \right]
 =
 \begin{cases}
  1,&d_1=d_2,\\
  0,&d_1\neq d_2.
 \end{cases}
\end{equation}
Thus Hermite polynomials of different degrees are uncorrelated
when evaluated at a standard Gaussian.

For a symmetric order-$d$ tensor $\mathsf A$, its entries are unchanged
when the indices are permuted.
We group these entries by the number of times each coordinate
appears.  Write $\bm{\alpha}=(\alpha_1,\ldots,\alpha_n)$ for these
counts, and let $A_{\bm{\alpha}}$ denote their common value.
For example, for an order-three tensor on $\R^2$,
$\bm{\alpha}=(2,1)$ groups the entries
$A_{112}=A_{121}=A_{211}$.
Write $|\bm{\alpha}|=\sum_j\alpha_j$ and
$\bm{\alpha}!=\prod_j\alpha_j!$.
There are $d!/\bm{\alpha}!$ entries in the group indexed by $\bm{\alpha}$.

\begin{proposition}[Hermite Expansion and Normalization; {\cite[Proposition~11.33]{ODonnell14}}]
\label{prop:hermite-expansion-normalization}
Let $\mathsf A$ be a symmetric order-$d$ tensor.
For every $\bm{x}\in\R^n$,
\[
 \HermitePairing{d}{\bm{x}}{\mathsf A}
 =
 \sum_{|\bm{\alpha}|=d}
 \sqrt{\frac{d!}{\bm{\alpha}!}}\,A_{\bm{\alpha}}
 \prod_{j=1}^n
 \frac{\He_{\alpha_j}(x_j)}{\sqrt{\alpha_j!}}.
\]
For $\bm{G}\sim N(\bm{0},I_n)$, the product form of Hermite orthogonality
shows that the products
in this sum are orthonormal.
Consequently,
\begin{equation}\label{eq:hermite-isometry}
 \E_{\bm{G}}[\HermitePairing{d}{\bm{G}}{\mathsf A}^2]
 =
 \sum_{|\bm{\alpha}|=d}\frac{d!}{\bm{\alpha}!}A_{\bm{\alpha}}^2
 =\norm{\mathsf A}_F^2.
\end{equation}
\end{proposition}

We give the derivation in \Cref{subsec:app-preliminaries-proofs:hermite-polynomials}.

\subsection{Spectral Perturbation Bounds}

We use Weyl's eigenvalue perturbation inequality in the following form.
\begin{proposition}[Weyl's Eigenvalue Perturbation Inequality; {\cite[Corollary~III.2.6]{Bhatia97}}]
\label{prop:weyl-perturbation}
Let $H,\widetilde H\in\R^{p\times p}$ be symmetric, with eigenvalues
in nonincreasing order.  Then, for every $j\in[p]$,
\[
 \abs{\lambda_j(\widetilde H)-\lambda_j(H)}
 \leq\norm{\widetilde H-H}_{\op}\,.
\]
\end{proposition}

\subsection{Heavy-Tail Mean Estimation}

\begin{lemma}[Coordinatewise Trimmed Mean;
 {\cite[Theorem~1.1.1]{OliveiraOrensteinRico25}}]
\label{lem:coordinatewise-trimmed-mean}
Let $\bm W=(W_1,\ldots,W_p)$ be a random vector satisfying
\[
 \E_{\bm W}[W_j^2]\leq\mathsf v
 \qquad(j\in[p]),
\]
where $\mathsf v>0$.
For $0<\varepsilon\leq\sqrt{\mathsf v}$ and
$\delta\in(0,\nicefrac13)$, there is an estimator
$\textnormal{\textsc{TrimmedMean}}(\bm W;\mathsf v,\varepsilon,\delta)$
that returns $\widehat{\bm\mu}\in\R^p$ satisfying
\[
 \max_{j\in[p]}|\widehat\mu_j-\E_{\bm W}[W_j]|
 \leq\varepsilon
\]
with probability at least $1-\delta$.
Its sample complexity is
\[
 O\!\left(
 \frac{\mathsf v}{\varepsilon^2}
 \log\!\left(\frac p\delta\right)
 \right),
\]
and its running time is
\[
 O\!\left(
 \frac{p\mathsf v}{\varepsilon^2}
 \log\!\left(\frac p\delta\right)
 \right).
\]
\end{lemma}

\section{Structure of Restricted Hermite Moments}
\label{sec:hermite:main}

In this section, we show that the restricted Hermite moments of the latent maximum corresponding to unrecovered regressors
are close to weighted sums of tensor powers of said regressors.
As we will see, the even and odd moments have the same rank-one components, and the odd
moments allow us to recover the regressor directions.
We also establish the bounds needed to estimate and project these
moments and to apply the observable-operator learning procedure.

Fix a non-empty set $Q\subseteq[k]$ of remaining regressor indices.
Recall that $Y_i=\ip{\Regressor_i}{\Covariate}+b_i+\xi_i$ is the latent response of component
$i$ and $Y_{\max,Q}=\max_{i\in Q}Y_i$ is the remaining maximum.
Write $\sigma_i^2=1+\norm{\Regressor_i}_2^2$ for the variance of $Y_i$ and
$\NormalizedRegressor_i=\Regressor_i/\norm{\Regressor_i}_2$ for its normalized regressor.

\paragraph{Restricted Hermite Moments.}
We use the normalized Hermite tensors from \Cref{subsec:preliminaries:hermite-polynomials}.
Let $\Covariate'\sim N(\bm{0},I_n)$ be independent of the covariate and
latent noises, and put
\[
 \widetilde{\Covariate}\coloneqq \Covariate+\iu \Covariate'.
\]
Our algorithm uses this complexified covariate to evaluate
the Hermite features.

We consider the restricted Hermite moments
\begin{align}
 \RestrictedMomentTensor{d}{Q}
 &\coloneqq
 \E_{\Covariate,\NoiseVector_Q}[\HermiteTensor{d}{\Covariate}\TailIndicator{Y_{\max,Q}}]
 \notag\\
 &=
 \frac1{\sqrt{d!}}\,
 \E_{\Covariate,\NoiseVector_Q,\Covariate'}[
  \RealPart(\TensorPower{\widetilde{\Covariate}}{d})\TailIndicator{Y_{\max,Q}}
 ]\,.
\end{align}
We note that compared to conditional moments,
these restricted moments are unnormalized and each expectation is restricted to the tail slice of $Y_{\max,Q}$.

\paragraph{Weights and Scaled Regressors.}
As we will soon see,
the even restricted moments approximately decompose as a weighted sum of rank-one tensors.
The odd moments contain the same components, with an additional factor in the weights.
We now define these quantities.
For $\sigma>0$, let
\[
 \phi_{\sigma^2}(t)
 \coloneqq\frac{1}{\sqrt{2\pi}\sigma}e^{-t^2/(2\sigma^2)}
\]
denote the density of $N(0,\sigma^2)$.
Recall $\He_d(\cdot)$ denotes the $d$-th Hermite polynomial in 1-D.
For every integer $d\geq0$, define the restricted moment coefficient
\begin{equation}
 \MomentCoefficient{d}{b,\sigma^2}\coloneqq
 \int_T^{T+1}
 \phi_{\sigma^2}(t-b)
 \frac{\He_d((t-b)/\sigma)}{\sigma^{d}\sqrt{d!}} \odif{t}.
\end{equation}
For $i\in[k]$ and $1\leq r\leq s$, define the degree-$2r$ weights
\begin{equation}\label{eq:component-weights}
 \ComponentWeight{2r}{i}\coloneqq \MomentCoefficient{2r}{b_i,\sigma_i^2}\norm{\Regressor_i}_2^{2r}.
\end{equation}
At level $s$, define the scaled regressors
\begin{equation}\label{eq:scaled-regressors}
 \ScaledRegressor_i\coloneqq
 \frac{\MomentCoefficient{2s+1}{b_i,\sigma_i^2}}{\MomentCoefficient{2s}{b_i,\sigma_i^2}}\,\Regressor_i.
\end{equation}

We are now ready to state the main result of this section.

\begin{theorem}[Restricted Moment Structure]\label{thm:hermite:restricted-moment-structure}
Consider the Gaussian self-selection model
(\Cref{def:gaussian-self-selection}) under the boundedness and
uncoveredness assumptions (\Cref{ass:boundedness-uncoveredness}).
Use the hyperparameters $K,T,s$ from \Cref{hyp:algorithm}.
The following holds for every non-empty $Q\subseteq[k]$.

\begin{enumerate}[label=(\roman*),ref=\thetheorem(\roman*),leftmargin=2.2em]
\item\label[theorem]{thm:hermite:structure-moment-decomposition} For every $1\leq r\leq s$ and $\bm{\zeta}\in\R^n$,
\begin{align}
 \Mat_r(\RestrictedMomentTensor{2r}{Q})
 &=
 \sum_{i\in Q}\ComponentWeight{2r}{i}
 \TensorPower{\NormalizedRegressor_i}{r}(\TensorPower{\NormalizedRegressor_i}{r})^\top+E_r\,,
 \\
 \Mat_s(\RestrictedMomentTensor{2s+1}{Q}\mathbin{\lrcorner}\bm{\zeta})
 &=
 \sum_{i\in Q}\ComponentWeight{2s}{i}
 \ip{\bm{\zeta}}{\ScaledRegressor_i}\TensorPower{\NormalizedRegressor_i}{s}(\TensorPower{\NormalizedRegressor_i}{s})^\top+E(\bm{\zeta})\,.
\end{align}
The errors satisfy
\begin{equation}
\begin{aligned}
 \norm{E_r}_F
 &\leq e^{8s\log K-T^2/(2K)}\sum_{i\in Q}\ComponentWeight{2r}{i}\,,\\
 \norm{E(\bm{\zeta})}_F
 &\leq e^{8s\log K-T^2/(2K)}
       \norm{\bm{\zeta}}_2
       \sum_{i\in Q}\ComponentWeight{2s}{i}\,.
\end{aligned}
\end{equation}

\item\label[theorem]{thm:hermite:structure-component-weights}
For every $i, j\in Q$ and $1\leq r\leq s$, the component weights satisfy
\begin{align}
 e^{-T^2}
 \leq\Prb(T\leq Y_i\leq T+1)
 &\leq\ComponentWeight{2r}{i}
 \leq e^{4r\log K-T^2/(3K)}\,,
 \\
 K^{-(s-r)}\frac{\ComponentWeight{2s}{i}}{\ComponentWeight{2s}{j}}
 &\leq\frac{\ComponentWeight{2r}{i}}{\ComponentWeight{2r}{j}}
 \leq K^{s-r}\frac{\ComponentWeight{2s}{i}}{\ComponentWeight{2s}{j}}\,.
\end{align}

\item\label[theorem]{thm:hermite:structure-approximate-orthogonality} The final lifted vectors are nearly orthogonal:
\begin{equation}
 \abs{\ip{\TensorPower{\NormalizedRegressor_i}{s}}{\TensorPower{\NormalizedRegressor_j}{s}}}\leq e^{-s/K}
 \qquad(i\neq j).
\end{equation}

\item\label[theorem]{thm:hermite:structure-scaled-regressors}
For every $i, j\in Q$, $\ScaledRegressor_i$ is a positive multiple of $\Regressor_i$, and
\begin{align}
 \frac{T}{2\sqrt{2s+1}}
 \frac{\norm{\Regressor_i}_2}{\sigma_i^2}
 &\leq\norm{\ScaledRegressor_i}_2
 \leq
 \frac{T+1+B}{\sqrt{2s+1}}
 \frac{\norm{\Regressor_i}_2}{\sigma_i^2}\,,
  \\
 \norm{\ScaledRegressor_i-\ScaledRegressor_j}_2^2
 &\geq\frac1K
 \max\left\{ \norm{\ScaledRegressor_i}_2^2,\norm{\ScaledRegressor_j}_2^2 \right\} \qquad(i\ne j)\,.
\end{align}
\end{enumerate}
\end{theorem}

\Cref{thm:hermite:structure-moment-decomposition} approximately decomposes the Hermite moments as weighted rank-one tensor sums of the
regressors with errors bounded relative to the total component weight.
In \Cref{sec:tensor-decomposition:main}, we estimate projections of these
sums and use the rank-one matrices shared by the final even moment
and the odd contractions to recover a block of normalized regressors.

\Cref{thm:hermite:structure-component-weights} bounds the component weights on the Gaussian
tail scale and compares them with the slice probabilities.
In \Cref{sec:peeling:main}, we use the slice-probability comparison
to bound the peeling bias relative to the total component weight.
In \Cref{subsec:tensor-decomposition:iterative-projections}, the weight bounds set the estimation
accuracy, and the comparison across degrees lets us bound the effect
of the intermediate projections on the final even moment.

\Cref{thm:hermite:structure-approximate-orthogonality} shows that the final lifted vectors are
nearly orthogonal. In \Cref{subsubsec:tensor-decomposition:operator-construction-proof}, we use
this property to compare the eigenvalues of the even moment's
rank-one sum with the component weights and choose a block on
which to rescale the odd moments.

\Cref{thm:hermite:structure-scaled-regressors} shows that each $\ScaledRegressor_i$
is a positive multiple of $\Regressor_i$, with multiplier of order
$T/(\sigma_i^2\sqrt{2s+1})$, and that distinct scaled regressors
remain separated. In \Cref{subsubsec:tensor-decomposition:direction-recovery-proof}, we use
their separation to recover the common eigenvectors from a random
contraction. The lower bound on their norms controls the error
when we normalize the estimated scaled regressors.

The rest of this section is dedicated to proving \Cref{thm:hermite:restricted-moment-structure}.

\subsection{Decomposition of the Restricted Hermite Moments}

To prove \Cref{thm:hermite:structure-moment-decomposition}, we compare the restricted
moment of the remaining maximum with the sum of the restricted moments
of the individual responses. \Cref{subsubsec:hermite:single-response-moments} computes
each summand. \Cref{subsubsec:hermite:moment-decomposition-error} gives the probability and
second-moment bounds needed to control the difference.

\subsubsection{Single-Response Moments}
\label{subsubsec:hermite:single-response-moments}
For the contribution of component $i$, we first condition on
$Y_i=t$ and then integrate over the slice $[T,T+1]$.
The following identity gives the conditional Hermite moment.

\begin{restatable}[Conditional Hermite Identity]{proposition}{conditionalhermiteidentity}
\label{prop:hermite:conditional-hermite-identity}
Let $\bm{G}\in\R^n$ and $Y\in\R$ be centered jointly Gaussian,
with $\bm{G}\sim N(\bm{0},I_n)$ and $\operatorname{Var}(Y)=\sigma^2>0$.
Put $\bm{v}=\E_{\bm{G},Y}[\bm{G}Y]$.
For every integer $d\geq0$ and $t\in\R$,
\[
 \E_{\bm{G},Y}[\HermiteTensor{d}{\bm{G}}\mid Y=t]
 =\frac{\He_d(t/\sigma)}{\sigma^d\sqrt{d!}}
   \TensorPower{\bm{v}}{d}.
\]
\end{restatable}

The proof in \Cref{subsec:app-hermite:single-response-moments} applies the
scalar Gaussian noise identity (\Cref{prop:hermite-noise-operator})
to one-dimensional projections
of $\bm{G}$.
In our model, $\E_{\Covariate,\xi_i}[\Covariate(Y_i-b_i)]=\Regressor_i$.
Apply \Cref{prop:hermite:conditional-hermite-identity} to
$(\Covariate,Y_i-b_i)$, and integrate its conditional identity at
$t-b_i$ against $\phi_{\sigma_i^2}(t-b_i)$ over $[T,T+1]$.
This gives the following corollary.

\begin{corollary}[Single-Response Hermite Identity]
\label{cor:hermite:single-response-hermite-identity}
In the Gaussian self-selection model (\Cref{def:gaussian-self-selection}),
for every $i\in[k]$ and integer $d\geq0$,
\begin{equation}
 \E_{\Covariate,\xi_i}[
   \HermiteTensor{d}{\Covariate}\TailIndicator{Y_i}]
 =\MomentCoefficient{d}{b_i,\sigma_i^2}\TensorPower{\Regressor_i}{d}.
\end{equation}
\end{corollary}

\subsubsection{Moment Decomposition Error}
\label{subsubsec:hermite:moment-decomposition-error}
The single-response sum counts $Y_i$ whenever it lies in the slice,
even if another remaining response is larger. We first bound the
probability that another response exceeds $t$ given $Y_i=t$, then
bound the conditional second moment of the Hermite feature.
Conditional on $Y_i=t$, the uncoveredness assumption ensures that the
mean of each other response is below $t$ by at least a constant fraction
of $t$, while the boundedness assumption controls its conditional variance.

\begin{lemma}[A Large Response Is Rarely Exceeded]
\label{lem:hermite:large-response-exceedance}
Consider the Gaussian self-selection model
(\Cref{def:gaussian-self-selection}) under the boundedness and
uncoveredness assumptions (\Cref{ass:boundedness-uncoveredness}).
Use the hyperparameters $K,T,s$ from \Cref{hyp:algorithm}.
For every non-empty $Q\subseteq[k]$, $i\in Q$, and $t\in[T,T+1]$,
\[
 \Prb(Y_{\max,Q}>t\mid Y_i=t)\leq k e^{-2t^2/K}.
\]
\end{lemma}

\begin{proof}
By the choice of \Cref{hyp:algorithm},
\begin{equation}\label{eq:hermite:exceedance-k-choice}
 K\geq\frac{16(1+B^2)^3}{\Delta^4},
 \qquad
 T\geq\frac{4B(1+B^2)}{\Delta^2}.
\end{equation}
Conditional on $Y_i=t$, another response must exceed $t$
for $Y_{\max,Q}$ to exceed $t$.
Fix $j\in Q\setminus\{i\}$ and write
$c_{ij}\coloneqq\ip{\Regressor_i}{\Regressor_j}+\rho$.
Independence of the centered noise vector and $\Covariate$ gives
\[
 \E_{\Covariate,\xi_i,\xi_j}[(Y_i-b_i)(Y_j-b_j)]
 =\Regressor_i^\top\E_{\Covariate}[\Covariate\Covariate^\top]\Regressor_j+\rho
 =c_{ij}.
\]
The conditional Gaussian formula (\Cref{prop:gaussian-conditioning})
gives
\[
 Y_j\mid Y_i=t
 \sim N\!\left(
 b_j+\frac{c_{ij}}{\sigma_i^2}(t-b_i),\,
 \sigma_j^2-\frac{c_{ij}^2}{\sigma_i^2}
 \right)
 \eqqcolon N(\mu,\tau^2).
\]
By the Gaussian tail bound (\Cref{prop:gaussian-tail}),
\begin{align*}
 t-\mu
 &=\frac{(\sigma_i^2-c_{ij})t+c_{ij}b_i-\sigma_i^2b_j}{\sigma_i^2}\\
 &\geq\frac{\Delta^2t}{\sigma_i^2}-2B
 \alignnote{boundedness and uncoveredness (\Cref{ass:boundedness-uncoveredness})}\\
 &\geq\frac{\Delta^2t}{2\sigma_i^2}>0,
 \alignnote{$t\geq T$, by \eqref{eq:hermite:exceedance-k-choice}}\\
 \frac{(t-\mu)^2}{2\tau^2}
 &\geq\frac{\Delta^4t^2}{8\sigma_i^4\tau^2}\\
 &\geq\frac{\Delta^4t^2}{8(1+B^2)^3}
 \alignnote{$\tau^2\leq\sigma_j^2$; $\sigma_i^2,\sigma_j^2\leq1+B^2$}\\
 &\geq\frac{2t^2}{K}.
 \alignnote{choice of $K$~\eqref{eq:hermite:exceedance-k-choice}}
\end{align*}
It follows that
\begin{align*}
 \Prb(Y_j>t\mid Y_i=t)
 &\leq\exp\!\left(-\frac{(t-\mu)^2}{2\tau^2}\right)
 \leq e^{-2t^2/K},\\
 \Prb(Y_{\max,Q}>t\mid Y_i=t)
 &\leq\sum_{j\in Q\setminus\{i\}}
       \Prb(Y_j>t\mid Y_i=t)
 \alignnote{union bound}\\
 &\leq(|Q|-1)e^{-2t^2/K}
 \leq k e^{-2t^2/K}.\qedhere
\end{align*}
\end{proof}

To bound the contribution of responses that enter the slice
without attaining the maximum, we combine
\Cref{lem:hermite:large-response-exceedance} with a conditional second-moment
bound for the Hermite feature by Cauchy--Schwarz.

\begin{restatable}[Conditional Hermite Second Moment]{proposition}{conditionalhermitesecondmoment}
\label{prop:hermite:conditional-hermite-second-moment}
Let $\bm{G}\in\R^n$ and $Y\in\R$ be centered jointly Gaussian,
with $\bm{G}\sim N(\bm{0},I_n)$ and $\operatorname{Var}(Y)=\sigma^2>0$.
Put $\bm{v}=\E_{\bm{G},Y}[\bm{G}Y]$.
For every integer $d\geq0$, symmetric order-$d$ tensor $\mathsf A$,
$t\in\R$, and $0<\theta<1/2$,
\[
 \left(
  \E_{\bm{G},Y}[\HermitePairing{d}{\bm{G}}{\mathsf A}^2\mid Y=t]
 \right)^{1/2}
 \leq
 \frac{\theta^{-d/2}}{\sqrt{1-2\theta}}
 \exp\!\left(
  \frac{\theta\norm {\bm{v}}_2^2t^2}{2\sigma^4(1-2\theta)}
 \right)\norm{\mathsf A}_F.
\]
\end{restatable}

We prove the proposition in \Cref{subsec:app-hermite:moment-decomposition-error}.
Taking $(\bm{G},Y)=(\Covariate,Y_i-b_i)$ gives the following corollary.

\begin{corollary}[Second Conditional Hermite Moment]
\label{cor:hermite:single-response-second-moment}
Consider the Gaussian self-selection model
(\Cref{def:gaussian-self-selection}) under the boundedness and
uncoveredness assumptions (\Cref{ass:boundedness-uncoveredness}).
Let $K\geq2$, $T\geq B+1$, and let $s\geq1$ be an integer.
For every $i\in[k]$, $1\leq d\leq2s+1$, $t\in[T,T+1]$, and
symmetric order-$d$ tensor $\mathsf A$ with $\norm{\mathsf A}_F\leq1$,
\[
 \left(\E_{\Covariate,\xi_i}[\HermitePairing{d}{\Covariate}{\mathsf A}^2\mid Y_i=t]\right)^{1/2}
 \leq e^{7s\log K+T^2/K^4}.
\]
\end{corollary}

\begin{proof}
Apply \Cref{prop:hermite:conditional-hermite-second-moment}
with $(\bm{G},Y)=(\Covariate,Y_i-b_i)$ and $\theta=K^{-4}$, conditioning at $t-b_i$.
Here $\bm{v}=\Regressor_i$, and
\[
 \frac{\norm{\Regressor_i}_2}{\sigma_i^2}
 =\frac{\norm{\Regressor_i}_2}{1+\norm{\Regressor_i}_2^2}
 \leq\frac12.
\]
Therefore
\begin{align*}
 \left(
  \E_{\Covariate,\xi_i}[\HermitePairing{d}{\Covariate}{\mathsf A}^2\mid Y_i=t]
 \right)^{1/2}
 &\leq
 \frac{K^{2d}}{\sqrt{1-2K^{-4}}}
 \exp\!\left(\frac{(t-b_i)^2}{8K^4(1-2K^{-4})}\right)
 \alignnote{$\norm{\mathsf A}_F\leq1$}\\
 &\leq\sqrt2\,K^{2d}e^{(t-b_i)^2/(4K^4)}
 \alignnote{$K^{-4}\leq1/16$}\\
 &\leq e^{7s\log K+T^2/K^4}.
 \alignnote{$d\leq3s$, $\sqrt2\leq K^s$, $|t-b_i|\leq2T$}
\end{align*}
\end{proof}

\subsection{Bounds on the Component Weights}

To prove \Cref{thm:hermite:structure-component-weights}, we first establish positivity
and ratio bounds for the restricted moment coefficients. We use these
bounds to compare each component weight with its slice probability,
then to compare the relative weights of two components across degrees.

\subsubsection{Restricted Moment Coefficients}
The weights depend on the even coefficients, while the scaled regressors depend
on the ratio of two consecutive coefficients. We first bound these ratios for
arbitrary degrees.
\begin{lemma}[Restricted Moment Coefficient Estimates]
\label{lem:hermite:restricted-moment-coefficients}
Use the choice of $K$ from \Cref{hyp:algorithm}.
Let $s\geq1$ be an integer, and suppose $|b|\leq B$, $T\geq3B$,
and $(T-B)^2\geq Ks$.
For every $\sigma^2\in[1+\Delta^2,1+B^2]$,
the coefficients $\MomentCoefficient{d}{b,\sigma^2}$ are positive for
$0\leq d\leq2s+2$.
Moreover, for every pair of integers $0\leq {d_1}\leq {d_2}\leq2s+2$,
\begin{equation}\label{eq:hermite:coefficient-degree-ratio}
 \sqrt{\frac{{d_1}!}{{d_2}!}}
 \left(\frac{T}{2\sigma^2}\right)^{{d_2}-{d_1}}
 \leq\frac{\MomentCoefficient{d_2}{b,\sigma^2}}{\MomentCoefficient{d_1}{b,\sigma^2}}
 \leq
 \sqrt{\frac{{d_1}!}{{d_2}!}}
 \left(\frac{T+1+B}{\sigma^2}\right)^{{d_2}-{d_1}}.
\end{equation}
\end{lemma}

\begin{proof}
The requirement on $K$ gives $K\geq16(1+B^2)$.
Consequently,
\[
 (T-B)^2\geq Ks\geq16s\sigma^2.
\]
Fix $x\geq4\sqrt s$.
We use the standard Hermite recurrence~\cite[Table~18.9.1]{DLMF},
\[
 \He_{j+1}(x)=x\He_j(x)-j\He_{j-1}(x).
\]
Starting from $\He_0(x)=1$ and $\He_1(x)=x$, we show by induction
that the ratios $q_j=\He_j(x)/\He_{j-1}(x)$ are well-defined and
lie in $[3x/4,x]$ for $1\leq j\leq2s+2$.
Indeed, $q_1=x$, and for $1\leq j<2s+2$ the recurrence gives
\[
 \frac{3x}{4}\leq x-\frac{4j}{3x}
 \leq q_{j+1}=x-\frac{j}{q_j}\leq x,
\]
since $j\leq2s+1\leq3s$ and $x^2\geq16s$.
Multiplying these ratios gives, for $1\leq d\leq2s+2$,
\[
 (3x/4)^d\leq\He_d(x)=\prod_{j=1}^d q_j\leq x^d.
\]
For $t\in[T,T+1]$, the shifted argument satisfies
$(t-b)/\sigma\geq4\sqrt s$. Thus the integrand defining
$\MomentCoefficient{d}{b,\sigma^2}$ is positive for $0\leq d\leq2s+2$.

We next compare coefficients of different degrees.
Fix $0\leq d_1\leq d_2\leq2s+2$.
For $t\in[T,T+1]$, the recurrence bounds give
\[
 \frac{\He_{d_2}((t-b)/\sigma)}{\He_{d_1}((t-b)/\sigma)}
 =\prod_{j=d_1+1}^{d_2}
   \frac{\He_j((t-b)/\sigma)}{\He_{j-1}((t-b)/\sigma)}
 \in\left[
 \left(\frac{3(t-b)}{4\sigma}\right)^{d_2-d_1},
 \left(\frac{t-b}{\sigma}\right)^{d_2-d_1}
 \right].
\]
When $d_1=d_2$, the product is empty and equals one.
The hypotheses imply $3(t-b)/4\geq3(T-B)/4\geq T/2$
and $t-b\leq T+1+B$.
Including the normalization factors therefore gives
\[
 \sqrt{\frac{d_1!}{d_2!}}
 \left(\frac{T}{2\sigma^2}\right)^{d_2-d_1}
 \leq
 \frac{\sigma^{-d_2}\He_{d_2}((t-b)/\sigma)/\sqrt{d_2!}}
      {\sigma^{-d_1}\He_{d_1}((t-b)/\sigma)/\sqrt{d_1!}}
 \leq
 \sqrt{\frac{d_1!}{d_2!}}
 \left(\frac{T+1+B}{\sigma^2}\right)^{d_2-d_1}.
\]
Multiply by the positive function
$\phi_{\sigma^2}(t-b)\sigma^{-d_1}\He_{d_1}((t-b)/\sigma)/\sqrt{d_1!}$
and integrate over $[T,T+1]$. This yields
\[
 \sqrt{\frac{d_1!}{d_2!}}
 \left(\frac{T}{2\sigma^2}\right)^{d_2-d_1}\MomentCoefficient{d_1}{b,\sigma^2}
 \leq\MomentCoefficient{d_2}{b,\sigma^2}
 \leq
 \sqrt{\frac{d_1!}{d_2!}}
 \left(\frac{T+1+B}{\sigma^2}\right)^{d_2-d_1}\MomentCoefficient{d_1}{b,\sigma^2}.
\]
Dividing by $\MomentCoefficient{d_1}{b,\sigma^2}>0$ proves
\eqref{eq:hermite:coefficient-degree-ratio}.
\end{proof}

\subsubsection{Slice Probabilities and Weights}
The degree-zero coefficient is the slice probability. Comparing it with the
degree-$2r$ coefficient gives the weight bounds in \Cref{thm:hermite:structure-component-weights}
and the comparison needed for the error bound in \Cref{thm:hermite:structure-moment-decomposition}.

\begin{lemma}[Even-Moment Weight Bounds]
\label{lem:hermite:even-moment-weights}
Consider the Gaussian self-selection model
(\Cref{def:gaussian-self-selection}) under the boundedness and
uncoveredness assumptions (\Cref{ass:boundedness-uncoveredness}).
Use the choice of $K$ from \Cref{hyp:algorithm}.
Let $s\geq1$ be an integer, and suppose $T\geq5(B+1)$
and $(T-B)^2\geq Ks$.
For every $i\in[k]$ and $1\leq r\leq s$,
\begin{equation}\label{eq:hermite:weight-slice-lower}
 \Prb(T\leq Y_i\leq T+1)
 \leq\ComponentWeight{2r}{i}
 \leq e^{4r\log K+T^2/(2K^4)}
       \Prb(T\leq Y_i\leq T+1).
\end{equation}
Consequently,
\begin{equation}\label{eq:hermite:weight-slice-upper}
 e^{-T^2}\leq\ComponentWeight{2r}{i}
 \leq e^{4r\log K-T^2/(3K)}.
\end{equation}
\end{lemma}

\begin{proof}
Fix $i\in[k]$ and $1\leq r\leq s$.
We first compare the weight with the slice probability.
Since $\He_0=1$,
\[
 \MomentCoefficient{0}{b_i,\sigma_i^2}
 =\Prb(T\leq Y_i\leq T+1)>0.
\]
Consequently,
\begin{align*}
 \frac{\ComponentWeight{2r}{i}}{\MomentCoefficient{0}{b_i,\sigma_i^2}}
 &=\frac{\MomentCoefficient{2r}{b_i,\sigma_i^2}}{\MomentCoefficient{0}{b_i,\sigma_i^2}}
   \norm{\Regressor_i}_2^{2r}\\
 &\geq
 \frac{(T\norm{\Regressor_i}_2/(2\sigma_i^2))^{2r}}
      {\sqrt{(2r)!}}
 \alignnote{by \Cref{lem:hermite:restricted-moment-coefficients}}\\
 &\geq
 \left(\frac{T^2\Delta^2}{8r(1+B^2)^2}\right)^r
 \alignnote{$(2r)!\leq(2r)^{2r}$;
  $\norm{\Regressor_i}_2\geq\Delta$;
  $\sigma_i^2\leq1+B^2$}\\
 &\geq\left(\frac{2T^2}{Kr}\right)^r
 \alignnote{choice of $K$: $K\Delta^2\geq16(1+B^2)^2$}\\
 &\geq1.
 \alignnote{$T^2\geq Ks$, $r\leq s$}
\end{align*}
Multiplying by the positive coefficient $\MomentCoefficient{0}{b_i,\sigma_i^2}$
proves the lower comparison in \eqref{eq:hermite:weight-slice-lower}.
For the upper comparison, use
$2\norm{\Regressor_i}_2\leq1+\norm{\Regressor_i}_2^2=\sigma_i^2$
and $T\geq B+1$.  Then
\begin{align*}
 \frac{\ComponentWeight{2r}{i}}{\MomentCoefficient{0}{b_i,\sigma_i^2}}
 &\leq
 \frac{\bigl((T+1+B)\norm{\Regressor_i}_2/\sigma_i^2\bigr)^{2r}}
      {\sqrt{(2r)!}}
 \alignnote{by \Cref{lem:hermite:restricted-moment-coefficients}}\\
 &\leq\frac{T^{2r}}{\sqrt{(2r)!}}
 \alignnote{$2\norm{\Regressor_i}_2\leq\sigma_i^2$, $T+1+B\leq2T$}\\
 &=K^{4r}
   \left(\frac{(T^2/K^4)^{2r}}{(2r)!}\right)^{1/2}\\
 &\leq K^{4r}
   \left(\sum_{\ell=0}^{\infty}
         \frac{(T^2/K^4)^\ell}{\ell!}\right)^{1/2}
 \alignnote{all terms are nonnegative}\\
 &=e^{4r\log K+T^2/(2K^4)}.
\end{align*}
Multiplying by $\MomentCoefficient{0}{b_i,\sigma_i^2}>0$ proves the upper comparison
in \eqref{eq:hermite:weight-slice-lower}.

It remains to prove \eqref{eq:hermite:weight-slice-upper}.
The hypotheses give $K\geq128$ and $T^2\geq K$, so $T\geq5$.
For the lower bound,
\begin{align*}
 \Prb(T\leq Y_i\leq T+1)
 &\geq\phi_{\sigma_i^2}(T+1-b_i)
 \alignnote{the density decreases on $[T,T+1]$}\\
 &\geq
 \exp\!\left(-\frac{(T+1+B)^2}{2}
             -\frac12\log(2\pi K)\right)
 \alignnote{$1\leq\sigma_i^2\leq K$}\\
 &\geq e^{-3T^2/4-K/4}
 \alignnote{$T+1+B\leq6T/5$, $K\geq128$}\\
 &\geq e^{-T^2}.
 \alignnote{$T^2\geq K$}
\end{align*}
Combining this with \eqref{eq:hermite:weight-slice-lower} gives the
lower bound on $\ComponentWeight{2r}{i}$.
For the upper bound,
\begin{align*}
 \ComponentWeight{2r}{i}
 &\leq e^{4r\log K+T^2/(2K^4)}\phi_{\sigma_i^2}(T-b_i)
 \alignnote{by \eqref{eq:hermite:weight-slice-lower}; decreasing density}\\
 &\leq e^{4r\log K+T^2/(2K^4)-T^2/(2K)}
 \alignnote{$T-B\geq T/2$, $K\geq4\sigma_i^2$}\\
 &\leq e^{4r\log K-T^2/(3K)}.
 \alignnote{$K^3\geq3$}
\end{align*}
Together, these bounds prove the claim.
\end{proof}

\subsubsection{Weights Across Degrees}
We have bounded each weight at a fixed degree. For the intermediate
projections, we also need to know how the ratio of two component weights
changes between levels $r$ and $s$.

\begin{lemma}[Weights Across Degrees]
\label{lem:hermite:weights-across-degrees}
Consider the Gaussian self-selection model
(\Cref{def:gaussian-self-selection}) under the boundedness and
uncoveredness assumptions (\Cref{ass:boundedness-uncoveredness}).
Use the hyperparameters $K,T,s$ from \Cref{hyp:algorithm}.
For every $i,j\in[k]$ and $1\leq r\leq s$,
\[
 \abs{
 \log\frac{\ComponentWeight{2r}{j}/\ComponentWeight{2r}{i}}
          {\ComponentWeight{2s}{j}/\ComponentWeight{2s}{i}}
 }
 \leq(s-r)\log K.
\]
\end{lemma}

\begin{proof}
By the choice of \Cref{hyp:algorithm},
\begin{equation}\label{eq:hermite:cross-degree-k-choice}
 K\geq\frac{16(1+B^2)^3}{\Delta^4}.
\end{equation}
The same choice gives $T\geq2(B+1)$.
Fix $i,j\in[k]$ and $1\leq r\leq s$.
The coefficients are positive by \Cref{lem:hermite:restricted-moment-coefficients}.
We compare the two components directly:
\begin{align}
 \frac{\ComponentWeight{2r}{j}/\ComponentWeight{2r}{i}}
      {\ComponentWeight{2s}{j}/\ComponentWeight{2s}{i}}
 &=
 \frac{\MomentCoefficient{2s}{b_i,\sigma_i^2}/\MomentCoefficient{2r}{b_i,\sigma_i^2}}
      {\MomentCoefficient{2s}{b_j,\sigma_j^2}/\MomentCoefficient{2r}{b_j,\sigma_j^2}}
 \left(\frac{\norm{\Regressor_i}_2}{\norm{\Regressor_j}_2}\right)^{2(s-r)}
 \notag\\
 &\leq
 \frac{
  \sqrt{\frac{(2r)!}{(2s)!}}
  \bigl((T+1+B)\norm{\Regressor_i}_2/\sigma_i^2\bigr)^{2(s-r)}
 }{
  \sqrt{\frac{(2r)!}{(2s)!}}
  \bigl(T\norm{\Regressor_j}_2/(2\sigma_j^2)\bigr)^{2(s-r)}
 }
 \alignnote{by \Cref{lem:hermite:restricted-moment-coefficients}}
 \notag\\
 &=
 \left(
  \frac{2(T+1+B)}{T}
  \frac{\norm{\Regressor_i}_2}{\norm{\Regressor_j}_2}
  \frac{\sigma_j^2}{\sigma_i^2}
 \right)^{2(s-r)}
 \notag\\
 &\leq
 \left(\frac{3B(1+B^2)}{\Delta}\right)^{2(s-r)}
 \alignnote{by \Cref{ass:boundedness-uncoveredness}; $T\geq2(B+1)$}
 \notag\\
 &\leq
 \left(\frac{16(1+B^2)^3}{\Delta^4}\right)^{s-r}
 \alignnote{$9B^2\Delta^2\leq16(1+B^2)$}
 \notag\\
 &\leq K^{s-r}.
 \alignnote{choice of $K$~\eqref{eq:hermite:cross-degree-k-choice}}
\end{align}
Interchanging $i$ and $j$ gives the same bound for the reciprocal.
Taking logarithms proves the claim.
\end{proof}

\subsection{Proof of \texorpdfstring{\Cref*{thm:hermite:restricted-moment-structure}}{Theorem 5.1}}
We first combine the single-response identity and the error estimates
from \Cref{subsubsec:hermite:single-response-moments,subsubsec:hermite:moment-decomposition-error} to obtain
the even and odd decompositions. The weight comparisons give
\Cref{thm:hermite:structure-component-weights}. We finish by applying uncoveredness to the
lifted vectors and the coefficient ratio bounds to the scaled regressors.

\begin{proof}[Proof of \cref{thm:hermite:restricted-moment-structure}]
We first prove \Cref{thm:hermite:structure-moment-decomposition} by identifying the rank-one terms
and bounding the contribution
of latent responses that enter the slice without attaining the maximum.

Let $I_Q$ be the maximizing index within $Q$, using the global tie rule.
Decomposing the slice by this index gives
\[
 \TailIndicator{Y_{\max,Q}}
 =\sum_{i\in Q}\1\{I_Q=i,\ T\leq Y_i\leq T+1\}.
\]
For $1\leq d\leq2s+1$, define
\[
 \MomentRemainder{d}^Q\coloneqq\sum_{i\in Q}\E_{\Covariate,\NoiseVector_Q}\!\left[
 \HermiteTensor{d}{\Covariate}\1\{T\leq Y_i\leq T+1,\ I_Q\ne i\}\right].
\]
Since $\1\{I_Q=i\}=1-\1\{I_Q\ne i\}$, we obtain
\begin{align}
 \RestrictedMomentTensor{d}{Q}
 &=\sum_{i\in Q}\E_{\Covariate,\NoiseVector_Q}\!\left[
  \HermiteTensor{d}{\Covariate}\TailIndicator{Y_i}\1\{I_Q=i\}
 \right]\notag\\
 &=\sum_{i\in Q}\E_{\Covariate,\xi_i}[
  \HermiteTensor{d}{\Covariate}\TailIndicator{Y_i}]-\MomentRemainder{d}^Q\notag\\
 &=\sum_{i\in Q}\MomentCoefficient{d}{b_i,\sigma_i^2}\TensorPower{\Regressor_i}{d}-\MomentRemainder{d}^Q.
 \label{eq:hermite:winner-decomposition}
\end{align}
The last equality applies \Cref{cor:hermite:single-response-hermite-identity} to each latent response.
Thus $\MomentRemainder{d}^Q$ collects the contributions of responses that enter the
slice without attaining the maximum.

We now bound the Frobenius norm of the error tensor.
Fix a symmetric order-$d$ tensor $\mathsf A$ with $\norm{\mathsf A}_F\leq1$.
Conditional on $Y_i=t$, ties have probability zero, so
$\Prb(I_Q\ne i\mid Y_i=t)=\Prb(Y_{\max,Q}>t\mid Y_i=t)$.
For $t\in[T,T+1]$, \Cref{lem:hermite:large-response-exceedance} bounds the probability
that another response exceeds $t$ by $k e^{-2t^2/K}$ conditional on
$Y_i=t$. Under the same conditioning, \Cref{cor:hermite:single-response-second-moment}
bounds the $L^2$ norm of the Hermite pairing by
$e^{7s\log K+T^2/K^4}$.
Conditioning on each latent response gives
\begin{align*}
 \abs{\ip{\MomentRemainder{d}^Q}{\mathsf A}}
 &=
 \abs{\sum_{i\in Q}\int_T^{T+1}
  \E_{\Covariate,\NoiseVector_Q}\!\left[
   \HermitePairing{d}{\Covariate}{\mathsf A}\1\{I_Q\ne i\}
   \mid Y_i=t
  \right]\phi_{\sigma_i^2}(t-b_i) \odif{t}}\\
 &\leq\sum_{i\in Q}\int_T^{T+1}
  \E_{\Covariate,\NoiseVector_Q}\!\left[
   \abs{\HermitePairing{d}{\Covariate}{\mathsf A}}\1\{I_Q\ne i\}
   \mid Y_i=t
  \right]\phi_{\sigma_i^2}(t-b_i) \odif{t}
 \alignnote{triangle inequality}\\
 &\leq\sum_{i\in Q}\int_T^{T+1}
  \left(\E_{\Covariate,\xi_i}[
   \HermitePairing{d}{\Covariate}{\mathsf A}^2\mid Y_i=t
  ]\right)^{1/2}
  \Prb(I_Q\ne i\mid Y_i=t)^{1/2}
  \phi_{\sigma_i^2}(t-b_i) \odif{t}
 \alignnote{conditional Cauchy--Schwarz}\\
 &\leq e^{7s\log K+\frac12\log k+T^2/K^4}
  \sum_{i\in Q}\int_T^{T+1}
   e^{-t^2/K}\phi_{\sigma_i^2}(t-b_i) \odif{t}
 \alignnote{by \Cref{cor:hermite:single-response-second-moment,lem:hermite:large-response-exceedance}}\\
 &\leq e^{8s\log K+T^2/K^4}
  \sum_{i\in Q}\int_T^{T+1}
   e^{-t^2/K}\phi_{\sigma_i^2}(t-b_i) \odif{t}
 \alignnote{choice of $s$: $\log k\leq s$; $\log K\geq1/2$}\\
 &\leq e^{8s\log K}
  \sum_{i\in Q}\int_T^{T+1}
   e^{-t^2/(2K)}\phi_{\sigma_i^2}(t-b_i) \odif{t}
 \alignnote{$t\geq T$, $K^3\geq2$}\\
 &\leq e^{8s\log K-T^2/(2K)}
  \sum_{i\in Q}\Prb(T\leq Y_i\leq T+1).
 \alignnote{$t\geq T$}
\end{align*}

Since $\MomentRemainder{d}^Q$ is symmetric, taking the supremum over
symmetric $\mathsf A$ with $\norm{\mathsf A}_F\leq1$ bounds its
Frobenius norm. By \Cref{lem:hermite:even-moment-weights}, each slice probability
is at most $\ComponentWeight{2r}{i}$ for every $1\leq r\leq s$. Hence
\begin{equation}\label{eq:hermite:relative-tensor-error}
 \norm{\MomentRemainder{d}^Q}_F
 \leq e^{8s\log K-T^2/(2K)}
      \sum_{i\in Q}\ComponentWeight{2r}{i}
 \qquad(1\leq d\leq2s+1,\ 1\leq r\leq s).
\end{equation}

We now apply flattening to obtain the even moment matrices.
For $1\leq r\leq s$, recall that
$\Regressor_i=\norm{\Regressor_i}_2\NormalizedRegressor_i$, which gives
\[
 \Mat_r(\TensorPower{\Regressor_i}{2r})
 =\norm{\Regressor_i}_2^{2r}\TensorPower{\NormalizedRegressor_i}{r}(\TensorPower{\NormalizedRegressor_i}{r})^\top.
\]
Applying $\Mat_r$ to the winner decomposition
\eqref{eq:hermite:winner-decomposition}
with $d=2r$ therefore gives
\begin{align*}
 \Mat_r(\RestrictedMomentTensor{2r}{Q})
 &=\sum_{i\in Q}\MomentCoefficient{2r}{b_i,\sigma_i^2}\norm{\Regressor_i}_2^{2r}
   \TensorPower{\NormalizedRegressor_i}{r}(\TensorPower{\NormalizedRegressor_i}{r})^\top-\Mat_r(\MomentRemainder{2r}^Q)\\
 &=\sum_{i\in Q}\ComponentWeight{2r}{i}
   \TensorPower{\NormalizedRegressor_i}{r}(\TensorPower{\NormalizedRegressor_i}{r})^\top+E_r.
\end{align*}
Here $E_r=-\Mat_r(\MomentRemainder{2r}^Q)$.
The last equality restates the weight definition
in \eqref{eq:component-weights}.
Flattening preserves the Frobenius norm, so
\eqref{eq:hermite:relative-tensor-error} with $d=2r$
gives the even error bound in \Cref{thm:hermite:structure-moment-decomposition}.

For the odd moment, fix $\bm{\zeta}\in\R^n$.
Contracting one index gives
\[
 \TensorPower{\Regressor_i}{(2s+1)}\mathbin{\lrcorner}\bm{\zeta}
 =\ip{\bm{\zeta}}{\Regressor_i}\TensorPower{\Regressor_i}{2s}.
\]
Contracting one index against $\bm{\zeta}$ and then applying $\Mat_s$
to the winner decomposition
\eqref{eq:hermite:winner-decomposition}
with $d=2s+1$ gives
\begin{align*}
 &\Mat_s(\RestrictedMomentTensor{2s+1}{Q}\mathbin{\lrcorner}\bm{\zeta})\\
 &\quad=\sum_{i\in Q}\MomentCoefficient{2s+1}{b_i,\sigma_i^2}
   \ip{\bm{\zeta}}{\Regressor_i}\norm{\Regressor_i}_2^{2s}
   \TensorPower{\NormalizedRegressor_i}{s}(\TensorPower{\NormalizedRegressor_i}{s})^\top
   -\Mat_s(\MomentRemainder{2s+1}^Q\mathbin{\lrcorner}\bm{\zeta})\\
 &\quad=\sum_{i\in Q}\ComponentWeight{2s}{i}
   \ip{\bm{\zeta}}{\ScaledRegressor_i}\TensorPower{\NormalizedRegressor_i}{s}(\TensorPower{\NormalizedRegressor_i}{s})^\top+E(\bm{\zeta}).
\end{align*}
Here $E(\bm{\zeta})=-\Mat_s(\MomentRemainder{2s+1}^Q\mathbin{\lrcorner}\bm{\zeta})$.
The last equality restates the definitions in
\eqref{eq:component-weights} and
\eqref{eq:scaled-regressors}, since
\[
 \MomentCoefficient{2s+1}{b_i,\sigma_i^2}\ip{\bm{\zeta}}{\Regressor_i}\norm{\Regressor_i}_2^{2s}
 =\MomentCoefficient{2s}{b_i,\sigma_i^2}\norm{\Regressor_i}_2^{2s}
   \ip{\bm{\zeta}}{\frac{\MomentCoefficient{2s+1}{b_i,\sigma_i^2}}{\MomentCoefficient{2s}{b_i,\sigma_i^2}}\Regressor_i}
 =\ComponentWeight{2s}{i}\ip{\bm{\zeta}}{\ScaledRegressor_i}.
\]
By Cauchy--Schwarz in the contracted index,
\[
 \norm{E(\bm{\zeta})}_F
 =\norm{\MomentRemainder{2s+1}^Q\mathbin{\lrcorner}\bm{\zeta}}_F
 \leq\norm{\MomentRemainder{2s+1}^Q}_F\norm{\bm{\zeta}}_2.
\]
Taking $d=2s+1$ and $r=s$ in
\eqref{eq:hermite:relative-tensor-error} gives the
odd error bound and proves \Cref{thm:hermite:structure-moment-decomposition}.

\emph{Proof of \Cref{thm:hermite:structure-component-weights}.}
The proof of \Cref{lem:hermite:even-moment-weights} bounds the slice probability
below by $e^{-T^2}$ and above by $\ComponentWeight{2r}{i}$, and bounds
the weight above by $e^{4r\log K-T^2/(3K)}$.
By \Cref{lem:hermite:weights-across-degrees}, the log ratio of two component weights
changes by at most $(s-r)\log K$ between levels $r$ and $s$.
Exponentiating gives the ratio bounds and proves \Cref{thm:hermite:structure-component-weights}.

\emph{Proof of \Cref{thm:hermite:structure-approximate-orthogonality}.}
By \Cref{ass:boundedness-uncoveredness}, for distinct $i,j$,
\[
 \abs{\ip{\NormalizedRegressor_i}{\NormalizedRegressor_j}}
 \leq1-\frac{\Delta^2}{B^2}.
\]
Taking tensor powers gives
\[
 \abs{\ip{\TensorPower{\NormalizedRegressor_i}{s}}{\TensorPower{\NormalizedRegressor_j}{s}}}
 =\abs{\ip{\NormalizedRegressor_i}{\NormalizedRegressor_j}}^s
 \leq\left(1-\frac{\Delta^2}{B^2}\right)^s
 \leq e^{-s\Delta^2/B^2}
 \leq e^{-s/K}.
\]
The last inequality uses $K\Delta^2/B^2\geq1$, a consequence
of the choice of \Cref{hyp:algorithm}.
This proves \Cref{thm:hermite:structure-approximate-orthogonality}.

\emph{Proof of \Cref{thm:hermite:structure-scaled-regressors}.}
By \Cref{lem:hermite:restricted-moment-coefficients},
\[
 0<\frac{T}{2\sigma_i^2\sqrt{2s+1}}
 \leq\frac{\MomentCoefficient{2s+1}{b_i,\sigma_i^2}}
          {\MomentCoefficient{2s}{b_i,\sigma_i^2}}
 \leq\frac{T+1+B}{\sigma_i^2\sqrt{2s+1}}.
\]
Multiplying by $\norm{\Regressor_i}_2$ gives the norm bounds in
\Cref{thm:hermite:structure-scaled-regressors}.
In particular, $\ScaledRegressor_i/\norm{\ScaledRegressor_i}_2=\NormalizedRegressor_i$.

To prove separation, project $\ScaledRegressor_i-\ScaledRegressor_j$ onto $\NormalizedRegressor_j^\perp$.
Since $\ScaledRegressor_j$ lies in the direction $\NormalizedRegressor_j$,
\begin{align*}
 \norm{\ScaledRegressor_i-\ScaledRegressor_j}_2^2
 &\geq\norm{(I_n-\NormalizedRegressor_j\NormalizedRegressor_j^\top)(\ScaledRegressor_i-\ScaledRegressor_j)}_2^2
 \alignnote{orthogonal projection}\\
 &=\norm{\ScaledRegressor_i}_2^2
   \left(1-\ip{\NormalizedRegressor_i}{\NormalizedRegressor_j}^2\right)\\
 &\geq\norm{\ScaledRegressor_i}_2^2
   \left(1-(1-\Delta^2/B^2)^2\right)
 \alignnote{by \Cref{ass:boundedness-uncoveredness}}\\
 &\geq\frac{\Delta^2}{B^2}\norm{\ScaledRegressor_i}_2^2
 \alignnote{$0<\Delta^2/B^2\leq1$}\\
 &\geq\frac1K\norm{\ScaledRegressor_i}_2^2.
 \alignnote{choice of $K$ (\Cref{hyp:algorithm}): $K\Delta^2\geq B^2$}
\end{align*}
Interchanging $i$ and $j$ gives the separation bound in
\Cref{thm:hermite:structure-scaled-regressors} and completes the proof.
\end{proof}

\section{Implicit Tensor Decomposition}
\label{sec:tensor-decomposition:main}

\Cref{sec:hermite:main} shows that the restricted Hermite
moments of the maximum response among the unrecovered regressors are
close to weighted sums of tensor powers of those regressors.
At degree $s=\Theta(\log n)$, the lifted vectors
$\TensorPower{\NormalizedRegressor_i}{s}$ are nearly orthogonal.
In this section, we present a tensor decomposition algorithm for these implicit high-order tensors.
More specifically, given a set $\widehat{\mathcal R}\subseteq\R^n\times\R$ of recovered
regressors and biases, with a bijection
$\match:\widehat{\mathcal R}\to R\subsetneq[k]$,
we will recover unit vectors parallel to some subset of the remaining regressors $\Regressor_i, i\in Q\coloneqq [k]\setminus R$.

For \textsc{ImplicitTensorDecomposition}, it suffices that the recovered pairs
have sum error at most $e^{-K^6s}$.

\begin{theorem}[Implicit Tensor Decomposition]\label{thm:tensor-decomposition:implicit-tensor-decomposition}
Consider the Gaussian self-selection model
(\Cref{def:gaussian-self-selection}) under the boundedness and
uncoveredness assumptions (\Cref{ass:boundedness-uncoveredness}).
Suppose the recovered list
$\widehat{\mathcal R}=((\widehat{\Regressor}_i,\widehat b_i))_{i\in R}$
is indexed by its matching to $R\subsetneq[k]$ and satisfies
\[
 \norm{\widehat{\Regressor}_i-\Regressor_i}_2+|\widehat b_i-b_i|
 \leq e^{-K^6s}\qquad(i\in R).
\]
Put $Q=[k]\setminus R$.
Given fresh samples independent of $\widehat{\mathcal R}$ and
a confidence parameter $\delta\in(0,\nicefrac13)$,
\mbox{\textnormal{\textsc{ImplicitTensorDecomposition}}}
(\cref{alg:tensor-decomposition:implicit-tensor-decomposition}) returns unit vectors
$\widehat{\NormalizedRegressor}_1,\ldots,\widehat{\NormalizedRegressor}_{i_\star}$, where
$1\leq i_\star\leq |Q|$, such that, with probability at least
$1-\delta$, there is an
injection $\pi:[i_\star]\to Q$
satisfying
\begin{equation}
 \max_{\ell\in[i_\star]}
 \norm{\widehat{\NormalizedRegressor}_\ell-
       \frac{\Regressor_{\pi(\ell)}}{\norm{\Regressor_{\pi(\ell)}}_2}
  }_2
 \leq e^{-s/(2K)}\leq\frac{\Delta^2}{16B^2}.
\end{equation}
Its sample complexity is
\[
 \Ot\!\left(
 (nk)^{O(K^{10})}
 \log(\nicefrac1\delta)
 \right),
\]
and its running time is
\[
 (nk)^{\poly(K)}
 \log(\nicefrac1\delta).
\]
\end{theorem}

The rest of this section is dedicated to proving \Cref{thm:tensor-decomposition:implicit-tensor-decomposition}.

\subsection{Proof Roadmap}
\label{subsec:tensor-decomposition:roadmap}
Our implicit tensor decomposition algorithm consists of three key subroutines.
The first algorithm is an iterative projection algorithm that produces projections of high-order Hermite moments without having to explicitly form the tensors.
Given the projected moments, we construct a family of matrices that encode the regressors within its eigenvalues.
Finally, we decode the regressors using an eigendecomposition.

\paragraph{Iterative Projections (\Cref{subsec:tensor-decomposition:iterative-projections}).}
Given the recovered list $\widehat{\mathcal R}$, fresh samples,
a sampling accuracy, and a failure probability,
\textsc{EstimateProjectedMoments} returns a projected even-moment estimate
$\widehat A$ and projected odd moment estimates
$\widehat M_1,\ldots,\widehat M_n$.
We quantify what it means to estimate below.

For the remaining set $Q$, write
\begin{equation}\label{eq:population-moment-matrices}
\begin{aligned}
 \EvenMomentMatrix{r}{Q}
 &\coloneqq
 \sum_{i\in Q}\ComponentWeight{2r}{i}
 \TensorPower{\NormalizedRegressor_i}{r}
 (\TensorPower{\NormalizedRegressor_i}{r})^\top,
 \qquad 1\leq r\leq s,\\
 \OddMomentMatrix{s}{Q}{\bm{\zeta}}
 &\coloneqq
 \sum_{i\in Q}\ComponentWeight{2s}{i}
 \ip{\bm{\zeta}}{\ScaledRegressor_i}
 \TensorPower{\NormalizedRegressor_i}{s}
 (\TensorPower{\NormalizedRegressor_i}{s})^\top,
 \qquad \bm{\zeta}\in\R^n.
\end{aligned}
\end{equation}

\begin{definition}[Projected Moment Estimate]
\label{def:projected-moment-estimate}
Fix a non-empty $Q\subseteq[k]$ and the target moment matrices
$\EvenMomentMatrix{s}{Q}$ and $\OddMomentMatrix{s}{Q}{\cdot}$ from
\eqref{eq:population-moment-matrices}.
Let $ \ProjectionMap{s}:(\R^n)^{\otimes s}\to\R^{p_s}$
have orthonormal rows.
For $\varepsilon\geq0$, we call symmetric matrices
$\widehat{A},\widehat{M}_1,\ldots,\widehat{M}_n
 \in\R^{p_s\times p_s}$
an \emph{$\varepsilon$-estimate of the projected moment matrices for $Q$ relative to $ \ProjectionMap{s}$} if the following bounds hold:
\begin{enumerate}[label=(\roman*),ref=\thedefinition(\roman*),leftmargin=2.2em]
\item\label[definition]{def:projected-moment-estimate-even}
The even estimate satisfies
\begin{equation}
 \norm{\widehat{A}- \ProjectionMap{s} \EvenMomentMatrix{s}{Q} \ProjectionMap{s}^\top}_F
 \leq\varepsilon\tr(\EvenMomentMatrix{s}{Q}).
\end{equation}

\item\label[definition]{def:projected-moment-estimate-odd}
For every $\bm{\zeta}\in\R^n$, the odd estimates satisfy
\begin{equation}
 \norm{\sum_{\ell=1}^n \zeta_\ell\widehat{M}_\ell
       - \ProjectionMap{s} \OddMomentMatrix{s}{Q}{\bm{\zeta}} \ProjectionMap{s}^\top}_F
 \leq\varepsilon\tr(\EvenMomentMatrix{s}{Q})\norm{\bm{\zeta}}_2.
\end{equation}
\end{enumerate}
\end{definition}

We are now ready to state the guarantee of the first subroutine.
\begin{theorem}[Projected-Moment Accuracy]
\label{thm:tensor-decomposition:projected-moment-accuracy}
Consider the Gaussian self-selection model
(\Cref{def:gaussian-self-selection}) under the boundedness and
uncoveredness assumptions (\Cref{ass:boundedness-uncoveredness}), with $k\geq2$
and the hyperparameters from \Cref{hyp:algorithm}.
Fix a recovered list
$\widehat{\mathcal R}=((\widehat{\Regressor}_i,\widehat b_i))_{i\in R}$,
indexed by its matching to $R\subsetneq[k]$, such that
\[
 \norm{\widehat{\Regressor}_i-\Regressor_i}_2+|\widehat b_i-b_i|
 \leq e^{-K^6s}\qquad(i\in R),
\]
and put $Q=[k]\setminus R$.

For $\delta\in(0,\nicefrac13)$, run
\textnormal{\textsc{EstimateProjectedMoments}}
(\Cref{alg:tensor-decomposition:estimate-projected-moments}) with sampling accuracy $e^{-K^6s}$,
failure probability $\delta$, and fresh observations
independent of $\widehat{\mathcal R}$.
With probability at least $1-\delta$, the returned matrices $\widehat A, \widehat M_1, \dots, \widehat M_n$
form an $e^{-K^2s}$-estimate of the projected moment matrices for $Q$
relative to the constructed map $\ProjectionMap{s}$
(\Cref{def:projected-moment-estimate}), and
\[
 \norm{(I-\ProjectionMap{s}^\top\ProjectionMap{s})
       (\EvenMomentMatrix{s}{Q})^{1/2}}_{\op}^2
 \leq2K^{s-1}e^{-K^2s}
       \tr(\EvenMomentMatrix{s}{Q}).
\]

Its sample complexity is
\[
 \Ot\!\left((nk)^{O(K^{10})}
            \log(\nicefrac1\delta)\right),
\]
and its running time is
\[
 (nk)^{\poly(K)}\log(\nicefrac1\delta).
\]
\end{theorem}

In order for \Cref{thm:tensor-decomposition:projected-moment-accuracy} to compute these projected moments in polynomial time,
it must avoid explicitly materializing high-order tensors.
This is accomplished by iteratively building the projections degree by degree.

\paragraph{Observable-Operator Construction (\Cref{subsubsec:tensor-decomposition:operator-construction-proof}).}
Given $\widehat A$ and $\widehat M_1,\ldots,\widehat M_n$,
\textsc{ConstructObservableOperators} returns a family of matrices $\widehat D_1,\ldots,\widehat D_n$ that encode a subset of the unrecovered regressors.
We formalize the precise spectral condition below.

\begin{restatable}[Approximate Observable-Operator Family]{definition}{commoneigenbasisdefinition}
\label{def:approximate-observable-operator-family}
Fix a non-empty set $J$ and vectors $(\ScaledRegressor_i)_{i\in J}$ in $\R^n$.
For $\varepsilon_{\mathrm{fam}}\geq0$, we call symmetric matrices
$\widehat{D}_1,\ldots,\widehat{D}_n\in\R^{|J|\times|J|}$
an $\varepsilon_{\mathrm{fam}}$-approximate observable-operator family
for $(\ScaledRegressor_i)_{i\in J}$ if there is an orthonormal basis
$(\bm{z}_i)_{i\in J}$ of $\R^{|J|}$ such that, for every $\bm{\zeta}\in\R^n$,
\begin{equation}
 \norm{\sum_{\ell=1}^n \zeta_\ell\widehat{D}_\ell
   -\sum_{i\in J}\ip{\bm{\zeta}}{\ScaledRegressor_i}\bm{z}_i\bm{z}_i^\top}_F
 \leq\varepsilon_{\mathrm{fam}}\norm{\bm{\zeta}}_2.
\end{equation}
\end{restatable}

We now state the guarantees of the second subroutine.

\begin{restatable}[Observable-Operator Construction]{theorem}{selectnormalizelemma}
\label{thm:tensor-decomposition:observable-operator-construction}
Fix a non-empty $Q\subseteq[k]$.
Consider the rank-one sums $\EvenMomentMatrix{s}{Q},\OddMomentMatrix{s}{Q}{\bm{\zeta}}$ from
\eqref{eq:population-moment-matrices}.
Let $ \ProjectionMap{s}\in\R^{p_s\times n^s}$ have orthonormal rows for some
$1\leq p_s\leq|Q|$ and
$\varepsilon_{\mathrm{mom}},\varepsilon_{\mathrm{proj}},
\varepsilon_\perp\geq0$.
Suppose the following conditions hold.
\begin{enumerate}[label=(\alph*),ref=\thetheorem(\alph*)]

\item
For every $i\in Q$, $\norm{\ScaledRegressor_i}_2\leq K^3$.

\item\label[assumption]{ass:tensor-decomposition:approximate-orthogonality}
The lifted vectors satisfy
\begin{equation}
 \norm{
  \bigl(\ip{\TensorPower{\NormalizedRegressor_i}{s}}
            {\TensorPower{\NormalizedRegressor_j}{s}}\bigr)_{i,j\in Q}-I
 }_{\op}\leq\varepsilon_\perp.
\end{equation}

\item
The matrices $\widehat A,\widehat M_1,\ldots,\widehat M_n$
form an $\varepsilon_{\mathrm{mom}}$-estimate of the projected
moment matrices for $Q$ relative to $\ProjectionMap{s}$
(\Cref{def:projected-moment-estimate}).

\item\label[assumption]{ass:tensor-decomposition:projection-error}
The projection error satisfies
\begin{equation}
 \norm{(I-\ProjectionMap{s}^\top\ProjectionMap{s})
       (\EvenMomentMatrix{s}{Q})^{1/2}}_{\op}^2
 \leq\varepsilon_{\mathrm{proj}}\tr(\EvenMomentMatrix{s}{Q}).
\end{equation}

\item\label[assumption]{ass:tensor-decomposition:small-errors}
The errors satisfy $\varepsilon_{\mathrm{mom}}+\varepsilon_{\mathrm{proj}}
+\varepsilon_\perp\leq\frac{1}{16k^2}$.
\end{enumerate}
\sloppy
Then \textnormal{\textsc{ConstructObservableOperators}}
(\Cref{alg:tensor-decomposition:construct-observable-operators}), given these matrices, returns a
\mbox{$16K^3k^3(
 \varepsilon_{\mathrm{mom}}
 +\varepsilon_{\mathrm{proj}}
 +\varepsilon_\perp)$-approximate} observable-operator family
(\Cref{def:approximate-observable-operator-family}) for $(\ScaledRegressor_i)_{i\in J}$
indexed by some $J\subseteq Q$.
Moreover, the routine runs in $\poly(n,k)$ time.
\end{restatable}

\Cref{thm:tensor-decomposition:observable-operator-construction} is the main bottleneck which prevents our algorithm from recovering all regressors at once.
This is because the construction roughly requires the eigenvalues of $\widehat A$ to be on the same relative scale.
However, due to the irreducible bias, we have no control over the smallest eigenvalues and must truncate to the top eigenspace instead.

\paragraph{Recovering Normalized Regressors (\Cref{subsubsec:tensor-decomposition:direction-recovery-proof}).}
Given an approximate observable-operator family
$\widehat D_1,\ldots,\widehat D_n$ of dimension $i_\star\geq2$
and a failure probability, \textsc{RecoverDirections} returns
$i_\star$ unit vectors estimating the normalized regressors in
the block up to permutation.

\begin{restatable}[Normalized-Regressor Recovery]{theorem}{diagonalizereadlemma}
\label{thm:tensor-decomposition:normalized-regressor-recovery}
Let $J\subseteq[k]$ have $p\geq2$ elements, and suppose the vectors
$(\ScaledRegressor_i)_{i\in J}$ in $\R^n$ satisfy
\[
 K^{-3}\leq\norm{\ScaledRegressor_i}_2\leq K^3\quad(i\in J),
 \qquad \norm{\ScaledRegressor_i-\ScaledRegressor_{i'}}_2\geq K^{-3}\quad(i\ne i').
\]
Fix $\varepsilon_{\mathrm{fam}}\geq0$.  Suppose the symmetric matrices
$\widehat{D}_1,\ldots,\widehat{D}_n\in\R^{p\times p}$
form an $\varepsilon_{\mathrm{fam}}$-approximate observable-operator family
for $(\ScaledRegressor_i)_{i\in J}$
(\Cref{def:approximate-observable-operator-family}).
Given these matrices and $\delta\in(0,\nicefrac13)$,
\textnormal{\textsc{RecoverDirections}}
(\cref{alg:tensor-decomposition:recover-directions})
returns unit vectors $\widehat{\NormalizedRegressor}_1,\ldots,\widehat{\NormalizedRegressor}_p$ such
that, with probability at least $1-\delta$, there is a bijection
$\pi:[p]\to J$ satisfying
\begin{equation}
 \max_{i\in[p]}
 \norm{\widehat{\NormalizedRegressor}_i-\frac{\ScaledRegressor_{\pi(i)}}{\norm{\ScaledRegressor_{\pi(i)}}_2}}_2
 \leq32K^6\sqrt n\,p^2\varepsilon_{\mathrm{fam}}.
\end{equation}
Moreover, the routine runs in time
$\poly(n,k)\log(\nicefrac1\delta)$.
\end{restatable}

Given an observable-operator family, the natural approach is to perform an eigendecomposition.
Assuming a simple spectrum, this allows us to extract the eigenvectors and then the eigenvalues.
We will use random contractions to ensure a sufficient eigenvalue gap.

\paragraph{Putting it Together (\Cref{subsec:tensor-decomposition:decomposition-proof}).}
Once we have analyzed the subroutines, it remains only to put them together and verify that their assumptions hold.
The rest of this section presents the proofs of the results above.

\subsection{Iterative Projections}
\label{subsec:tensor-decomposition:iterative-projections}

Let $Q\subseteq [k]$ index the set of unrecovered regressors.
The first step in the roadmap (\Cref{subsec:tensor-decomposition:roadmap}) to analyzing our implicit tensor decomposition algorithm is an iterative projection procedure that computes estimates of the projections of restricted even and odd Hermite moments $A_Q^{(s)}, M_Q^{(s)}(b)$ (\Cref{thm:tensor-decomposition:projected-moment-accuracy}).
This allows us to avoid explicitly materializing these high-order tensors.
In this section, we present an iterative projection procedure inspired by \citet{LL24,DK25}.
Starting from the trivial degree-$1$ projection map $C_1 = I_n$,
we will show how to implicitly compute the degree-$(r+1)$ projection map $C_{r+1}$ as well as the projected moment $\ProjectionMap{r+1}\EvenMomentMatrix{r+1}{Q}\ProjectionMap{r+1}^\top$,
given the degree-$r$ projection map and projected moment.
We will store an implicit representation of $C_s$ and the projected odd moments can then be obtained given this representation.

We first describe the implicit representation in \Cref{subsubsec:tensor-decomposition:projection-maps}.
Our algorithm is easiest to state assuming a moment estimation oracle \textsc{EstimateMomentsByPeeling} (\Cref{alg:peeling:estimate-moments-by-peeling}) that takes as input projected samples $(\widetilde{\Covariate}, C_r\widetilde\Covariate^{\otimes r}, Y_{\max})$ and estimates the expectation of the projected moment.
In \Cref{subsubsec:tensor-decomposition:moment-estimation},
we describe the specifics of the oracle needed for this section and defer further details and its proof to \Cref{sec:peeling:main}.
As our projection map is computed iteratively,
we analyze the cumulative error in \Cref{subsubsec:tensor-decomposition:cumulative-projection-error}.
In \Cref{subsubsec:tensor-decomposition:projected-moment-proof}, we combine all the ingredients to prove the guarantees of our iterative projection procedure (\Cref{thm:tensor-decomposition:projected-moment-accuracy}).

\subsubsection{Representation of Projection Maps}\label{subsubsec:tensor-decomposition:projection-maps}
Fix a recovered list $\widehat{\mathcal R}$ matched bijectively
to $R\subsetneq[k]$, and put $Q=[k]\setminus R$.
Recall that $\NormalizedRegressor_i=\Regressor_i/\norm{\Regressor_i}_2$, and
$\EvenMomentMatrix{r}{Q}=\sum_{i\in Q}\ComponentWeight{2r}{i}
\TensorPower{\NormalizedRegressor_i}{r}(\TensorPower{\NormalizedRegressor_i}{r})^\top$.
We would like to construct a projection map $\ProjectionMap{s}\in\R^{p_s\times n^s}$ so that $\Range(\ProjectionMap{s}^\top)$ approximately contains the lifted vectors $\TensorPower{\NormalizedRegressor_i}{s}$ of components whose weights $\ComponentWeight{2s}{i}$ are sufficiently large relative to $\tr(\EvenMomentMatrix{s}{Q})$.
We build $\ProjectionMap{s}$ one tensor order at a time, starting with $\ProjectionMap{1}=I_n$ and $p_1=n$.
For $2\leq r\leq s$, put
$p_r=\min\{|Q|,p_{r-1}n\}$.
Since $\TensorPower{\NormalizedRegressor_i}{r}=\TensorPower{\NormalizedRegressor_i}{(r-1)}\otimes \NormalizedRegressor_i$,
\[
 \TensorPower{\NormalizedRegressor_i}{(r-1)}\in\Range(\ProjectionMap{r-1}^\top)
 \quad\Longrightarrow\quad
 \TensorPower{\NormalizedRegressor_i}{r}\in\Range(\ProjectionMap{r-1}^\top\otimes I_n).
\]
For $r\geq 2$, assume that we have used \textsc{EstimateMomentsByPeeling} (\Cref{alg:peeling:estimate-moments-by-peeling}) to produce an intermediate estimate $\widehat A_r$ of the projected target moment matrix
\begin{equation}\label{eq:tensor-decomposition:projected-even-target}
 (\ProjectionMap{r-1}\otimes I_n)
 \EvenMomentMatrix{r}{Q}
 (\ProjectionMap{r-1}\otimes I_n)^\top.
\end{equation}
Note that the case of $C_1=I_n$ is trivial.
This target matrix has rank at most $|Q|$.
Let the rows of $V_r$ be the $p_r$ leading
orthonormal eigenvectors of $\widehat A_r$, and set
\begin{equation}\label{eq:tensor-decomposition:projection-recursion}
 \ProjectionMap{r}\coloneqq V_r(\ProjectionMap{r-1}\otimes I_n).
\end{equation}
Our algorithm can proceed to the next iteration.
Note that we only need to store $V_2,\ldots,V_s$.

Remark that we only used the even moments to determine the projection map.
For the odd moment $M_Q^{(s)}(\bm{\zeta})$, it calls \textsc{EstimateMomentsByPeeling}
(\Cref{alg:peeling:estimate-moments-by-peeling}) once again to estimate all coordinate
contractions
\begin{equation}\label{eq:tensor-decomposition:projected-odd-target}
\ProjectionMap{s}\OddMomentMatrix{s}{Q}{\bm{e}_\ell}\ProjectionMap{s}^\top
\end{equation}
for $\ell\in[n]$.
Since the odd moment is linear in the contraction direction,
the matrices $\widehat M_1,\ldots,\widehat M_n$ determine the family
\begin{equation}
 \widehat M(\bm{\zeta})\coloneqq\sum_{\ell=1}^n \zeta_\ell\widehat M_\ell.
\end{equation}

\subsubsection{Mean Estimation Oracle}\label{subsubsec:tensor-decomposition:moment-estimation}
We now describe how to obtain estimates of \eqref{eq:tensor-decomposition:projected-even-target} and \eqref{eq:tensor-decomposition:projected-odd-target}
using \textsc{EstimateMomentsByPeeling} (\Cref{alg:peeling:estimate-moments-by-peeling}).
To form $\widehat A_{r+1}$, we use the factors $V_2,\ldots,V_r$.
Once $V_s$ is available, the same recursion gives the final odd features.
For this purpose,
we use the fact that the flattened Hermite moments can be written as expectations of outer products of the complexified covariates.
For each observation, draw an independent $\Covariate'\sim N(0,I_n)$
and put $\widetilde{\Covariate}=\Covariate+\iu \Covariate'$.
By the representation of the projection map $C_r$ \eqref{eq:tensor-decomposition:projection-recursion},
the projected powers
satisfy the implicit tensor recursion
\begin{equation}
\begin{aligned}
 \ProjectionMap{1}\widetilde{\Covariate}
 &=\widetilde{\Covariate},\\
 \ProjectionMap{r+1}\TensorPower{\widetilde{\Covariate}}{(r+1)}
 &=V_{r+1}\left[
  \bigl(\ProjectionMap{r}\TensorPower{\widetilde{\Covariate}}{r}\bigr)
  \otimes\widetilde{\Covariate}
 \right].
\end{aligned}
\end{equation}
Thus given $V_2,\ldots,V_r$, we can evaluate $C_r \widetilde \Covariate^{\otimes r}$
without forming $\TensorPower{\widetilde{\Covariate}}{r}$ or $\ProjectionMap{r}$.
The recursion has at most $s$ levels and uses tensor products
of vectors and matrix-vector products in dimension at most $n^2$.
Thus forming each projected sample takes time polynomial in $s$ and $n$.

For each fresh observation, the even call at order $r+1$ uses
$(\ProjectionMap{r}\TensorPower{\widetilde{\Covariate}}{r})
 \otimes\widetilde{\Covariate}$, and the final odd call uses
$\ProjectionMap{s}\TensorPower{\widetilde{\Covariate}}{s}$.
For $1\leq r<s$, the following even feature is used to
estimate $\widehat A_{r+1}$, which determines $V_{r+1}$.
\begin{align}
 \Psi_{r+1}
 &\coloneqq\frac{\RealPart[
 \bigl((\ProjectionMap{r}\TensorPower{\widetilde{\Covariate}}{r})
       \otimes\widetilde{\Covariate}\bigr)
 \bigl((\ProjectionMap{r}\TensorPower{\widetilde{\Covariate}}{r})
       \otimes\widetilde{\Covariate}\bigr)^\top]}
      {\sqrt{(2r+2)!}},\\
 \E_{\Covariate'}[\Psi_{r+1}\mid \Covariate]
 &=(\ProjectionMap{r}\otimes I_n)
 \Mat_{r+1}(\mathsf H_{2r+2}(\Covariate))
 (\ProjectionMap{r}\otimes I_n)^\top,
\\
 \E_{\Covariate,\NoiseVector_Q,\Covariate'}
 [\Psi_{r+1}\TailIndicator{Y_{\max,Q}}]
 &\approx(\ProjectionMap{r}\otimes I_n)\EvenMomentMatrix{r+1}{Q}
         (\ProjectionMap{r}\otimes I_n)^\top.
\label{eq:tensor-decomposition:even-moment-estimation-target}
\end{align}
At the final level, the odd features used to estimate
$(\widehat M_\ell)_{\ell\in[n]}$ are
\begin{align}
 \Psi_{\odd,\ell}
 &\coloneqq\frac{\RealPart[
 (\ProjectionMap{s}\TensorPower{\widetilde{\Covariate}}{s})
 (\ProjectionMap{s}\TensorPower{\widetilde{\Covariate}}{s})^\top
 \widetilde X_\ell]}
      {\sqrt{(2s+1)!}},\qquad \ell\in[n],
 \\
 \E_{\Covariate'}[\Psi_{\odd,\ell}\mid \Covariate]
 &=\ProjectionMap{s}
 \Mat_s(\HermiteContraction{2s+1}{\Covariate}{\bm{e}_\ell})
 \ProjectionMap{s}^\top,
\\
 \E_{\Covariate,\NoiseVector_Q,\Covariate'}
 [\Psi_{\odd,\ell}\TailIndicator{Y_{\max,Q}}]
 &\approx\ProjectionMap{s}\OddMomentMatrix{s}{Q}{\bm e_\ell}
         \ProjectionMap{s}^\top.
\label{eq:tensor-decomposition:odd-moment-estimation-target}
\end{align}
We call \textsc{EstimateMomentsByPeeling}
(\Cref{alg:peeling:estimate-moments-by-peeling}) to estimate \eqref{eq:tensor-decomposition:even-moment-estimation-target} for $1\leq r<s$ and \eqref{eq:tensor-decomposition:odd-moment-estimation-target}.
Since $Y_{\max,Q}$ is not observed, the subroutine uses a weight function
$\widehat\omega_{\widehat{\mathcal R}}(\Covariate,Y_{\max})$
as a proxy of the indicator. 
We give its estimation guarantee in
\Cref{sec:peeling:main}.

\subsubsection{Cumulative Projection Error}\label{subsubsec:tensor-decomposition:cumulative-projection-error}
We now present the pseudocode of our iterative projection algorithm (\Cref{alg:tensor-decomposition:estimate-projected-moments}) and bound the cumulative projection error in terms of the accuracy of the intermediate moment estimates.

\begin{lemma}[Recursive Projection Error]
\label{lem:tensor-decomposition:recursive-projection-error}
Consider the Gaussian self-selection model
(\Cref{def:gaussian-self-selection}) under the boundedness and
uncoveredness assumptions (\Cref{ass:boundedness-uncoveredness}).
Use the hyperparameters $K,T,s$ from \Cref{hyp:algorithm}.
Suppose $Q=[k]\setminus R\ne\varnothing$ indexes the unrecovered regressors.
Form $ \ProjectionMap{1},\ldots, \ProjectionMap{s}$ by the intermediate
eigenvector steps of \textnormal{\textsc{EstimateProjectedMoments}}
(\cref{alg:tensor-decomposition:estimate-projected-moments}).

For $2\leq r\leq s$, let $\gamma_r$ be the Frobenius error
\[
  \gamma_r \coloneqq \norm{\widehat{A}_r
  -( \ProjectionMap{r-1}\otimes I_n)
   \EvenMomentMatrix{r}{Q}( \ProjectionMap{r-1}\otimes I_n)^\top}_F\,.
\]
For $1\leq t\leq s$, put
$\Pi_t=(\ProjectionMap{t}^\top\ProjectionMap{t})\otimes I_{n^{s-t}}$.
\begin{enumerate}[label=(\roman*),ref=\thelemma(\roman*),leftmargin=2.2em]
\item\label[lemma]{lem:tensor-decomposition:projection-error-levelwise}
For every $1\leq t\leq s$,
\begin{equation}
 \norm{(I-\Pi_t)(\EvenMomentMatrix{s}{Q})^{1/2}}_{\op}^2
 \leq2\sum_{r=2}^t\gamma_r
       \max_{i\in Q}\frac{\ComponentWeight{2s}{i}}{\ComponentWeight{2r}{i}}.
\end{equation}

\item\label[lemma]{lem:tensor-decomposition:projection-error-relative}
Fix $\varepsilon_{\mathrm{mean}}>0$. If
$\gamma_r\leq\varepsilon_{\mathrm{mean}}\tr(\EvenMomentMatrix{r}{Q})$
for every $2\leq r\leq s$, then, for every $1\leq t\leq s$,
\begin{equation}
 \norm{(I-\Pi_t)(\EvenMomentMatrix{s}{Q})^{1/2}}_{\op}^2
 \leq2\varepsilon_{\mathrm{mean}}K^{s-1}\tr(\EvenMomentMatrix{s}{Q}).
\end{equation}
\end{enumerate}
\end{lemma}

\begin{algorithm}[H]
\caption[\textsc{EstimateProjectedMoments}]{\textsc{EstimateProjectedMoments}}
\label{alg:tensor-decomposition:estimate-projected-moments}
\KwIn{The recovered list $\widehat{\mathcal R}$;
sampling accuracy $\varepsilon_{\mathrm{est}}\in(0,1)$;
failure probability $\delta\in(0,\nicefrac13)$}
\KwOut{Symmetric $p_s\times p_s$ matrices
$\widehat A,\widehat M_1,\ldots,\widehat M_n$}
$p_1\gets n$\;
\For{$r=1,\ldots,s-1$}{
  \tcp{\Cref{alg:peeling:estimate-moments-by-peeling,thm:peeling:hermite-moment-estimation-even}}
  $\widehat A_{r+1}\gets$
  \PeeledMomentMean{$\widehat{\mathcal R}, \{V_2, \dots, V_r\}, \Psi_{r+1},
  \varepsilon_{\mathrm{est}},\delta/s$}\;
  $p_{r+1}\gets
  \min\{k-|\widehat{\mathcal R}|,p_rn\}$\;
  Compute the eigenvalues
  $\nu_1\geq\cdots\geq\nu_{p_rn}$ of $\widehat A_{r+1}$
  and corresponding orthonormal eigenvectors
  $\bm{v}_1,\ldots,\bm{v}_{p_rn}$ using a fixed deterministic convention\;
  $V_{r+1}\gets
  \begin{bmatrix}
   \bm{v}_1^\top\\
   \vdots\\
   \bm{v}_{p_{r+1}}^\top
  \end{bmatrix}$\;
}
$\widehat A\gets V_s\widehat A_sV_s^\top$\;
\tcp{\Cref{alg:peeling:estimate-moments-by-peeling,thm:peeling:hermite-moment-estimation-odd}}
$(\widehat M_1,\ldots,\widehat M_n)\gets$
\PeeledMomentMean{$\widehat{\mathcal R}, \{V_2, \dots, V_s\},
(\Psi_{\odd,\ell})_{\ell\in[n]},
\varepsilon_{\mathrm{est}},\delta/s$}\;
\Return $(\widehat A,\widehat M_1,\ldots,\widehat M_n)$\;
\end{algorithm}

We are now ready to prove \Cref{lem:tensor-decomposition:recursive-projection-error}.

\begin{proof}
To prove \Cref{lem:tensor-decomposition:projection-error-levelwise}, we bound the error at
each projection, transfer the bounds to order $s$, and add them
using the nested projection maps. For \Cref{lem:tensor-decomposition:projection-error-relative},
we compare the component weights at orders $r$ and $s$.

\emph{Error at one projection.}
Fix $2\leq r\leq s$.  We compare the leading eigenspace of
$\widehat A_r$ with the projected target moment matrix
\[
 L_r\coloneqq(\ProjectionMap{r-1}\otimes I_n)
 \EvenMomentMatrix{r}{Q}(\ProjectionMap{r-1}\otimes I_n)^\top.
\]
Let $P_r=V_r^\top V_r$, where the rows of $V_r$ are the transposes of the $p_r$
leading orthonormal eigenvectors of $\widehat A_r$.
Since the component weights are positive, $L_r$ is positive
semidefinite with rank at most $|Q|$.
By the dimension rule in \Cref{alg:tensor-decomposition:estimate-projected-moments},
$p_r=\min\{|Q|,p_{r-1}n\}$.
If $P_r=I$, the following one-level bound is immediate.
Otherwise $p_r=|Q|$ and $\lambda_{p_r+1}(L_r)=0$.
For $\bm{x}\in\Range(I-P_r)$, Weyl's eigenvalue perturbation
inequality (\Cref{prop:weyl-perturbation}) gives
\begin{align}
 \bm{x}^\top L_r\bm{x}
 &\leq \bm{x}^\top\widehat A_r\bm{x}+\gamma_r\norm{\bm{x}}_2^2
 \alignnote{$\norm{\widehat A_r-L_r}_{\op}\leq\gamma_r$}\notag\\
 &\leq(\lambda_{p_r+1}(\widehat A_r)+\gamma_r)\norm{\bm{x}}_2^2
 \alignnote{$\bm{x}\in\Range(I-P_r)$}\notag\\
 &\leq2\gamma_r\norm{\bm{x}}_2^2
 \alignnote{by \Cref{prop:weyl-perturbation}}\notag\\
 0\preceq(I-P_r)L_r(I-P_r)
 &\preceq2\gamma_r(I-P_r).
 \alignnote{all $\bm{x}\in\Range(I-P_r)$}\label{eq:tensor-decomposition:discarded-covariance}
\end{align}
The one-level bound \eqref{eq:tensor-decomposition:discarded-covariance} applies for
$2\leq r\leq s$.
We still need to compare these bounds at order $s$.

\emph{Projection maps at order $s$.}
We now express the successive projection maps at order $s$.
Their nested ranges will let us telescope the projection maps.
Since $\ProjectionMap{r-1}$ has orthonormal rows,
\begin{align}
 (\ProjectionMap{r-1}\otimes I_n)
 (\ProjectionMap{r-1}\otimes I_n)^\top
 &=(\ProjectionMap{r-1}\ProjectionMap{r-1}^\top)\otimes I_n\notag\\
 &=I_{p_{r-1}n}.
 \label{eq:tensor-decomposition:projection-coordinate-isometry}
\end{align}
The recursion $\ProjectionMap{r}=V_r(\ProjectionMap{r-1}\otimes I_n)$
and $P_r=V_r^\top V_r$ therefore give
\begin{align}
 \Pi_{r-1}
 &=(\ProjectionMap{r-1}^\top\ProjectionMap{r-1})
   \otimes I_{n^{s-r+1}}\notag\\
 &=\bigl[(\ProjectionMap{r-1}\otimes I_n)^\top
   (\ProjectionMap{r-1}\otimes I_n)\bigr]
   \otimes I_{n^{s-r}}\notag\\
 \Pi_r
 &=(\ProjectionMap{r}^\top\ProjectionMap{r})
   \otimes I_{n^{s-r}}\notag\\
 &=\bigl[(\ProjectionMap{r-1}\otimes I_n)^\top P_r
   (\ProjectionMap{r-1}\otimes I_n)\bigr]
   \otimes I_{n^{s-r}}\notag\\
 \Pi_{r-1}\Pi_r
 &=\bigl[(\ProjectionMap{r-1}\otimes I_n)^\top I_{p_{r-1}n}P_r
   (\ProjectionMap{r-1}\otimes I_n)\bigr]
   \otimes I_{n^{s-r}}
 \alignnote{by \eqref{eq:tensor-decomposition:projection-coordinate-isometry}}\notag\\
 &=\Pi_r
 \label{eq:tensor-decomposition:projection-nesting}\\
 D_r\coloneqq\Pi_{r-1}-\Pi_r
 &=\bigl[(\ProjectionMap{r-1}\otimes I_n)^\top
       (I_{p_{r-1}n}-P_r)
       (\ProjectionMap{r-1}\otimes I_n)\bigr]
       \otimes I_{n^{s-r}}.
 \label{eq:tensor-decomposition:projection-level-projector}
\end{align}
Since the $\Pi_r$ are orthogonal projectors, \eqref{eq:tensor-decomposition:projection-nesting}
gives, for $2\leq r\leq s$ and $1\leq t\leq s$,
\begin{align}
 D_r^2
 &=(\Pi_{r-1}-\Pi_r)^2\notag\\
 &=\Pi_{r-1}-2\Pi_r+\Pi_r=D_r
 \alignnote{by \eqref{eq:tensor-decomposition:projection-nesting}}\notag\\
 I-\Pi_t
 &=\sum_{r=2}^tD_r.
 \alignnote{$\Pi_1=I$; telescoping}
 \label{eq:tensor-decomposition:projection-telescoping}
\end{align}
Thus each $D_r$ is an orthogonal projector, with
$\norm{D_r}_{\op}\leq1$.
We next bound $D_r\EvenMomentMatrix{s}{Q}D_r$.

\emph{Error at order $s$.}
We first express \eqref{eq:tensor-decomposition:discarded-covariance} at order $s$.
By \eqref{eq:tensor-decomposition:projection-level-projector} and the definition
of $L_r$,
\begin{align}
 D_r(\EvenMomentMatrix{r}{Q}\otimes I_{n^{s-r}})D_r
 &=\bigl[
   (\ProjectionMap{r-1}\otimes I_n)^\top
   (I-P_r)L_r(I-P_r)
   (\ProjectionMap{r-1}\otimes I_n)
   \bigr]\otimes I_{n^{s-r}}\notag\\
 &\preceq2\gamma_rD_r.
 \alignnote{by \eqref{eq:tensor-decomposition:discarded-covariance}}
 \label{eq:tensor-decomposition:projection-lifted-covariance}
\end{align}
Since the $\NormalizedRegressor_i$ are unit vectors,
$\TensorPower{(\NormalizedRegressor_i \NormalizedRegressor_i^\top)}{(s-r)}\preceq I_{n^{s-r}}$.
With zeroth tensor powers interpreted as $1$, we obtain
\begin{align*}
 \EvenMomentMatrix{s}{Q}
 &=\sum_{i\in Q}\ComponentWeight{2s}{i}
   \TensorPower{\NormalizedRegressor_i}{r}(\TensorPower{\NormalizedRegressor_i}{r})^\top
   \otimes\TensorPower{(\NormalizedRegressor_i \NormalizedRegressor_i^\top)}{(s-r)}\\
 &\preceq\sum_{i\in Q}\ComponentWeight{2s}{i}
   \TensorPower{\NormalizedRegressor_i}{r}(\TensorPower{\NormalizedRegressor_i}{r})^\top\otimes I_{n^{s-r}}
 \alignnote{unit $\NormalizedRegressor_i$}\\
 &\preceq
   \left(\max_{i\in Q}
   \frac{\ComponentWeight{2s}{i}}{\ComponentWeight{2r}{i}}\right)
   (\EvenMomentMatrix{r}{Q}\otimes I_{n^{s-r}}).
 \alignnote{positive component weights}
\end{align*}
Conjugating by $D_r$, taking operator norms, and using
$\norm{D_r}_{\op}\leq1$ gives
\begin{align}
 \norm{D_r\EvenMomentMatrix{s}{Q}D_r}_{\op}
 &\leq
   \left(\max_{i\in Q}
   \frac{\ComponentWeight{2s}{i}}{\ComponentWeight{2r}{i}}\right)
   \norm{D_r(\EvenMomentMatrix{r}{Q}\otimes I_{n^{s-r}})D_r}_{\op}
 \alignnote{positive semidefinite order}\notag\\
 &\leq2\gamma_r
   \max_{i\in Q}\frac{\ComponentWeight{2s}{i}}{\ComponentWeight{2r}{i}}.
 \alignnote{by \eqref{eq:tensor-decomposition:projection-lifted-covariance}}
 \label{eq:tensor-decomposition:projection-one-level-loss}
\end{align}
It remains to add these bounds over the levels.

\emph{Accumulating projection error.}
We now combine \eqref{eq:tensor-decomposition:projection-one-level-loss}
using the telescoping identity
\eqref{eq:tensor-decomposition:projection-telescoping}.
The Gram-matrix norm identity (\Cref{prop:gram-matrix-identities}) gives
\begin{align*}
 \norm{(I-\Pi_t)(\EvenMomentMatrix{s}{Q})^{1/2}}_{\op}^2
 &=\norm{(\EvenMomentMatrix{s}{Q})^{1/2}(I-\Pi_t)
         (\EvenMomentMatrix{s}{Q})^{1/2}}_{\op}
 \alignnote{by \Cref{prop:gram-matrix-identities}}\\
 &=\norm{\sum_{r=2}^t
   (\EvenMomentMatrix{s}{Q})^{1/2}D_r(\EvenMomentMatrix{s}{Q})^{1/2}}_{\op}
 \alignnote{by \eqref{eq:tensor-decomposition:projection-telescoping}}\\
 &\leq\sum_{r=2}^t
   \norm{(\EvenMomentMatrix{s}{Q})^{1/2}D_r
         (\EvenMomentMatrix{s}{Q})^{1/2}}_{\op}
 \alignnote{triangle inequality}\\
 &=\sum_{r=2}^t\norm{D_r\EvenMomentMatrix{s}{Q}D_r}_{\op}
 \alignnote{by \Cref{prop:gram-matrix-identities}}\\
 &\leq2\sum_{r=2}^t\gamma_r
   \max_{i\in Q}\frac{\ComponentWeight{2s}{i}}{\ComponentWeight{2r}{i}}.
 \alignnote{by \eqref{eq:tensor-decomposition:projection-one-level-loss}}
\end{align*}
This proves \Cref{lem:tensor-decomposition:projection-error-levelwise}.
To prove \Cref{lem:tensor-decomposition:projection-error-relative}, we compare the trace
at order $r$ with the trace at order $s$. The component-weight
comparison (\Cref{thm:hermite:structure-component-weights}) gives, for $i,j\in Q$,
$\ComponentWeight{2s}{i}\ComponentWeight{2r}{j}
 \leq K^{s-r}\ComponentWeight{2r}{i}\ComponentWeight{2s}{j}$.
Summing over $j$ and using the positive weights yields
\begin{equation}\label{eq:tensor-decomposition:projection-weight-comparison}
 \max_{i\in Q}\frac{\ComponentWeight{2s}{i}}{\ComponentWeight{2r}{i}}
 \leq K^{s-r}
       \frac{\tr(\EvenMomentMatrix{s}{Q})}{\tr(\EvenMomentMatrix{r}{Q})}.
\end{equation}
Substituting the assumed bounds on $\gamma_r$ and
\eqref{eq:tensor-decomposition:projection-weight-comparison} in \Cref{lem:tensor-decomposition:projection-error-levelwise} gives
\begin{align*}
 \norm{(I-\Pi_t)(\EvenMomentMatrix{s}{Q})^{1/2}}_{\op}^2
 &\leq2\varepsilon_{\mathrm{mean}}\tr(\EvenMomentMatrix{s}{Q})
       \sum_{r=2}^t K^{s-r}
 \alignnote{by \eqref{eq:tensor-decomposition:projection-weight-comparison}}\\
 &\leq2\varepsilon_{\mathrm{mean}}K^{s-1}\tr(\EvenMomentMatrix{s}{Q}).
 \alignnote{geometric sum, $K\geq2$}
\end{align*}
This concludes the proof.
\end{proof}

\subsubsection{Proof of \texorpdfstring{\Cref*{thm:tensor-decomposition:projected-moment-accuracy}}{Theorem 6.3}}
\label{subsubsec:tensor-decomposition:projected-moment-proof}
We are now ready to prove the guarantees of our iterative projection algorithm (\Cref{alg:tensor-decomposition:estimate-projected-moments}).
The proof calls on the estimation guarantees for \textsc{EstimateMomentsByPeeling}
(\Cref{alg:peeling:estimate-moments-by-peeling}) that are presented in \Cref{sec:peeling:main} as well as \Cref{lem:tensor-decomposition:recursive-projection-error}.

\begin{proof}[Proof of \Cref{thm:tensor-decomposition:projected-moment-accuracy}]
For the estimation calls to \textsc{EstimateMomentsByPeeling}
(\Cref{alg:peeling:estimate-moments-by-peeling}) for even moments at order $2\leq r\leq s$, put
$U=\ProjectionMap{r-1}\otimes I_n$. 
The even moment estimation guarantees in \Cref{thm:peeling:hermite-moment-estimation-even} and the
restricted-moment decomposition in \Cref{thm:hermite:structure-moment-decomposition} give
\begin{align*}
   \gamma_r 
   &\coloneqq \norm{\widehat{A}_r
  -( \ProjectionMap{r-1}\otimes I_n)
   \EvenMomentMatrix{r}{Q}( \ProjectionMap{r-1}\otimes I_n)^\top}_F \\
 &\leq
 \underbrace{\norm{\widehat A_r-U\Mat_r(\RestrictedMomentTensor{2r}{Q})
                       U^\top}_F}_{\text{peeling estimate}}
 +\underbrace{\norm{UE_rU^\top}_F}_{\text{restricted-moment remainder}}
 \alignnote{by \Cref{thm:hermite:structure-moment-decomposition}}\\[6pt]
 &\leq
 \underbrace{\left(1+\sqrt{k/\pi}\right)e^{-K^6s}
       +kK^{4s}e^{-T^2/(2K)}\tr(\EvenMomentMatrix{r}{Q})}
       _{\text{by \Cref{thm:peeling:hermite-moment-estimation-even}}}\\[6pt]
 &\qquad+
 \underbrace{K^{8s}e^{-T^2/(2K)}\tr(\EvenMomentMatrix{r}{Q})}
       _{\text{by \Cref{thm:hermite:structure-moment-decomposition} and $\norm{U}_{\op}=1$}}.
\end{align*}
By \Cref{thm:hermite:structure-component-weights}, $1\leq e^{T^2}\tr(\EvenMomentMatrix{r}{Q})$.
Applying this bound to the additive term and simplifying the bound on $\gamma_r$ using the choice of hyperparameters (\Cref{prop:app-tensor-decomposition:iterative-projection-hyperparameter-bounds}) gives
\[
 \gamma_r\leq e^{-K^2s}\tr(\EvenMomentMatrix{r}{Q}).
\]
Applying the odd moment estimation guarantees from \Cref{thm:peeling:hermite-moment-estimation-odd} to the final call to \textsc{EstimateMomentsByPeeling}
(\Cref{alg:peeling:estimate-moments-by-peeling}) for the odd moment and repeating the
calculation at $r=s$ gives, simultaneously for every $\bm{\zeta}\in\R^n$,
\[
 \norm{\widehat M(\bm{\zeta})-\ProjectionMap{s}\OddMomentMatrix{s}{Q}{\bm{\zeta}}
                    \ProjectionMap{s}^\top}_F
 \leq e^{-K^2s}\tr(\EvenMomentMatrix{s}{Q})\norm{\bm{\zeta}}_2.
\]

Each call draws independent samples and satisfies its bound with probability
at least $1-\nicefrac\delta{s}$. 
By a union bound over
the $s-1$ even calls and the one odd call, all bounds hold with probability
at least $1-\delta$. 
We proceed conditioning on the successful estimation.

For the final even estimate, take $U=\ProjectionMap{s-1}\otimes I_n$.
Then
\begin{align*}
 \norm{\widehat A-\ProjectionMap{s}\EvenMomentMatrix{s}{Q}
                    \ProjectionMap{s}^\top}_F
 &=\norm{V_s\widehat A_sV_s^\top
        -V_sU\EvenMomentMatrix{s}{Q}U^\top V_s^\top}_F
 \alignnote{$\widehat A=V_s\widehat A_sV_s^\top$, $\ProjectionMap{s}=V_sU$}\\
 &=\norm{V_s\bigl(\widehat A_s-U\EvenMomentMatrix{s}{Q}U^\top\bigr)
         V_s^\top}_F\\
 &\leq\norm{\widehat A_s-U\EvenMomentMatrix{s}{Q}U^\top}_F
 \alignnote{$\norm{V_s}_{\op}=1$}\\
 &\leq e^{-K^2s}\tr(\EvenMomentMatrix{s}{Q}).
 \alignnote{even-call bound at $r=s$}
\end{align*}
Together with the odd-moment bound for every $\bm{\zeta}\in\R^n$,
this shows that the returned matrices $\widehat A,\widehat M_1,\ldots,\widehat M_n$
form an $e^{-K^2s}$-estimate of the projected moment matrices for $Q$
relative to the constructed map $\ProjectionMap{s}$
(\Cref{def:projected-moment-estimate}).
Moreover, by \Cref{lem:tensor-decomposition:projection-error-relative}, it follows that
\[
 \norm{(I-\ProjectionMap{s}^\top\ProjectionMap{s})
       (\EvenMomentMatrix{s}{Q})^{1/2}}_{\op}^2
 \leq 2K^{s-1}e^{-K^2s}\tr(\EvenMomentMatrix{s}{Q}).
\]

\emph{Sample complexity and running time.}
By \Cref{thm:peeling:hermite-moment-estimation}, calling \textsc{EstimateMomentsByPeeling}
(\Cref{alg:peeling:estimate-moments-by-peeling}) with $p$ projected coordinates and $m$
matrices uses
\[
 \Ot\!\left(p^2m\,e^{2K^6s}100^s\log\!\left(\frac{s}{\delta}\right)\right)
\]
projected samples and takes time at most the sample count times
$\poly(n,k,p,m)$. The first even matrix has dimension $n^2$,
and each later even matrix has dimension at most $nk$. The odd call estimates $n$ matrices of dimension at most
$k$. 
After simplifying by the choice of hyperparameters (\Cref{prop:app-tensor-decomposition:iterative-projection-hyperparameter-bounds}),
this amounts to a total sample complexity of
\begin{equation}\label{eq:tensor-decomposition:projection-explicit-samples}
 \Ot\!\left((nk)^{4K^{10}+12K^4}\log(\nicefrac1\delta)\right).
\end{equation}
Including the construction of the projected samples and the
eigendecompositions gives total running time
$(nk)^{\poly(K)}\log(\nicefrac1\delta)$.
\end{proof}

\subsection{Observable-Operator Learning}
\label{subsec:tensor-decomposition:observable-operator-learning}

In \Cref{subsec:tensor-decomposition:iterative-projections}, we constructed a projection map
and estimated the even and odd moments in a lower-dimensional space.
In this subsection, we use these estimates to recover a block of normalized regressors.
On a high-level, we will construct an observable-operator family (\Cref{def:approximate-observable-operator-family}) whose eigenvalues encode projections $\ip{\bm{\zeta}}{\ScaledRegressor_i}$ of the scaled regressors.
Following \citet{AHK12}, we diagonalize a random contraction
and use its eigenvectors to read the coordinates of the scaled regressors.

\subsubsection{Proof of
\texorpdfstring{\Cref*{thm:tensor-decomposition:observable-operator-construction}}{Theorem 6.5}}
\label{subsubsec:tensor-decomposition:operator-construction-proof}
We first prove \Cref{thm:tensor-decomposition:observable-operator-construction} by constructing an observable-operator family for a subset of the remaining
scaled regressors.
Let $Q\subseteq [k]$ index the unrecovered regressors and $\widehat A, \widehat M_1, \dots, \widehat M_n$ denote the estimates to the projected moments returned by \Cref{alg:tensor-decomposition:estimate-projected-moments}.
Roughly speaking, we have
\begin{align*}
  \widehat A &\approx \sum_{i\in Q} \ComponentWeight{2s}{i} (C_s \NormalizedRegressor_i^{\otimes s}) (C_s\NormalizedRegressor_i^{\otimes s})^\top\,, \\
  \sum_\ell \zeta_\ell \widehat M_\ell
  &\approx \sum_{i\in Q} \ComponentWeight{2s}{i} \ip{\bm{\zeta}}{\ScaledRegressor_i} (C_s \NormalizedRegressor_i^{\otimes s}) (C_s\NormalizedRegressor_i^{\otimes s})^\top\,.
\end{align*}
We would like to rescale the odd moments by the even moment to cancel out the weights $\ComponentWeight{2s}{i}$
\[
  \widehat A^{-\frac12 } \left( \sum_\ell \zeta_\ell \widehat M_\ell \right) \widehat A^{-\frac12}
  \approx \sum_{i\in Q} \ip{\bm{\zeta}}{\ScaledRegressor_i} (C_s \NormalizedRegressor_i^{\otimes s}) (C_s\NormalizedRegressor_i^{\otimes s})^\top\,.
\]
For this to be a well-conditioned operation,
we would need the eigenvalues of $\widehat A$ to be bounded away from $0$.
Unfortunately, due to the irreducible bias that we swept under the rug and the fact that we only have approximate orthonormality,
we cannot guarantee this is the case.
Hence we first truncate $\widehat A$ to its leading eigenspace up to the largest gap between eigenvalues.
Then we apply the rescaling procedure using the truncated even moment and construct an observable-operator family for the scaled regressors corresponding to the truncated eigenvectors.
The pseudocode for this procedure can be found in \Cref{alg:tensor-decomposition:construct-observable-operators}.

\begin{algorithm}[H]
\caption[\textsc{ConstructObservableOperators}]{\textsc{ConstructObservableOperators}}\label{alg:tensor-decomposition:construct-observable-operators}
\KwIn{Symmetric matrices $\widehat A,\widehat M_1,\ldots,\widehat M_n \in\R^{p_s\times p_s}$}
\KwOut{A dimension $i_\star$ and matrices
$\widehat D_1,\ldots,\widehat D_n$}
Compute the eigenvalues $\widehat\lambda_1\geq\cdots\geq
\widehat\lambda_{p_s}$ and corresponding orthonormal eigenvectors
$\widehat{\bm{v}}_1,\ldots,\widehat{\bm{v}}_{p_s}$ of $\widehat A$\;
$\widehat\lambda_h\gets\max\{\widehat\lambda_h,0\}$ for
$h\in[p_s]$, and $\widehat\lambda_{p_s+1}\gets0$\;
$i_\star\gets\arg\max_{1\leq h\leq p_s}
 (\widehat\lambda_h-\widehat\lambda_{h+1})$\;
$\widehat V\gets
 \begin{bmatrix}
  \widehat{\bm{v}}_1^\top\\
  \vdots\\
  \widehat{\bm{v}}_{i_\star}^\top
 \end{bmatrix}$\;
$\widehat\Lambda\gets
  \Diag(\widehat\lambda_1,\ldots,\widehat\lambda_{i_\star})$\;
\For{$\ell=1,\ldots,n$}{
  $\widehat D_\ell\gets\widehat\Lambda^{-1/2}
  \widehat V\widehat M_\ell\widehat V^\top
  \widehat\Lambda^{-1/2}$\;
}
\Return $(i_\star,\widehat D_1,\ldots,\widehat D_n)$\;
\end{algorithm}

Before proving \Cref{thm:tensor-decomposition:observable-operator-construction},
we state a useful spectral inequality whose proof is deferred to \Cref{subsec:app-tensor-decomposition:observable-operator-learning}.
\begin{restatable}[Largest-Gap Eigenvalue Bounds]{proposition}{largestgapeigenvaluebounds}
\label{prop:tensor-decomposition:largest-gap-eigenvalue-bounds}
Let $1\leq p\leq q$, let $W\in\R^{p\times q}$, and let
$\Lambda\in\R^{q\times q}$ be nonzero and positive semidefinite.
Suppose the symmetric matrix $A\in\R^{p\times p}$
satisfies
\begin{align*}
 \norm{A-WW^\top}_{\op}
 \leq\frac{\norm{\Lambda}_{\op}}{8q}\,,\qquad
 \norm{W^\top W-\Lambda}_{\op}
 \leq\frac{\norm{\Lambda}_{\op}}{8q}\,.
\end{align*}
Write the eigenvalues in nonincreasing order and put
$
  \widehat\lambda_h
 =\max\{\lambda_h(A),0\}, h\in[p]
$,
where we understand $\widehat\lambda_{p+1}=0$ and
$\lambda_{q+1}(\Lambda)=0$.
Then for the largest eigengap indexed by
$ 
  i_\star\in\arg\max_{1\leq h\leq p}
 (\widehat\lambda_h-\widehat\lambda_{h+1})
$,
we have
\begin{align*}
 \widehat\lambda_{i_\star}
 \geq\frac{\norm{\Lambda}_{\op}}{2q}\,,\qquad
 \widehat\lambda_{i_\star}
 -\lambda_{i_\star+1}(\Lambda)
 \geq\frac{\norm{\Lambda}_{\op}}{2q}\,.
\end{align*}
\end{restatable}

We are now ready to analyze \Cref{alg:tensor-decomposition:construct-observable-operators}.
The proof proceeds by decomposing the error into three interpretable terms.
To avoid cumbersome calculations that distract from our understanding,
we defer the analysis of the ones that are a consequence of the near orthogonality of the rank-one tensors to the appendix.
\begin{proof}[Proof of \cref{thm:tensor-decomposition:observable-operator-construction}]
We first describe the error decomposition assuming that
$\widehat\Lambda$ is positive definite.
We verify positivity in the perturbation bounds below.
Since $|Q|\leq k$, the error assumption also holds with
$|Q|$ in place of $k$.

Order the components by decreasing $\ComponentWeight{2s}{i}$.
Define the factor matrix
\[
 H\coloneqq
 \left[\sqrt{\ComponentWeight{2s}{i}}\,\TensorPower{\NormalizedRegressor_i}{s}\right]_{i\in Q}\,,
 \qquad
 \Lambda\coloneqq\Diag_{i\in Q}(\ComponentWeight{2s}{i})\,,
\]
so that
\begin{align*}
 \EvenMomentMatrix{s}{Q}=HH^\top\,,\qquad
 \OddMomentMatrix{s}{Q}{\bm{\zeta}}=H\Diag_{i\in Q}(\ip{\bm{\zeta}}{\ScaledRegressor_i})H^\top\,.
\end{align*}
Now, let $\widehat V^\top\coloneq [\widehat{\bm{v}}_1, \dots, \widehat{\bm{v}}_{i_\star}]$ be the matrix of top $i_\star$ eigenvectors of $\widehat A$,
as chosen in \Cref{alg:tensor-decomposition:construct-observable-operators}.
Define the ideal projected, truncated, and rescaled factor matrix
\begin{equation}\label{eq:tensor-decomposition:rescaled-components}
 F=[\bm{f}_i]_{i\in Q}
 =\widehat\Lambda^{-1/2}\widehat V
 \ProjectionMap{s} H\,,\qquad
 \bm f_i \coloneqq \sqrt{\ComponentWeight{2s}{i}} \widehat\Lambda^{-1/2}\widehat V (\ProjectionMap{s} \NormalizedRegressor_i^{\otimes s})\,.
\end{equation}
Note that by \Cref{thm:hermite:structure-component-weights}, $\ComponentWeight{2s}{i} > 0$ so that $\Lambda\neq 0$ and its inverse square root is well-defined.

Let {$J\subseteq Q$} index the first $i_\star$ components
in the fixed order, and write $F_J$ and $F_{Q\setminus J}$
for the corresponding column submatrices.
Expanding into outer products gives
\begin{align*}
 F\Diag_{i\in Q}(\ip{\bm{\zeta}}{\ScaledRegressor_i})F^\top
 &=\sum_{i\in Q}\ip{\bm{\zeta}}{\ScaledRegressor_i}\bm{f}_i\bm{f}_i^\top\\
 &=\sum_{i\in J}\ip{\bm{\zeta}}{\ScaledRegressor_i}\bm{f}_i\bm{f}_i^\top
   +\sum_{i\in Q\setminus J}\ip{\bm{\zeta}}{\ScaledRegressor_i}\bm{f}_i\bm{f}_i^\top\\
 &=F_J\Diag_{i\in J}(\ip{\bm{\zeta}}{\ScaledRegressor_i})F_J^\top
 +
 F_{Q\setminus J}\Diag_{i\in Q\setminus J}(\ip{\bm{\zeta}}{\ScaledRegressor_i})
 F_{Q\setminus J}^\top\,.
\end{align*}
We will construct an orthogonal matrix $Z=[\bm{z}_i]_{i\in J}$.
For $\widehat D(\bm{\zeta})=\sum_{\ell=1}^n \zeta_\ell\widehat D_\ell$,
the corresponding error decomposes as
\begin{align}
 &\widehat D(\bm{\zeta})-\sum_{i\in J}\ip{\bm{\zeta}}{\ScaledRegressor_i}\bm{z}_i\bm{z}_i^\top
 \notag\\
 &\quad=
 \underbrace{
 \widehat D(\bm{\zeta})-F\Diag_{i\in Q}(\ip{\bm{\zeta}}{\ScaledRegressor_i})F^\top
 }_{\text{odd estimation error}}
 \notag\\
 &\qquad+
 \underbrace{
 F_{Q\setminus J}
 \Diag_{i\in Q\setminus J}(\ip{\bm{\zeta}}{\ScaledRegressor_i})
 F_{Q\setminus J}^\top
 }_{\text{spectral truncation error}}
 \notag\\
 &\qquad+
 \underbrace{
 F_J\Diag_{i\in J}(\ip{\bm{\zeta}}{\ScaledRegressor_i})F_J^\top
 -\sum_{i\in J}\ip{\bm{\zeta}}{\ScaledRegressor_i}\bm{z}_i\bm{z}_i^\top
 }_{\text{approximate orthonormality}}.
 \label{eq:tensor-decomposition:operator-family-error-decomposition}
\end{align}
By the triangle inequality, it suffices to bound the
Frobenius norm of each term.

The first term comes from estimating the projected odd moment.
Rescaling can blow up this estimation error by at most a factor
of $1/\widehat\lambda_{i_\star}$ in Frobenius norm.
We therefore need a lower bound on the smallest eigenvalue used
for rescaling.

For the second term, suppose the projected lifted vectors are
orthonormal and the even moment is known exactly. Then the rows
of $\widehat V$ span the projected directions indexed by $J$, so
$\widehat V\ProjectionMap{s}\TensorPower{\NormalizedRegressor_i}{s}=0$
for every $i\in Q\setminus J$. For the estimated even moment, we
apply Sylvester's equation (\Cref{prop:sylvester-equation}) and
the eigenvalue gap to bound the contribution from these components.

In the same exact case, rescaling cancels the factor
$\sqrt{\ComponentWeight{2s}{i}}$ in each column of $F_J$.
The resulting vectors already form an orthonormal basis, so
the third term vanishes. When the even moment is accurate on
the block and little weight enters from outside $J$, only a
small correction is needed to make the rescaled vectors orthonormal.

\emph{Perturbation bounds.}
We first compare the estimated even moment with the projected
moment matrix, and the projected component matrix with
the component weights.
These bounds will verify the hypotheses of
\Cref{prop:tensor-decomposition:largest-gap-eigenvalue-bounds}.
By the even-moment estimation assumption
(\Cref{def:projected-moment-estimate-even}),
\begin{align}
 \norm{\widehat A-\ProjectionMap{s} \EvenMomentMatrix{s}{Q}\ProjectionMap{s}^\top}_F
 &\leq\varepsilon_{\mathrm{mom}}\tr(\EvenMomentMatrix{s}{Q})
 \alignnote{by \Cref{def:projected-moment-estimate-even}}\notag\\
 &\leq|Q|\varepsilon_{\mathrm{mom}}\norm{\Lambda}_{\op}.
 \alignnote{$\tr(\EvenMomentMatrix{s}{Q})=\tr(\Lambda)\leq|Q|\norm{\Lambda}_{\op}$}
 \label{eq:tensor-decomposition:projected-even-error}
\end{align}

The second comparison separates approximate orthogonality
from the change caused by projection.
The approximate-orthogonality hypothesis
(\Cref{ass:tensor-decomposition:approximate-orthogonality}) gives
\begin{align}
 \norm{H^\top H-\Lambda}_{\op}
 &=\norm{
 \Lambda^{1/2}
 \left[
 \bigl(\ip{\TensorPower{\NormalizedRegressor_i}{s}}{\TensorPower{\NormalizedRegressor_j}{s}}\bigr)_{i,j\in Q}-I
 \right]
 \Lambda^{1/2}}_{\op}\notag\\
 &\leq\varepsilon_\perp\norm{\Lambda}_{\op}.
 \alignnote{by \Cref{ass:tensor-decomposition:approximate-orthogonality}}
 \label{eq:tensor-decomposition:weighted-component-orthogonality}
\end{align}
Using the orthonormality of the rows of $\ProjectionMap{s}$
and the norm identity (\Cref{prop:gram-matrix-identities}), we obtain
\begin{align}
 \norm{H^\top\ProjectionMap{s}^\top\ProjectionMap{s} H-\Lambda}_{\op}
 &\leq
 \underbrace{\norm{H^\top H-\Lambda}_{\op}}
 _{\text{approximate orthogonality}}\notag\\
 &\quad+
 \underbrace{\norm{H^\top(I-\ProjectionMap{s}^\top\ProjectionMap{s})H}_{\op}}
 _{\text{projection error}}
 \alignnote{triangle inequality}\notag\\
 &\leq\varepsilon_\perp\norm{\Lambda}_{\op}
 +\norm{(I-\ProjectionMap{s}^\top\ProjectionMap{s})H}_{\op}^2
 \alignnote{by \eqref{eq:tensor-decomposition:weighted-component-orthogonality}
 and \Cref{prop:gram-matrix-identities}}\notag\\
 &=\varepsilon_\perp\norm{\Lambda}_{\op}\notag\\
 &\quad+\norm{(I-\ProjectionMap{s}^\top\ProjectionMap{s})
       (\EvenMomentMatrix{s}{Q})^{1/2}}_{\op}^2
 \alignnote{by \Cref{prop:gram-matrix-identities}}\notag\\
 &\leq\varepsilon_\perp\norm{\Lambda}_{\op}
 +\varepsilon_{\mathrm{proj}}\tr(\EvenMomentMatrix{s}{Q})
 \alignnote{by \Cref{ass:tensor-decomposition:projection-error}}\notag\\
 &\leq
 \left(|Q|\varepsilon_{\mathrm{proj}}+\varepsilon_\perp\right)
 \norm{\Lambda}_{\op}.
 \label{eq:tensor-decomposition:projected-weight-error}
\end{align}

The bounds \eqref{eq:tensor-decomposition:projected-even-error} and
\eqref{eq:tensor-decomposition:projected-weight-error} verify the
hypotheses of \Cref{prop:tensor-decomposition:largest-gap-eigenvalue-bounds}, since
\begin{align*}
 |Q|\varepsilon_{\mathrm{mom}}
 +|Q|\varepsilon_{\mathrm{proj}}+\varepsilon_\perp
 &\leq |Q|\left(
 \varepsilon_{\mathrm{mom}}
 +\varepsilon_{\mathrm{proj}}+\varepsilon_\perp
 \right)\\
 &\leq\frac1{16|Q|}.
 \alignnote{by \Cref{ass:tensor-decomposition:small-errors}, $|Q|\leq k$}
\end{align*}
Each coefficient is therefore at most $1/(8|Q|)$, and
$\norm{\cdot}_{\op}\leq\norm{\cdot}_F$ converts the even
estimation bound to operator norm.
Applying \Cref{prop:tensor-decomposition:largest-gap-eigenvalue-bounds} with
$W=\ProjectionMap{s} H$ gives
\begin{equation}\label{eq:tensor-decomposition:selected-eigenvalue-lower-bound}
 \widehat\lambda_{i_\star}
 \geq\frac{\norm{\Lambda}_{\op}}{2|Q|}.
\end{equation}
This bound shows that $\widehat\Lambda^{-1/2}$ is defined
and controls its norm.

\emph{Odd estimation error.}
We first bound the odd estimation error in
\eqref{eq:tensor-decomposition:operator-family-error-decomposition}.
Using the definitions of $\widehat D(\bm{\zeta})$ and $F$ and the
odd-moment estimation assumption (\Cref{def:projected-moment-estimate-odd}),
we obtain
\begin{align}
 &\norm{\widehat D(\bm{\zeta})
       -F\Diag_{i\in Q}(\ip{\bm{\zeta}}{\ScaledRegressor_i})F^\top}_F
 \notag\\
 &\quad=
 \norm{\widehat\Lambda^{-1/2}\widehat V
 \left(
 \sum_{\ell=1}^n \zeta_\ell\widehat M_\ell
 -\ProjectionMap{s} \OddMomentMatrix{s}{Q}{\bm{\zeta}}\ProjectionMap{s}^\top
 \right)
 \widehat V^\top\widehat\Lambda^{-1/2}}_F
 \alignnote{by \Cref{alg:tensor-decomposition:construct-observable-operators} and \eqref{eq:tensor-decomposition:rescaled-components}}
 \notag\\
 &\quad\leq
 \frac{
 \norm{\sum_{\ell=1}^n \zeta_\ell\widehat M_\ell
 -\ProjectionMap{s} \OddMomentMatrix{s}{Q}{\bm{\zeta}}\ProjectionMap{s}^\top}_F}
 {\widehat\lambda_{i_\star}}
 \notag\\
 &\quad\leq
 \frac{\varepsilon_{\mathrm{mom}}\tr(\EvenMomentMatrix{s}{Q})}
      {\widehat\lambda_{i_\star}}\norm{\bm{\zeta}}_2
 \alignnote{by assumption (\Cref{def:projected-moment-estimate-odd})}
 \notag\\
 &\quad\leq2|Q|^2\varepsilon_{\mathrm{mom}}\norm{\bm{\zeta}}_2.
 \alignnote{by \eqref{eq:tensor-decomposition:selected-eigenvalue-lower-bound}, $\tr(\Lambda)\leq|Q|\norm{\Lambda}_{\op}$}
 \label{eq:tensor-decomposition:rescaled-odd-error}
\end{align}

This bounds the first term in the error decomposition \eqref{eq:tensor-decomposition:operator-family-error-decomposition}.
For the last two terms, take $Z$ to be an orthogonal polar factor
of $F_J$.
We defer the bounds for the spectral truncation error (\Cref{lem:app-tensor-decomposition:truncation-error}) and the approximate orthonormality term \Cref{lem:app-tensor-decomposition:orthonormality-error} to the appendix.

\emph{Putting it together.}
Applying the triangle inequality to
\eqref{eq:tensor-decomposition:operator-family-error-decomposition} and substituting
the three bounds gives
\begin{align}
 &\norm{\widehat D(\bm{\zeta})
       -\sum_{i\in J}\ip{\bm{\zeta}}{\ScaledRegressor_i}\bm{z}_i\bm{z}_i^\top}_F\notag\\
 &\quad\leq\Bigl(
 \underbrace{2|Q|^2\varepsilon_{\mathrm{mom}}}
 _{\text{odd estimation error}}
 \alignnote{by \eqref{eq:tensor-decomposition:rescaled-odd-error}}\notag\\
 &\qquad+
 \underbrace{K^3\norm{F_{Q\setminus J}}_F^2}
 _{\text{spectral truncation error}}
 \alignnote{by \Cref{lem:app-tensor-decomposition:truncation-error}}\notag\\
 &\qquad+
 \underbrace{3K^3\bigl(2|Q|^2\varepsilon_{\mathrm{mom}}
           +\norm{F_{Q\setminus J}}_F^2\bigr)}
 _{\text{approximate orthonormality}}
 \Bigr)\norm{\bm{\zeta}}_2
 \alignnote{by \Cref{lem:app-tensor-decomposition:orthonormality-error}}\notag\\
 &\quad\leq
 \left(
 8K^3|Q|^2\varepsilon_{\mathrm{mom}}
 +4K^3\norm{F_{Q\setminus J}}_F^2
 \right)\norm{\bm{\zeta}}_2
 \alignnote{$K\geq1$}
 \\
 &\quad\leq
 \Bigl(
 8K^3|Q|^2\varepsilon_{\mathrm{mom}}
 +4K^3|Q|^3(
 \varepsilon_{\mathrm{mom}}
 +\varepsilon_{\mathrm{proj}}
 +\varepsilon_\perp)
 \Bigr)\norm{\bm{\zeta}}_2
 \alignnote{by \Cref{lem:app-tensor-decomposition:truncation-outside-components}}\notag\\
 &\quad\leq16K^3|Q|^3
 \left(
 \varepsilon_{\mathrm{mom}}
 +\varepsilon_{\mathrm{proj}}
 +\varepsilon_\perp
 \right)\norm{\bm{\zeta}}_2.
 \alignnote{$|Q|\geq1$}
 \label{eq:tensor-decomposition:operator-family-error-bound}
\end{align}
The same basis works for every $\bm{\zeta}\in\R^n$.
Since $|Q|\leq k$, \eqref{eq:tensor-decomposition:operator-family-error-bound} gives the claimed
family approximation in the sense of \Cref{def:approximate-observable-operator-family}.
All matrices used by the routine have dimension at most
$p_s\leq|Q|\leq k$, so it runs in $\poly(n,k)$ time.
\end{proof}

\subsubsection{Proof of
\texorpdfstring{\Cref*{thm:tensor-decomposition:normalized-regressor-recovery}}{Theorem 6.6}}
\label{subsubsec:tensor-decomposition:direction-recovery-proof}

\Cref{thm:tensor-decomposition:observable-operator-construction} gives an approximate observable-operator family
for a block $(\ScaledRegressor_i)_{i\in J}$ that approximates
\[
  D(\bm{\zeta})=\sum_{i\in J}\ip{\bm{\zeta}}{\ScaledRegressor_i}\bm{z}_i\bm{z}_i^\top.
\]
Here, we prove \Cref{thm:tensor-decomposition:normalized-regressor-recovery},
by demonstrating how to recover the normalized regressors given this family.

Observe that the eigenvalues of $D(\bm{e}_j)$ are
the $j$-th coordinates of the vectors $\ScaledRegressor_i$.
A first attempt might be to compute the eigenvalues
for each $D(\bm{e}_j)$ and assemble the vectors.
However, an eigendecomposition outputs each coordinate as an unordered
list, and does not specify the mapping to the particular regressor.

In order to match the coordinates to regressors,
we must also recover the common eigenvectors $\bm{z}_i$'s.
This can be accomplished via an eigendecomposition of $D(\bm{\zeta})$,
assuming its spectrum is simple.
We can ensure a quantitative separation between the eigenvalues
with constant probability by picking a random unit vector $\bm{\zeta}$.
We then recover the coordinates through
\[
 \bm{z}_i^\top D(\bm{e}_j)\bm{z}_i=(\ScaledRegressor_i)_j
 \qquad(i\in J,\ j\in[n]).
\]
The pseudocode for this subroutine can be found in \Cref{alg:tensor-decomposition:recover-directions}.

\begin{algorithm}[H]
\caption[\textsc{RecoverDirections}]{\textsc{RecoverDirections}}
\label{alg:tensor-decomposition:recover-directions}
\KwIn{Symmetric matrices $(\widehat D_j)_{j\in[n]}$
of common dimension $p\geq2$, and confidence parameter $\delta$}
\KwOut{Unit vectors $\widehat{\NormalizedRegressor}_1,\ldots,\widehat{\NormalizedRegressor}_p$}
$L\gets\log_2(\nicefrac1\delta)$\;
Draw independent
$\bm{\zeta}_1,\ldots,\bm{\zeta}_L\sim\operatorname{Unif}(\mathbb S^{n-1})$\;
\For{$\iota=1,\ldots,L$}{
  $\widehat D(\bm{\zeta}_\iota)\gets
   \sum_{j=1}^n(\bm{\zeta}_\iota)_j\widehat D_j$\;
  Compute the eigenvalues
  $\widehat\nu_1^{(\iota)}\geq\cdots\geq\widehat\nu_p^{(\iota)}$
  and corresponding orthonormal eigenvectors
  $\widehat{\bm{z}}_1^{(\iota)},\ldots,\widehat{\bm{z}}_p^{(\iota)}$
  of $\widehat D(\bm{\zeta}_\iota)$\;
}
$\iota_\star\gets
 \arg\max_{1\leq\iota\leq L}
 \min_{1\leq i<p}
 (\widehat\nu_i^{(\iota)}-\widehat\nu_{i+1}^{(\iota)})$\;
For $i\in[p]$ and $j\in[n]$, set
$(\widehat{\ScaledRegressor}_i)_j\gets
 (\widehat{\bm{z}}_i^{(\iota_\star)})^\top
 \widehat D_j\widehat{\bm{z}}_i^{(\iota_\star)}$\;
$\widehat{\NormalizedRegressor}_i\gets\widehat{\ScaledRegressor}_i/\norm{\widehat{\ScaledRegressor}_i}_2$
for $i\in[p]$\;
\Return $(\widehat{\NormalizedRegressor}_i)_{i\in[p]}$\;
\end{algorithm}

We now prove that $O(\log(\nicefrac1\delta))$
random unit directions suffice to recover every normalized vector $\ScaledRegressor_i/\norm{\ScaledRegressor_i}_2$ with bounded error.
\begin{proof}[Proof of \cref{thm:tensor-decomposition:normalized-regressor-recovery}]
We first show that our algorithm finds an approximate common
eigenbasis for the family.
We then show that evaluating the coordinate quadratic forms
in this basis and normalizing the resulting vectors gives
the claimed error guarantees.

Distances between unit vectors are at most $2$, so it
suffices to prove the claim when
\begin{equation}\label{eq:tensor-decomposition:direction-recovery-small-error}
 \varepsilon_{\mathrm{fam}}
 <\frac1{16K^6\sqrt n\,p^2}.
\end{equation}

\emph{Recovering the common eigenvectors.}
Our strategy is to apply a Davis--Kahan type eigenvector bound
(\Cref{prop:davis-kahan-eigenvector}) to control the eigenvector
error by the operator-norm error divided by an eigenvalue gap.
For unit contractions, the operator-norm error is bounded above by $\varepsilon_{\mathrm{fam}}$ by assumption,
so we first establish the required eigenvalue gap.

Let $(\bm{z}_i)_{i\in J}$ be the orthonormal basis supplied by assumption (\Cref{def:approximate-observable-operator-family}),
and put
\[
 D(\bm{\zeta})=\sum_{i\in J}\ip{\bm{\zeta}}{\ScaledRegressor_i}\bm{z}_i\bm{z}_i^\top.
\]
Its eigenvalues are $(\ip{\bm{\zeta}}{\ScaledRegressor_i})_{i\in J}$.
For a uniformly random unit vector $\bm{\zeta}\sim\mathbb S^{n-1}$,
spherical anti-concentration
(\Cref{prop:spherical-anticoncentration}) and a union bound give
\begin{align*}
 &\Prb_{\bm{\zeta}}\!\left(
 \min_{i\neq i'\in J}
 \abs{\ip{\bm{\zeta}}{\ScaledRegressor_i-\ScaledRegressor_{i'}}}
 \leq\frac1{K^3\sqrt n\,p^2}
 \right)\\
 &\quad\leq
 \sum_{\substack{i,i'\in J\\i<i'}}
 \frac1{K^3p^2\norm{\ScaledRegressor_i-\ScaledRegressor_{i'}}_2}
 \alignnote{union bound and \Cref{prop:spherical-anticoncentration}}\\
 &\quad\leq\frac{\binom p2}{p^2}<\frac12.
 \alignnote{$\norm{\ScaledRegressor_i-\ScaledRegressor_{i'}}_2\geq K^{-3}$}
\end{align*}
The probability that none of the $L = \log_2(\nicefrac1\delta)$ independent trials
has this separation is at most
$
 2^{-L}\leq\delta.
$
For the rest of the proof, suppose that at least one trial
has this separation.

We now show that the eigendecomposition chosen by our
algorithm approximates the common eigenbasis.
By assumption, the family approximation bound and $\norm{\bm{\zeta}_\iota}_2=1$ give
\begin{equation}\label{eq:tensor-decomposition:random-contraction-error}
 \norm{\widehat D(\bm{\zeta}_\iota)-D(\bm{\zeta}_\iota)}_{\op}
 \leq\norm{\widehat D(\bm{\zeta}_\iota)-D(\bm{\zeta}_\iota)}_F
 \leq\varepsilon_{\mathrm{fam}}
 \qquad(\iota\in[L]).
\end{equation}
By Weyl's inequality
(\Cref{prop:weyl-perturbation}), each estimated
consecutive eigengap differs from the corresponding true
eigengap by at most $2\varepsilon_{\mathrm{fam}}$.
Since our algorithm maximizes the minimum estimated gap,
\begin{align}
 \min_{1\leq i<p}
 \bigl(\widehat\nu_i^{(\iota_\star)}
       -\widehat\nu_{i+1}^{(\iota_\star)}\bigr)
 &\geq\frac1{K^3\sqrt n\,p^2}
       -2\varepsilon_{\mathrm{fam}}
 \alignnote{separated trial and choice of $\iota_\star$}\notag\\
 &\geq\frac1{K^3\sqrt n\,p^2}
       -\frac1{8K^6\sqrt n\,p^2}
 \alignnote{by \eqref{eq:tensor-decomposition:direction-recovery-small-error}}\notag\\
 &\geq\frac3{4K^3\sqrt n\,p^2}.
 \alignnote{$K\geq1$}
 \label{eq:tensor-decomposition:estimated-contraction-gap}
\end{align}
As this is at least
$3/(4K^3\sqrt n\,p^2)>2\varepsilon_{\mathrm{fam}}$
by \eqref{eq:tensor-decomposition:direction-recovery-small-error},
Weyl's inequality shows that $\widehat D(\bm{\zeta}_{\iota_\star})$
also has simple spectrum.
Matching the true eigenvectors in eigenvalue order gives
a bijection $\pi:[p]\to J$.
Write $\widehat{\bm{z}}_i=\widehat{\bm{z}}_i^{(\iota_\star)}$ for the
eigenvectors used by the algorithm.
With $\widehat D(\bm{\zeta}_{\iota_\star})$ as the reference matrix,
the Davis--Kahan eigenvector bound gives, after choosing
the signs of the true eigenvectors,
\begin{align}
 \norm{\widehat{\bm{z}}_i-\bm{z}_{\pi(i)}}_2
 &\leq
 \frac{
  2\sqrt2\,\norm{\widehat D(\bm{\zeta}_{\iota_\star})
                    -D(\bm{\zeta}_{\iota_\star})}_{\op}
 }{
  \min_{1\leq i'<p}
  \bigl(\widehat\nu_{i'}^{(\iota_\star)}
       -\widehat\nu_{i'+1}^{(\iota_\star)}\bigr)
 }
 \alignnote{by \Cref{prop:davis-kahan-eigenvector}}\notag\\
 &\leq
 2\sqrt2\,\varepsilon_{\mathrm{fam}}
 \frac{4K^3\sqrt n\,p^2}{3}
 \alignnote{by \eqref{eq:tensor-decomposition:random-contraction-error}
 and \eqref{eq:tensor-decomposition:estimated-contraction-gap}}\notag\\
 &\leq4K^3\sqrt n\,p^2\varepsilon_{\mathrm{fam}}
 \alignnote{$2\sqrt2\leq3$}
 \label{eq:tensor-decomposition:contracted-eigenvector-error}\\
 &<\frac1{4K^3}.
 \alignnote{by \eqref{eq:tensor-decomposition:direction-recovery-small-error}}
 \label{eq:tensor-decomposition:eigenvector-distance-smallness}
\end{align}

\emph{Recovering the normalized regressors.}
We now use the approximate common eigenvectors to estimate
the scaled regressors.
Since $\bm{z}_{\pi(i)}$ is an eigenvector of every $D(\bm{\zeta})$
and both vectors are unit, the coordinate construction gives
\begin{align}
 \ip{\bm{\zeta}}{\widehat{\ScaledRegressor}_i-\ScaledRegressor_{\pi(i)}}
 &=\widehat{\bm{z}}_i^\top
   \bigl(\widehat D(\bm{\zeta})-D(\bm{\zeta})\bigr)\widehat{\bm{z}}_i\notag\\
 &\quad+(\widehat{\bm{z}}_i-\bm{z}_{\pi(i)})^\top
 \bigl(D(\bm{\zeta})-\ip{\bm{\zeta}}{\ScaledRegressor_{\pi(i)}}I\bigr)
 (\widehat{\bm{z}}_i-\bm{z}_{\pi(i)}).
 \label{eq:tensor-decomposition:quadratic-readout-error-decomposition}
\end{align}
For every unit vector $\bm{\zeta}$,
\begin{align}
 \norm{D(\bm{\zeta})-\ip{\bm{\zeta}}{\ScaledRegressor_{\pi(i)}}I}_{\op}
 &=\max_{j\in J}\abs{\ip{\bm{\zeta}}{\ScaledRegressor_j-\ScaledRegressor_{\pi(i)}}}
 \alignnote{common orthonormal eigenbasis}\notag\\
 &\leq\max_{j\in J}
 \bigl(\norm{\ScaledRegressor_j}_2+\norm{\ScaledRegressor_{\pi(i)}}_2\bigr)
 \alignnote{Cauchy--Schwarz and triangle inequality}\notag\\
 &\leq2K^3.
 \alignnote{$\norm{\ScaledRegressor_j}_2\leq K^3$}
 \label{eq:tensor-decomposition:centered-coordinate-operator-bound}
\end{align}
Taking the supremum over unit $\bm{\zeta}$ in
\eqref{eq:tensor-decomposition:quadratic-readout-error-decomposition}
and using the family approximation assumption
(\Cref{def:approximate-observable-operator-family}) gives
\begin{align}
 \norm{\widehat{\ScaledRegressor}_i-\ScaledRegressor_{\pi(i)}}_2
 &=\sup_{\norm{\bm{\zeta}}_2=1}
   \abs{\ip{\bm{\zeta}}{\widehat{\ScaledRegressor}_i-\ScaledRegressor_{\pi(i)}}}\notag\\
 &\leq\varepsilon_{\mathrm{fam}}
       +2K^3\norm{\widehat{\bm{z}}_i-\bm{z}_{\pi(i)}}_2^2
 \alignnote{by \eqref{eq:tensor-decomposition:centered-coordinate-operator-bound}}\notag\\
 &\leq\varepsilon_{\mathrm{fam}}
       +\frac12\norm{\widehat{\bm{z}}_i-\bm{z}_{\pi(i)}}_2
 \alignnote{by \eqref{eq:tensor-decomposition:eigenvector-distance-smallness}}\notag\\
 &\leq
 \bigl(1+2K^3\sqrt n\,p^2\bigr)\varepsilon_{\mathrm{fam}}
 \alignnote{by \eqref{eq:tensor-decomposition:contracted-eigenvector-error}}\notag\\
 &\leq3K^3\sqrt n\,p^2\varepsilon_{\mathrm{fam}}.
 \alignnote{$K,n,p\geq1$}
 \label{eq:tensor-decomposition:scaled-regressor-error}
\end{align}
By \eqref{eq:tensor-decomposition:scaled-regressor-error} and
\eqref{eq:tensor-decomposition:direction-recovery-small-error},
\[
 \norm{\widehat{\ScaledRegressor}_i-\ScaledRegressor_{\pi(i)}}_2
 <\frac3{16K^3}
 <K^{-3}\leq\norm{\ScaledRegressor_{\pi(i)}}_2.
\]
Thus every $\widehat{\ScaledRegressor}_i$ is nonzero.
Since $\widehat{\ScaledRegressor}_i=\norm{\widehat{\ScaledRegressor}_i}_2\widehat{\NormalizedRegressor}_i$
and $\norm{\widehat{\NormalizedRegressor}_i}_2=1$, normalization gives
\begin{align}
 &\norm{\widehat{\NormalizedRegressor}_i-
       \frac{\ScaledRegressor_{\pi(i)}}{\norm{\ScaledRegressor_{\pi(i)}}_2}}_2\notag\\
 &\quad=
 \Bigl\|
 \frac{\widehat{\ScaledRegressor}_i-\ScaledRegressor_{\pi(i)}}{\norm{\ScaledRegressor_{\pi(i)}}_2}
 +\left(1-\frac{\norm{\widehat{\ScaledRegressor}_i}_2}
                   {\norm{\ScaledRegressor_{\pi(i)}}_2}\right)\widehat{\NormalizedRegressor}_i
 \Bigr\|_2\notag\\
 &\quad\leq
 \frac{\norm{\widehat{\ScaledRegressor}_i-\ScaledRegressor_{\pi(i)}}_2
       +\abs{\norm{\ScaledRegressor_{\pi(i)}}_2-\norm{\widehat{\ScaledRegressor}_i}_2}}
      {\norm{\ScaledRegressor_{\pi(i)}}_2}
 \alignnote{triangle inequality}\notag\\
 &\quad\leq
 \frac{2\norm{\widehat{\ScaledRegressor}_i-\ScaledRegressor_{\pi(i)}}_2}
      {\norm{\ScaledRegressor_{\pi(i)}}_2}
 \alignnote{reverse triangle inequality}\notag\\
 &\quad\leq6K^6\sqrt n\,p^2\varepsilon_{\mathrm{fam}}.
 \alignnote{by \eqref{eq:tensor-decomposition:scaled-regressor-error};
 $\norm{\ScaledRegressor_{\pi(i)}}_2\geq K^{-3}$}
\end{align}
This gives the claimed error bound for normalized vectors.
Each of the $L=O(\log(\nicefrac1\delta))$ trials forms an
$p\times p$ matrix and computes its eigendecomposition.
Since $p\leq k$, these computations and the coordinate
quadratic forms take $\poly(n,k)\log(\nicefrac1\delta)$ time.
\end{proof}

\subsection{Proof of \texorpdfstring{\Cref*{thm:tensor-decomposition:implicit-tensor-decomposition}}{Theorem}}
\label{subsec:tensor-decomposition:decomposition-proof}

Having analyzed each of the subroutines in our implicit tensor decomposition algorithm (\Cref{alg:tensor-decomposition:implicit-tensor-decomposition}),
we are now ready to put them together and prove its guarantees.
The pseudocode can be found in \Cref{alg:tensor-decomposition:implicit-tensor-decomposition}.

\begin{algorithm}[H]
\caption{\textsc{ImplicitTensorDecomposition}}\label{alg:tensor-decomposition:implicit-tensor-decomposition}
\KwIn{The recovered list $\widehat{\mathcal R}$ and failure probability
$\delta$}
\KwOut{A non-empty list of unit vectors $(\widehat{\NormalizedRegressor}_1,\ldots,\widehat{\NormalizedRegressor}_{i_\star})$ giving rough estimates of normalized regressors}
\tcp{\Cref{alg:tensor-decomposition:estimate-projected-moments,thm:tensor-decomposition:projected-moment-accuracy}}
$(\widehat{A},\widehat{M}_1,\ldots,\widehat{M}_n)\leftarrow$
\ProjectedMomentFamily{$\widehat{\mathcal R},
e^{-K^6s},\delta/2$}\;
\tcp{\Cref{alg:tensor-decomposition:construct-observable-operators,thm:tensor-decomposition:observable-operator-construction}}
$(i_\star,\widehat{D}_1,\ldots,\widehat{D}_n)\leftarrow$
\SelectAndNormalize{$\widehat{A},(\widehat{M}_\ell)_{\ell\in[n]}$}\;
\If{$i_\star=1$}{
  $\widehat{\ScaledRegressor}_1\gets((\widehat{D}_\ell)_{11})_{\ell\in[n]}$\;
  \Return $(\widehat{\ScaledRegressor}_1/\norm{\widehat{\ScaledRegressor}_1}_2)$\;
}
\tcp{\Cref{alg:tensor-decomposition:recover-directions,thm:tensor-decomposition:normalized-regressor-recovery}}
$(\widehat{\NormalizedRegressor}_1,\ldots,\widehat{\NormalizedRegressor}_{i_\star})\leftarrow$
\DiagonalizeAndRead{$(\widehat{D}_\ell)_{\ell\in[n]},
\delta/2$}\;
\Return $(\widehat{\NormalizedRegressor}_1,\ldots,\widehat{\NormalizedRegressor}_{i_\star})$\;
\end{algorithm}

We use \Cref{thm:tensor-decomposition:projected-moment-accuracy} to estimate the projected
moments and \Cref{thm:tensor-decomposition:observable-operator-construction} to construct an observable-operator
family for a non-empty block. For a singleton block we read the coordinates
directly. For a larger block we apply \Cref{thm:tensor-decomposition:normalized-regressor-recovery}
to recover the normalized regressors.

\begin{proof}[Proof of \cref{thm:tensor-decomposition:implicit-tensor-decomposition}]
We verify that the conditions to the respective results are satisfied.

\emph{Projected moments and projection error.}
In order to apply \Cref{thm:tensor-decomposition:projected-moment-accuracy}, we must check
the accuracy of the recovered list and the inputs to
\textsc{EstimateProjectedMoments}.

By hypothesis, the recovered regressors have estimation error at most
$e^{-K^6s}$.
Hence \Cref{thm:tensor-decomposition:projected-moment-accuracy} applies.
Our algorithm calls \textsc{EstimateProjectedMoments} with
 sampling accuracy $e^{-K^6s}$ and failure probability
$\nicefrac\delta2$, using fresh observations.
As a result, with probability at least $1-\delta/2$, the following holds,
where $\varepsilon_{\mathrm{mean}}=e^{-K^2s}$:
\begin{align*}
 \norm{\widehat A-\ProjectionMap{s} \EvenMomentMatrix{s}{Q}
                    \ProjectionMap{s}^\top}_F
 &\leq\varepsilon_{\mathrm{mean}}\tr(\EvenMomentMatrix{s}{Q}),\\
 \norm{\sum_{\ell=1}^n \zeta_\ell\widehat M_\ell
       -\ProjectionMap{s} \OddMomentMatrix{s}{Q}{\bm{\zeta}}\ProjectionMap{s}^\top}_F
 &\leq\varepsilon_{\mathrm{mean}}
       \tr(\EvenMomentMatrix{s}{Q})\norm{\bm{\zeta}}_2,\\
 \norm{(I-\ProjectionMap{s}^\top\ProjectionMap{s})
       (\EvenMomentMatrix{s}{Q})^{1/2}}_{\op}^2
 &\leq2\varepsilon_{\mathrm{mean}}K^{s-1}
       \tr(\EvenMomentMatrix{s}{Q}).
\end{align*}
The odd-moment bound holds simultaneously for every $\bm{\zeta}\in\R^n$.
We proceed assuming fixed estimates that satisfy these bounds.

\emph{Constructing the observable-operator family.}
In order to apply \Cref{thm:tensor-decomposition:observable-operator-construction}, we must check the
moment and projection errors, approximate orthogonality of the lifted
vectors, the upper bound on the scaled-regressor norms, and the total
error bound. 

By construction, $\ProjectionMap{s}\in\R^{p_s\times n^s}$
has orthonormal rows and $1\leq p_s\leq |Q|$. The preceding estimates give
\begin{align*}
 \varepsilon_{\mathrm{mom}}
 &=\varepsilon_{\mathrm{mean}},\\
 \varepsilon_{\mathrm{proj}}
 &=2K^{s-1}\varepsilon_{\mathrm{mean}}.
\end{align*}
The approximate-orthogonality bound
(\Cref{thm:hermite:structure-approximate-orthogonality}) gives the required estimate with
$\varepsilon_\perp=k e^{-s/K}$:
\begin{align*}
 \norm{
 \bigl(\ip{\TensorPower{\NormalizedRegressor_i}{s}}{\TensorPower{\NormalizedRegressor_j}{s}}\bigr)_{i,j\in Q}-I
 }_{\op}
 &\leq\max_{i\in Q}\sum_{j\in Q\setminus\{i\}}
 \abs{\ip{\TensorPower{\NormalizedRegressor_i}{s}}{\TensorPower{\NormalizedRegressor_j}{s}}}
 \alignnote{maximum absolute row sum}\\
 &\leq(|Q|-1)e^{-s/K}\leq k e^{-s/K}=\varepsilon_\perp.
 \alignnote{by \Cref{thm:hermite:structure-approximate-orthogonality}}
\end{align*}
By mundane calculations deferred to the appendix
(\Cref{prop:app-hermite:scaled-regressor-geometry}),
\begin{equation}\label{eq:tensor-decomposition:scaled-regressor-geometry}
\begin{aligned}
 K^{-3}\leq\norm{\ScaledRegressor_i}_2&\leq K^3
 \qquad(i\in Q),\\
 \norm{\ScaledRegressor_i-\ScaledRegressor_j}_2&\geq K^{-3}
 \qquad(i,j\in Q,\ i\ne j).
\end{aligned}
\end{equation}
Similarly, by further calculations deferred to the appendix
(\Cref{prop:app-tensor-decomposition:observable-operator-accuracy}),
\[
 \varepsilon_{\mathrm{mom}}
 +\varepsilon_{\mathrm{proj}}
 +\varepsilon_\perp
 \leq\frac{e^{-s/(2K)}}{512K^9\sqrt n\,k^5}
 \leq\frac1{16k^2}.
\]
All hypotheses of \Cref{thm:tensor-decomposition:observable-operator-construction} are now
satisfied.
Applying the theorem gives a non-empty $J\subseteq Q$ and
returned matrices of dimension $i_\star=|J|$.
There is an orthonormal basis $(\bm{z}_i)_{i\in J}$ of
$\R^{i_\star}$ such that, for every $\bm{\zeta}\in\R^n$,
\begin{equation}\label{eq:tensor-decomposition:operator-family-approximation}
 \norm{\sum_{\ell=1}^n \zeta_\ell\widehat D_\ell
       -\sum_{i\in J}\ip{\bm{\zeta}}{\ScaledRegressor_i}\bm{z}_i\bm{z}_i^\top}_F
 \leq\varepsilon_{\mathrm{fam}}\norm{\bm{\zeta}}_2,
\end{equation}
where
\[
 \varepsilon_{\mathrm{fam}}
 \coloneqq16K^3k^3
 \left(\varepsilon_{\mathrm{mom}}
       +\varepsilon_{\mathrm{proj}}
       +\varepsilon_\perp\right).
\]

\emph{Recovering the normalized regressors.}
In order to apply \Cref{thm:tensor-decomposition:normalized-regressor-recovery} when $|J|\geq2$,
we must check the family approximation bound and the norm and
separation bounds for the scaled regressors. These are given by
\eqref{eq:tensor-decomposition:operator-family-approximation} and \eqref{eq:tensor-decomposition:scaled-regressor-geometry},
respectively. To obtain the direction accuracy in \Cref{thm:tensor-decomposition:implicit-tensor-decomposition},
it remains to bound the resulting error. By the choice of $K$ and
the hyperparameter calculation in \Cref{prop:app-tensor-decomposition:observable-operator-accuracy},
\begin{equation}\label{eq:tensor-decomposition:operator-family-accuracy}
 32K^6\sqrt n\,k^2\varepsilon_{\mathrm{fam}}
 \leq e^{-s/(2K)}
 \leq\frac{\Delta^2}{16B^2}.
\end{equation}
If $J=\{j\}$, the scalar coordinate matrices give
$\widehat{\ScaledRegressor}_1=((\widehat D_\ell)_{11})_{\ell\in[n]}$
with
$\norm{\widehat{\ScaledRegressor}_1-\ScaledRegressor_j}_2\leq\varepsilon_{\mathrm{fam}}
<K^{-3}\leq\norm{\ScaledRegressor_j}_2$
by \eqref{eq:tensor-decomposition:operator-family-approximation},
\eqref{eq:tensor-decomposition:scaled-regressor-geometry}, and \eqref{eq:tensor-decomposition:operator-family-accuracy}.
Normalizing then gives an estimate of the normalized regressor $\NormalizedRegressor_j$ with error at most
$2K^3\varepsilon_{\mathrm{fam}}$.
For $|J|\geq2$, applying \Cref{thm:tensor-decomposition:normalized-regressor-recovery} with failure
probability $\nicefrac\delta2$ gives unit vectors matched bijectively to the
normalized regressors indexed by $J$, with error at most
$32K^6\sqrt n\,k^2\varepsilon_{\mathrm{fam}}$.
Since $\ScaledRegressor_i/\norm{\ScaledRegressor_i}_2=\Regressor_i/\norm{\Regressor_i}_2$ by
\Cref{thm:hermite:structure-scaled-regressors}, both cases give a bijection
$\pi:[i_\star]\to J$ satisfying
\begin{align*}
 \max_{i\in[i_\star]}
 \norm{\widehat{\NormalizedRegressor}_i-
       \frac{\Regressor_{\pi(i)}}{\norm{\Regressor_{\pi(i)}}_2}}_2
 &\leq32K^6\sqrt n\,k^2\varepsilon_{\mathrm{fam}}
 \alignnote{\Cref{thm:tensor-decomposition:normalized-regressor-recovery} and singleton bound}\\
 &\leq e^{-s/(2K)}
 \leq\frac{\Delta^2}{16B^2}.
 \alignnote{by \eqref{eq:tensor-decomposition:operator-family-accuracy}}
\end{align*}
By a union bound over all failure events above, the total failure
probability is at most $\delta$.

\emph{Sample complexity and running time.}
By \Cref{thm:tensor-decomposition:projected-moment-accuracy},
\textsc{EstimateProjectedMoments} uses
\[
 \Ot\!\left((nk)^{O(K^{10})}
            \log(\nicefrac1\delta)\right)
\]
samples and runs in time
\[
 (nk)^{\poly(K)}\log(\nicefrac1\delta).
\]
The observable-operator routines use no further observations.
By \Cref{thm:tensor-decomposition:observable-operator-construction,thm:tensor-decomposition:normalized-regressor-recovery},
their total running time is
$\poly(n,k)\log(\nicefrac1\delta)$.
Combining these bounds proves the sample and time complexity
guarantees.
\end{proof}

\section{Restricted Least-Squares}
\label{sec:refinement:main}

In this section, we show how to refine the estimates of normalized regressors returned by \textsc{ImplicitTensorDecomposition} (\Cref{alg:tensor-decomposition:implicit-tensor-decomposition}) from \Cref{sec:tensor-decomposition:main} into accurate estimates of the full regressors and biases.
Our goal is to prove the result below.

\begin{theorem}[Restricted Least-Squares Refinement]\label{thm:refinement:restricted-least-squares}
Fix $\varepsilon\in(0,1)$, $\delta\in(0,\nicefrac13)$,
and $i\in[k]$.
Suppose the unit vector $\widehat{\NormalizedRegressor}$ satisfies
\[
 \norm{\widehat{\NormalizedRegressor}-\frac{\Regressor_i}{\norm{\Regressor_i}_2}}_2
 \leq\frac{\Delta^2}{16B^2}.
\]
Then \textsc{RestrictedLeastSquares} (\Cref{alg:refinement:restricted-least-squares}) returns estimates
$\widehat{\Regressor}$ and $\widehat b$ satisfying
\begin{equation}
 \norm{\widehat{\Regressor}-\Regressor_i}_2+|\widehat b-b_i|\leq\varepsilon
\end{equation}
with probability at least $1-\delta$.
Its sample complexity is
\[
 \Ot\!\left(
 \left(\frac{nk}{\varepsilon}\right)^{5K^2}
 \log(\nicefrac1\delta)
 \right)
\]
and its running time is
\[
 \Ot\!\left(
 \left(\frac{nk}{\varepsilon}\right)^{\poly(K)}
 \log(\nicefrac1\delta)
 \right)\,.
\]
\end{theorem}

Let $\RefinementEvent{\widehat{\NormalizedRegressor}}
=\{\Gamma\leq\ip{\widehat{\NormalizedRegressor}}{\Covariate}\leq\Gamma+1\}$
be the large-projection event, with $\Gamma = \Theta(\sqrt{\log(\nicefrac{n}{\eps})})$ chosen as in \Cref{alg:refinement:restricted-least-squares}.
When $\widehat{\NormalizedRegressor}$ is close to
$\Regressor_i/\norm{\Regressor_i}_2$, the response $Y_i$ is likely
to attain the maximum on this event. We therefore estimate the
least-squares solution restricted to this event. Its normal equations
are $\RestrictedSecondMoment{\widehat{\NormalizedRegressor}}
\begin{bsmallmatrix}\bm w\\b\end{bsmallmatrix}
=\RestrictedResponse{\widehat{\NormalizedRegressor}}$, where
\begin{equation}
\begin{aligned}
 \RestrictedSecondMoment{\widehat{\NormalizedRegressor}}
 &\coloneqq\E_{\Covariate}\!\left[
 \begin{bmatrix}\Covariate\\1\end{bmatrix}
 \begin{bmatrix}\Covariate^\top&1\end{bmatrix}
 \1_{\RefinementEvent{\widehat{\NormalizedRegressor}}}
 \right],\\
 \RestrictedResponse{\widehat{\NormalizedRegressor}}
 &\coloneqq\E_{\Covariate,\NoiseVector}\!\left[
 \begin{bmatrix}\Covariate\\1\end{bmatrix}Y_{\max}
 \1_{\RefinementEvent{\widehat{\NormalizedRegressor}}}
 \right].
\end{aligned}
\end{equation}
Our algorithm computes
$\RestrictedSecondMoment{\widehat{\NormalizedRegressor}}$
in closed form and applies its inverse to a trimmed-mean estimate of
$\RestrictedResponse{\widehat{\NormalizedRegressor}}$.
Its pseudocode can be found in \Cref{alg:refinement:restricted-least-squares}.

\begin{algorithm}[H]
\caption[\textsc{RestrictedLeastSquares}]{\textsc{RestrictedLeastSquares}}\label{alg:refinement:restricted-least-squares}
\KwIn{A unit vector $\widehat{\NormalizedRegressor}$, accuracy $\varepsilon\in(0,1)$,
and failure probability $\delta\in(0,\nicefrac13)$}
\KwOut{Estimates $\widehat{\Regressor}\in\R^n$ and $\widehat b\in\R$}
${\Gamma}\gets
 2K\sqrt{\log\!\left(\frac{n}{\varepsilon}\right)}$\;
$\Pr[\RefinementEvent{\widehat{\NormalizedRegressor}}]\gets\PhiN({\Gamma}+1)-\PhiN({\Gamma})$\;
$\nu_1\gets\phi(\Gamma)-\phi(\Gamma+1)$\;
$\nu_2\gets \Pr[\RefinementEvent{\widehat{\NormalizedRegressor}}]+{\Gamma}\phi({\Gamma})
 -({\Gamma}+1)\phi({\Gamma}+1)$\;
Compute
\[
 \RestrictedSecondMoment{\widehat{\NormalizedRegressor}}
 =\begin{pmatrix}
 \nu_2\widehat{\NormalizedRegressor}\widehat{\NormalizedRegressor}^\top
 +\Pr[\RefinementEvent{\widehat{\NormalizedRegressor}}]
  (I-\widehat{\NormalizedRegressor}\widehat{\NormalizedRegressor}^\top)
 &\nu_1\widehat{\NormalizedRegressor}\\
 \nu_1\widehat{\NormalizedRegressor}^\top
 &\Pr[\RefinementEvent{\widehat{\NormalizedRegressor}}]
 \end{pmatrix}
\]

$\mathsf v_{\refine}\gets96B^2\log k$\tcp*[r]{coordinate second-moment bound}
$\varepsilon'\gets
 \frac{\varepsilon\,
 \lambda_{\min}(\RestrictedSecondMoment{\widehat{\NormalizedRegressor}})}
 {2\sqrt{2(n+1)}}$\;
\tcp{\Cref{lem:coordinatewise-trimmed-mean}}
$\widehat{\bm c}_{\widehat{\NormalizedRegressor}}\gets$
\TrimmedMean{$
 \begin{bsmallmatrix}\Covariate\\1\end{bsmallmatrix}Y_{\max}\1_{\RefinementEvent{\widehat{\NormalizedRegressor}}};
 \mathsf v_{\refine},\varepsilon',\delta
$}\;
$\begin{bsmallmatrix}\widehat{\Regressor}\\\widehat b\end{bsmallmatrix}\gets\RestrictedSecondMoment{\widehat{\NormalizedRegressor}}^{-1}\widehat{\bm c}_{\widehat{\NormalizedRegressor}}$\;
\Return $(\widehat{\Regressor},\widehat b)$\;
\end{algorithm}

The rest of this section is dedicated to proving \Cref{thm:refinement:restricted-least-squares}.

\subsection{Proof of \texorpdfstring{\Cref*{thm:refinement:restricted-least-squares}}{Theorem}}

Before analyzing \Cref{alg:refinement:restricted-least-squares},
we bound the conditional second moments of
$\ip{\bm v}{\Covariate}+v_{n+1}$ and $Y_{\max}-Y_i$ to control the bias
from using the maximum response.
We also bound the second moment
of each coordinate of
$\begin{bsmallmatrix}\Covariate\\1\end{bsmallmatrix}
 Y_{\max}\1_{\RefinementEvent{\widehat{\NormalizedRegressor}}}$
to control the sampling error.

\begin{lemma}[Moment Bounds for Restricted Least Squares]
\label{lem:refinement:moment-bounds}
Consider the Gaussian self-selection model
(\Cref{def:gaussian-self-selection}) under the boundedness
assumption (\Cref{ass:boundedness-uncoveredness}).
Use $K$ from \Cref{hyp:algorithm}, and let ${\Gamma}\geq8B^2/\Delta^2$.
Fix a unit vector $\widehat{\NormalizedRegressor}$ and $i\in[k]$ such that
\[
 \ip{\Regressor_i-\Regressor_j}{\widehat{\NormalizedRegressor}}
 \geq\frac{3\Delta^2}{4B}
 \qquad(j\ne i).
\]
Then the following bounds hold.
\begin{enumerate}[label=(\roman*),ref=\thelemma(\roman*),leftmargin=2.2em]
\item\label[lemma]{lem:refinement:moment-bounds-covariate}
For every $\bm v\in\R^n$ and $v_{n+1}\in\R$ satisfying
$\norm{\bm v}_2^2+|v_{n+1}|^2=1$,
\[
 \E_{\Covariate}\!\left[
 \left(
 \begin{bmatrix}\bm v^\top&v_{n+1}\end{bmatrix}
 \begin{bmatrix}\Covariate\\1\end{bmatrix}
 \right)^2
 \,\middle|\,\RefinementEvent{\widehat{\NormalizedRegressor}}
 \right]
 \leq5\Gamma^2.
\]

\item\label[lemma]{lem:refinement:moment-bounds-response-difference}
The response difference satisfies
\[
 \E_{\Covariate,\NoiseVector}[(Y_{\max}-Y_i)^2
 \mid\RefinementEvent{\widehat{\NormalizedRegressor}}]
 \leq6B^2k e^{-2\Gamma^2/K}.
\]

\item\label[lemma]{lem:refinement:moment-bounds-coordinate}
For every $j\in[n]$,
\[
 \E_{\Covariate,\NoiseVector}[
 (X_jY_{\max}\1_{\RefinementEvent{\widehat{\NormalizedRegressor}}})^2]
 \leq96B^2\log k.
\]

\item\label[lemma]{lem:refinement:moment-bounds-intercept}
The response coordinate satisfies
\[
 \E_{\Covariate,\NoiseVector}[
 (Y_{\max}\1_{\RefinementEvent{\widehat{\NormalizedRegressor}}})^2]
 \leq96B^2\log k.
\]
\end{enumerate}
\end{lemma}

\begin{proof}
We first bound the conditional covariate and response moments,
then bound the coordinate second moments used by
\textsc{TrimmedMean}.
By the choice of \Cref{hyp:algorithm},
\begin{equation}\label{eq:refinement:response-gap-k-choice}
 K\Delta^4
 \geq16(1+B^2)^3
 \geq16B^2(4B^2+2).
\end{equation}
The last inequality uses
$(1+B^2)^3-B^2(4B^2+2)=B^4(B^2-1)+B^2+1\geq0$.

\emph{Conditional covariate moment.}
Fix $\bm v\in\R^n$ and $v_{n+1}\in\R$ satisfying
$\norm{\bm v}_2^2+|v_{n+1}|^2=1$.
For $y\in[\Gamma,\Gamma+1]$, the conditional Gaussian formula
(\Cref{prop:gaussian-conditioning}) gives
\begin{align*}
 \E_{\Covariate}[\ip{\bm v}{\Covariate}^2
 \mid\ip{\widehat{\NormalizedRegressor}}{\Covariate}=y]
 &=y^2\ip{\bm v}{\widehat{\NormalizedRegressor}}^2
 +\norm{(I-\widehat{\NormalizedRegressor}
            \widehat{\NormalizedRegressor}^\top)\bm v}_2^2\\
 &=\norm{\bm v}_2^2
 +(y^2-1)\ip{\bm v}{\widehat{\NormalizedRegressor}}^2\\
 &\leq y^2\norm{\bm v}_2^2
 \alignnote{Cauchy--Schwarz and $y\geq1$}\\
 &\leq(\Gamma+1)^2\norm{\bm v}_2^2
 \leq4\Gamma^2\norm{\bm v}_2^2.
\end{align*}
Averaging over $y$ conditional on
$\RefinementEvent{\widehat{\NormalizedRegressor}}$
and applying the triangle inequality gives
\begin{align*}
 &\left(
 \E_{\Covariate}\!\left[
 (\ip{\bm v}{\Covariate}+v_{n+1})^2
 \,\middle|\,\RefinementEvent{\widehat{\NormalizedRegressor}}\right]
 \right)^{1/2}\\
 &\quad\leq
 \left(
 \E_{\Covariate}\!\left[
 \ip{\bm v}{\Covariate}^2
 \,\middle|\,\RefinementEvent{\widehat{\NormalizedRegressor}}\right]
 \right)^{1/2}
 +|v_{n+1}|\\
 &\quad\leq2\Gamma\norm{\bm v}_2+|v_{n+1}|\\
 &\quad\leq
 \sqrt{4\Gamma^2+1}
 \sqrt{\norm{\bm v}_2^2+|v_{n+1}|^2}
 \alignnote{Cauchy--Schwarz}\\
 &\quad=\sqrt{4\Gamma^2+1}
 \leq\sqrt5\,\Gamma.
 \alignnote{normalization and $\Gamma\geq1$}
\end{align*}
Squaring proves \Cref{lem:refinement:moment-bounds-covariate}.

\emph{Conditional response difference.}
The truncated Gaussian second-moment bound
(\Cref{prop:truncated-gaussian-second-moment}) gives
\[
 \E_G[G^2\1_{\{G>0\}}]
 \leq\operatorname{Var}_G(G)\Prb_G(G>0)
\]
for every nondegenerate Gaussian $G$ with $\E_G[G]\leq0$.
We apply this bound to the conditional response differences.

Fix $j\ne i$ and $y\in[{\Gamma},{\Gamma}+1]$.
Applying the conditional Gaussian formula
(\Cref{prop:gaussian-conditioning}) to
\[
 Y_j-Y_i-(b_j-b_i)
 =\ip{\Regressor_j-\Regressor_i}{\Covariate}+\xi_j-\xi_i
\]
gives
\begin{equation}
 Y_j-Y_i\mid\ip{\widehat{\NormalizedRegressor}}{\Covariate}=y
 \sim N\!\left(
 y\ip{\Regressor_j-\Regressor_i}{\widehat{\NormalizedRegressor}}+b_j-b_i,
 \norm{(I-\widehat{\NormalizedRegressor}\widehat{\NormalizedRegressor}^\top)(\Regressor_j-\Regressor_i)}_2^2+2(1-\rho)
 \right).
\end{equation}
Here the noises are independent of $\Covariate$, and
$\operatorname{Var}(\xi_j-\xi_i)=2(1-\rho)\leq2$.
The assumed directional gap and boundedness give
\begin{align}
 \E_{\Covariate,\NoiseVector}[Y_j-Y_i\mid\ip{\widehat{\NormalizedRegressor}}{\Covariate}=y]
 &=y\ip{\Regressor_j-\Regressor_i}{\widehat{\NormalizedRegressor}}+b_j-b_i
 \notag\\
 &\leq-\frac{3y\Delta^2}{4B}+2B
 \notag\\
 &\leq-\frac{y\Delta^2}{2B},
 \alignnote{$y\geq\Gamma\geq8B^2/\Delta^2$}
 \label{eq:refinement:response-difference-mean}\\
 \operatorname{Var}_{\Covariate,\NoiseVector}
 (Y_j-Y_i\mid\ip{\widehat{\NormalizedRegressor}}{\Covariate}=y)
 &\leq\norm{\Regressor_j-\Regressor_i}_2^2+2
 \notag\\
 &\leq4B^2+2.
 \label{eq:refinement:response-difference-variance}
\end{align}
If the conditional variance is zero, then, almost surely
conditional on $\ip{\widehat{\NormalizedRegressor}}{\Covariate}=y$,
\begin{align*}
 &Y_j-Y_i
 =y\ip{\Regressor_j-\Regressor_i}{\widehat{\NormalizedRegressor}}+b_j-b_i
 \leq-\frac{y\Delta^2}{2B}<0,
 \alignnote{by \eqref{eq:refinement:response-difference-mean}}\\
 &\E_{\Covariate,\NoiseVector}\!\left[
 (Y_j-Y_i)^2\1_{\{Y_j>Y_i\}}
 \,\middle|\,\ip{\widehat{\NormalizedRegressor}}{\Covariate}=y
 \right]=0.
\end{align*}
We apply the nondegenerate Gaussian bounds below to the remaining terms.
Since $Y_{\max}=\max_{j\in[k]}Y_j$,
\[
 (Y_{\max}-Y_i)^2
 \leq\sum_{j\ne i}
 (Y_j-Y_i)^2\1_{\{Y_j>Y_i\}}.
\]
The truncated Gaussian second-moment bound and a Gaussian
tail bound (\Cref{prop:gaussian-tail}) therefore give
\begin{align*}
 &\E_{\Covariate,\NoiseVector}[(Y_{\max}-Y_i)^2\mid\ip{\widehat{\NormalizedRegressor}}{\Covariate}=y]\\
 &\quad\leq\sum_{j\ne i}
 \E_{\Covariate,\NoiseVector}\left[
 (Y_j-Y_i)^2\1_{\{Y_j>Y_i\}}
 \,\middle|\,\ip{\widehat{\NormalizedRegressor}}{\Covariate}=y
 \right]\\
 &\quad\leq\sum_{j\ne i}
 \operatorname{Var}_{\Covariate,\NoiseVector}
 (Y_j-Y_i\mid\ip{\widehat{\NormalizedRegressor}}{\Covariate}=y)\,
 \Prb_{\Covariate,\NoiseVector}(Y_j>Y_i\mid\ip{\widehat{\NormalizedRegressor}}{\Covariate}=y)
 \alignnote{by \Cref{prop:truncated-gaussian-second-moment}}\\
 &\quad\leq(4B^2+2)\sum_{j\ne i}
 \Prb_{\Covariate,\NoiseVector}(Y_j>Y_i\mid\ip{\widehat{\NormalizedRegressor}}{\Covariate}=y)
 \alignnote{by \eqref{eq:refinement:response-difference-variance}}\\
 &\quad\leq(4B^2+2)k
 \exp\!\left(
 -\frac{y^2\Delta^4}{8B^2(4B^2+2)}
 \right)
 \alignnote{Gaussian tail bound; by \eqref{eq:refinement:response-difference-mean},
 \eqref{eq:refinement:response-difference-variance}}\\
 &\quad\leq(4B^2+2)k e^{-2y^2/K}
 \alignnote{choice of $K$~\eqref{eq:refinement:response-gap-k-choice}}\\
 &\quad\leq6B^2k e^{-2{\Gamma}^2/K}.
 \alignnote{$B\geq1$; $y\geq{\Gamma}$}
\end{align*}
Averaging over $y\in[{\Gamma},{\Gamma}+1]$ conditional on
$\RefinementEvent{\widehat{\NormalizedRegressor}}$ proves the second bound.

\emph{Coordinate second moment.}
It remains to control the second moment of each coordinate of
$\begin{bsmallmatrix}\Covariate\\1\end{bsmallmatrix}
Y_{\max}\1_{\RefinementEvent{\widehat{\NormalizedRegressor}}}$ used by
\textsc{TrimmedMean}.
The responses $Y_i$ are Gaussian, with means $b_i$.
The Gaussian maximum bound gives
\begin{align*}
 \left(\E_{\Covariate,\NoiseVector}
  \max_{i\in[k]}|Y_i|^4\right)^{1/4}
 &\leq\max_{i\in[k]}|\E_{\Covariate,\NoiseVector}[Y_i]|
 +3\sqrt{\log k}
 \left(\max_{i\in[k]}\operatorname{Var}_{\Covariate,\NoiseVector}(Y_i)\right)^{1/2}
 \alignnote{by \Cref{prop:gaussian-maximum-fourth-moment}}\\
 &\leq B+3\sqrt{1+B^2}\sqrt{\log k}
 \alignnote{boundedness (\Cref{ass:boundedness-uncoveredness})}\\
 &\leq6B\sqrt{\log k}.
 \alignnote{$B\geq1$, $k\geq2$}
\end{align*}
Since $|Y_{\max}|\leq\max_{i\in[k]}|Y_i|$,
Cauchy--Schwarz therefore gives
\begin{align*}
 \E_{\Covariate,\NoiseVector}[(X_jY_{\max}\1_{\RefinementEvent{\widehat{\NormalizedRegressor}}})^2]
 &\leq\sqrt{\E_{\Covariate}[X_j^4]}
 \left(\E_{\Covariate,\NoiseVector}\max_{i\in[k]}|Y_i|^4\right)^{1/2}\\
 &\leq36\sqrt3 B^2\log k
 \alignnote{standard Gaussian fourth moment}\\
 &\leq96B^2\log k.
\end{align*}
The same fourth-moment bound gives
\[
 \E_{\Covariate,\NoiseVector}[(Y_{\max}\1_{\RefinementEvent{\widehat{\NormalizedRegressor}}})^2]
 \leq36B^2\log k\leq96B^2\log k.
\]
This proves the remaining bounds.
\end{proof}

{We next lower bound the restricted second moment to control
how much solving the normal equations can amplify the estimation error.

\begin{lemma}[Lower Bound on the Restricted Second Moment]
\label{lem:refinement:second-moment-lower-bound}
Let $\Covariate\sim N(\bm0,I_n)$, let
$\widehat{\NormalizedRegressor}\in\R^n$ be a unit vector,
and let $\Gamma\geq1$.
For every $\bm v\in\R^n$ and $v_{n+1}\in\R$ satisfying
$\norm{\bm v}_2^2+|v_{n+1}|^2=1$,
\[
 \begin{bmatrix}\bm v\\v_{n+1}\end{bmatrix}^{\!\top}
 \RestrictedSecondMoment{\widehat{\NormalizedRegressor}}
 \begin{bmatrix}\bm v\\v_{n+1}\end{bmatrix}
 \geq
 \frac{\Pr[\RefinementEvent{\widehat{\NormalizedRegressor}}]}
      {1024\Gamma^4}.
\]
\end{lemma}

\begin{proof}
Put $p=\Pr[\RefinementEvent{\widehat{\NormalizedRegressor}}]$,
and fix $\bm v,v_{n+1}$ satisfying
$\norm{\bm v}_2^2+|v_{n+1}|^2=1$.
The conditional Gaussian formula
(\Cref{prop:gaussian-conditioning}) gives
\begin{align*}
 \E_{\Covariate}\!\left[
 \ip{\bm v}{\Covariate}
 \,\middle|\,
 \ip{\widehat{\NormalizedRegressor}}{\Covariate}=y\right]
 &=y\ip{\bm v}{\widehat{\NormalizedRegressor}},\\
 \operatorname{Var}_{\Covariate}\!\left(
 \ip{\bm v}{\Covariate}
 \,\middle|\,
 \ip{\widehat{\NormalizedRegressor}}{\Covariate}=y\right)
 &=\norm{\bm v}_2^2
   -\ip{\bm v}{\widehat{\NormalizedRegressor}}^2.
\end{align*}
The second-moment identity therefore gives
\[
 \E_{\Covariate}\!\left[
 (\ip{\bm v}{\Covariate}+v_{n+1})^2
 \,\middle|\,
 \ip{\widehat{\NormalizedRegressor}}{\Covariate}=y\right]\\
 =
 \norm{\bm v}_2^2-\ip{\bm v}{\widehat{\NormalizedRegressor}}^2
 +\left(y\ip{\bm v}{\widehat{\NormalizedRegressor}}+v_{n+1}\right)^2.
\]
Starting from the definition, we obtain
\begin{align}
 \begin{bmatrix}\bm v\\v_{n+1}\end{bmatrix}^{\!\top}
 \RestrictedSecondMoment{\widehat{\NormalizedRegressor}}
 \begin{bmatrix}\bm v\\v_{n+1}\end{bmatrix}
 &=\E_{\Covariate}\!\left[
 (\ip{\bm v}{\Covariate}+v_{n+1})^2
 \1_{\RefinementEvent{\widehat{\NormalizedRegressor}}}\right]
 \notag\\
 &=\int_\Gamma^{\Gamma+1}
 \E_{\Covariate}\!\left[
 (\ip{\bm v}{\Covariate}+v_{n+1})^2
 \,\middle|\,
 \ip{\widehat{\NormalizedRegressor}}{\Covariate}=y\right]
 \phi(y)\odif{y}
 \notag\\
 &=p\left(
 \norm{\bm v}_2^2-\ip{\bm v}{\widehat{\NormalizedRegressor}}^2
 \right)
 +
 \int_\Gamma^{\Gamma+1}
 \left(y\ip{\bm v}{\widehat{\NormalizedRegressor}}+v_{n+1}\right)^2
 \phi(y)\odif{y}.
 \label{eq:refinement:restricted-quadratic-form}
\end{align}

\emph{Case I: $\norm{\bm v}_2\geq1/(4\Gamma)$.}
We lower bound the integral on $[\Gamma,\Gamma+\Gamma^{-1}]$.
Since
\[
 p\leq\int_\Gamma^\infty\phi(y)\odif{y}
 \leq\frac1\Gamma\int_\Gamma^\infty y\phi(y)\odif{y}
 =\frac{\phi(\Gamma)}\Gamma,
\]
for $y\in[\Gamma,\Gamma+\Gamma^{-1}]$ we have
\begin{equation}\label{eq:refinement:density-lower-bound}
 \phi(y)
 =\phi(\Gamma)e^{-(y^2-\Gamma^2)/2}
 \geq e^{-3/2}\phi(\Gamma)
 \geq p\Gamma e^{-3/2}.
\end{equation}
Centering the interval at its midpoint via a change of variables gives
\begin{align*}
 &\int_\Gamma^{\Gamma+\Gamma^{-1}}
 \left(y\ip{\bm v}{\widehat{\NormalizedRegressor}}+v_{n+1}\right)^2
 \odif{y}\\
 &\quad=\int_{-1/(2\Gamma)}^{1/(2\Gamma)}
 \left(
 t\ip{\bm v}{\widehat{\NormalizedRegressor}}
 +v_{n+1}
 +\left(\Gamma+\frac1{2\Gamma}\right)
  \ip{\bm v}{\widehat{\NormalizedRegressor}}
 \right)^2\odif{t}\\
 &\quad=\frac1\Gamma
 \left(
 v_{n+1}
 +\left(\Gamma+\frac1{2\Gamma}\right)
  \ip{\bm v}{\widehat{\NormalizedRegressor}}
 \right)^2
 +\frac{\ip{\bm v}{\widehat{\NormalizedRegressor}}^2}{12\Gamma^3}\\
 &\quad\geq
 \frac{\ip{\bm v}{\widehat{\NormalizedRegressor}}^2}{12\Gamma^3}.
\end{align*}
Substituting into \eqref{eq:refinement:restricted-quadratic-form}
and using \eqref{eq:refinement:density-lower-bound}, we obtain
\begin{align*}
 &\begin{bmatrix}\bm v\\v_{n+1}\end{bmatrix}^{\!\top}
 \RestrictedSecondMoment{\widehat{\NormalizedRegressor}}
 \begin{bmatrix}\bm v\\v_{n+1}\end{bmatrix}\\
 &\quad\geq
 \frac{p e^{-3/2}}{12\Gamma^2}\left(\norm{\bm v}_2^2
       -\ip{\bm v}{\widehat{\NormalizedRegressor}}^2\right)
 +\frac{pe^{-3/2}}{12\Gamma^2}
  \ip{\bm v}{\widehat{\NormalizedRegressor}}^2
  \alignnote{$\norm{\bm v}_2^2 - \ip{\bm v}{\widehat\NormalizedRegressor}\geq 0$}\\
 &\quad\geq\frac{p}{64\Gamma^2}\norm{\bm v}_2^2
 \alignnote{$12e^{3/2}\leq64$, $\Gamma\geq1$}\\
 &\quad\geq\frac{p}{1024\Gamma^4}.
 \alignnote{$\norm{\bm v}_2\geq1/(4\Gamma)$}
\end{align*}

\emph{Case II: $\norm{\bm v}_2<1/(4\Gamma)$.}
The normalization gives
\[
 |v_{n+1}|
 =\sqrt{1-\norm{\bm v}_2^2}
 \geq\sqrt{1-\frac1{16\Gamma^2}}
 \geq\frac34.
\]
For every $y\in[\Gamma,\Gamma+1]$,
\begin{align*}
 \left|y\ip{\bm v}{\widehat{\NormalizedRegressor}}+v_{n+1}\right|
 &\geq |v_{n+1}|
       -y\left|\ip{\bm v}{\widehat{\NormalizedRegressor}}\right|\\
 &\geq\frac34-(\Gamma+1)\norm{\bm v}_2\\
 &\geq\frac34-\frac{\Gamma+1}{4\Gamma}
 \geq\frac14.
\end{align*}
By \eqref{eq:refinement:restricted-quadratic-form},
\begin{align*}
 &\begin{bmatrix}\bm v\\v_{n+1}\end{bmatrix}^{\!\top}
 \RestrictedSecondMoment{\widehat{\NormalizedRegressor}}
 \begin{bmatrix}\bm v\\v_{n+1}\end{bmatrix}\\
 &\quad\geq
 \int_\Gamma^{\Gamma+1}
 \left(y\ip{\bm v}{\widehat{\NormalizedRegressor}}+v_{n+1}\right)^2
 \phi(y)\odif{y}\\
 &\quad\geq\frac{p}{16}
 \geq\frac{p}{1024\Gamma^4}.
\end{align*}
This proves the lemma.
\end{proof}
}

We are now ready to prove the guarantees of our \textsc{RestrictedLeastSquares} algorithm.
\begin{proof}[Proof of \cref{thm:refinement:restricted-least-squares}]
We first compute the restricted second moment and use
\Cref{lem:refinement:second-moment-lower-bound} to bound the operator norm of its inverse.
This controls blowup in the bias from replacing $Y_i$ by $Y_{\max}$ and the sampling error in the restricted response when solving the normal equations.

{\emph{Restricted second moment.}
To compute the second restricted moment in closed form,
we first evaluate the following scalar quantities:
\begin{align}
 \Pr[\RefinementEvent{\widehat{\NormalizedRegressor}}]
 &=\PhiN({\Gamma}+1)-\PhiN({\Gamma}),
 \notag\\
 \nu_1
 &\coloneqq\int_\Gamma^{\Gamma+1}y\phi(y)\odif{y}
 =\phi(\Gamma)-\phi(\Gamma+1),
 \alignnote{$\phi'(y)=-y\phi(y)$}\notag\\
 \nu_2
 &\coloneqq\int_{\Gamma}^{{\Gamma}+1}y^2\phi(y) \odif{y}
 \notag\\
 &=-\int_{\Gamma}^{{\Gamma}+1}y\phi'(y) \odif{y}
 \alignnote{$\phi'(y)=-y\phi(y)$}\notag\\
 &=-[y\phi(y)]_{\Gamma}^{{\Gamma}+1}
   +\int_{\Gamma}^{{\Gamma}+1}\phi(y) \odif{y}
 \alignnote{integration by parts}\notag\\
 &=\Pr[\RefinementEvent{\widehat{\NormalizedRegressor}}]
   +{\Gamma}\phi({\Gamma})-({\Gamma}+1)\phi({\Gamma}+1).
\end{align}
Within this proof, put $p=\Pr[\RefinementEvent{\widehat{\NormalizedRegressor}}]$.
Starting from the definition, we have
\begin{align}
 \RestrictedSecondMoment{\widehat{\NormalizedRegressor}}
 &=\E_{\Covariate}\!\left[
 \begin{bmatrix}\Covariate\\1\end{bmatrix}
 \begin{bmatrix}\Covariate^\top&1\end{bmatrix}
 \1_{\RefinementEvent{\widehat{\NormalizedRegressor}}}\right]
 \notag\\
 &=\int_\Gamma^{\Gamma+1}
 \E_{\Covariate}\!\left[
 \begin{bmatrix}\Covariate\\1\end{bmatrix}
 \begin{bmatrix}\Covariate^\top&1\end{bmatrix}
 \,\middle|\,
 \ip{\widehat{\NormalizedRegressor}}{\Covariate}=y
 \right]\phi(y)\odif{y}
 \notag\\
 &=\int_\Gamma^{\Gamma+1}
 \begin{pmatrix}
 y^2\widehat{\NormalizedRegressor}\widehat{\NormalizedRegressor}^\top
 +I-\widehat{\NormalizedRegressor}\widehat{\NormalizedRegressor}^\top
 &y\widehat{\NormalizedRegressor}\\
 y\widehat{\NormalizedRegressor}^\top&1
 \end{pmatrix}\phi(y)\odif{y}
 \alignnote{by \Cref{prop:gaussian-conditioning}}\notag\\
 &=\begin{pmatrix}
 \nu_2\widehat{\NormalizedRegressor}\widehat{\NormalizedRegressor}^\top
 +p(I-\widehat{\NormalizedRegressor}\widehat{\NormalizedRegressor}^\top)
 &\nu_1\widehat{\NormalizedRegressor}\\
 \nu_1\widehat{\NormalizedRegressor}^\top&p
 \end{pmatrix}.
\end{align}
Moreover, by \Cref{lem:refinement:second-moment-lower-bound},
\begin{equation}\label{eq:refinement:inverse-norm-bound}
 \norm{\RestrictedSecondMoment{\widehat{\NormalizedRegressor}}^{-1}}_{\op}
 \leq\frac{1024\Gamma^4}{p}.
\end{equation}

}\emph{Error decomposition.}
Since the noise $\xi_i$ is centered and independent of $\Covariate$,
\[
 \E_{\Covariate,\xi_i}\!\left[
 \begin{bmatrix}\Covariate\\1\end{bmatrix}Y_i
 \1_{\RefinementEvent{\widehat{\NormalizedRegressor}}}\right]
 =\RestrictedSecondMoment{\widehat{\NormalizedRegressor}}
 \begin{bmatrix}\Regressor_i\\b_i\end{bmatrix}.
\]
Thus the bias in the restricted response is
\begin{equation}
\begin{aligned}
 \RefinementBias{\widehat{\NormalizedRegressor}}
 &\coloneqq\RestrictedResponse{\widehat{\NormalizedRegressor}}
 -\RestrictedSecondMoment{\widehat{\NormalizedRegressor}}
 \begin{bmatrix}\Regressor_i\\b_i\end{bmatrix}\\
 &=\E_{\Covariate,\NoiseVector}\!\left[
 \begin{bmatrix}\Covariate\\1\end{bmatrix}(Y_{\max}-Y_i)
 \1_{\RefinementEvent{\widehat{\NormalizedRegressor}}}\right].
\end{aligned}
\end{equation}
The integrand vanishes whenever $Y_{\max}=Y_i$.
Since
$\begin{bsmallmatrix}\widehat{\Regressor}\\\widehat b\end{bsmallmatrix}
=\RestrictedSecondMoment{\widehat{\NormalizedRegressor}}^{-1}
\widehat{\bm c}_{\widehat{\NormalizedRegressor}}$,
adding and subtracting $\RestrictedResponse{\widehat{\NormalizedRegressor}}$ gives
\begin{align}
 &\norm{\widehat{\Regressor}-\Regressor_i}_2+|\widehat b-b_i|
 \notag\\
 &\quad\leq\sqrt2\norm{
 \begin{bmatrix}\widehat{\Regressor}-\Regressor_i\\\widehat b-b_i\end{bmatrix}}_2
 \alignnote{Cauchy--Schwarz}\notag\\
 &\quad=\sqrt2\norm{\RestrictedSecondMoment{\widehat{\NormalizedRegressor}}^{-1}
 \bigl(\RefinementBias{\widehat{\NormalizedRegressor}}
 +\widehat{\bm c}_{\widehat{\NormalizedRegressor}}
 -\RestrictedResponse{\widehat{\NormalizedRegressor}}\bigr)}_2
 \notag\\
 &\quad\leq
 \underbrace{\sqrt2\norm{\RestrictedSecondMoment{\widehat{\NormalizedRegressor}}^{-1}}_{\op}
 \norm{\RefinementBias{\widehat{\NormalizedRegressor}}}_2}_{\text{bias}}
 \notag\\
 &\qquad+
 \underbrace{\sqrt2\norm{\RestrictedSecondMoment{\widehat{\NormalizedRegressor}}^{-1}}_{\op}
 \norm{\widehat{\bm c}_{\widehat{\NormalizedRegressor}}
 -\RestrictedResponse{\widehat{\NormalizedRegressor}}}_2}_{\text{sampling error}}.
 \alignnote{triangle inequality}
 \label{eq:refinement:error-decomposition}
\end{align}
We will show that each term is at most $\nicefrac{\varepsilon}2$.

\emph{Conditional mean gap.}
We first verify the hypothesis of
\Cref{lem:refinement:moment-bounds}.
Recall $\NormalizedRegressor_i=\Regressor_i/\norm{\Regressor_i}_2$.
Boundedness and uncoveredness (\Cref{ass:boundedness-uncoveredness}) give
\begin{align}
 \ip{\Regressor_i-\Regressor_j}{\widehat{\NormalizedRegressor}}
 &=\ip{\Regressor_i-\Regressor_j}{\NormalizedRegressor_i}
   +\ip{\Regressor_i-\Regressor_j}{\widehat{\NormalizedRegressor}-\NormalizedRegressor_i}
 \notag\\
 &\geq\ip{\Regressor_i-\Regressor_j}{\NormalizedRegressor_i}
       -\norm{\Regressor_i-\Regressor_j}_2\norm{\widehat{\NormalizedRegressor}-\NormalizedRegressor_i}_2
 \alignnote{Cauchy--Schwarz}\notag\\
 &\geq\frac{\Delta^2}{B}
       -2B\frac{\Delta^2}{16B^2}
 \alignnote{by \Cref{ass:boundedness-uncoveredness}; estimation assumption}\notag\\
 &=\frac{7\Delta^2}{8B}
 \geq\frac{3\Delta^2}{4B}
 \qquad(j\ne i).
\end{align}
The choice of $\Gamma$ also gives
$\Gamma\geq K\geq8B^2/\Delta^2$.

\emph{Bias.}
We apply Cauchy--Schwarz to the conditional expectation.
Fix $\bm v\in\R^n$ and $v_{n+1}\in\R$ satisfying
$\norm{\bm v}_2^2+|v_{n+1}|^2=1$.
Then
\begin{align*}
 &\frac1p\left|
 \begin{bmatrix}\bm v\\v_{n+1}\end{bmatrix}^{\!\top}
 \RefinementBias{\widehat{\NormalizedRegressor}}\right|\\
 &\quad=\left|\E_{\Covariate,\NoiseVector}\!\left[
 (\ip{\bm v}{\Covariate}+v_{n+1})(Y_{\max}-Y_i)
 \,\middle|\,\RefinementEvent{\widehat{\NormalizedRegressor}}\right]\right|
 \alignnote{definition of $\RefinementBias{\widehat{\NormalizedRegressor}}$}\\
 &\quad\leq\sqrt5\,\Gamma
 \left(\E_{\Covariate,\NoiseVector}[(Y_{\max}-Y_i)^2
 \mid\RefinementEvent{\widehat{\NormalizedRegressor}}]\right)^{1/2}
 \alignnote{by \Cref{lem:refinement:moment-bounds-covariate}}\\
 &\quad\leq\sqrt{30}\,B\Gamma\sqrt k\,e^{-\Gamma^2/K}.
 \alignnote{by \Cref{lem:refinement:moment-bounds-response-difference}}
\end{align*}
Taking the supremum over
$\norm{\bm v}_2^2+|v_{n+1}|^2=1$
and multiplying by the bound on
$\norm{\RestrictedSecondMoment{\widehat{\NormalizedRegressor}}^{-1}}_{\op}$
in \eqref{eq:refinement:inverse-norm-bound} gives
\begin{align}
 \sqrt2\norm{\RestrictedSecondMoment{\widehat{\NormalizedRegressor}}^{-1}}_{\op}
 \norm{\RefinementBias{\widehat{\NormalizedRegressor}}}_2
 &\leq1024\sqrt{60}\,B\Gamma^5\sqrt k\,e^{-\Gamma^2/K}
 \alignnote{by \eqref{eq:refinement:inverse-norm-bound}}\notag\\
 &\leq2^{13}B\Gamma^5\sqrt k\,e^{-\Gamma^2/K}
 \alignnote{$\sqrt{60}\leq8$}\notag\\
 &\leq\frac\varepsilon2.
\end{align}
The last inequality follows by simplifying using $K$
from \Cref{hyp:algorithm} and the choice of $\Gamma$
(\Cref{prop:app-refinement:refinement-bounds}).

\emph{Sampling error.}
The bias is now at most $\nicefrac\varepsilon2$.
To control the sampling term in
\eqref{eq:refinement:error-decomposition}, we apply the trimmed-mean
guarantee (\Cref{lem:coordinatewise-trimmed-mean}) to
$\begin{bsmallmatrix}\Covariate\\1\end{bsmallmatrix}
Y_{\max}\1_{\RefinementEvent{\widehat{\NormalizedRegressor}}}$ with coordinate accuracy
\[
 \varepsilon'
 =\frac{\varepsilon\,
 \lambda_{\min}(\RestrictedSecondMoment{\widehat{\NormalizedRegressor}})}
 {2\sqrt{2(n+1)}}.
\]
By
\Cref{lem:refinement:moment-bounds-coordinate,lem:refinement:moment-bounds-intercept},
the second moment of each of these $n+1$ coordinates is at most
$96B^2\log k=\mathsf v_{\refine}$.
The requested accuracy is admissible because
\[
 0<\varepsilon'\leq\frac12\leq\sqrt{\mathsf v_{\refine}},
\]
since
$0<\lambda_{\min}(\RestrictedSecondMoment{\widehat{\NormalizedRegressor}})
\leq p\leq1$, by evaluating the quadratic form at
$\begin{bsmallmatrix}\bm0\\1\end{bsmallmatrix}$.
Thus, with probability at least $1-\delta$,
\begin{align}
 \norm{\widehat{\bm c}_{\widehat{\NormalizedRegressor}}
 -\RestrictedResponse{\widehat{\NormalizedRegressor}}}_2
 &\leq\sqrt{n+1}\,\varepsilon',
 \notag\\
 \sqrt2\norm{\RestrictedSecondMoment{\widehat{\NormalizedRegressor}}^{-1}}_{\op}
 \norm{\widehat{\bm c}_{\widehat{\NormalizedRegressor}}
 -\RestrictedResponse{\widehat{\NormalizedRegressor}}}_2
 &\leq\sqrt{2(n+1)}
 \norm{\RestrictedSecondMoment{\widehat{\NormalizedRegressor}}^{-1}}_{\op}
 \varepsilon'
 \notag\\
 &=\frac\varepsilon2.
\end{align}
Combining these bounds in \eqref{eq:refinement:error-decomposition}
proves $\norm{\widehat{\Regressor}-\Regressor_i}_2+|\widehat b-b_i|\leq\varepsilon$.

\emph{Sample complexity and running time.}
To bound the sample complexity, we integrate the density bound
\eqref{eq:refinement:density-lower-bound} over
$[\Gamma,\Gamma+\Gamma^{-1}]$. This gives
\begin{align}
 \Pr[\RefinementEvent{\widehat{\NormalizedRegressor}}]
 &=\int_{\Gamma}^{{\Gamma}+1}\phi(y) \odif{y}
 \notag\\
 &\geq\int_0^{{\Gamma}^{-1}}\phi({\Gamma}+h) \odif{h}
 \alignnote{${\Gamma}^{-1}\leq1$}\notag\\
 &\geq\frac{e^{-3/2}\phi({\Gamma})}{{\Gamma}}
 \notag\\
 &\geq\frac{e^{-{\Gamma}^2/2}}{16{\Gamma}}.
 \alignnote{$e^{3/2}\sqrt{2\pi}\leq16$}
 \label{eq:refinement:large-projection-probability}
\end{align}
Substituting the inverse bound \eqref{eq:refinement:inverse-norm-bound} and
\eqref{eq:refinement:large-projection-probability} into the trimmed-mean sample bound
(\Cref{lem:coordinatewise-trimmed-mean}) gives
\[
 O\!\left(
 \frac{nB^2\log k}{\varepsilon^2}
 {\Gamma}^{10}e^{{\Gamma}^2}
 \log\!\left(\frac n\delta\right)
 \right)
\]
observations.
By simplifying using $K$ from \Cref{hyp:algorithm}
and the choice of $\Gamma$ (\Cref{prop:app-refinement:refinement-bounds}),
the trimmed-mean estimate uses
\[
 \Ot\!\left(
 \left(\frac{nk}{\varepsilon}\right)^{5K^2}
 \log(\nicefrac1\delta)
 \right)
\]
observations.
All operations run in polynomial time with respect to the number of samples so the algorithm terminates in 
\[
 \Ot\!\left(
 \left(\frac{nk}{\varepsilon}\right)^{\poly(K)}
 \log(\nicefrac1\delta)
 \right)
\]
time.
\end{proof}

\section{Moment Peeling \& Estimation}
\label{sec:peeling:main}

A key subroutine of our implicit tensor decomposition algorithm from \Cref{sec:tensor-decomposition:main} is a mean estimation oracle for the projected Hermite moments of the unrecovered regressors.
For a set $Q$ of unrecovered regressors, recall that
$Y_{\max,Q}=\max_{i\in Q}Y_i$ and that its degree-$d$
restricted Hermite moment is
\[
 \RestrictedMomentTensor{d}{Q}
 =\E_{\Covariate,\NoiseVector_Q}\!\left[
 \HermiteTensor{d}{\Covariate}\TailIndicator{Y_{\max,Q}}\right].
\]
In this section, we implement this oracle given samples of the form $(\widetilde{\Covariate}, U\widetilde \Covariate^{\otimes r}, Y_{\max})$ where $U$ is some projection matrix.
Our goal is to prove the following result,
which quantifies the estimation error relative to the trace $\tr(\EvenMomentMatrix{r}{Q}) = \sum_{i\in Q} \lambda_i^{(2r)}$.

\begin{theorem}[Hermite Moment Estimation by Peeling]
\label{thm:peeling:hermite-moment-estimation}
Consider the Gaussian self-selection model
(\Cref{def:gaussian-self-selection}) under the boundedness and
uncoveredness assumptions (\Cref{ass:boundedness-uncoveredness}).
Use the hyperparameters $K,T,s$ from
\Cref{hyp:algorithm}.

Let $2\leq r\leq s$, and fix a projection map
$U:\R^{n^r}\to\R^p$, where $U\in\R^{p\times n^r}$
has orthonormal rows.
Assume we are given fresh samples
$(\widetilde{\Covariate},U\TensorPower{\widetilde{\Covariate}}{r},Y_{\max})$
where $\widetilde{\Covariate}=\Covariate+\iu\Covariate'$.

Fix a sampling accuracy $\varepsilon_{\mathrm{est}}\in(0,1)$
and a failure probability $\delta\in(0,\nicefrac13)$.
Let $\widehat{\mathcal R}$ be a list of slope and bias estimates matched to some
$R\subsetneq[k]$, and put $Q=[k]\setminus R$.
Index the estimates as $(\widehat{\Regressor}_i,\widehat b_i)_{i\in R}$
according to this matching.
If the estimates satisfy
\[
 \norm{\widehat{\Regressor}_i-\Regressor_i}_2+|\widehat b_i-b_i|
 \leq e^{-K^6s}
 \qquad(i\in R)\,,
\]
then the following guarantees hold.
\begin{enumerate}[label=(\roman*),ref=\thetheorem(\roman*),leftmargin=2.2em]
\item\label[theorem]{thm:peeling:hermite-moment-estimation-even}
For the projected even moment at order $r$,
\textnormal{\textsc{EstimateMomentsByPeeling}} (\cref{alg:peeling:estimate-moments-by-peeling}),
given sampling accuracy $\varepsilon_{\mathrm{est}}$ and
confidence parameter $\delta$, returns a symmetric matrix
$\widehat{M}\in\R^{p\times p}$ satisfying, with probability at
least $1-\delta$,
\begin{equation}
\begin{aligned}
 \norm{\widehat{M}-U\Mat_r(\RestrictedMomentTensor{2r}{Q})U^\top}_F
 &\leq\varepsilon_{\mathrm{est}}
  + \frac{\sqrt{k}}{\sqrt\pi}e^{-K^6s}
  + e^{4s\log K-T^2/(2K)}
         k\tr(\EvenMomentMatrix{r}{Q})\,.
\end{aligned}
\end{equation}

\item\label[theorem]{thm:peeling:hermite-moment-estimation-odd}
For the projected odd moment at $r=s$,
\textnormal{\textsc{EstimateMomentsByPeeling}} returns symmetric matrices
$\widehat{M}_1,\ldots,\widehat{M}_n\in\R^{p\times p}$ satisfying,
with probability at least $1-\delta$, simultaneously for every
$\bm{\zeta}\in\R^n$,
\begin{equation}
\begin{aligned}
 \norm{
  \sum_{\ell=1}^n \zeta_\ell\widehat{M}_\ell
  -U\Mat_s(\RestrictedMomentTensor{2s+1}{Q}\mathbin{\lrcorner}\bm{\zeta})U^\top
 }_F
 &\leq
 \left(
  \varepsilon_{\mathrm{est}}
  +\frac{\sqrt{k}}{\sqrt\pi}e^{-K^6s}
  +e^{4s\log K-T^2/(2K)}k\tr(\EvenMomentMatrix{s}{Q})
 \right)\norm{\bm{\zeta}}_2\,.
\end{aligned}
\end{equation}
\end{enumerate}
Each call requires at most
\[
 \Ot\!\left(
 \frac{p^2n}{\varepsilon_{\mathrm{est}}^2}
 100^s\log(\nicefrac1\delta)
 \right)
\]
projected samples.
Moreover, its running time is at most
$
  \poly(n,k,p, \nicefrac1{\varepsilon_{\mathrm{est}}}, \log(\nicefrac1\delta))
$.
\end{theorem}

As we will see in the roadmap (\Cref{subsec:peeling:roadmap}),
the estimation error can be decomposed into three terms.
The rest of this section aims to bound each of these terms which can then be combined into a proof of \Cref{thm:peeling:hermite-moment-estimation}.

\subsection{Proof Roadmap}
\label{subsec:peeling:roadmap}

Our estimation algorithm is easiest to present and analyze for expectations of the form
\[
  \E_{\Covariate,\NoiseVector,\Covariate'}
  [\Psi(\Covariate, \Covariate')\,\TailIndicator{Y_{\max,Q}}]\,,
\]
where $\Psi\in \R^{p\times p}$ is a random matrix.
We will then apply our analysis with $\Psi$ being projections of the even and odd Hermite moments.

One key challenge is that we do not directly observe $Y_{\max, Q}$ and hence we must use a proxy $\widehat\omega_{\widehat{\mathcal R}}(\Covariate,Y_{\max})$ for the indicator variable $\TailIndicator{Y_{\max, Q}}$.
We define this proxy in \Cref{subsubsec:peeling:signed-weight}.
Then, we discuss the specific instantiations of $\Psi$ in \Cref{subsubsec:peeling:projected-hermite-features}.
We then decompose the estimation error into three terms in \Cref{subsubsec:peeling:error-decomposition} and state our bounds for each term.

{\subsubsection{Signed Peeling Weight}
\label{subsubsec:peeling:signed-weight}}
We now define the proxy of the indicator $\TailIndicator{Y_{\max, Q}}$.
Let $\widehat{\mathcal R}$ be the list of recovered slope and bias estimates corresponding to indices $R\subsetneq[k]$ and put $Q=[k]\setminus R$.

\paragraph{Ideal Peeling Weight.}
We first define an idealized peeling weight $\omega_R(\Covariate,Y_{\max})$ assuming we know the exact slopes and biases $(\Regressor_i,b_i), i\in R$.
Recall for an index set $I\subseteq[k]$,
we write $Y_{\max,I}=\max_{i\in I}Y_i$ to denote its maximum response.
Define the conditional CDF of the maximum response conditioned on the covariate,
\begin{equation}
 F_I(t\mid\bm x)
 \coloneqq
 \Prb_{\NoiseVector_I}(Y_{\max,I}\leq t\mid\Covariate=\bm x)\,.
\end{equation}

Recall the peeling identity from \Cref{subsubsec:overview:moment-peeling}.
Since $Y_{\max}=\max\{Y_{\max,R},Y_{\max,Q}\}$,
\begin{align*}
 \TailIndicator{Y_{\max}}-\TailIndicator{Y_{\max,R}}
 &=\TailIndicator{Y_{\max,Q}}\\
 &\quad-\TailIndicator{\min\{Y_{\max,R},Y_{\max,Q}\}}.
\end{align*}
We would like to define the proxy $\omega_R(\Covariate, Y_{\max})$ to be the LHS,
since it is equal to the indicator we want and a term we will show to be negligible.
However, we do not directly observe $Y_{\max, R}$.
Our saving grace is that we only need equality under expectation so that we can define the proxy to match the LHS in conditional expectation conditioned on the covariate $\Covariate$.
\begin{align}
 \omega_R(\Covariate,Y_{\max})
 &\coloneqq
 \TailIndicator{Y_{\max}}
 -\bigl(F_R(T+1\mid\Covariate)-F_R(T\mid\Covariate)\bigr)\,,
 \label{eq:true-peeling-weight} \\
 \E_{\NoiseVector}[\omega_R(\Covariate,Y_{\max})\mid \Covariate]
 &= \E_{\NoiseVector}[\TailIndicator{Y_{\max}}-\TailIndicator{Y_{\max,R}} \mid \Covariate]\,.
 \notag
\end{align}

\paragraph{Peeling Weight.} 
In reality, we only have a list $\widehat{\mathcal R}$ of estimates to the exact slopes and biases $(\Regressor_i,b_i), i\in R$.
Since we only need said regressors to evaluate the conditional CDF,
we can simply plug in our estimates instead.
Write $\widehat{\mathcal R}=((\widehat{\Regressor}_i,\widehat b_i))_{i\in R}$,
where $R\subseteq [k]$ corresponds to the index of the recovered regressors.
Let $\bm\eta_R$ be an independent copy of $\NoiseVector_R$. Define
\begin{equation}
 \widehat F_{\widehat{\mathcal R}}(t\mid\bm x)
 \coloneqq\Prb_{\bm\eta_R}\!\left[
 \max_{i\in R}\bigl(\ip{\widehat{\Regressor}_i}{\bm x}+\widehat b_i+\eta_i\bigr)
 \leq t\right].
\end{equation}
The noise law depends only on $|R|$ and $\rho$, so our algorithm
can estimate $\widehat{F}_{\widehat{\mathcal R}}(t\mid \bm{x})$ using simulated noise $\bm \eta_R$ without knowing $R$.
We define the signed peeling weight to be
\begin{equation}\label{eq:estimated-peeling-weight}
 \widehat{\omega}_{\widehat{\mathcal R}}(\Covariate,Y_{\max})
 \coloneqq
 \TailIndicator{Y_{\max}}
 -\bigl(
 \widehat{F}_{\widehat{\mathcal R}}(T+1\mid \Covariate)
 -\widehat{F}_{\widehat{\mathcal R}}(T\mid \Covariate)
 \bigr)
 \in [-1,1]\,.
\end{equation}

{\subsubsection{Projected Hermite Features}
\label{subsubsec:peeling:projected-hermite-features}}
We now discuss the specific projected Hermite features $\Psi$ whose restricted expectations we wish to estimate.
For $1\leq r\leq s$, fix a projection map
$U:\R^{n^r}\to\R^p$ with orthonormal rows,
so $UU^\top=I_{p}$.
Recall that part of our implicit tensor decomposition algorithm presented in \Cref{subsec:tensor-decomposition:iterative-projections} is an iterative projection procedure
that recursively computes the projected moments and projection maps degree by degree.
Extending the construction by one tensor degree requires an even
moment estimate to determine the next selected subspace.

Define the functions used by \textsc{EstimateMomentsByPeeling} as
\begin{equation}\label{eq:projected-even-hermite-feature}
 \Psi_{\even}^{U}(\Covariate,\Covariate')
 \coloneqq
 \frac{\RealPart\!\left[
  (U\TensorPower{\widetilde{\Covariate}}{r})
  (U\TensorPower{\widetilde{\Covariate}}{r})^\top
 \right]}{\sqrt{(2r)!}}
\end{equation}
and
\begin{equation}\label{eq:projected-odd-hermite-feature}
 \Psi_{\odd,\ell}^{U}(\Covariate,\Covariate')
 \coloneqq
 \frac{\RealPart\!\left[
  (U\TensorPower{\widetilde{\Covariate}}{r})
  (U\TensorPower{\widetilde{\Covariate}}{r})^\top
  \widetilde X_\ell
 \right]}{\sqrt{(2r+1)!}},
 \qquad \ell\in[n].
\end{equation}
These definitions apply to every $r\geq1$ and give real symmetric
$p\times p$ matrices.
Here $\widetilde \Covariate = \Covariate + \iu \Covariate'$ denotes the complexified covariate.
Averaging over the complex component $\Covariate'$ 
gives
\begin{equation}
\begin{aligned}
 \E_{\Covariate'}[\Psi_{\even}^{U}(\Covariate,\Covariate')\mid \Covariate]
 &=U\Mat_r(\HermiteTensor{2r}{\Covariate})U^\top,\\
 \E_{\Covariate'}[\Psi_{\odd,\ell}^{U}(\Covariate,\Covariate')\mid \Covariate]
 &=U
   \Mat_r(\HermiteContraction{2r+1}{\Covariate}{\bm{e}_\ell})U^\top.
\end{aligned}
\end{equation}

Our algorithm estimates the weighted expectations below. In this display,
the projected Hermite features are evaluated at $(\Covariate,\Covariate')$.
\begin{equation}
\begin{aligned}
 \E_{\Covariate,\NoiseVector,\Covariate'}\!\left[
   \Psi_{\even}^{U}
   \red{\widehat\omega_{\widehat{\mathcal R}}(\Covariate,Y_{\max})}\right]
 &\approx\E_{\Covariate,\NoiseVector_Q,\Covariate'}\!\left[
   \Psi_{\even}^{U}\red{\TailIndicator{Y_{\max,Q}}}\right]\\
 &=U\Mat_r(\RestrictedMomentTensor{2r}{Q})U^\top,\\[3pt]
 \E_{\Covariate,\NoiseVector,\Covariate'}\!\left[
   \Psi_{\odd,\ell}^{U}
   \red{\widehat\omega_{\widehat{\mathcal R}}(\Covariate,Y_{\max})}\right]
 &\approx\E_{\Covariate,\NoiseVector_Q,\Covariate'}\!\left[
   \Psi_{\odd,\ell}^{U}\red{\TailIndicator{Y_{\max,Q}}}\right]\\
 &=U\Mat_r(\RestrictedMomentTensor{2r+1}{Q}
                   \mathbin{\lrcorner}\bm e_\ell)U^\top
 \qquad(\ell\in[n]).
\end{aligned}
\end{equation}

We use the even feature for $2\leq r\leq s$ and the odd coordinate
features at $r=s$.

For each observation $(\Covariate,Y_{\max})$, we use the same
auxiliary noise vector $\bm\eta_R$ for every feature
$\Psi_\ell(\Covariate,\Covariate')$, $\ell\in[m]$,
evaluated on that observation.
We draw $\bm\eta_R$ independently of the samples from the self-selection model.

\begin{algorithm}[H]
\caption{\textsc{EstimateMomentsByPeeling}}\label{alg:peeling:estimate-moments-by-peeling}
\KwIn{$\widehat{\mathcal R}$; a fresh sample block; a fixed family
$(\Psi_\ell)_{\ell\in[m]}$ of symmetric $p\times p$ matrix features
evaluated from the supplied samples;
accuracy $\varepsilon_{\mathrm{est}}\in(0,1)$;
failure probability $\delta\in(0,\nicefrac13)$}
\KwOut{Symmetric matrices $\widehat{M}_1,\ldots,\widehat{M}_m$}
For each observation, draw $\bm\eta_R$ and compute
$\widetilde Y\gets\max_{i\in R}
 \bigl(\ip{\widehat{\Regressor}_i}{\Covariate}+\widehat b_i+\eta_i\bigr)$\;
\tcp{\Cref{lem:coordinatewise-trimmed-mean}}
$(\widehat{M}_1,\ldots,\widehat{M}_m)\gets$
\TrimmedMean{$
 \bigl((\TailIndicator{Y_{\max}}-\TailIndicator{\widetilde Y})
       \Psi_\ell(\Covariate,\Covariate')\bigr)_{\ell\in[m]};
 100^{s+2},\frac{\varepsilon_{\mathrm{est}}}{p\sqrt m},\delta
$}\;
\Return $(\widehat{M}_1,\ldots,\widehat{M}_m)$\;
\end{algorithm}

{\subsubsection{Error Decomposition}
\label{subsubsec:peeling:error-decomposition}}

We now decompose the estimation error.
Let $\Psi:\R^n\times\R^n\to\R^{p\times p}$ be measurable
with $\E_{\Covariate,\Covariate'}[\norm{\Psi(\Covariate,\Covariate')}_F]<\infty$,
and let $\omega_R(\Covariate,Y_{\max})$ and $\widehat\omega_{\widehat{\mathcal R}}(\Covariate,Y_{\max})$ denote the true and estimated peeling weights from \eqref{eq:true-peeling-weight} and \eqref{eq:estimated-peeling-weight}.
For a matrix estimate $\widehat{M}$, the error decomposes as follows.
\begin{equation}\label{eq:peeling:error-decomposition}
\begin{aligned}
 &{\widehat{M}}-\E_{\Covariate,\NoiseVector_Q,\Covariate'}[\Psi(\Covariate,\Covariate'){\TailIndicator{Y_{\max,Q}}}]\\
 &\quad=
 \underbrace{
  {\widehat{M}}-{\E}_{\Covariate,Y_{\max},\Covariate'}[
   \Psi(\Covariate,\Covariate')\widehat{\omega}_{\widehat{\mathcal R}}(\Covariate,Y_{\max})]
 }_{\text{sampling error}}\\
 &\qquad+
 \underbrace{
  \E_{\Covariate,Y_{\max},\Covariate'}\!\left[
   \Psi(\Covariate,\Covariate')\bigl(
    \red{\widehat{\omega}_{\widehat{\mathcal R}}(\Covariate,Y_{\max})}-\red{\omega_R(\Covariate,Y_{\max})}
   \bigr)\right]
 }_{\text{error from recovered regressors}}\\
 &\qquad+
 \underbrace{
   \E_{\Covariate,Y_{\max},\Covariate'}[\Psi(\Covariate,\Covariate')\red{\omega_R(\Covariate,Y_{\max})}]
   -\E_{\Covariate,\NoiseVector_Q,\Covariate'}[\Psi(\Covariate,\Covariate')\red{\TailIndicator{Y_{\max,Q}}}]
 }_{\text{bias due to signed peeling weight}}.
\end{aligned}
\end{equation}

The three terms in each bound are the sampling error, the error
from the recovered regressors, and the bias due to the signed
peeling weight, respectively.
We now state the bounds for each term.

\paragraph{Sampling error.}
The first term in the error decomposition \eqref{eq:peeling:error-decomposition} is the sampling error in estimating the weighted
expectation formed from the recovered list. We bound this error
for the even estimate and for every contraction of the odd estimates.

\begin{lemma}[Sampling Error for Projected Hermite Moments]
\label{lem:peeling:projected-moment-sampling}
Fix integers $2\leq r\leq s$, a matrix $U\in\R^{p\times n^r}$
with orthonormal rows, and a recovered list $\widehat{\mathcal R}$.
Run \textnormal{\textsc{EstimateMomentsByPeeling}} (\cref{alg:peeling:estimate-moments-by-peeling})
with accuracy $\varepsilon_{\mathrm{est}}\in(0,1)$ and failure
probability $\delta\in(0,\nicefrac13)$. Each call uses fresh samples
from the Gaussian self-selection model (\Cref{def:gaussian-self-selection}),
independent of $U$ and $\widehat{\mathcal R}$, and an independent
standard Gaussian vector $\Covariate'$ for each sample.
\begin{enumerate}[label=(\roman*),ref=\thelemma(\roman*),leftmargin=2.2em]
\item\label[lemma]{lem:peeling:projected-moment-sampling-even}
Given the even feature $\Psi_{\even}^U$, the algorithm returns
$\widehat M$ satisfying, with probability at least $1-\delta$,
\[
 \norm{\widehat M-\E_{\Covariate,Y_{\max},\Covariate'}\!\left[
 \widehat\omega_{\widehat{\mathcal R}}(\Covariate,Y_{\max})
 \Psi_{\even}^U(\Covariate,\Covariate')\right]}_F
 \leq\varepsilon_{\mathrm{est}}.
\]
\item\label[lemma]{lem:peeling:projected-moment-sampling-odd}
At $r=s$, given the odd features $(\Psi_{\odd,\ell}^U)_{\ell\in[n]}$,
the algorithm returns $\widehat M_1,\ldots,\widehat M_n$ satisfying,
with probability at least $1-\delta$, simultaneously for every
$\bm{\zeta}\in\R^n$,
\begin{equation}
\begin{aligned}
 &\norm{\sum_{\ell=1}^n \zeta_\ell\left(\widehat M_\ell
 -\E_{\Covariate,Y_{\max},\Covariate'}\!\left[
 \widehat\omega_{\widehat{\mathcal R}}(\Covariate,Y_{\max})
 \Psi_{\odd,\ell}^U(\Covariate,\Covariate')\right]\right)}_F\\
 &\qquad\leq\varepsilon_{\mathrm{est}}\norm{\bm{\zeta}}_2.
\end{aligned}
\end{equation}
\end{enumerate}
\end{lemma}

We prove \Cref{lem:peeling:projected-moment-sampling} in \Cref{subsec:peeling:sampling-error}.
The proof applies a general mean estimation guarantee under an entrywise
second-moment bound, which we verify for the weighted even and odd
Hermite features.

\paragraph{Recovered-regressor error.}
The second term in the error decomposition \eqref{eq:peeling:error-decomposition} measures the change in the weighted expectation
when the true recovered regressors are replaced by their estimates.

For a scalar function $\psi$, write
$\norm{\psi}_{L^2(\Covariate,\Covariate')}
=(\E_{\Covariate,\Covariate'}[
 |\psi(\Covariate,\Covariate')|^2])^{1/2}$.
\begin{lemma}[Peeling-Weight Stability]
\label{lem:peeling:weight-stability}
Consider the Gaussian self-selection model
(\Cref{def:gaussian-self-selection}).
Fix a list $\widehat{\mathcal R}$ of pairs in $\R^n\times\R$,
a set $R\subsetneq[k]$, and a bijection
$\match:\widehat{\mathcal R}\to R$.
Index the pairs as $(\widehat{\Regressor}_i,\widehat b_i)_{i\in R}$
according to this matching.
Let $\Covariate'\sim N(\bm{0},I_n)$ be independent of the covariate and noises.
For every measurable function $\psi:\R^n\times\R^n\to\R$
with $\norm{\psi}_{L^2(\Covariate,\Covariate')}<\infty$,
\begin{equation}
\begin{aligned}
 &\abs{\E_{\Covariate,Y_{\max},\Covariate'}\!\left[
  \psi(\Covariate,\Covariate')\bigl(
   \red{\widehat{\omega}_{\widehat{\mathcal R}}(\Covariate,Y_{\max})}-\red{\omega_R(\Covariate,Y_{\max})}
  \bigr)\right]}\\
 &\quad\leq\frac1{\sqrt\pi}\norm{\psi}_{L^2(\Covariate,\Covariate')}
 \left(
  \sum_{i\in R}
  \bigl(\norm{\widehat{\Regressor}_i-\Regressor_i}_2+|\widehat b_i-b_i|\bigr)^2
 \right)^{1/2}.
\end{aligned}
\end{equation}
\end{lemma}

We prove \Cref{lem:peeling:weight-stability} in \Cref{subsec:peeling:recovered-regressor-error} using a Gaussian translation bound.

\paragraph{Peeling bias.}
Even with the true recovered regressors, the signed weight can
differ from the slice indicator of the unrecovered maximum.
The third term in the error decomposition \eqref{eq:peeling:error-decomposition} captures this bias.
We identify this bias in the lemma below and bound its magnitude.

\begin{lemma}[Peeling Bias Bound]
\label{lem:peeling:bias-bound}
Consider the Gaussian self-selection model
(\Cref{def:gaussian-self-selection}).
Fix $R\subsetneq[k]$, put $Q=[k]\setminus R$, and let $\Covariate'$
be a random element independent of the covariate and noises.
\begin{enumerate}[label=(\roman*),ref=\thelemma(\roman*),leftmargin=2.2em]
\item\label[lemma]{lem:peeling:bias-identity}
For every integrable scalar- or matrix-valued function
$\psi(\Covariate,\Covariate')$,
\begin{equation}
\begin{aligned}
 &\E_{\Covariate,\NoiseVector,\Covariate'}\!\left[
  \psi(\Covariate,\Covariate')\red{\omega_R(\Covariate,Y_{\max})}\right]
 -\E_{\Covariate,\NoiseVector,\Covariate'}\!\left[
  \psi(\Covariate,\Covariate')\red{\TailIndicator{Y_{\max,Q}}}\right]\\
 &\qquad=
 -\E_{\Covariate,\NoiseVector,\Covariate'}\!\left[
  \psi(\Covariate,\Covariate')
  \TailIndicator{\min\{Y_{\max,R},Y_{\max,Q}\}}\right].
\end{aligned}
\end{equation}
\item\label[lemma]{lem:peeling:cross-block-bias-bound}
Suppose the boundedness and uncoveredness assumptions
(\Cref{ass:boundedness-uncoveredness}) hold, and use the choice of $K$ from
\Cref{hyp:algorithm}. Let $s\geq1$ be an integer and suppose
$T\geq4BK$ and $(T-B)^2\geq Ks$. If $\Covariate'\sim N(\bm{0},I_n)$ and
$\psi(\Covariate,\Covariate')$ is a real polynomial of degree at
most $2s+1$, then for every $1\leq r\leq s$,
\[
\begin{aligned}
 &\left|\E_{\Covariate,\NoiseVector,\Covariate'}\!\left[
  \psi(\Covariate,\Covariate')
  \TailIndicator{\min\{Y_{\max,R},Y_{\max,Q}\}}\right]\right|\\
 &\qquad\leq k e^{4s\log K-T^2/(2K)}
 \norm{\psi}_{L^2(\Covariate,\Covariate')}\tr(\EvenMomentMatrix{r}{Q}).
\end{aligned}
\]
\end{enumerate}
\end{lemma}

The identity in \Cref{lem:peeling:bias-identity} follows by conditioning
on $\Covariate$ in the peeling identity above. 
We prove \Cref{lem:peeling:cross-block-bias-bound} in \Cref{subsec:peeling:peeling-bias} by applying H\"older's inequality and then bounding the resulting factors using Gaussian hypercontractivity and standard tail bounds.

Finally, we combine all the ingredients from the roadmap to prove \Cref{thm:peeling:hermite-moment-estimation} in \Cref{subsec:peeling:moment-estimation-proof}.
The rest of this section is dedicated to the proofs of the statements above.

\subsection{Sampling Error}
\label{subsec:peeling:sampling-error}

To prove \Cref{lem:peeling:projected-moment-sampling}, we first bound the second
moment of every weighted matrix entry, then apply a general mean
estimation guarantee. Since the signed peeling weight has absolute
value at most one, it suffices to bound the entries of the projected
Hermite features.

Recall that
$\norm{\psi}_{L^2(\Covariate,\Covariate')}
=(\E_{\Covariate,\Covariate'}[
 |\psi(\Covariate,\Covariate')|^2])^{1/2}$.
\begin{restatable}[Projected-Feature Second Moments]{proposition}{projectedfeaturesecondmoments}
\label{prop:peeling:projected-feature-second-moments}
Let $r\geq1$, and fix a projection map
$U:\R^{n^r}\to\R^p$, where $U\in\R^{p\times n^r}$
has orthonormal rows.
Let $\Covariate,\Covariate'\sim N(\bm{0},I_n)$ be independent.
For the projected features $\Psi_{\even}^U$ and
$(\Psi_{\odd,\ell}^U)_{\ell\in[n]}$ defined in
\eqref{eq:projected-even-hermite-feature} and \eqref{eq:projected-odd-hermite-feature}, a matrix $H\in\R^{p\times p}$
with $\norm{H}_F\leq1$, and any $\bm{\zeta}\in\R^n$, write
\begin{align*}
 \psi(\Covariate,\Covariate')
 &\coloneqq\tr\!\left(H^\top\Psi_{\even}^{U}(\Covariate,\Covariate')\right),\\
 \psi_{\bm{\zeta}}(\Covariate,\Covariate')
 &\coloneqq\tr\!\left(
 H^\top\sum_{\ell=1}^n \zeta_\ell\Psi_{\odd,\ell}^{U}(\Covariate,\Covariate')
 \right).
\end{align*}
These are real polynomials in $(\Covariate,\Covariate')\in\R^{2n}$ of degrees at
most $2r$ and $2r+1$, respectively, and
\[
 \E_{\Covariate,\Covariate'}[\psi(\Covariate,\Covariate')^2]\leq100^{r+2},
 \qquad
 \E_{\Covariate,\Covariate'}[\psi_{\bm{\zeta}}(\Covariate,\Covariate')^2]\leq100^{r+2}\norm{\bm{\zeta}}_2^2.
\]
Moreover, their conditional means satisfy
\begin{equation}
\begin{aligned}
 \left\|\E_{\Covariate'}[\psi(\Covariate,\Covariate')\mid \Covariate]\right\|_{L^2(\Covariate)}
 &\leq1,\\
 \left\|\E_{\Covariate'}[\psi_{\bm{\zeta}}(\Covariate,\Covariate')\mid \Covariate]\right\|_{L^2(\Covariate)}
 &\leq\norm{\bm{\zeta}}_2.
\end{aligned}
\end{equation}
\end{restatable}
The proof of \Cref{prop:peeling:projected-feature-second-moments} follows from a standard application of
Wick's identity (\Cref{prop:wick-identity}) and Hermite orthogonality (\Cref{prop:hermite-expansion-normalization}).
As such, we defer the details to \Cref{subsec:app-peeling:second-moment-bounds}.

We use the following mean estimation bound for a family of
matrix-valued functions. Its proof applies the coordinatewise
trimmed mean to all entries on the same sample block and is deferred
to \Cref{subsec:app-peeling:sampling-error}.

\begin{restatable}[Mean Estimation with the Peeling Weight]{lemma}{peelingmeanestimation}
\label{lem:peeling:weighted-mean-estimation}
Let $\varepsilon_{\mathrm{est}}\in(0,1)$ and
$\delta\in(0,\nicefrac13)$.
Fix a recovered list $\widehat{\mathcal R}$ and measurable symmetric
matrix-valued functions
$\Psi_\ell:\R^n\times\R^n\to\R^{p\times p}$, $\ell\in[m]$,
where $p,m\geq1$.
Let $(\Covariate,Y_{\max})$ be a fresh sample from the Gaussian self-selection
model (\Cref{def:gaussian-self-selection}), and let
$\Covariate'\sim N(\bm{0},I_n)$ be independent of $(\Covariate,Y_{\max})$.
For the signed peeling weight $\widehat{\omega}_{\widehat{\mathcal R}}(\Covariate,Y_{\max})$ \eqref{eq:estimated-peeling-weight} formed from $\widehat{\mathcal R}$, suppose
\[
 \E_{\Covariate,\Covariate'}\!\left[
  (\Psi_\ell(\Covariate,\Covariate'))_{ac}^{\,2}\right]
 \leq100^{s+2}\qquad(\ell\in[m],\ a,c\in[p]).
\]
Given the recovered list, fresh samples, and these functions,
\mbox{\textnormal{\textsc{EstimateMomentsByPeeling}}} (\cref{alg:peeling:estimate-moments-by-peeling}) returns
symmetric matrices $\widehat{M}_1,\ldots,\widehat{M}_m$ such that,
with probability at least $1-\delta$, simultaneously for every
$\bm{\zeta}\in\R^m$,
\[
 \norm{\sum_{\ell=1}^m \zeta_\ell\left(
  {\widehat{M}_\ell}-{\E}_{\Covariate,Y_{\max},\Covariate'}\!\left[
  \widehat{\omega}_{\widehat{\mathcal R}}(\Covariate,Y_{\max})
  \Psi_\ell(\Covariate,\Covariate')\right]\right)}_F
 \leq \varepsilon_{\mathrm{est}}\norm{\bm{\zeta}}_2.
\]
Its sample complexity is
\[
 \Ot\!\left(
 \frac{p^2m}{\varepsilon_{\mathrm{est}}^2}
 100^s\log(\nicefrac1\delta)
 \right)
\]
and its running time is $\poly(n, k, p, m, 100^s, \nicefrac1{\eps_{\mathrm{est}}}, \log(\nicefrac1\delta))$.
\end{restatable}

Armed with the necessary second moment bound (\Cref{prop:peeling:projected-feature-second-moments}) to apply \Cref{lem:peeling:weighted-mean-estimation},
we are now ready to analyze the sampling error. 
\begin{proof}[Proof of \cref{lem:peeling:projected-moment-sampling}]
Take $H=\bm e_a\bm e_c^\top$ in \Cref{prop:peeling:projected-feature-second-moments}.
The feature bound gives
\begin{align*}
 \E_{\Covariate,\Covariate'}[
  (\Psi_{\even}^{U}(\Covariate,\Covariate'))_{ac}^2]
 &\leq100^{r+2}\leq100^{s+2}.
 \alignnote{by \Cref{prop:peeling:projected-feature-second-moments}, $r\leq s$}
\end{align*}
Applying \Cref{lem:peeling:weighted-mean-estimation} with $m=1$ gives
\Cref{lem:peeling:projected-moment-sampling-even}.

For the odd call, take $r=s$ and $\bm{\zeta}=\bm e_\ell$ in the same
proposition. This gives the required second-moment bound for every
entry of every odd coordinate matrix. Applying
\Cref{lem:peeling:weighted-mean-estimation} with $m=n$ gives
\Cref{lem:peeling:projected-moment-sampling-odd}, simultaneously for every contraction.
\end{proof}

\subsection{Recovered-Regressor Error}
\label{subsec:peeling:recovered-regressor-error}

Having controlled the error due to sampling,
we turn to the proof of \cref{lem:peeling:weight-stability},
which quantifies the error from replacing the true regressors in the ideal peeling weight with their estimates.

First, we state a useful translation bound of the Gaussian measure without proof.
\begin{proposition}[Gaussian Translation Bound;
 {\cite[Section~2]{DevroyeMehrabianReddad23}}]
\label{prop:gaussian-translation}
Let $\bm{G}\sim N(\bm{0},I_p)$, where $p\geq1$.
For every $\bm{\mu},\bm{\nu}\in\R^p$ and measurable $\mathcal E\subseteq\R^p$,
\[
 \abs{\Prb(\bm{G}+\bm{\mu}\in\mathcal E)-\Prb(\bm{G}+\bm{\nu}\in\mathcal E)}
 \leq\frac{\norm{\bm{\mu}-\bm{\nu}}_2}{\sqrt{2\pi}}.
\]
\end{proposition}

This is all we need to control the error due to approximate regressor recovery.
\begin{proof}[Proof of \cref{lem:peeling:weight-stability}]
Let $\widehat{\mathcal R}$ be a list of recovered slope and bias estimates corresponding to regressors indexed by $R$ and $\psi$ be a function with
$\norm{\psi}_{L^2(\Covariate,\Covariate')}<\infty$.
We wish to bound
\[
 \abs{\E_{\Covariate,Y_{\max},\Covariate'}\!\left[
  \psi(\Covariate,\Covariate')\bigl(
  \widehat\omega_{\widehat{\mathcal R}}(\Covariate,Y_{\max})
  -\omega_R(\Covariate,Y_{\max})\bigr)\right]}.
\]
We first fix the covariate and compare the probabilities of the same
slice under the true and estimated Gaussian means. We then integrate
this comparison against $\psi$ and apply Cauchy--Schwarz.

If the recovered list is empty, the two peeling weights coincide.
Otherwise, index the estimates by corresponding indices $R\subseteq [k]$ of the regressors and fix $\bm x\in\R^n$.
Put
\[
 \mu_i=\ip{\Regressor_i}{\bm x}+b_i,\qquad
 \widehat\mu_i=\ip{\widehat{\Regressor}_i}{\bm x}+\widehat b_i\qquad(i\in R).
\]
Let $\bm\eta_R$ be an independent copy of $\NoiseVector_R$ and set
\[
 \mathcal E=\{\bm z\in\R^{|R|}:T\leq\max_{i\in R}z_i\leq T+1\}.
\]
Expanding the two peeling weights gives, for every $y\in\R$,
\begin{align}
 \widehat\omega_{\widehat{\mathcal R}}(\bm x,y)-\omega_R(\bm x,y)
 &=\TailIndicator{y}
   -\Prb_{\bm\eta_R}(\bm\eta_R+\widehat{\bm\mu}\in\mathcal E)
   \notag\\
 &\quad-\TailIndicator{y}
   +\Prb_{\bm\eta_R}(\bm\eta_R+\bm\mu\in\mathcal E)
   \notag\\
 &=\Prb_{\bm\eta_R}(\bm\eta_R+\bm\mu\in\mathcal E)
   -\Prb_{\bm\eta_R}(\bm\eta_R+\widehat{\bm\mu}\in\mathcal E).
\end{align}
We use $\eta_i=\sqrt\rho\,g+\sqrt{1-\rho}\,g_i$, with
$g,(g_i)_{i\in R}$ independent standard Gaussians. We condition on
a Gaussian part whose variance is at least $1/2$.

\emph{Case I ($\rho\leq1/2$).} Condition on $g$.
The two response vectors then have conditional laws
\begin{gather*}
 N\!\left((\mu_i+\sqrt\rho\,g)_{i\in R},(1-\rho)I_{|R|}\right),\\
 N\!\left((\widehat\mu_i+\sqrt\rho\,g)_{i\in R},
          (1-\rho)I_{|R|}\right).
\end{gather*}
Their mean vectors differ by $\widehat{\bm\mu}-\bm\mu$.
Applying the Gaussian translation bound to $\mathcal E$ gives
\begin{align*}
 &\abs{\Prb_{\bm\eta_R}(\bm\eta_R+\widehat{\bm\mu}\in\mathcal E\mid g)
       -\Prb_{\bm\eta_R}(\bm\eta_R+\bm\mu\in\mathcal E\mid g)}\\
 &\quad\leq
 \frac{\norm{\widehat{\bm\mu}-\bm\mu}_2}{\sqrt{2\pi(1-\rho)}}
 \alignnote{by \Cref{prop:gaussian-translation}, after scaling}\\
 &\quad\leq
 \frac{\norm{\widehat{\bm\mu}-\bm\mu}_2}{\sqrt\pi}.
 \alignnote{$1-\rho\geq1/2$}
\end{align*}
Averaging over $g$ proves the same bound without conditioning.

\emph{Case II ($\rho\geq1/2$).} Condition on $(g_i)_{i\in R}$.
Put
\[
 m=\max_{i\in R}\{\mu_i+\sqrt{1-\rho}\,g_i\},\qquad
 \widehat m=\max_{i\in R}\{\widehat\mu_i+\sqrt{1-\rho}\,g_i\}.
\]
The two recovered maxima are then $m+\sqrt\rho\,g$ and
$\widehat m+\sqrt\rho\,g$. Moreover,
\[
 |\widehat m-m|\leq\max_{i\in R}|\widehat\mu_i-\mu_i|
 \leq\norm{\widehat{\bm\mu}-\bm\mu}_2.
\]
Applying the scalar Gaussian translation bound to $[T,T+1]$ gives
\begin{align*}
 &\abs{\Prb_g(\widehat m+\sqrt\rho\,g\in[T,T+1])
       -\Prb_g(m+\sqrt\rho\,g\in[T,T+1])}\\
 &\quad\leq\frac{|\widehat m-m|}{\sqrt{2\pi\rho}}
 \alignnote{by \Cref{prop:gaussian-translation}, after scaling}\\
 &\quad\leq
 \frac{\norm{\widehat{\bm\mu}-\bm\mu}_2}{\sqrt\pi}.
 \alignnote{$\rho\geq1/2$}
\end{align*}
Averaging over $(g_i)_{i\in R}$ proves the same bound without conditioning.
Thus, for every $\bm x\in\R^n$ and $y\in\R$,
\begin{equation}\label{eq:peeling:pointwise-weight-perturbation}
 \abs{\widehat\omega_{\widehat{\mathcal R}}(\bm x,y)-\omega_R(\bm x,y)}
 \leq\frac1{\sqrt\pi}
 \left(\sum_{i\in R}
 \bigl(\ip{\widehat{\Regressor}_i-\Regressor_i}{\bm x}+\widehat b_i-b_i\bigr)^2\right)^{1/2}.
\end{equation}

We now integrate this pointwise bound against $\psi$.
The difference of the weights does not depend on $Y_{\max}$, so
\begin{align*}
 &\abs{\E_{\Covariate,Y_{\max},\Covariate'}\!\left[
  \psi(\Covariate,\Covariate')\bigl(
   \widehat{\omega}_{\widehat{\mathcal R}}(\Covariate,Y_{\max})
   -\omega_R(\Covariate,Y_{\max})\bigr)\right]}\\
 &\quad\leq\frac1{\sqrt\pi}
 \E_{\Covariate,\Covariate'}\!\left[
  |\psi(\Covariate,\Covariate')|
  \left(
   \sum_{i\in R}
   \bigl(\ip{\widehat{\Regressor}_i-\Regressor_i}{\Covariate}+\widehat b_i-b_i\bigr)^2
  \right)^{1/2}\right]
 \alignnote{by \eqref{eq:peeling:pointwise-weight-perturbation}}\\
 &\quad\leq\frac1{\sqrt\pi}\norm{\psi}_{L^2(\Covariate,\Covariate')}
 \left(
  \E_{\Covariate}\sum_{i\in R}
  \bigl(\ip{\widehat{\Regressor}_i-\Regressor_i}{\Covariate}+\widehat b_i-b_i\bigr)^2
 \right)^{1/2}
 \alignnote{Cauchy--Schwarz}\\
 &\quad=\frac1{\sqrt\pi}\norm{\psi}_{L^2(\Covariate,\Covariate')}
 \left(
  \sum_{i\in R}
  \left(\norm{\widehat{\Regressor}_i-\Regressor_i}_2^2+|\widehat b_i-b_i|^2\right)
 \right)^{1/2}
 \alignnote{$\Covariate\sim N(\bm0,I_n)$}\\
 &\quad\leq\frac1{\sqrt\pi}\norm{\psi}_{L^2(\Covariate,\Covariate')}
 \left(
  \sum_{i\in R}
  \bigl(\norm{\widehat{\Regressor}_i-\Regressor_i}_2+|\widehat b_i-b_i|\bigr)^2
 \right)^{1/2}.
\end{align*}
This proves the lemma.
\end{proof}

\subsection{Peeling Bias}
\label{subsec:peeling:peeling-bias}

Finally, we prove \Cref{lem:peeling:bias-bound}, which gives a bound on the bias resulting from using the peeling weight $\omega_R(\Covariate, Y_{\max})$ as a proxy for the indicator $\TailIndicator{Y_{\max, Q}}$.
Our main tool is this pairwise bound stated below whose proof follows from an application of Gaussian hypercontractivity and standard tail bounds to the factors obtained via H\"older's inequality.
As such, the proof is deferred to \Cref{subsec:app-peeling:peeling-bias}.
\begin{restatable}[Pairwise Exceedance Bound]{lemma}{pairwiseexceedancebound}
\label{lem:peeling:pairwise-exceedance}
Consider the Gaussian self-selection model
(\Cref{def:gaussian-self-selection}) under the boundedness and
uncoveredness assumptions (\Cref{ass:boundedness-uncoveredness}).
Use the choice of $K$ from \Cref{hyp:algorithm}.
Let $s\geq1$ be an integer and suppose $T\geq4BK$ and $(T-B)^2\geq Ks$.
Let $(\Covariate,\NoiseVector)$ follow the Gaussian self-selection
model (\Cref{def:gaussian-self-selection}). Let $\Covariate'\sim N(\bm{0},I_n)$ be independent of
the covariate and noises.
Fix a real polynomial $\psi(\Covariate,\Covariate')$ of degree at most $2s+1$.
For every pair of distinct indices $i,j\in[k]$ and every
$1\leq r\leq s$,
\begin{equation}
\begin{aligned}
 &\E_{\Covariate,\xi_i,\xi_j,\Covariate'}\!\left[
  |\psi(\Covariate,\Covariate')|\1\{Y_i>T,Y_j>T\}\right]\\
 &\quad\leq e^{4s\log K-T^2/(2K)}
 \norm{\psi}_{L^2(\Covariate,\Covariate')}\ComponentWeight{2r}{j}.
\end{aligned}
\end{equation}
\end{restatable}

We are now ready to bound the bias due to the peeling weights.
\begin{proof}[Proof of \cref{lem:peeling:bias-bound}]
For \Cref{lem:peeling:bias-identity}, fix $R\subsetneq[k]$, put
$Q=[k]\setminus R$, and let $\psi(\Covariate,\Covariate')$ be integrable,
where $\Covariate'$ is independent of the covariate and noises.
We have
\begin{align*}
 &\E_{\Covariate,\NoiseVector,\Covariate'}\!\left[
  \psi\bigl(\red{\omega_R(\Covariate,Y_{\max})}
            -\red{\TailIndicator{Y_{\max,Q}}}\bigr)\right]\\
 &\quad=\E_{\Covariate,\NoiseVector,\Covariate'}\!\left[
  \psi\left(\begin{aligned}
   &\TailIndicator{Y_{\max}}-\TailIndicator{Y_{\max,R}}\\
   &\quad-\TailIndicator{Y_{\max,Q}}
  \end{aligned}\right)\right]
 \alignnote{conditioning on $\Covariate,\Covariate'$}\\
 &\quad=-\E_{\Covariate,\NoiseVector,\Covariate'}\!\left[
  \psi\,\TailIndicator{\min\{Y_{\max,R},Y_{\max,Q}\}}\right].
 \alignnote{the max--min identity}
\end{align*}
This proves the identity, entrywise for a matrix feature.

For \Cref{lem:peeling:cross-block-bias-bound}, we will bound a superset of the indicator event.
Indeed, observe that for $\TailIndicator{\min\{Y_{\max,R},Y_{\max,Q}\}}=1$,
it must be that some $i\in R$ and $j\in Q$ satisfy $Y_i > T$ and some $Y_j \geq Y_i > T$,
or vice versa. 
Hence
\[
  \TailIndicator{\min\{Y_{\max,R},Y_{\max,Q}\}}
  \leq \sum_{i\in R, j\in Q} \1\{Y_i>T, Y_j>T\}\,.
\]
Suppose $\psi$ is a real polynomial
of degree at most $2s+1$, $\Covariate'$ is standard Gaussian, and
$T\geq4BK$ and $(T-B)^2\geq Ks$. Under the boundedness and uncoveredness assumptions,
\Cref{lem:peeling:pairwise-exceedance} applies so that
\begin{align*}
 &\abs{\E_{\Covariate,\NoiseVector,\Covariate'}\!\left[
  \psi(\Covariate,\Covariate')
  \TailIndicator{\min\{Y_{\max,R},Y_{\max,Q}\}}\right]}\\
 &\quad\leq\sum_{i\in R,\,j\in Q}
 \E_{\Covariate,\xi_i,\xi_j,\Covariate'}\!\left[
  |\psi(\Covariate,\Covariate')|\1\{Y_i>T,Y_j>T\}\right]
 \alignnote{union bound}\\
 &\quad\leq e^{4s\log K-T^2/(2K)}
 \norm{\psi}_{L^2(\Covariate,\Covariate')}
 |R|\sum_{j\in Q}\ComponentWeight{2r}{j}
 \alignnote{by \Cref{lem:peeling:pairwise-exceedance}}\\
 &\quad\leq k e^{4s\log K-T^2/(2K)}
 \norm{\psi}_{L^2(\Covariate,\Covariate')}\tr(\EvenMomentMatrix{r}{Q}).
 \alignnote{$|R|\leq k$}
\end{align*}
This proves the bound.
\end{proof}

\subsection{Proof of \texorpdfstring{\Cref*{thm:peeling:hermite-moment-estimation}}{Theorem 8.1}}
\label{subsec:peeling:moment-estimation-proof}
We have bounded the sampling error, the error from the recovered
regressors, and the peeling bias. We combine these bounds for the even
estimate and then for every contraction of the odd estimate.

\begin{proof}[Proof of \cref{thm:peeling:hermite-moment-estimation}]
Fix a projection map $U$ with orthonormal rows and a list
$\widehat{\mathcal R}$ of recovered slope and bias estimates corresponding to regressors indexed by $R$.
Each pair satisfies
$\norm{\widehat{\Regressor}_i-\Regressor_i}_2+|\widehat b_i-b_i|\leq e^{-K^6s}$
for $i\in R$. We first prove the even bound
at order $r$, then apply the same argument to each contraction of
the odd moment at $r=s$.

\emph{Target expectations.}
We first express the two target matrices as expectations of the
projected features restricted to $T\leq Y_{\max,Q}\leq T+1$.
By the definition of the restricted moments and independence of $\Covariate'$
from the covariate and latent noises, for every $\bm{\zeta}\in\R^n$,
\begin{equation}\label{eq:peeling:restricted-moment-targets}
\begin{aligned}
 \E_{\Covariate,\NoiseVector_Q,\Covariate'}[\Psi_{\even}^{U}(\Covariate,\Covariate')\TailIndicator{Y_{\max,Q}}]
 &=U\Mat_r(\RestrictedMomentTensor{2r}{Q})U^\top,\\
 \E_{\Covariate,\NoiseVector_Q,\Covariate'}\!\left[
  \left(\sum_{\ell=1}^n \zeta_\ell\Psi_{\odd,\ell}^{U}(\Covariate,\Covariate')\right)
  \TailIndicator{Y_{\max,Q}}
 \right]
 &=U
   \Mat_r(\RestrictedMomentTensor{2r+1}{Q}\mathbin{\lrcorner}\bm{\zeta})U^\top.
\end{aligned}
\end{equation}

\emph{Even-moment error.}
We first compare $\widehat M$ with the expectation formed from
the estimated peeling weight. We then replace the estimated weight
by the true weight and the true weight by the slice indicator.
By \eqref{eq:peeling:restricted-moment-targets}, the triangle inequality gives
\begin{align}
 &\norm{\widehat{M}-U\Mat_r(\RestrictedMomentTensor{2r}{Q})U^\top}_F
 \notag\\
 &\quad\leq
 \norm{\widehat{M}-
  \E_{\Covariate,Y_{\max},\Covariate'}[
   \widehat{\omega}_{\widehat{\mathcal R}}(\Covariate,Y_{\max})
   \Psi_{\even}^{U}(\Covariate,\Covariate')]
 }_F
 \label{eq:peeling:even-moment-sampling}\\
 &\qquad+
 \norm{\E_{\Covariate,Y_{\max},\Covariate'}\!\left[
  \bigl(\red{\widehat{\omega}_{\widehat{\mathcal R}}(\Covariate,Y_{\max})}
        -\red{\omega_R(\Covariate,Y_{\max})}\bigr)
  \Psi_{\even}^{U}(\Covariate,\Covariate')
 \right]}_F
 \label{eq:peeling:even-moment-regressor-error}\\
 &\qquad+
 \left\|\begin{aligned}
  &\E_{\Covariate,Y_{\max},\Covariate'}[
    \red{\omega_R(\Covariate,Y_{\max})}\Psi_{\even}^{U}(\Covariate,\Covariate')]\\
  &\quad-\E_{\Covariate,\NoiseVector_Q,\Covariate'}[
    \Psi_{\even}^{U}(\Covariate,\Covariate')\red{\TailIndicator{Y_{\max,Q}}}]
 \end{aligned}\right\|_F.
 \label{eq:peeling:even-moment-bias}
\end{align}

\emph{Sampling error in \eqref{eq:peeling:even-moment-sampling}.}
The first term is the sampling error bounded in
\Cref{lem:peeling:projected-moment-sampling-even}. Thus it is at most
$\varepsilon_{\mathrm{est}}$ with probability at least $1-\delta$.

\emph{Recovered-regressor error in
\eqref{eq:peeling:even-moment-regressor-error}.}
The remaining two errors are matrix-valued, while the
recovered-regressor and peeling-bias bounds that we proved apply to scalar features.
For an error matrix $D$, we use the following Frobenius duality fact
\[
 \norm{D}_F=\sup_{\norm{H}_F\leq1}\abs{\tr(H^\top D)}.
\]
Thus, we fix $H\in\R^{p\times p}$ with $\norm{H}_F\leq1$
and bound the pairing uniformly over $H$.
For the even feature, put
\[
 \bar{\psi}(\Covariate)\coloneqq
 \E_{\Covariate'}[\tr(H^\top\Psi_{\even}^{U}(\Covariate,\Covariate'))\mid \Covariate].
\]
By \Cref{prop:peeling:projected-feature-second-moments},
\[
 \deg\bar{\psi}\leq2r\leq2s,
 \qquad
 \norm{\bar{\psi}}_{L^2(\Covariate)}\leq1.
\]
The signed peeling weights are independent of $\Covariate'$ conditional on
$(\Covariate,Y_{\max})$, so the bound from the recovered estimates
(\Cref{lem:peeling:weight-stability}) gives
\begin{align*}
 &\abs{\E_{\Covariate,Y_{\max},\Covariate'}\!\left[
  \tr(H^\top\Psi_{\even}^{U}(\Covariate,\Covariate'))
  \bigl(\widehat{\omega}_{\widehat{\mathcal R}}(\Covariate,Y_{\max})
        -\omega_R(\Covariate,Y_{\max})\bigr)
 \right]}\\
 &\quad=\abs{\E_{\Covariate,Y_{\max}}\!\left[
  \bar{\psi}(\Covariate)
  \bigl(\widehat{\omega}_{\widehat{\mathcal R}}(\Covariate,Y_{\max})
        -\omega_R(\Covariate,Y_{\max})\bigr)
 \right]}
 \alignnote{averaging over $\Covariate'$}\\
 &\quad\leq\frac{\norm{\bar{\psi}}_{L^2(\Covariate)}}{\sqrt\pi}
 \left(
  \sum_{i\in R}
  \bigl(\norm{\widehat{\Regressor}_i-\Regressor_i}_2+|\widehat b_i-b_i|\bigr)^2
 \right)^{1/2}
 \alignnote{by \Cref{lem:peeling:weight-stability}}\\
 &\quad\leq\frac{\sqrt{k}}{\sqrt\pi}e^{-K^6s}.
 \alignnote{the hypotheses; $\norm{\bar{\psi}}_{L^2(\Covariate)}\leq1$, $|R|\leq k$}
\end{align*}
Taking the supremum over $H$ gives the bound on
\eqref{eq:peeling:even-moment-regressor-error}.

\emph{Signed-weight bias in \eqref{eq:peeling:even-moment-bias}.}
It remains to bound the bias from using the true peeling weight
in place of the slice indicator. 
We will again use Frobenius duality.
Let $H$ be a unit Frobenius norm matrix and $\bar\psi$ be as above.
\begin{align*}
 &\abs{\begin{aligned}
  &\E_{\Covariate,Y_{\max}}[\bar{\psi}(\Covariate)\omega_R(\Covariate,Y_{\max})]\\
  &\quad-\E_{\Covariate,\NoiseVector_Q}[\bar{\psi}(\Covariate)\TailIndicator{Y_{\max,Q}}]
 \end{aligned}}\\
 &\quad\leq
 k e^{4s\log K-T^2/(2K)}
 \norm{\bar{\psi}}_{L^2(\Covariate)}\tr(\EvenMomentMatrix{r}{Q})
 \alignnote{by \Cref{lem:peeling:bias-identity,lem:peeling:cross-block-bias-bound}}\\
 &\quad\leq k e^{4s\log K-T^2/(2K)}\tr(\EvenMomentMatrix{r}{Q}).
 \alignnote{$\norm{\bar{\psi}}_{L^2(\Covariate)}\leq1$}
\end{align*}
Taking the same supremum bounds \eqref{eq:peeling:even-moment-bias}.
Adding the three bounds gives, with probability at least $1-\delta$,
\begin{align*}
 \norm{\widehat{M}-U\Mat_r(\RestrictedMomentTensor{2r}{Q})U^\top}_F
 &\leq\varepsilon_{\mathrm{est}}
  +\frac{\sqrt{k}}{\sqrt\pi}e^{-K^6s}
 +k e^{4s\log K-T^2/(2K)}\tr(\EvenMomentMatrix{r}{Q})\,.
\end{align*}
This proves \Cref{thm:peeling:hermite-moment-estimation-even}.

\emph{Odd-moment error.}
At $r=s$, we apply the even-case decomposition to each contraction
of the odd moment. For $\bm{\zeta}\in\R^n$, the estimator is
$\sum_{\ell=1}^n \zeta_\ell\widehat M_\ell$, and the second identity in
\eqref{eq:peeling:restricted-moment-targets} gives its estimation target.

By \Cref{lem:peeling:projected-moment-sampling-odd}, the sampling error is at most
$\varepsilon_{\mathrm{est}}\norm{\bm{\zeta}}_2$ simultaneously for
every $\bm{\zeta}\in\R^n$, with probability at least $1-\delta$.

For the two deterministic errors, fix $\bm{\zeta}\in\R^n$. As in the
even case, we pair each error matrix with an arbitrary
$H\in\R^{p\times p}$ satisfying $\norm{H}_F\leq1$. The scalar
feature for this contraction is
\[
 \bar{\psi}_{\bm{\zeta}}(\Covariate)\coloneqq
 \E_{\Covariate'}\!\left[
  \tr\!\left(H^\top
   \sum_{\ell=1}^n \zeta_\ell\Psi_{\odd,\ell}^{U}(\Covariate,\Covariate')\right)
  \mid \Covariate
 \right].
\]
By \Cref{prop:peeling:projected-feature-second-moments},
\[
 \deg\bar{\psi}_{\bm{\zeta}}\leq2s+1,
 \qquad
 \norm{\bar{\psi}_{\bm{\zeta}}}_{L^2(\Covariate)}\leq\norm{\bm{\zeta}}_2.
\]

For the recovered-regressor error, \Cref{lem:peeling:weight-stability} yields the following bound
\[
 \frac{\norm{\bar\psi_{\bm{\zeta}}}_{L^2(\Covariate)}}{\sqrt\pi}
 \left(\sum_{i\in R}
 \bigl(\norm{\widehat{\Regressor}_i-\Regressor_i}_2+|\widehat b_i-b_i|\bigr)^2\right)^{1/2}
 \leq \frac{\sqrt{k}}{\sqrt\pi}e^{-K^6s}\norm{\bm{\zeta}}_2\,.
\]

For the peeling bias, $\deg\bar\psi_{\bm{\zeta}}\leq2s+1$, so
\Cref{lem:peeling:bias-identity,lem:peeling:cross-block-bias-bound} yield the following bound
\[
 k e^{4s\log K-T^2/(2K)}\tr(\EvenMomentMatrix{s}{Q})
 \norm{\bar\psi_{\bm{\zeta}}}_{L^2(\Covariate)}
 \leq k e^{4s\log K-T^2/(2K)}
 \tr(\EvenMomentMatrix{s}{Q})\norm{\bm{\zeta}}_2.
\]

Both bounds hold for every $H$ with $\norm{H}_F\leq1$.
Taking the supremum in each and summing all error bounds gives
\begin{align*}
 \left\|
  \sum_{\ell=1}^n \zeta_\ell\widehat{M}_\ell
  -U\Mat_s(\RestrictedMomentTensor{2s+1}{Q}\mathbin{\lrcorner}\bm{\zeta})U^\top
 \right\|_F
 &\leq\left(
  \varepsilon_{\mathrm{est}}
  +\frac{\sqrt{k}}{\sqrt\pi}e^{-K^6s}
  +k e^{4s\log K-T^2/(2K)}\tr(\EvenMomentMatrix{s}{Q})
 \right)\norm{\bm{\zeta}}_2.
\end{align*}
Thus the estimator satisfies \Cref{thm:peeling:hermite-moment-estimation-odd}
for every $\bm{\zeta}\in\R^n$, with probability at least $1-\delta$.

\emph{Sample complexity and running time.}
The even call estimates $m=1$ matrix, and the odd call estimates
$m=n$ matrices on the same sample block. By \Cref{lem:peeling:weighted-mean-estimation},
each call uses
\[
 \Ot\!\left(\frac{p^2m}{\varepsilon_{\mathrm{est}}^2}
 100^s\log(\nicefrac1\delta)\right)
\]
projected samples. Evaluating the features and weights and computing the
trimmed means takes the sample count times $\poly(n,k,p,m)$
operations. Since $s=K^4\log n$, this gives the stated running time.
\end{proof}

This gives the moment estimation oracle used at each recursive step of \Cref{subsec:tensor-decomposition:iterative-projections}.

\section{Proof of \texorpdfstring{\Cref*{thm:intro:polynomial-time-recovery}}{Theorem}}\label{sec:recovery:main}
We have shown how to obtain rough estimates for a subset of normalized regressors through implicit tensor decomposition (\Cref{thm:tensor-decomposition:implicit-tensor-decomposition})
and refine these unit vectors into accurate estimates of the corresponding regressors and biases using a variant of OLS (\Cref{thm:refinement:restricted-least-squares}).
We now combine these guarantees to prove \Cref{thm:intro:polynomial-time-recovery}.
Our full algorithm simply repeats the two aforementioned subroutines until all regressors are recovered.
Its pseudocode can be found in \Cref{alg:recovery:recover-all-regressors}.

\begin{algorithm}[H]
\caption[\textsc{RecoverAllRegressors}]{\textsc{RecoverAllRegressors}}\label{alg:recovery:recover-all-regressors}
\KwIn{Known $B,\Delta,\rho$, $k$, $\varepsilon$, $\delta$, and
sample access to $(\Covariate,Y_{\max})$}
\KwOut{$k$ pairs in $\R^n\times\R$}
Use fresh observations in every call\;
Choose $K,s,T,\varepsilon_{\mathrm{int}},{\Gamma}$
as in \Cref{hyp:algorithm}\;
$\widehat{\mathcal R}\gets()$\;
\While{$|\widehat{\mathcal R}|<k$}{
  \tcp{\Cref{alg:tensor-decomposition:implicit-tensor-decomposition,thm:tensor-decomposition:implicit-tensor-decomposition}}
  $\widehat{\mathcal D}\gets$
  \ImplicitTensorDecomposition{$\widehat{\mathcal R},\frac{\delta}{2k}$}\;
  \ForEach{$\widehat{\NormalizedRegressor}\in\widehat{\mathcal D}$}{
    \tcp{\Cref{alg:refinement:restricted-least-squares,thm:refinement:restricted-least-squares}}
    $(\widehat{\Regressor},\widehat b)\gets$
    \Refine{$\widehat{\NormalizedRegressor},\varepsilon_{\mathrm{int}},
            \frac{\delta}{2k}$}\;
    Append $(\widehat{\Regressor},\widehat b)$ to $\widehat{\mathcal R}$\;
  }
}
\Return $\widehat{\mathcal R}$\;
\end{algorithm}

We now analyze the guarantees of \Cref{alg:recovery:recover-all-regressors},
which yields a proof of our main theorem.
\begin{proof}[Proof of \cref{thm:intro:polynomial-time-recovery}]
We first show that each round appends estimates of previously
unrecovered regressors with the accuracy needed for the next
decomposition call.
We then count the calls to bound the failure probability,
sample complexity, and running time.

\emph{Accuracy after each round.}
We first prove correctness when the decomposition and refinement
calls satisfy their guarantees.
At the beginning of each round,
assume we recovered a subset of regressors indexed by
$R\subseteq[k]$, with estimates $\widehat{\mathcal R}$
up to sum error $\varepsilon_{\mathrm{int}}$:
\[
 \norm{\widehat{\Regressor}_i-\Regressor_i}_2+|\widehat b_i-b_i|
 \leq\varepsilon_{\mathrm{int}}\qquad(i\in R).
\]
Here the estimates $(\widehat{\Regressor}_i,\widehat b_i)$ are indexed
according to their matching with $R$.
This vacuously holds for the initial empty list.

Fix a current list satisfying this property, with $R\ne[k]$.
By the choice of $\varepsilon_{\mathrm{int}}$ in
\Cref{hyp:algorithm},
$\varepsilon_{\mathrm{int}}\leq e^{-K^6s}$.
Thus $\widehat{\mathcal R}$ satisfies the input accuracy requirement
of \textsc{ImplicitTensorDecomposition} (\Cref{thm:tensor-decomposition:implicit-tensor-decomposition}).
By \Cref{thm:tensor-decomposition:implicit-tensor-decomposition}, \textsc{ImplicitTensorDecomposition} returns a list $\widehat{\mathcal D}$ of unit vectors matched to some normalized regressors indexed by a non-empty subset $J\subseteq[k]\setminus R$, with errors at most $\Delta^2/(16B^2)$.
Execute \textsc{RestrictedLeastSquares} on each $\widehat{\NormalizedRegressor}\in \widehat{\mathcal D}$ with requested accuracy $\varepsilon_{\mathrm{int}}$.
\Cref{thm:refinement:restricted-least-squares} guarantees that the subroutine returns $\varepsilon_{\mathrm{int}}$-accurate estimates $(\widehat{\Regressor}_i,\widehat b_i)$ in sum error for regressors indexed by $J$.

\emph{Completion and failure probability.}
Each round adds at least one previously unmatched regressor,
so there are at most $k$ calls to \textsc{ImplicitTensorDecomposition}.
Since each regressor is refined once, there are also
at most $k$ calls to \textsc{RestrictedLeastSquares}.
At termination, the matching gives a permutation $\pi$ of $[k]$ such that
\[
 \max_{i\in[k]}\left(
 \norm{\widehat{\Regressor}_i-\Regressor_{\pi(i)}}_2
 +|\widehat b_i-b_{\pi(i)}|\right)
 \leq\varepsilon_{\mathrm{int}}\leq\varepsilon.
\]
By a union bound over these at most $2k$ calls, each using fresh
observations and failure probability $\delta/(2k)$,
the algorithm
succeeds with probability at least $1-\delta$.

\emph{Sample complexity and running time.}
There are at most $k$ calls to each routine.
Each call uses failure probability $\delta/(2k)$, and each refinement
call uses accuracy $\varepsilon_{\mathrm{int}}$.
By the choice of \Cref{hyp:algorithm},
$s=K^4\log n$ and
\[
 \frac1{\varepsilon_{\mathrm{int}}}
 =\max\left\{\frac1\varepsilon,e^{K^6s}\right\}
 \leq\left(\frac n\varepsilon\right)^{\poly(K)}.
\]
By the sample complexity and running time guarantees of
\Cref{thm:tensor-decomposition:implicit-tensor-decomposition,thm:refinement:restricted-least-squares},
the total sample complexity and running time are both at most
\[
\begin{aligned}
 &\Ot\!\left(
 k\left[(nk)^{\poly(K)}
 +\left(\frac{nk}{\varepsilon_{\mathrm{int}}}\right)^{\poly(K)}
 \right]\log\!\left(\frac{2k}{\delta}\right)\right)\\
 &\qquad=
 \Ot\!\left(
 \left(\frac{nk}{\varepsilon}\right)^{\poly(K)}
 \log(\nicefrac1\delta)\right).
\end{aligned}
\]
\end{proof}

\section*{Acknowledgements}
We thank Anay Mehrotra and Manolis Zampetakis for feedback on an initial draft of this work.
Felix Zhou acknowledges the support of the Natural Sciences and Engineering Research Council of Canada (NSERC).
Part of this work was completed during an internship at Google Research.

\section*{AI Methodology}
We used \texttt{gpt-5.6-sol} extensively for brainstorming, calculations, and literature search
as well as \texttt{gpt-6-astra} and \texttt{claude-fable-5.1} to check for correctness and refine the presentation of the proofs.
The author verified the correctness and originality of all content, including references.

\printbibliography

\clearpage

\appendix
\crefalias{section}{appendix}
\crefalias{subsection}{subappendix}
\crefalias{subsubsection}{subsubappendix}

\section{Additional Preliminaries}

\subsection{Gaussian Moments and Elementary Bounds}

\begin{proposition}[Moment-Generating Function of a Gaussian Quadratic Form;
 {\cite[Theorem~5.2b]{RencherSchaalje08}}]
\label{prop:gaussian-quadratic-mgf}
Let $G_0\sim N(0,1)$, $m\in\R$, and $\tau\geq0$.
For every $a\in\R$ with $2a\tau^2<1$,
\[
 \E_{G_0}\!\left[e^{a(m+\tau G_0)^2}\right]
 =
 \frac{\exp\!\left(\dfrac{am^2}{1-2a\tau^2}\right)}
      {\sqrt{1-2a\tau^2}}.
\]
\end{proposition}

\begin{proposition}[Gaussian Noise Operator on Hermite Polynomials;
 {\cite[Definition~11.12 and Proposition~11.37]{ODonnell14}}]
\label{prop:hermite-noise-operator}
For a polynomial $f:\R\to\R$ and $c\in[-1,1]$, the Gaussian noise
operator is given by $(U_c f)(x)=\E_g[f(cx+\sqrt{1-c^2}\,g)]$,
where $g\sim N(0,1)$.
For every integer $d\geq0$, it satisfies $U_c\He_d=c^d\He_d$.
Thus, if $U,V$ are jointly standard Gaussian with correlation $c$,
then for every $x\in\R$,
\begin{equation}
 \E_{U,V}[\He_d(U)\mid V=x]=c^d\He_d(x).
\end{equation}
\end{proposition}

\begin{proposition}[Gaussian Hypercontractivity for Polynomials;
 {\cite[Theorem~11.23 and Proposition~11.37]{ODonnell14}}]
\label{prop:gaussian-hypercontractivity}
Fix integers $n\geq1$ and $d\geq0$, and a real polynomial
$P:\R^n\to\R$ of total degree at most $d$.
Let $\bm{G}\sim N(\bm{0},I_n)$.
For every $2\leq p<\infty$,
\[
 \norm{P}_{L^p(\bm{G})}
 \leq(p-1)^{d/2}\norm{P}_{L^2(\bm{G})},
 \qquad
 \norm{P}_{L^p(\bm{G})}
 \coloneqq\bigl(\E_{\bm{G}}[|P(\bm{G})|^p]\bigr)^{1/p}.
\]
\end{proposition}

The following is a Gaussian fourth-moment specialization of the
standard bounds for sub-Gaussian maxima; see
\cite[Proposition~2.5.2 and Exercise~2.5.10]{Vershynin18}.
For completeness, we provide a proof below with explicit numerical constants.

\begin{proposition}[Fourth Moment of a Gaussian Maximum]
\label{prop:gaussian-maximum-fourth-moment}
Let $\bm{G}=(G_1,\ldots,G_n)$ be a Gaussian vector, where $n\geq2$.
Then
\[
 \left(
 \E_{\bm{G}}\!\left[\max_{j\in[n]}|G_j|^4\right]
 \right)^{1/4}
 \leq\max_{j\in[n]}|\E_{\bm{G}}[G_j]|
 +3\sqrt{\log n}
 \left(\max_{j\in[n]}\operatorname{Var}_{\bm{G}}(G_j)\right)^{1/2}.
\]
\end{proposition}

\begin{proof}
We first bound the maximum after subtracting the means.
Put $p=6\log n>4$.
Each centered variable $G_j-\E_{\bm{G}}[G_j]$ has the law of a linear
polynomial in a standard Gaussian, with $L^2$ norm
$\operatorname{Var}_{\bm{G}}(G_j)^{1/2}$.
Gaussian hypercontractivity
(\Cref{prop:gaussian-hypercontractivity}) therefore gives
\begin{align*}
 &\left(
 \E_{\bm{G}}\max_{j\in[n]}|G_j-\E_{\bm{G}}[G_j]|^4
 \right)^{1/4}\\
 &\quad\leq\left(
 \E_{\bm{G}}\max_{j\in[n]}|G_j-\E_{\bm{G}}[G_j]|^p
 \right)^{1/p}
 \alignnote{$p>4$}\\
 &\quad\leq\left(
 \sum_{j\in[n]}\E_{\bm{G}}|G_j-\E_{\bm{G}}[G_j]|^p
 \right)^{1/p}\\
 &\quad\leq n^{1/p}\sqrt{p-1}
 \left(\max_{j\in[n]}\operatorname{Var}_{\bm{G}}(G_j)\right)^{1/2}
 \alignnote{by \Cref{prop:gaussian-hypercontractivity}}\\
 &\quad\leq e^{1/6}\sqrt{6\log n}
 \left(\max_{j\in[n]}\operatorname{Var}_{\bm{G}}(G_j)\right)^{1/2}
 \alignnote{$n^{1/p}=e^{1/6}$; $p=6\log n$}\\
 &\quad\leq3\sqrt{\log n}
 \left(\max_{j\in[n]}\operatorname{Var}_{\bm{G}}(G_j)\right)^{1/2}.
\end{align*}
Applying the triangle inequality inside the maximum and then
Minkowski's inequality gives
\[
 \left(\E_{\bm{G}}\max_{j\in[n]}|G_j|^4\right)^{1/4}
 \leq\max_{j\in[n]}|\E_{\bm{G}}[G_j]|
 +\left(
 \E_{\bm{G}}\max_{j\in[n]}|G_j-\E_{\bm{G}}[G_j]|^4
 \right)^{1/4}.
\]
Combining these bounds proves the claim.
\end{proof}

\begin{proposition}[Truncated Gaussian Second Moment;
 {\cite[Theorem~5.8(6),(8)]{Hansen22}}]
\label{prop:truncated-gaussian-second-moment}
Let $G$ be a nondegenerate Gaussian with $\E_G[G]\leq0$.
Then
\[
 \E_G[G^2\1_{\{G>0\}}]
 \leq\operatorname{Var}_G(G)\Prb_G(G>0).
\]
\end{proposition}

\begin{proposition}[Spherical Anti-Concentration; {\cite[Lemma~4]{DKM09}}]
\label{prop:spherical-anticoncentration}
Let $\bm{\zeta}$ be uniform on the unit sphere $\mathbb S^{n-1}$.
For every nonzero $\bm{v}\in\R^n$ and $0<t<1$,
\[
 \Prb_{\bm{\zeta}}\!\left(
  \abs{\ip {\bm{\zeta}} {\bm{v}}}\leq\frac{t\norm {\bm{v}}_2}{\sqrt n}
 \right)\leq t.
\]
\end{proposition}

\begin{proposition}[Wick's Identity~\cite{Isserlis18}]
\label{prop:wick-identity}
Let $\bm{G}=(G_1,\ldots,G_n)$ be a centered jointly Gaussian vector,
where $n\geq1$.
If $n$ is odd, then $\E_{\bm{G}}[\prod_{a=1}^nG_a]=0$.
If $n=2q$, then
\[
 \E_{\bm{G}}\!\left[\prod_{a=1}^{2q}G_a\right]
 =\sum_{P\in\mathcal P_2(2q)}
   \prod_{\{a,b\}\in P}\E_{\bm{G}}[G_aG_{b}],
\]
where $\mathcal P_2(2q)$ is the set of partitions of $[2q]$
into unordered pairs.
\end{proposition}

By multilinearity, the same formula holds for complex linear forms
in a real Gaussian vector.  Conjugated forms are treated as separate
factors in the product.

\begin{proposition}[Joint Gaussian Tail Bound]
\label{prop:joint-gaussian-tail}
Fix a centered jointly Gaussian pair $(G_1,G_2)$ with positive
variances $\sigma_1^2,\sigma_2^2$ and covariance $c$ satisfying
\[
 D\coloneqq\sigma_1^2\sigma_2^2-c^2>0,
 \qquad
 c\leq\min\{\sigma_1^2,\sigma_2^2\}.
\]
For every $t>0$,
\[
 \Prb(G_1>t,G_2>t)
 \leq\exp\!\left(
  -\frac{t^2(\sigma_1^2+\sigma_2^2-2c)}{2D}
 \right).
\]
\end{proposition}

\begin{proof}
Fix $t>0$ and put
\[
 H\coloneqq
 (\sigma_2^2-c)G_1+(\sigma_1^2-c)G_2.
\]
The coefficients are nonnegative by the hypothesis on $c$.
They cannot both vanish, since that would give
$\sigma_1^2=\sigma_2^2=c$ and $D=0$.
Consequently, $\sigma_1^2+\sigma_2^2-2c>0$, and
$G_1>t,G_2>t$ implies
$H>t(\sigma_1^2+\sigma_2^2-2c)$.
The covariance and variance calculations give
\begin{align*}
 \operatorname{Cov}(G_1,H)
 &=(\sigma_2^2-c)\sigma_1^2+(\sigma_1^2-c)c=D,\\
 \operatorname{Cov}(G_2,H)
 &=(\sigma_2^2-c)c+(\sigma_1^2-c)\sigma_2^2=D,\\
 \operatorname{Var}(H)
 &=(\sigma_2^2-c)\operatorname{Cov}(G_1,H)
   +(\sigma_1^2-c)\operatorname{Cov}(G_2,H)\\
 &=D(\sigma_1^2+\sigma_2^2-2c)>0.
\end{align*}
The Gaussian tail bound (\Cref{prop:gaussian-tail}) gives
\begin{align*}
 \Prb(G_1>t,G_2>t)
 &\leq\Prb\!\left(H>t(\sigma_1^2+\sigma_2^2-2c)\right)\\
 &\leq\exp\!\left(
  -\frac{t^2(\sigma_1^2+\sigma_2^2-2c)^2}
        {2\operatorname{Var}(H)}
 \right)
 \alignnote{Gaussian tail bound}\\
 &=\exp\!\left(
  -\frac{t^2(\sigma_1^2+\sigma_2^2-2c)}{2D}
 \right).
\end{align*}
This proves the claim.
\end{proof}

\subsection{Linear Algebra}

\begin{proposition}[Gram-Matrix Identities; Follows from {\cite[Section~I.2]{Bhatia97}}]
\label{prop:gram-matrix-identities}
For every real matrix $B$, the matrices $B^\top B$ and
$BB^\top$ have the same nonzero eigenvalues, with multiplicity.
Moreover,
\[
 \norm{B^\top B}_{\op}
 =\norm{BB^\top}_{\op}
 =\norm {B}_{\op}^2\,.
\]
\end{proposition}

\begin{proposition}[Polar-Factor Bound; {\cite[Theorem~1.1 and Corollary~2.3]{Higham86}}]
\label{prop:polar-factor-bounds}
Let $B\in\R^{p\times p}$, and let $Z$ be an orthogonal
polar factor of $B$.  Writing $\sigma_h(B)$ for the
singular values, we have
\[
 \norm{B-Z}_F^2
 =\sum_{h=1}^p(\sigma_h(B)-1)^2
 \leq\norm{BB^\top-I}_F^2.
\]
\end{proposition}

\begin{proposition}[Sylvester's Equation;
 {\cite[Theorem~VII.2.8]{Bhatia97}}]
\label{prop:sylvester-equation}
Let $A\in\R^{p\times p}$ and $B\in\R^{q\times q}$ be symmetric.
Suppose $|\lambda_h(A)-\lambda_i(B)|\geq g>0$
for every $h\in[p]$ and $i\in[q]$.
If $AX-XB=C$, then
\[
 \norm{X}_F\leq\frac{\norm{C}_F}{g}.
\]
\end{proposition}

We use the Davis--Kahan sine theorem in the following form.
\begin{proposition}[Davis--Kahan Eigenvector Bound;
 {\cite[Corollary~1]{YuWangSamworth2015}}]
\label{prop:davis-kahan-eigenvector}
Let $H,\widetilde H\in\R^{p\times p}$ be symmetric, with eigenvalues
in nonincreasing order.
Fix $j\in[p]$ and $g>0$ such that
\[
 \abs{\lambda_j(H)-\lambda_h(H)}\geq g
 \qquad(h\in[p]\setminus\{j\}).
\]
Let $\bm{x},\widetilde{\bm{x}}$ be unit eigenvectors corresponding to
$\lambda_j(H)$ and $\lambda_j(\widetilde H)$, respectively.
The sign of $\widetilde{\bm{x}}$ can be chosen so that
\[
 \norm{\widetilde{\bm{x}}-\bm{x}}_2
 \leq\frac{2\sqrt2\,\norm{\widetilde H-H}_{\op}}{g}.
\]
\end{proposition}

\section{Deferred Proofs for \texorpdfstring{\Cref{sec:preliminaries:main}}{Section 4}}

\subsection{Hermite Polynomials}
\label{subsec:app-preliminaries-proofs:hermite-polynomials}

\begin{proof}[Proof of \cref{prop:hermite-expansion-normalization}]
We first verify orthogonality for products of normalized Hermite
polynomials over independent Gaussian coordinates.
Let $\bm{G}=(G_1,\ldots,G_n)\sim N(\bm{0},I_n)$, and let
$a_1,\ldots,a_n,b_1,\ldots,b_n$ be nonnegative integers.
Then
\begin{equation}
\begin{aligned}
 &\E_{\bm{G}}\!\left[
  \left(\prod_{j=1}^n
   \frac{\He_{a_j}(G_j)}{\sqrt{a_j!}}\right)
  \left(\prod_{j=1}^n
   \frac{\He_{b_j}(G_j)}{\sqrt{b_j!}}\right)
 \right]\\
 &\quad=
 \prod_{j=1}^n
 \E_{G_j}\!\left[
  \frac{\He_{a_j}(G_j)\He_{b_j}(G_j)}
       {\sqrt{a_j!b_j!}}
 \right]\\
 &\quad=
 \begin{cases}
  1,&a_j=b_j\text{ for every }j,\\
  0,&\text{otherwise}.
 \end{cases}
\end{aligned}
\end{equation}
The first equality uses independence.
The second applies Hermite orthogonality
(\Cref{subsec:preliminaries:hermite-polynomials}) in each coordinate.

For a symmetric tensor $\mathsf A$, use the coordinate counts
$\bm{\alpha}$ from \Cref{subsec:preliminaries:hermite-polynomials}.
The complex Gaussian definition and independence of the coordinates
of $\Covariate'$ give, for indices $j_1,\ldots,j_d$ with counts $\bm{\alpha}$,
\begin{align*}
 [\HermiteTensor{d}{\bm{x}}]_{j_1,\ldots,j_d}
 &=\frac1{\sqrt{d!}}
   \RealPart\!\left(
    \E_{\Covariate'}\left[\prod_{\ell=1}^d(x_{j_\ell}+\iu X'_{j_\ell})\right]
   \right)\\
 &=\frac1{\sqrt{d!}}
   \RealPart\!\left(
    \prod_{j=1}^n\E_{X'_j}[(x_j+\iu X'_j)^{\alpha_j}]
   \right)
 \alignnote{independence}\\
 &=\frac1{\sqrt{d!}}\prod_{j=1}^n\He_{\alpha_j}(x_j).
 \alignnote{definition of $\He_{\alpha_j}$}
\end{align*}
The last product is real.

Grouping the tensor entries gives the expansion in
\Cref{prop:hermite-expansion-normalization}.
Product orthogonality then gives the stated isometry.
\end{proof}

\section{Deferred Proofs for \texorpdfstring{\Cref{sec:hermite:main}}{Section 5}}

\subsection{Weights and Scaled Regressors}

We specialize the norm and separation bounds of \Cref{thm:hermite:structure-scaled-regressors}
to the hyperparameter choice used in \textsc{ImplicitTensorDecomposition}.

\begin{proposition}[Geometry of the Scaled Regressors]
\label{prop:app-hermite:scaled-regressor-geometry}
Under the hypotheses and notation of \Cref{thm:hermite:restricted-moment-structure},
the scaled regressors satisfy
\begin{equation}\label{eq:app-hermite:scaled-regressor-geometry}
 K^{-3}\leq\norm{\ScaledRegressor_i}_2\leq K^3\quad(i\in[k]),\qquad
 \norm{\ScaledRegressor_i-\ScaledRegressor_j}_2\geq K^{-3}\quad(i\neq j).
\end{equation}
\end{proposition}

\begin{proof}
By the choice of \Cref{hyp:algorithm},
\begin{align}
 K&\geq\frac{128B^2}{\Delta^2}
 \geq\max\left\{3B,\frac{\sqrt5(1+B^2)}{\Delta}\right\},
 \label{eq:app-hermite:scaled-geometry-k-choice}\\
 2s+1&\leq\frac{5T^2}{4K^4},
 \label{eq:app-hermite:scaled-geometry-degree-upper}\\
 2s+1&\geq\frac{T^2}{4K^4},
 \label{eq:app-hermite:scaled-geometry-degree-lower}\\
 T+1+B&\leq\frac{3T}{2}.
 \label{eq:app-hermite:scaled-geometry-threshold}
\end{align}
The first inequality follows from
$(1+B^2)^3\geq8B^2$ and $\Delta\leq1$.
The degree bounds use $T^2=4K^4s$ and $s\geq1$.
Finally, $T\geq2K^2\geq2(B+1)$ gives
\eqref{eq:app-hermite:scaled-geometry-threshold}.
The scale bounds in \Cref{thm:hermite:structure-scaled-regressors} give
\begin{align*}
 \norm{\ScaledRegressor_i}_2
 &\geq\frac{T\Delta}{2(1+B^2)\sqrt{2s+1}}
 \alignnote{by \Cref{thm:hermite:structure-scaled-regressors,ass:boundedness-uncoveredness}}\\
 &\geq\frac{K^2\Delta}{\sqrt5(1+B^2)}
 \alignnote{choice of $T,s$~\eqref{eq:app-hermite:scaled-geometry-degree-upper}}\\
 &\geq K,
 \alignnote{choice of $K$~\eqref{eq:app-hermite:scaled-geometry-k-choice}}\\
 \norm{\ScaledRegressor_i}_2
 &\leq\frac{(T+1+B)B}{\sqrt{2s+1}}
 \alignnote{by \Cref{thm:hermite:structure-scaled-regressors,ass:boundedness-uncoveredness}}\\
 &\leq3K^2B
 \alignnote{choice of $T,s$~\eqref{eq:app-hermite:scaled-geometry-threshold},
 \eqref{eq:app-hermite:scaled-geometry-degree-lower}}\\
 &\leq K^3.
 \alignnote{choice of $K$~\eqref{eq:app-hermite:scaled-geometry-k-choice}}
\end{align*}
For distinct $i,j$, the separation bound in
\Cref{thm:hermite:structure-scaled-regressors} now gives
\begin{align*}
 \norm{\ScaledRegressor_i-\ScaledRegressor_j}_2^2
 &\geq K^{-1}\max\{\norm{\ScaledRegressor_i}_2^2,\norm{\ScaledRegressor_j}_2^2\}
 \alignnote{by \Cref{thm:hermite:structure-scaled-regressors}}\\
 &\geq K\geq K^{-6}.
 \alignnote{$\norm{\ScaledRegressor_i}_2\geq K$ and $K\geq1$}
\end{align*}
Taking square roots, and using $K\geq K^{-3}$ for the norm lower
bound, proves \eqref{eq:app-hermite:scaled-regressor-geometry}.
\end{proof}

\subsection{Single-Response Moments}
\label{subsec:app-hermite:single-response-moments}

We prove the conditional Hermite identity in \Cref{prop:hermite:conditional-hermite-identity}
by applying the scalar Gaussian noise identity
(\Cref{prop:hermite-noise-operator}) to projections of the covariate.

\conditionalhermiteidentity*

\begin{proof}
We apply the Gaussian noise identity (\Cref{prop:hermite-noise-operator})
to scalar projections of $\bm{G}$.

The case $d=0$ is immediate.
Fix $d\geq1$ and $\bm{\theta}\ne0$.
The complex Gaussian definition of the Hermite tensor
(\Cref{subsec:preliminaries:hermite-polynomials}) gives, for every $\bm{x}\in\R^n$,
\begin{align*}
 \HermitePairing{d}{\bm{x}}{\TensorPower{\bm{\theta}}{d}}
 &=\frac1{\sqrt{d!}}
   \E_{\Covariate'}[\RealPart(\ip{\bm{\theta}}{\bm{x}+\iu \Covariate'}^d)]\\
 &=\frac{\norm{\bm{\theta}}_2^d}{\sqrt{d!}}
   \He_d\!\left(\frac{\ip{\bm{\theta}}{\bm{x}}}{\norm{\bm{\theta}}_2}\right).
 \alignnote{$\ip{\bm{\theta}}{\Covariate'}/\norm{\bm{\theta}}_2\sim N(0,1)$}
\end{align*}
The variables $\ip{\bm{\theta}}{\bm{G}}/\norm{\bm{\theta}}_2$ and $Y/\sigma$
are standard Gaussians with correlation
$\ip{\bm{\theta}}{\bm{v}}/(\sigma\norm{\bm{\theta}}_2)$.
Thus
\begin{align*}
 \E_{\bm{G},Y}\!\left[
  \HermitePairing{d}{\bm{G}}{\TensorPower{\bm{\theta}}{d}}\mid Y=t
 \right]
 &=\frac{\norm{\bm{\theta}}_2^d}{\sqrt{d!}}
   \E_{\bm{G},Y}\!\left[
    \He_d\!\left(\frac{\ip{\bm{\theta}}{\bm{G}}}{\norm{\bm{\theta}}_2}\right)
    \biggm|Y=t
   \right]\\
 &=\frac{\norm{\bm{\theta}}_2^d}{\sqrt{d!}}
   \left(\frac{\ip{\bm{\theta}}{\bm{v}}}{\sigma\norm{\bm{\theta}}_2}\right)^d
   \He_d(t/\sigma)
 \alignnote{by \Cref{prop:hermite-noise-operator}}\\
 &=\frac{\He_d(t/\sigma)}{\sigma^d\sqrt{d!}}
   \ip{\bm{\theta}}{\bm{v}}^d.
\end{align*}
Both tensors are symmetric, and the tensors $\TensorPower{\bm{\theta}}{d}$
span the symmetric order-$d$ tensor space.
This proves the identity.
\end{proof}

\subsection{Moment Decomposition Error}
\label{subsec:app-hermite:moment-decomposition-error}

We prove the conditional second-moment bound used in \Cref{subsubsec:hermite:moment-decomposition-error}.
Conditioning changes only one Gaussian coordinate.
Orthogonality in the other coordinates reduces the bound to
scalar moments, which we control by the complex Gaussian definition
and Jensen's inequality.

\conditionalhermitesecondmoment*

\begin{proof}
Choose orthonormal coordinates in which $\bm{v}=\norm {\bm{v}}_2\bm e_1$;
if $\bm{v}=0$, any orthonormal coordinates suffice.
Express $\bm{G}$ and $\mathsf A$ in these coordinates.
Rotational invariance of the auxiliary Gaussian in the definition
of $\mathsf H_d$ shows that the tensor pairing and Frobenius norm
are unchanged.
The conditional Gaussian formula
(\Cref{prop:gaussian-conditioning}) gives
\begin{align*}
 \bm{G}\mid Y=t
 &\sim N\!\left(\mu \bm e_1,\operatorname{Diag}(\tau^2,1,\ldots,1)\right),\\
 \mu&=\frac{\norm {\bm{v}}_2t}{\sigma^2},\\
 \tau^2&=1-\frac{\norm {\bm{v}}_2^2}{\sigma^2}\in[0,1].
\end{align*}
Under this conditional law, write $G_1=\mu+\tau G_0$,
where $G_0\sim N(0,1)$.
The remaining coordinates $G_2,\ldots,G_n$ are independent
standard Gaussians, independent of $G_0$.

We first use this independence to reduce the tensor second moment
to scalar moments.
Squaring the expansion in \Cref{prop:hermite-expansion-normalization} gives the
following double sum.  Orthogonality forces $\alpha_j=\beta_j$
for $j\geq2$; the common total degree then gives
$\alpha_1=\beta_1$.
\begin{align}
 \E_{\bm{G},Y}[\HermitePairing{d}{\bm{G}}{\mathsf A}^2\mid Y=t]
 &=
 \sum_{\substack{|\bm{\alpha}|=d\\|\bm{\beta}|=d}}
 \frac{d!\,A_{\bm{\alpha}} A_{\bm{\beta}}}{\sqrt{\bm{\alpha}!\bm{\beta}!}}\,
 \E_{G_0}\!\left[
  \frac{
   \He_{\alpha_1}(\mu+\tau G_0)
   \He_{\beta_1}(\mu+\tau G_0)}
  {\sqrt{\alpha_1!\beta_1!}}
 \right]
 \notag\\
 &\qquad\qquad{}\times
 \prod_{j=2}^n
 \E_{G_j}\!\left[
  \frac{\He_{\alpha_j}(G_j)\He_{\beta_j}(G_j)}
       {\sqrt{\alpha_j!\beta_j!}}
 \right]
 \notag\\
 &=
 \sum_{|\bm{\alpha}|=d}\frac{d!}{\bm{\alpha}!}A_{\bm{\alpha}}^2
 \E_{G_0}\!\left[
  \frac{\He_{\alpha_1}(\mu+\tau G_0)^2}{\alpha_1!}
 \right]
 \alignnote[.6\linewidth]{Hermite orthogonality~\eqref{eq:hermite-orthogonality};
  $|\bm{\alpha}|=|\bm{\beta}|=d$}
 \notag\\
 &\leq
 \norm{\mathsf A}_F^2
 \max_{0\leq\ell\leq d}
 \frac{\E_{G_0}[\He_\ell(\mu+\tau G_0)^2]}{\ell!}.
 \alignnote{by \eqref{eq:hermite-isometry}}
 \label{eq:app-hermite:conditional-hermite-scalar-reduction}
\end{align}
It remains to bound these scalar moments.
Let $g\sim N(0,1)$ be independent of $G_0$.
We use the complex Gaussian definition of $\He_\ell$ and
Jensen's inequality, followed by the Gaussian quadratic-form
moment-generating formula
(\Cref{prop:gaussian-quadratic-mgf}).
Since $\tau^2\leq1$ and $0<\theta<1/2$, both exponential moments
below are finite.
For $0\leq\ell\leq d$,
\begin{align}
 \frac{\E_{G_0}[\He_\ell(\mu+\tau G_0)^2]}{\ell!}
 &=\frac1{\ell!}
   \E_{G_0}\!\left[
    \left|\E_g[(\mu+\tau G_0+\iu g)^\ell]\right|^2
   \right]
 \alignnote{definition of $\He_\ell$}
 \notag\\
 &\leq
 \frac{\E_{G_0,g}[((\mu+\tau G_0)^2+g^2)^\ell]}{\ell!}
 \alignnote{Jensen's inequality}
 \notag\\
 &\leq
 \theta^{-\ell}
 \E_{G_0,g}\!\left[e^{\theta((\mu+\tau G_0)^2+g^2)}\right]
 \alignnote{$x^\ell/\ell!\leq\theta^{-\ell}e^{\theta x}$ for $x\geq0$}
 \notag\\
 &=\theta^{-\ell}
   \E_{G_0}[e^{\theta(\mu+\tau G_0)^2}]
   \E_g[e^{\theta g^2}]
 \alignnote{independence}
 \notag\\
 &=\frac{\theta^{-\ell}
   \exp\!\left(\frac{\theta\mu^2}{1-2\theta\tau^2}\right)}
  {\sqrt{(1-2\theta\tau^2)(1-2\theta)}}
 \alignnote{by \Cref{prop:gaussian-quadratic-mgf}}
 \notag\\
 &\leq\frac{\theta^{-\ell}}{1-2\theta}
   \exp\!\left(\frac{\theta\mu^2}{1-2\theta}\right).
 \alignnote{$\tau^2\leq1$}
 \label{eq:app-hermite:conditional-hermite-scalar-bound}
\end{align}
Combining \eqref{eq:app-hermite:conditional-hermite-scalar-reduction} and
\eqref{eq:app-hermite:conditional-hermite-scalar-bound} gives
\begin{align*}
 \left(
  \E_{\bm{G},Y}[\HermitePairing{d}{\bm{G}}{\mathsf A}^2\mid Y=t]
 \right)^{1/2}
 &\leq
 \frac{\theta^{-d/2}}{\sqrt{1-2\theta}}
 \exp\!\left(\frac{\theta\mu^2}{2(1-2\theta)}\right)\norm{\mathsf A}_F
 \alignnote{$\theta^{-\ell}\leq\theta^{-d}$}\\
 &=\frac{\theta^{-d/2}}{\sqrt{1-2\theta}}
 \exp\!\left(
  \frac{\theta\norm {\bm{v}}_2^2t^2}{2\sigma^4(1-2\theta)}
 \right)\norm{\mathsf A}_F.
 \alignnote{definition of $\mu$}
\end{align*}
This proves the claim.
\end{proof}

\section{Deferred Proofs for \texorpdfstring{\Cref{sec:tensor-decomposition:main}}{Section 6}}

\subsection{Iterative Projections}

\begin{proposition}[Hyperparameter Bounds for Iterative Projections]
\label{prop:app-tensor-decomposition:iterative-projection-hyperparameter-bounds}
Use the hyperparameters from \Cref{hyp:algorithm} and put
$\varepsilon_{\mathrm{est}}=e^{-K^6s}$. Then
\begin{align}
 &e^{T^2}\left(1+\sqrt{\frac{k}{\pi}}\right)
   \varepsilon_{\mathrm{est}}
  +\left(kK^{4s}+K^{8s}\right)e^{-T^2/(2K)}
 \leq e^{-K^2s},
 \label{eq:app-tensor-decomposition:projection-estimation-error-budget}\\
 &\frac{n^4+s n^2k^2}{\varepsilon_{\mathrm{est}}^2}\,100^s
 \leq(nk)^{4K^{10}+12K^4}.
 \label{eq:app-tensor-decomposition:projection-sample-factor}
\end{align}
\end{proposition}

\begin{proof}
By the choice of \Cref{hyp:algorithm},
\begin{align}
 K&\geq128,
 \\
 s&\geq K^4\log n\geq7,
 \\
 s&\leq2K^4\log n,
 \label{eq:app-tensor-decomposition:projection-order-upper-bound}\\
 T^2&=4K^4s,
 \label{eq:app-tensor-decomposition:projection-threshold-choice}\\
 \varepsilon_{\mathrm{est}}&=e^{-K^6s}.
 \label{eq:app-tensor-decomposition:projection-sampling-accuracy}
\end{align}

\emph{Error bound.}
We bound the left-hand side of
\eqref{eq:app-tensor-decomposition:projection-estimation-error-budget}.
The choice of $s$ and $k\leq n$ give
$k\leq e^{s/K^4}\leq K^s$ and $\sqrt{k}\leq e^s$.
For $K\geq128$,
\begin{align}
 K^6-4K^4-1&\geq(K^2-5)K^4\geq3K^2,
 \label{eq:app-tensor-decomposition:projection-accuracy-exponent}\\
 2K^3-8\log K&\geq2K^3-8K\geq3K^2.
 \alignnote{$\log K\leq K$}
 \label{eq:app-tensor-decomposition:projection-tail-exponent}
\end{align}
Since $1+\sqrt{k/\pi}\leq2e^s$ and
$kK^{4s}+K^{8s}\leq2K^{8s}$, substituting
$T^2=4K^4s$ and $\varepsilon_{\mathrm{est}}=e^{-K^6s}$ gives
\begin{align*}
 &e^{T^2}\left(1+\sqrt{\frac{k}{\pi}}\right)
   \varepsilon_{\mathrm{est}}
 +(kK^{4s}+K^{8s})e^{-T^2/(2K)}\\
 &\quad\leq
 2e^{-(K^6-4K^4-1)s}
 +2e^{-(2K^3-8\log K)s}\\
 &\quad\leq4e^{-3K^2s}
 \alignnote{by \eqref{eq:app-tensor-decomposition:projection-accuracy-exponent},
 \eqref{eq:app-tensor-decomposition:projection-tail-exponent}}\\
 &\quad\leq e^{-K^2s}.
 \alignnote{$K\geq128$, $s\geq1$}
\end{align*}
This proves
\eqref{eq:app-tensor-decomposition:projection-estimation-error-budget}.

\emph{Sample factor.}
We now bound the expression in
\eqref{eq:app-tensor-decomposition:projection-sample-factor}.
The choice of $s$ gives $s+1\leq2K^4n$, so
\begin{align*}
 \frac1{\varepsilon_{\mathrm{est}}^2}
 &=e^{2K^6s}
 \alignnote{choice of $\varepsilon_{\mathrm{est}}$~\eqref{eq:app-tensor-decomposition:projection-sampling-accuracy}}\\
 &\leq(nk)^{4K^{10}},
 \alignnote{choice of $s$~\eqref{eq:app-tensor-decomposition:projection-order-upper-bound}}\\
 100^s(n^4+s n^2k^2)
 &\leq(s+1)(nk)^{10K^4+4}
 \alignnote{choice of $s$~\eqref{eq:app-tensor-decomposition:projection-order-upper-bound}, $\log100\leq5$}\\
 &\leq2K^4(nk)^{10K^4+5}
 \alignnote{$s+1\leq2K^4n$, $n\leq nk$}\\
 &\leq(nk)^{12K^4}.
 \alignnote{$2K^4\leq(nk)^{K^4}$, $K^4\geq5$}
\end{align*}
Multiplying these bounds proves
\eqref{eq:app-tensor-decomposition:projection-sample-factor}.
\end{proof}

\subsection{Observable-Operator Learning}
\label{subsec:app-tensor-decomposition:observable-operator-learning}

We prove the largest-gap bound used to choose the block and then
\Cref{lem:app-tensor-decomposition:spectral-truncation-orthonormality}, which bounds the errors from
spectral truncation and approximate orthonormality.

\largestgapeigenvaluebounds*

\begin{proof}
We first compare the estimated eigenvalues with those of
$\Lambda$.
Pad the clipped estimated eigenvalues and the eigenvalues of
$WW^\top$ with zeros through index $q+1$.
The spectrum identity (\Cref{prop:gram-matrix-identities}) identifies
the padded eigenvalues of $WW^\top$ with those of
$W^\top W$, also padded by a zero at index $q+1$.
Clipping at zero cannot increase the error relative to
these nonnegative eigenvalues.
Weyl's inequality
(\Cref{prop:weyl-perturbation}) therefore gives
\begin{align}
 \abs{\widehat\lambda_h-\lambda_h(\Lambda)}
 &\leq
 \norm{A-WW^\top}_{\op}
 +\norm{W^\top W-\Lambda}_{\op}\notag\\
 &\leq
 \frac{\norm{\Lambda}_{\op}}{8q}
 +\frac{\norm{\Lambda}_{\op}}{8q}
 \alignnote{assumed perturbation bounds}\notag\\
 &=\frac{\norm{\Lambda}_{\op}}{4q}
 \qquad(1\leq h\leq q).
 \label{eq:app-tensor-decomposition:largest-gap-eigenvalue-comparison}
\end{align}
The same bound holds at $h=q+1$, where both eigenvalues
are zero.
In particular,
$\widehat\lambda_1\geq3\norm{\Lambda}_{\op}/4$.

We now use this comparison to bound the upper endpoint
of the largest gap and its separation from the next
true eigenvalue.
Telescoping gives
\begin{align*}
 \widehat\lambda_{i_\star}-\widehat\lambda_{i_\star+1}
 &\geq\frac1p\sum_{h=1}^p
       (\widehat\lambda_h-\widehat\lambda_{h+1})
 \alignnote{choice of $i_\star$}\\
 &=\frac{\widehat\lambda_1}{p}
 \geq\frac{3\norm{\Lambda}_{\op}}{4q},
 \alignnote{$p\leq q$}\\
 \widehat\lambda_{i_\star}
 &\geq
 \widehat\lambda_{i_\star}-\widehat\lambda_{i_\star+1}
 \geq\frac{\norm{\Lambda}_{\op}}{2q},\\
 \widehat\lambda_{i_\star}-\lambda_{i_\star+1}(\Lambda)
 &\geq
 \widehat\lambda_{i_\star}-\widehat\lambda_{i_\star+1}
 -\frac{\norm{\Lambda}_{\op}}{4q}
 \alignnote{by \eqref{eq:app-tensor-decomposition:largest-gap-eigenvalue-comparison}}\\
 &\geq\frac{\norm{\Lambda}_{\op}}{2q}.
\end{align*}
\end{proof}

\begin{lemma}[Spectral Truncation and Approximate Orthonormality]
\label{lem:app-tensor-decomposition:spectral-truncation-orthonormality}
Under the hypotheses of \Cref{thm:tensor-decomposition:observable-operator-construction}, let
$F=[\bm f_i]_{i\in Q}$ be the rescaled component matrix from
\eqref{eq:tensor-decomposition:rescaled-components}, and let $J$ index the
$i_\star$ largest component weights in the fixed order.

\begin{enumerate}[label=(\roman*),ref=\thelemma(\roman*),leftmargin=2.2em]
\item\label[lemma]{lem:app-tensor-decomposition:truncation-outside-components}
The components outside $J$ satisfy
\[
 \norm{F_{Q\setminus J}}_F^2
 \leq |Q|^3\bigl(\varepsilon_{\mathrm{mom}}
                 +\varepsilon_{\mathrm{proj}}+\varepsilon_\perp\bigr).
\]

\item\label[lemma]{lem:app-tensor-decomposition:truncation-error}
For every $\bm{\zeta}\in\R^n$, the contribution of components
outside $J$ to the rescaled odd moment satisfies
\[
 \norm{\sum_{i\in Q\setminus J}\ip{\bm{\zeta}}{\ScaledRegressor_i}\,
       \bm f_i\bm f_i^\top}_F
 \leq K^3\norm{F_{Q\setminus J}}_F^2\norm{\bm{\zeta}}_2.
\]

\item\label[lemma]{lem:app-tensor-decomposition:orthonormality-error}
Let $Z=[\bm z_i]_{i\in J}$ be an orthogonal polar factor
of $F_J$. For every $\bm{\zeta}\in\R^n$, replacing the vectors $\bm f_i$
in $J$ by $\bm z_i$ gives the error bound
\[
\begin{aligned}
 &\norm{\sum_{i\in J}\ip{\bm{\zeta}}{\ScaledRegressor_i}
       \bigl(\bm f_i\bm f_i^\top-\bm z_i\bm z_i^\top\bigr)}_F\\
 &\qquad\leq
 3K^3\bigl(2|Q|^2\varepsilon_{\mathrm{mom}}
          +\norm{F_{Q\setminus J}}_F^2\bigr)\norm{\bm{\zeta}}_2.
\end{aligned}
\]
\end{enumerate}
\end{lemma}

\begin{proof}
We first recall the perturbation estimates used to justify rescaling.
We then bound the components outside $J$ using Sylvester's equation
and use their bound to control the correction to an orthonormal basis.

Recall that the hypotheses give
$\norm{\ScaledRegressor_i}_2\leq K^3$ for every $i\in Q$ and
$\varepsilon_{\mathrm{mom}}+\varepsilon_{\mathrm{proj}}
+\varepsilon_\perp\leq1/(16k^2)$.

Let
\[
 H=\left[\sqrt{\ComponentWeight{2s}{i}}\,
          \TensorPower{\NormalizedRegressor_i}{s}\right]_{i\in Q},
 \qquad
 \Lambda=\Diag_{i\in Q}(\ComponentWeight{2s}{i}),
\]
with components ordered by decreasing weight.
Then $\EvenMomentMatrix{s}{Q}=HH^\top$ and
$\tr(\EvenMomentMatrix{s}{Q})=\tr(\Lambda)\leq|Q|\norm{\Lambda}_{\op}$.
The even-moment, projection, and approximate-orthogonality assumptions
give the perturbation bounds established in the proof of
\Cref{thm:tensor-decomposition:observable-operator-construction}. We recall them here in the form
used below.
\begin{align}
 \norm{\widehat A-\ProjectionMap{s}\EvenMomentMatrix{s}{Q}
                         \ProjectionMap{s}^\top}_F
 &\leq |Q|\varepsilon_{\mathrm{mom}}\norm{\Lambda}_{\op},
 \label{eq:app-tensor-decomposition:truncation-even-error}\\
 \norm{H^\top\ProjectionMap{s}^\top\ProjectionMap{s}H-\Lambda}_{\op}
 &\leq\bigl(|Q|\varepsilon_{\mathrm{proj}}+\varepsilon_\perp\bigr)
       \norm{\Lambda}_{\op}.
 \label{eq:app-tensor-decomposition:truncation-weight-error}
\end{align}
Since
\[
 |Q|\varepsilon_{\mathrm{mom}}+|Q|\varepsilon_{\mathrm{proj}}+\varepsilon_\perp
 \leq |Q|\bigl(\varepsilon_{\mathrm{mom}}
              +\varepsilon_{\mathrm{proj}}+\varepsilon_\perp\bigr)
 \leq\frac1{16|Q|},
\]
\Cref{prop:tensor-decomposition:largest-gap-eigenvalue-bounds}, with $W=\ProjectionMap{s}H$, gives
\begin{align}
 \widehat\lambda_{i_\star}
 &\geq\frac{\norm{\Lambda}_{\op}}{2|Q|},
 \label{eq:app-tensor-decomposition:truncation-eigenvalue-lower-bound}\\
 \widehat\lambda_{i_\star}-\lambda_{i_\star+1}(\Lambda)
 &\geq\frac{\norm{\Lambda}_{\op}}{2|Q|}.
 \label{eq:app-tensor-decomposition:truncation-outside-separation}
\end{align}
Thus $\widehat\Lambda$ is positive definite. We use
$F=\widehat\Lambda^{-1/2}\widehat V\ProjectionMap{s}H$ throughout.

The rows of $\widehat V$ are orthonormal eigenvectors
of $\widehat A$ with eigenvalues in $\widehat\Lambda$.
Thus $\widehat V\widehat A\widehat V^\top=\widehat\Lambda$,
and the even-moment estimation assumption
(\Cref{def:projected-moment-estimate-even}) gives
\begin{align}
 \norm{FF^\top-I}_F
 &=\norm{\widehat\Lambda^{-1/2}\widehat V
       (\ProjectionMap{s} \EvenMomentMatrix{s}{Q}\ProjectionMap{s}^\top-\widehat A)
       \widehat V^\top\widehat\Lambda^{-1/2}}_F\notag\\
 &\leq
 \norm{\widehat\Lambda^{-1/2}}_{\op}^2
 \norm{\widehat V}_{\op}^2
 \norm{\ProjectionMap{s} \EvenMomentMatrix{s}{Q}\ProjectionMap{s}^\top-\widehat A}_F
 \alignnote{Frobenius product inequality}\notag\\
 &=\frac{
 \norm{\ProjectionMap{s} \EvenMomentMatrix{s}{Q}\ProjectionMap{s}^\top-\widehat A}_F}
 {\widehat\lambda_{i_\star}}\notag\\
 &\leq
 \frac{\varepsilon_{\mathrm{mom}}\tr(\EvenMomentMatrix{s}{Q})}
      {\widehat\lambda_{i_\star}}
 \alignnote{by \Cref{def:projected-moment-estimate-even}}\notag\\
 &\leq2|Q|^2\varepsilon_{\mathrm{mom}}
 \alignnote{by \eqref{eq:app-tensor-decomposition:truncation-eigenvalue-lower-bound}}\notag\\
 &\leq\frac18.
 \alignnote{by \Cref{ass:tensor-decomposition:small-errors}, $|Q|\leq k$}
 \label{eq:app-tensor-decomposition:rescaled-even-error}
\end{align}
Consequently,
\begin{align}
 \norm{F}_{\op}
 &=\sqrt{\norm{FF^\top}_{\op}}\notag\\
 &\leq\sqrt{1+\norm{FF^\top-I}_F}\notag\\
 &\leq\sqrt{\frac98}<2.
 \alignnote{by \eqref{eq:app-tensor-decomposition:rescaled-even-error}}
 \label{eq:app-tensor-decomposition:rescaled-operator-norm}
\end{align}
For the remaining factor, the Gram-matrix norm identity
(\Cref{prop:gram-matrix-identities}) and the projected-weight bound give
\begin{align}
 \norm{\ProjectionMap{s}H}_{\op}^2
 &=\norm{H^\top\ProjectionMap{s}^\top\ProjectionMap{s}H}_{\op}\notag\\
 &\leq\norm{\Lambda}_{\op}
 +\norm{H^\top\ProjectionMap{s}^\top\ProjectionMap{s}H-\Lambda}_{\op}\notag\\
 &\leq2\norm{\Lambda}_{\op},
 \alignnote{by \eqref{eq:app-tensor-decomposition:truncation-weight-error}
 and \Cref{ass:tensor-decomposition:small-errors}}\notag\\
 \norm{\widehat\Lambda^{-1/2}}_{\op}
 \norm{\ProjectionMap{s}H}_{\op}
 &\leq\sqrt{
 \frac{2\norm{\Lambda}_{\op}}{\widehat\lambda_{i_\star}}}\notag\\
 &\leq2\sqrt{|Q|}.
 \alignnote{by \eqref{eq:app-tensor-decomposition:truncation-eigenvalue-lower-bound}}
 \label{eq:app-tensor-decomposition:truncation-multiplier-bound}
\end{align}

\emph{Components outside $J$.}
We now prove \Cref{lem:app-tensor-decomposition:truncation-outside-components}.
If $J=Q$, its left-hand side is zero. Otherwise, restricting
the residual to the columns outside $J$ gives the Sylvester equation
\[
 \widehat\Lambda F_{Q\setminus J}
 -F_{Q\setminus J}\Diag_{i\in Q\setminus J}(\ComponentWeight{2s}{i})
 =(\widehat\Lambda F-F\Lambda)_{Q\setminus J}.
\]
Since $J$ indexes the $i_\star$ largest weights, for
$1\leq h\leq i_\star$ and $i\in Q\setminus J$,
\begin{align}
 \widehat\lambda_h-\ComponentWeight{2s}{i}
 &\geq\widehat\lambda_{i_\star}
       -\lambda_{i_\star+1}(\Lambda)\notag\\
 &\geq\frac{\norm{\Lambda}_{\op}}{2|Q|}.
 \alignnote{by \eqref{eq:app-tensor-decomposition:truncation-outside-separation}}
\end{align}
Sylvester's equation (\Cref{prop:sylvester-equation}) therefore gives
\begin{align}
 \norm{F_{Q\setminus J}}_F
 &\leq\frac{2|Q|}{\norm{\Lambda}_{\op}}
       \norm{(\widehat\Lambda F-F\Lambda)_{Q\setminus J}}_F
 \alignnote{by \Cref{prop:sylvester-equation}}\notag\\
 &\leq\frac{2|Q|}{\norm{\Lambda}_{\op}}
       \norm{\widehat\Lambda F-F\Lambda}_F.
 \label{eq:app-tensor-decomposition:outside-components-sylvester-equation}
\end{align}

We now bound the right-hand side of the Sylvester equation using the
two perturbation errors.
Using $\widehat V\widehat A
=\widehat\Lambda\widehat V$ and
$HH^\top=\EvenMomentMatrix{s}{Q}$, we obtain
\begin{align}
 \widehat\Lambda F-F\Lambda
 &=\left(
 \widehat\Lambda F-FH^\top\ProjectionMap{s}^\top\ProjectionMap{s} H
 \right)
 +F(H^\top\ProjectionMap{s}^\top\ProjectionMap{s} H-\Lambda)
 \notag\\
 &=\widehat\Lambda^{-1/2}\widehat V
   (\widehat A-\ProjectionMap{s} \EvenMomentMatrix{s}{Q}\ProjectionMap{s}^\top)
   \ProjectionMap{s} H
 +F(H^\top\ProjectionMap{s}^\top\ProjectionMap{s} H-\Lambda).
 \alignnote{by \eqref{eq:tensor-decomposition:rescaled-components}}
 \label{eq:app-tensor-decomposition:component-sylvester-residual}
\end{align}
Since $F$ has $|Q|$ columns,
\begin{equation}\label{eq:app-tensor-decomposition:rescaled-frobenius-norm}
 \norm{F}_F\leq\sqrt{|Q|}\norm{F}_{\op}
 \leq2\sqrt{|Q|}.
\end{equation}
Using \eqref{eq:app-tensor-decomposition:component-sylvester-residual} and the
Frobenius product inequality gives
\begin{align}
 \norm{\widehat\Lambda F-F\Lambda}_F
 &\leq
 \norm{\widehat\Lambda^{-1/2}}_{\op}
 \norm{\ProjectionMap{s}H}_{\op}
 \norm{\widehat A-\ProjectionMap{s}\EvenMomentMatrix{s}{Q}\ProjectionMap{s}^\top}_F
 \notag\\
 &\quad+\norm{F}_F
 \norm{H^\top\ProjectionMap{s}^\top\ProjectionMap{s}H-\Lambda}_{\op}
 \alignnote{$\norm{\widehat V}_{\op}=1$}\notag\\
 &\leq2\sqrt{|Q|}\Bigl(
 \norm{\widehat A-\ProjectionMap{s}\EvenMomentMatrix{s}{Q}\ProjectionMap{s}^\top}_F
 \notag\\
 &\qquad+
 \norm{H^\top\ProjectionMap{s}^\top\ProjectionMap{s}H-\Lambda}_{\op}
 \Bigr)
 \alignnote{by \eqref{eq:app-tensor-decomposition:truncation-multiplier-bound},
 \eqref{eq:app-tensor-decomposition:rescaled-frobenius-norm}}\notag\\
 &\leq2|Q|^{3/2}\norm{\Lambda}_{\op}
 \left(
 \varepsilon_{\mathrm{mom}}
 +\varepsilon_{\mathrm{proj}}+\varepsilon_\perp
 \right).
 \alignnote{by \eqref{eq:app-tensor-decomposition:truncation-even-error},
 \eqref{eq:app-tensor-decomposition:truncation-weight-error}; $|Q|\geq1$}
 \label{eq:app-tensor-decomposition:outside-components-residual-bound}
\end{align}
Squaring \eqref{eq:app-tensor-decomposition:outside-components-sylvester-equation} gives
\begin{align}
 \norm{F_{Q\setminus J}}_F^2
 &\leq
 \frac{4|Q|^2}{\norm{\Lambda}_{\op}^2}
 \norm{\widehat\Lambda F-F\Lambda}_F^2
 \alignnote{by \eqref{eq:app-tensor-decomposition:outside-components-sylvester-equation}}\notag\\
 &\leq16|Q|^5
 \left(
 \varepsilon_{\mathrm{mom}}
 +\varepsilon_{\mathrm{proj}}+\varepsilon_\perp
 \right)^2
 \alignnote{by \eqref{eq:app-tensor-decomposition:outside-components-residual-bound}}\notag\\
 &\leq|Q|^3
 \left(
 \varepsilon_{\mathrm{mom}}
 +\varepsilon_{\mathrm{proj}}+\varepsilon_\perp
 \right).
 \alignnote{by \Cref{ass:tensor-decomposition:small-errors}, $|Q|\leq k$}
\end{align}

This proves \Cref{lem:app-tensor-decomposition:truncation-outside-components}, including the case $J=Q$.

\emph{Spectral truncation error.}
We prove \Cref{lem:app-tensor-decomposition:truncation-error} by bounding the
weighted sum of outer products of the columns outside $J$.
This sum equals
$F_{Q\setminus J}\Diag_{i\in Q\setminus J}(\ip{\bm{\zeta}}{\ScaledRegressor_i})
 F_{Q\setminus J}^\top$, so
\begin{align}
 &\norm{F_{Q\setminus J}
       \Diag_{i\in Q\setminus J}(\ip{\bm{\zeta}}{\ScaledRegressor_i})
       F_{Q\setminus J}^\top}_F\notag\\
 &\quad\leq
 \norm{F_{Q\setminus J}}_F
 \norm{\Diag_{i\in Q\setminus J}(\ip{\bm{\zeta}}{\ScaledRegressor_i})}_{\op}
 \norm{F_{Q\setminus J}}_{\op}
 \alignnote{Frobenius product inequality}\notag\\
 &\quad\leq K^3\norm{\bm{\zeta}}_2\norm{F_{Q\setminus J}}_F^2.
 \alignnote{Cauchy--Schwarz; $\norm{\ScaledRegressor_i}_2\leq K^3$}
\end{align}

This proves \Cref{lem:app-tensor-decomposition:truncation-error}.

\emph{Approximate orthonormality.}
We now prove \Cref{lem:app-tensor-decomposition:orthonormality-error}.
Choose an orthogonal polar factor
$Z=[\bm{z}_i]_{i\in J}$ of the square matrix $F_J$.
The orthogonal polar-factor bound
(\Cref{prop:polar-factor-bounds}) gives
\begin{align}
 \norm{F_J-Z}_F
 &\leq\norm{F_JF_J^\top-I}_F
 \alignnote{by \Cref{prop:polar-factor-bounds}}\notag\\
 &=\norm{(FF^\top-I)
       -F_{Q\setminus J}F_{Q\setminus J}^\top}_F\notag\\
 &\leq\norm{FF^\top-I}_F+\norm{F_{Q\setminus J}}_F^2
 \notag\\
 &\leq2|Q|^2\varepsilon_{\mathrm{mom}}
       +\norm{F_{Q\setminus J}}_F^2.
 \alignnote{by \eqref{eq:app-tensor-decomposition:rescaled-even-error}}
 \label{eq:app-tensor-decomposition:polar-factor-error}
\end{align}
We now bound the error from replacing $F_J$ by $Z$.
Since $\norm{F_J}_{\op}\leq\norm{F}_{\op}\leq2$
by \eqref{eq:app-tensor-decomposition:rescaled-operator-norm},
\begin{align}
 &\norm{F_J\Diag_{i\in J}(\ip{\bm{\zeta}}{\ScaledRegressor_i})F_J^\top
       -Z\Diag_{i\in J}(\ip{\bm{\zeta}}{\ScaledRegressor_i})Z^\top}_F\notag\\
 &\quad=\Bigl\|
 (F_J-Z)\Diag_{i\in J}(\ip{\bm{\zeta}}{\ScaledRegressor_i})F_J^\top\notag\\
 &\qquad+Z\Diag_{i\in J}(\ip{\bm{\zeta}}{\ScaledRegressor_i})(F_J-Z)^\top
 \Bigr\|_F\notag\\
 &\quad\leq
 K^3(\norm{F_J}_{\op}+1)\norm{F_J-Z}_F\norm{\bm{\zeta}}_2
 \alignnote{$Z$ orthogonal; $\norm{\ScaledRegressor_i}_2\leq K^3$}
 \\
 &\quad\leq
 3K^3\left(
 2|Q|^2\varepsilon_{\mathrm{mom}}
 +\norm{F_{Q\setminus J}}_F^2
 \right)\norm{\bm{\zeta}}_2.
 \alignnote{by \eqref{eq:app-tensor-decomposition:rescaled-operator-norm}
 and \eqref{eq:app-tensor-decomposition:polar-factor-error}}
\end{align}
This proves \Cref{lem:app-tensor-decomposition:orthonormality-error},
with the same orthonormal basis for every $\bm{\zeta}$.
\end{proof}

\subsection{Proof of \texorpdfstring{\Cref*{thm:tensor-decomposition:implicit-tensor-decomposition}}{Theorem 6.1}}

We verify that the moment, projection, and near-orthogonality bounds
give the direction accuracy required by \Cref{thm:tensor-decomposition:implicit-tensor-decomposition}.

\begin{proposition}[Observable-Operator Accuracy]
\label{prop:app-tensor-decomposition:observable-operator-accuracy}
Use the hyperparameters from \Cref{hyp:algorithm}, with $2\leq k\leq n$.
Fix $0\leq\varepsilon_{\mathrm{mean}}\leq e^{-K^2s}$, and put
\begin{equation}\label{eq:app-tensor-decomposition:operator-family-error}
 \varepsilon_{\mathrm{fam}}=16K^3k^3
 \left((1+2K^{s-1})\varepsilon_{\mathrm{mean}}
       +k e^{-s/K}\right).
\end{equation}
Then
\begin{equation}\label{eq:app-tensor-decomposition:operator-family-target}
 32K^6\sqrt n\,k^2\varepsilon_{\mathrm{fam}}
 \leq e^{-s/(2K)}\leq\frac{\Delta^2}{16B^2}.
\end{equation}
\end{proposition}

\begin{proof}
By the choice of \Cref{hyp:algorithm},
\begin{align}
 K&\geq128,
 \label{eq:app-tensor-decomposition:operator-accuracy-k-choice}\\
 s&\geq K^4\log n.
 \label{eq:app-tensor-decomposition:operator-accuracy-order-choice}
\end{align}
The assumed bound on $\varepsilon_{\mathrm{mean}}$ gives
\begin{align}
 (1+2K^{s-1})\varepsilon_{\mathrm{mean}}
 &\leq3e^{-(K^2-\log K)s}
 \alignnote{$\varepsilon_{\mathrm{mean}}\leq e^{-K^2s}$; $K\geq1$}\notag\\
 &\leq e^{-s/K}.
 \alignnote{$K^2-\log K-1\geq1/K$, $3\leq e^s$}
 \label{eq:app-tensor-decomposition:operator-family-estimation-term}
\end{align}
Consequently,
\begin{align}
 \varepsilon_{\mathrm{fam}}
 &\leq16K^3k^3(1+k)e^{-s/K}
 \alignnote{by \eqref{eq:app-tensor-decomposition:operator-family-error}
 and \eqref{eq:app-tensor-decomposition:operator-family-estimation-term}}\notag\\
 &\leq32K^3k^4e^{-s/K}.
 \alignnote{$k\geq1$}
 \label{eq:app-tensor-decomposition:operator-family-error-bound}
\end{align}
The remaining factor is absorbed by half of the exponential decay:
\begin{align}
 1024K^9\sqrt n\,k^6
 &\leq2^{10+9K}n^7
 \alignnote{$K\leq2^K$; $k\leq n$}\notag\\
 &\leq n^{17+9K}
 \alignnote{$n\geq2$}\notag\\
 &\leq n^{K^3/2}.
 \alignnote{choice of $K$~\eqref{eq:app-tensor-decomposition:operator-accuracy-k-choice}:
 $17+9K\leq K^3/2$}
 \label{eq:app-tensor-decomposition:operator-family-prefactor}
\end{align}
Combining these bounds gives
\begin{align*}
 32K^6\sqrt n\,k^2\varepsilon_{\mathrm{fam}}
 &\leq1024K^9\sqrt n\,k^6e^{-s/K}
 \alignnote{by \eqref{eq:app-tensor-decomposition:operator-family-error-bound}}\\
 &\leq e^{(K^3/2)\log n-s/K}
 \alignnote{by \eqref{eq:app-tensor-decomposition:operator-family-prefactor}}\\
 &\leq e^{-s/(2K)}.
 \alignnote{choice of $s$~\eqref{eq:app-tensor-decomposition:operator-accuracy-order-choice}}
\end{align*}
It remains to compare this error bound with the tolerance needed for refinement.
The choices in \Cref{hyp:algorithm} give $s\geq K^4\log2$ and
$K\Delta^2\geq128B^2$, and hence
\begin{equation}
 e^{-s/(2K)}\leq2^{-K^3/2}\leq K^{-1}
 \leq\frac{\Delta^2}{128B^2}
 \leq\frac{\Delta^2}{16B^2}.
\end{equation}
Here $2^{-K^3/2}\leq K^{-1}$ follows from
$K\leq2^K\leq2^{K^3/2}$ for $K\geq128$.
This proves \eqref{eq:app-tensor-decomposition:operator-family-target}.
\end{proof}

\section{Deferred Proofs for \texorpdfstring{\Cref{sec:refinement:main}}{Section 7}}

\begin{proposition}[Refinement Bounds]
\label{prop:app-refinement:refinement-bounds}
Use $K$ from \Cref{hyp:algorithm}, with $2\leq k\leq n$.
Fix $\varepsilon\in(0,1)$ and $\delta\in(0,\nicefrac13)$, and set
\[
 {\Gamma}=2K\sqrt{\log\!\left(\frac{n}{\varepsilon}\right)}.
\]
Then
\begin{align*}
 2^{13}B{\Gamma}^5\sqrt{k}\,e^{-{\Gamma}^2/K}
 &\leq\frac{\varepsilon}{2},\\
 \frac{n^2B^2\log k}{\varepsilon^2}
 {\Gamma}^{10}e^{{\Gamma}^2}
 \log\!\left(\frac n\delta\right)
 &=\Ot\!\left(
 \left(\frac{nk}{\varepsilon}\right)^{5K^2}
 \log(\nicefrac1\delta)
 \right).
\end{align*}
\end{proposition}

\begin{proof}
We first verify the bias bound, then bound the expression
governing the sample complexity and running time.
By the choice of \Cref{hyp:algorithm},
\begin{align}
 K&\geq128,
 \label{eq:app-refinement:k-choice}\\
 B&\leq K.
\end{align}
The choice of ${\Gamma}$ gives
\begin{equation}\label{eq:app-refinement:refinement-threshold-choice}
 {\Gamma}^2=4K^2\log\!\left(\frac{n}{\varepsilon}\right).
\end{equation}
For $K\geq128$, we have $K^6\leq2^{K/2}$.
Substituting the choice of ${\Gamma}$ gives
\begin{align*}
 2^{13}B{\Gamma}^5\sqrt{k}\,e^{-{\Gamma}^2/K}
 &=2^{18}BK^5\sqrt{k}\,
 \log^{5/2}\!\left(\frac{n}{\varepsilon}\right)
 \left(\frac{n}{\varepsilon}\right)^{-4K}
 \alignnote{choice of ${\Gamma}$~\eqref{eq:app-refinement:refinement-threshold-choice}}\\
 &=\frac{2^{18}BK^5\sqrt{k}}n
 \log^{5/2}\!\left(\frac{n}{\varepsilon}\right)
 \left(\frac{n}{\varepsilon}\right)^{1-4K}\varepsilon\\
 &\leq2^{18}K^6
 \left(\frac{n}{\varepsilon}\right)^{7/2-4K}\varepsilon
 \alignnote{$B\leq K$; $k\leq n$; $\log(n/\varepsilon)\leq n/\varepsilon$}\\
 &\leq2^{22-7K/2}\varepsilon
 \alignnote{$K^6\leq2^{K/2}$; $n/\varepsilon\geq2$}\\
 &\leq\frac{\varepsilon}{2}.
 \alignnote{choice of $K$~\eqref{eq:app-refinement:k-choice}}
\end{align*}

We now bound the expression governing the sample complexity
and running time.
The choice of ${\Gamma}$ and the ranges of $n,k,\varepsilon,\delta$ give
\begin{align}
 e^{{\Gamma}^2}
 &=\left(\frac{n}{\varepsilon}\right)^{4K^2}
 \leq\left(\frac{nk}{\varepsilon}\right)^{4K^2},
 \\
 {\Gamma}^2
 &\leq4K^2\log\!\left(\frac{nk}{\varepsilon}\right),
 \\
 \log\!\left(\frac n\delta\right)
 &\leq2\log\!\left(\frac{nk}{\varepsilon}\right)
       \log(\nicefrac1\delta).
 \alignnote{$\log(\nicefrac1\delta)\geq1$; $\log(\nicefrac{k}\varepsilon)\geq\log2$}
\end{align}
Since $K\leq2^K$ and $nk/\varepsilon\geq2$,
\[
 K^{12}\leq2^{12K}
 \leq\left(\frac{nk}{\varepsilon}\right)^{12K}.
\]
Using these bounds, $B\leq K$, $\log k\leq\log(nk/\varepsilon)$,
and $n^2/\varepsilon^2\leq(nk/\varepsilon)^2$, we obtain
\begin{align*}
 \frac{n^2B^2\log k}{\varepsilon^2}
 {\Gamma}^{10}e^{{\Gamma}^2}\log\!\left(\frac n\delta\right)
 &\leq2^{11}K^{12}
 \left(\frac{nk}{\varepsilon}\right)^{4K^2+2}
 \log^7\!\left(\frac{nk}{\varepsilon}\right)
 \log(\nicefrac1\delta)\\
 &\leq2^{11}
 \left(\frac{nk}{\varepsilon}\right)^{4K^2+12K+2}
 \log^7\!\left(\frac{nk}{\varepsilon}\right)
 \log(\nicefrac1\delta)
 \alignnote{$K^{12}\leq(nk/\varepsilon)^{12K}$}\\
 &=\Ot\!\left(
 \left(\frac{nk}{\varepsilon}\right)^{5K^2}
 \log(\nicefrac1\delta)
 \right).
 \alignnote{$4K^2+12K+2\leq5K^2$}
\end{align*}
This proves both bounds.
\end{proof}

\section{Deferred Proofs for \texorpdfstring{\Cref{sec:peeling:main}}{Section 8}}

We give the projected-feature second-moment calculation and the sampling
proof of \Cref{lem:peeling:weighted-mean-estimation}. We then prove the pairwise exceedance
bound used for the peeling bias.

\subsection{Second-Moment Bounds}
\label{subsec:app-peeling:second-moment-bounds}

\projectedfeaturesecondmoments*

\begin{proof}
Write $\widetilde{\Covariate}=\Covariate+\iu \Covariate'$.
Let $\mathsf A\in(\R^n)^{\otimes2r}$ be the tensor satisfying
$\Mat_r(\mathsf A)=U^\top H U$.
Expanding the matrix contractions gives
\[
 \psi(\Covariate,\Covariate')
 =\frac{\RealPart\ip{\mathsf A}{\TensorPower{\widetilde{\Covariate}}{2r}}}{\sqrt{(2r)!}},
 \qquad
 \psi_{\bm{\zeta}}(\Covariate,\Covariate')
 =\frac{\RealPart\ip{\mathsf A\otimes\bm{\zeta}}{\TensorPower{\widetilde{\Covariate}}{2r+1}}}
                    {\sqrt{(2r+1)!}}.
\]
Here the tensor pairing is extended bilinearly to complex tensors.
Since $UU^\top=I_p$ and flattening preserves the Frobenius norm,
\[
 \norm{\mathsf A}_F=\norm{U^\top H U}_F=\norm {H}_F\leq1,
 \qquad
 \norm{\mathsf A\otimes\bm{\zeta}}_F=\norm{\mathsf A}_F\norm {\bm{\zeta}}_2\leq\norm {\bm{\zeta}}_2.
\]
The polynomials $\psi(\Covariate,\Covariate')$ and $\psi_{\bm{\zeta}}(\Covariate,\Covariate')$ have degrees
at most $2r$ and $2r+1$, respectively.

We first bound the averages over $\Covariate'$.
The definition of the Hermite tensors gives
\[
 \E_{\Covariate'}[\psi(\Covariate,\Covariate')\mid \Covariate]
 =\HermitePairing{2r}{\Covariate}{\Sym(\mathsf A)},
 \qquad
 \E_{\Covariate'}[\psi_{\bm{\zeta}}(\Covariate,\Covariate')\mid \Covariate]
 =\HermitePairing{2r+1}{\Covariate}{\Sym(\mathsf A\otimes\bm{\zeta})}.
\]
Hermite orthogonality (\Cref{prop:hermite-expansion-normalization}) therefore gives
\begin{align*}
 \left\|\E_{\Covariate'}[\psi(\Covariate,\Covariate')\mid \Covariate]\right\|_{L^2(\Covariate)}
 &=\norm{\Sym(\mathsf A)}_F
 \leq\norm{\mathsf A}_F
 \leq1,\\
 \left\|\E_{\Covariate'}[\psi_{\bm{\zeta}}(\Covariate,\Covariate')\mid \Covariate]\right\|_{L^2(\Covariate)}
 &=\norm{\Sym(\mathsf A\otimes\bm{\zeta})}_F
 \leq\norm{\mathsf A\otimes\bm{\zeta}}_F
 \leq\norm {\bm{\zeta}}_2.
\end{align*}
Here averaged symmetrization is the orthogonal projection onto
symmetric tensors.

We now bound the sampled functions using Wick's identity
(\Cref{prop:wick-identity}).
The complex Gaussian vector $\widetilde{\Covariate}$ satisfies
\[
 \E_{\Covariate,\Covariate'}[\widetilde X_j\widetilde X_h]=0,
 \qquad
 \E_{\Covariate,\Covariate'}[\widetilde X_j\overline{\widetilde X_h}]
 =2\mathbf 1_{\{j=h\}}.
\]
For a real order-$d$ tensor $\mathsf G$, every nonzero Wick pairing in
$\E_{\Covariate,\Covariate'}\abs{\ip{\mathsf G}{\TensorPower{\widetilde{\Covariate}}{d}}}^2$ therefore matches
one of the $d$ unbarred factors with one of the $d$ barred factors.
Summing over the $d!$ such pairings gives
\[
 \E_{\Covariate,\Covariate'}\abs{\ip{\mathsf G}{\TensorPower{\widetilde{\Covariate}}{d}}}^2
 =2^d d!\,\ip{\mathsf G}{\Sym(\mathsf G)}
 =2^d d!\,\norm{\Sym(\mathsf G)}_F^2
 \leq2^d d!\,\norm{\mathsf G}_F^2.
\]
Since $(\RealPart z)^2\leq\abs {z}^2$, applying this identity
to $\mathsf A$ and $\mathsf A\otimes\bm{\zeta}$ and dividing by the respective factorials yields
\[
 \E_{\Covariate,\Covariate'}[\psi(\Covariate,\Covariate')^2]\leq2^{2r}\leq100^{r+2},
 \qquad
 \E_{\Covariate,\Covariate'}[\psi_{\bm{\zeta}}(\Covariate,\Covariate')^2]\leq2^{2r+1}\norm {\bm{\zeta}}_2^2
 \leq100^{r+2}\norm {\bm{\zeta}}_2^2.
\]
\end{proof}

\subsection{Sampling Error}
\label{subsec:app-peeling:sampling-error}

\peelingmeanestimation*

\begin{proof}
Index the recovered list by $R$.
For each observation, \Cref{alg:peeling:estimate-moments-by-peeling} uses the synthetic maximum
\[
 \widetilde Y
 =\max_{i\in R}
 \bigl(\ip{\widehat{\Regressor}_i}{\Covariate}
       +\widehat b_i+\eta_i\bigr),
\]
where $\bm\eta_R$ is an independent copy of $\NoiseVector_R$.
Averaging over this auxiliary noise gives
\begin{align*}
 &\E_{\bm\eta_R}\!\left[
  (\TailIndicator{Y_{\max}}-\TailIndicator{\widetilde Y})
  \Psi_\ell(\Covariate,\Covariate')
  \,\middle|\,\Covariate,Y_{\max},\Covariate'\right]\\
 &\quad=\widehat\omega_{\widehat{\mathcal R}}(\Covariate,Y_{\max})
       \Psi_\ell(\Covariate,\Covariate').
\end{align*}
The difference between indicators has absolute value at most one,
so the second moment of each weighted entry is at most $100^{s+2}$.
Apply \Cref{lem:coordinatewise-trimmed-mean} to all $mp^2$ weighted entries,
with second-moment bound $100^{s+2}$, coordinate accuracy
$\varepsilon_{\mathrm{est}}/(p\sqrt m)$, and failure probability $\delta$.
Here $\varepsilon_{\mathrm{est}}/(p\sqrt m)\leq1\leq\sqrt{100^{s+2}}$.
Thus, with probability at least $1-\delta$, every estimated entry
has error at most $\varepsilon_{\mathrm{est}}/(p\sqrt m)$.
Condition on this event and write
\[
 E_\ell\coloneqq\widehat M_\ell
 -\E_{\Covariate,Y_{\max},\Covariate'}\!\left[
  \widehat\omega_{\widehat{\mathcal R}}(\Covariate,Y_{\max})\Psi_\ell(\Covariate,\Covariate')\right].
\]
The matrices are symmetric because corresponding upper- and
lower-triangular entries have identical observation lists.
Consequently,
\[
 \norm{E_\ell}_F
 =\left(\sum_{a,c=1}^p|(E_\ell)_{ac}|^2\right)^{1/2}
 \leq p\,\frac{\varepsilon_{\mathrm{est}}}{p\sqrt m}
 =\frac{\varepsilon_{\mathrm{est}}}{\sqrt m}.
\]
For every $\bm{\zeta}\in\R^m$, the triangle inequality gives
\begin{align*}
 \norm{\sum_{\ell=1}^m \zeta_\ell E_\ell}_F
 &\leq\sum_{\ell=1}^m|\zeta_\ell|\norm{E_\ell}_F\\
 &\leq\norm {\bm{\zeta}}_2
       \left(\sum_{\ell=1}^m\norm{E_\ell}_F^2\right)^{1/2}
 \alignnote{Cauchy--Schwarz}\\
 &\leq\varepsilon_{\mathrm{est}}\norm {\bm{\zeta}}_2.
\end{align*}
This calculation holds for every $\bm{\zeta}$ on the same entrywise success
event, which proves the simultaneous guarantee.

Since $\log(mp^2/\delta)\leq2\log(mp/\delta)$, the number of samples is
\[
 O\!\left(
 \frac{p^2m}{\varepsilon_{\mathrm{est}}^2}
 100^s\log\!\left(\frac{mp}{\delta}\right)
 \right).
\]
On each observation, the synthetic maximum is computed once and the
difference between indicators is reused across the family.
Each input function is evaluated once.
By \Cref{lem:coordinatewise-trimmed-mean}, computing the coordinatewise mean
estimates takes $O(mp^2)$ time for each sample.
This proves the running-time claim.
\end{proof}

\subsection{Peeling Bias}
\label{subsec:app-peeling:peeling-bias}

\pairwiseexceedancebound*

\begin{proof}
By the choice of \Cref{hyp:algorithm},
\begin{align}
 K&\geq128,\\
 K&\geq\frac{16(1+B^2)^3}{\Delta^4}.
 \label{eq:app-peeling:joint-tail-k-choice}
\end{align}
Fix distinct $i,j\in[k]$ and $1\leq r\leq s$.
H\"older's inequality separates the polynomial from the
probability that both responses exceed $T$.
We use exponent $K$. Hypercontractivity bounds the polynomial
factor, while uncoveredness preserves additional Gaussian decay
in the probability factor. Thus
\begin{align}
 &\E_{\Covariate,\xi_i,\xi_j,\Covariate'}\!\left[
  |\psi(\Covariate,\Covariate')|\1\{Y_i>T,Y_j>T\}\right]\notag\\
 &\quad\leq
 \norm{\psi}_{L^K(\Covariate,\Covariate')}
 \Prb_{\Covariate,\xi_i,\xi_j}(Y_i>T,Y_j>T)^{1-1/K}.
 \alignnote{H\"older's inequality}
 \label{eq:app-peeling:holder-decomposition}
\end{align}
We first bound the polynomial factor.
By Gaussian hypercontractivity
(\Cref{prop:gaussian-hypercontractivity}) of the standard
Gaussian vector $(\Covariate,\Covariate')$,
\begin{align}
 \norm{\psi}_{L^K(\Covariate,\Covariate')}
 &\leq(K-1)^{(2s+1)/2}\norm{\psi}_{L^2(\Covariate,\Covariate')}
 \alignnote{by \Cref{prop:gaussian-hypercontractivity}}\notag\\
 &\leq K^{s+1/2}\norm{\psi}_{L^2(\Covariate,\Covariate')}.
\end{align}

We now bound the probability factor.
Put $Z_i=Y_i-b_i$ and $Z_j=Y_j-b_j$.
Since $|b_i|,|b_j|\leq B$,
\[
 \{Y_i>T,Y_j>T\}\subseteq\{Z_i>T-B,Z_j>T-B\}.
\]
For the centered responses, uncoveredness supplies additional
Gaussian decay that absorbs both the polynomial factors and the bias shifts.
The pair $(Z_i,Z_j)$ is centered jointly Gaussian, with
variances $\sigma_i^2=1+\norm{\Regressor_i}_2^2$ and
$\sigma_j^2=1+\norm{\Regressor_j}_2^2$, and covariance
$c_{ij}=\ip{\Regressor_i}{\Regressor_j}+\rho$.
Uncoveredness and $0\leq\rho\leq1$ give
\[
 |c_{ij}|\leq\min\{\sigma_i^2,\sigma_j^2\}-\Delta^2,
 \qquad
 D_{ij}\coloneqq\sigma_i^2\sigma_j^2-c_{ij}^2>0,
\]
and $\sigma_i^2-c_{ij},\sigma_j^2-c_{ij}\geq\Delta^2$.
Applying the joint Gaussian tail bound
(\Cref{prop:joint-gaussian-tail}) at $T-B>0$ and raising
the bound to the positive power $1-1/K$ gives
\begin{equation}\label{eq:app-peeling:joint-response-tail}
\begin{aligned}
 &\Prb_{\Covariate,\xi_i,\xi_j}(Y_i>T,Y_j>T)^{1-1/K}\\
 &\quad\leq\exp\!\left(
  -(T-B)^2\left(1-\frac1K\right)
  \frac{\sigma_i^2+\sigma_j^2-2c_{ij}}{2D_{ij}}
 \right).
\end{aligned}
\end{equation}

We next compare this exponent with the marginal tail of $Y_j$.
Using $\sigma_j^2\leq1+B^2$ and $D_{ij}\leq(1+B^2)^2$,
\begin{align}
 &\left(1-\frac1K\right)
 \frac{\sigma_i^2+\sigma_j^2-2c_{ij}}{2D_{ij}}\notag\\
 &\quad=
 \left(1-\frac1K\right)
 \left(
  \frac1{2\sigma_j^2}
  +\frac{(\sigma_j^2-c_{ij})^2}{2\sigma_j^2D_{ij}}
 \right)\notag\\
 &\quad\geq
 \left(1-\frac1K\right)
 \left(
  \frac1{2\sigma_j^2}
  +\frac{\Delta^4}{2(1+B^2)^3}
 \right)
 \alignnote{by \Cref{ass:boundedness-uncoveredness}}\notag\\
 &\quad\geq
 \left(1-\frac1K\right)
 \left(\frac1{2\sigma_j^2}+\frac8K\right)
 \alignnote{choice of $K$~\eqref{eq:app-peeling:joint-tail-k-choice}}\notag\\
 &\quad=
 \frac1{2\sigma_j^2}+\frac8K
 -\frac1{2K\sigma_j^2}-\frac8{K^2}\notag\\
 &\quad\geq\frac1{2\sigma_j^2}+\frac7K,
 \alignnote{$\sigma_j^2\geq1$, $K\geq16$}\notag\\
 &(T-B)^2\left(1-\frac1K\right)
 \frac{\sigma_i^2+\sigma_j^2-2c_{ij}}{2D_{ij}}
 -\frac{(T-b_j)^2}{2\sigma_j^2}\notag\\
 &\quad\geq\frac{7(T-B)^2}{K}
 -\frac{(T-b_j)^2-(T-B)^2}{2\sigma_j^2}\notag\\
 &\quad\geq\frac{7(T-B)^2}{K}-2BT
 \alignnote{$|b_j|\leq B$, $\sigma_j^2\geq1$}\notag\\
 &\quad\geq\left(\frac74-\frac12\right)\frac{T^2}{K}
 \alignnote{$T-B\geq T/2$, $T\geq4BK$}\notag\\
 &\quad\geq\frac{T^2}{K}.
 \label{eq:app-peeling:joint-tail-exponent-gap}
\end{align}
It remains to compare $e^{-(T-b_j)^2/(2\sigma_j^2)}$ with
$\ComponentWeight{2r}{j}$.
We use the Gaussian density on the initial interval
$[T,T+1/T]$ of the slice and then apply
\Cref{lem:hermite:even-moment-weights}.
The hypotheses give $T\geq5(B+1)$ and $(T-B)^2\geq Ks$,
as required by that lemma.
Since $T\geq2B$, $T\geq1$, and $\sigma_j^2\geq1$,
\begin{align}
 \Prb(T\leq Y_j\leq T+1)
 &\geq\int_T^{T+1/T}\phi_{\sigma_j^2}(u-b_j) \odif{u}\notag\\
 &\geq\frac1T\phi_{\sigma_j^2}(T+1/T-b_j)
 \alignnote{the density decreases}\notag\\
 &=\frac1T\phi_{\sigma_j^2}(T-b_j)
   \exp\!\left(-\frac{2(T-b_j)/T+T^{-2}}{2\sigma_j^2}\right)\notag\\
 &\geq\frac{e^{-2}}T\phi_{\sigma_j^2}(T-b_j).
 \alignnote{$|b_j|\leq B$, $T\geq2B$, $T\geq1$, $\sigma_j^2\geq1$}
 \label{eq:app-peeling:response-slice-mass}
\end{align}
Consequently,
\begin{align}
 &\Prb_{\Covariate,\xi_i,\xi_j}(Y_i>T,Y_j>T)^{1-1/K}
 \notag\\
 &\quad\leq e^{-T^2/K}e^{-(T-b_j)^2/(2\sigma_j^2)}
 \alignnote{by \eqref{eq:app-peeling:joint-response-tail} and \eqref{eq:app-peeling:joint-tail-exponent-gap}}\notag\\
 &\quad\leq
 e^{-T^2/K}e^2T\sqrt{2\pi\sigma_j^2}\,
 \Prb(T\leq Y_j\leq T+1)
 \alignnote{by \eqref{eq:app-peeling:response-slice-mass}}\notag\\
 &\quad\leq K^{3/2}T e^{-T^2/K}
 \Prb(T\leq Y_j\leq T+1)
 \alignnote{$\sqrt{2\pi\sigma_j^2}\leq K$, $e^2\leq\sqrt K$}\notag\\
 &\quad\leq K^2e^{-T^2/(2K)}
 \Prb(T\leq Y_j\leq T+1)
 \alignnote{$T/\sqrt K\leq e^{T^2/(2K)}$}\notag\\
 &\quad\leq K^2e^{-T^2/(2K)}\ComponentWeight{2r}{j}.
 \alignnote{by \Cref{lem:hermite:even-moment-weights}}
\end{align}
Substituting the two factor bounds into
\eqref{eq:app-peeling:holder-decomposition} gives
\begin{align*}
 \E_{\Covariate,\xi_i,\xi_j,\Covariate'}\!\left[
  |\psi(\Covariate,\Covariate')|\1\{Y_i>T,Y_j>T\}\right]
 &\leq K^{s+5/2}e^{-T^2/(2K)}
       \norm{\psi}_{L^2(\Covariate,\Covariate')}\ComponentWeight{2r}{j}\\
 &\leq e^{4s\log K-T^2/(2K)}
       \norm{\psi}_{L^2(\Covariate,\Covariate')}\ComponentWeight{2r}{j}.
 \alignnote{$s\geq1$}
\end{align*}

This proves \Cref{lem:peeling:pairwise-exceedance}.
\end{proof}

\end{document}